\documentstyle[epsf,bezier,12pt]{article}
\newcounter{ppp}
\newcounter{psi}
\newcounter{triple}
\newcounter{sigma}
\newcounter{pdtwo}
\newcounter{pdthree}
\newcounter{pdfour}
\newcounter{pdfive}

\newcounter{pdten}
\newcounter{pdeleven}
\newcounter{pdtwelve}

\newcounter{pdsixteen}

\newcounter{pdeighteen}

\newcounter{pdtwenty}
\newcounter{pdtwentyone}

\newcounter{pdtwentyfive}
\newcounter{pdtwentysix}

\newcounter{move}
\newcounter{adjust}

\newcommand{\area}{{\rm area}}
\newcommand{\Lab}{{\rm Lab}}
\newcommand{\yyy}{{\cal Y}}
\newcommand{\ttt}{{\cal T}}
\newcommand{\kkk}{{\cal K}}
\newcommand{\eee}{{\cal E}}
\newcommand{\aaa}{{\cal A}}

\newcommand{\bx}{{\bf X}(0)}
\newcommand{\ex}{{\bf E}(0)}
\newcommand{\fx}{{\bf F}(0)}
\newcommand{\bb}{{\cal B}}
\newcommand{\topp}{{\bf top}}
\newcommand{\bott}{{\bf bot}}
\newcommand{\vk}{van Kampen }

\newcommand{\Sym}{\hbox{Sym}}

\newcommand{\gam}{\Gamma}

\newcommand{\pleft}{p_{\hbox{left}}}
\newcommand{\pright}{p_{\hbox{right}}}
\newcommand{\ptop}{p_{\hbox{top}}}
\newcommand{\pbot}{p_{\hbox{bot}}}
\newcommand{\cee}{\alpha}
\newcommand{\dol}{\omega}
\newcommand{\iv}{^{-1}}

\newcommand{\bn}{{\bf N} }

\newcommand{\dti}{$\Delta$}

\newcommand{\rr}{{\cal R} }
\newcommand{\xxx}{{\cal X} }
\newcommand{\fff}{{\cal F} }

\newcommand{\ww}{{\cal W} }
\newcommand{\pp}{{\cal P} }
\newcommand{\qq}{{\cal Q} }
\newcommand{\sss}{{\cal S} }
\begin{document}

\title{Isoperimetric and Isodiametric Functions of Groups}
\author{Mark V. Sapir, Jean-Camille Birget, Eliyahu Rips\thanks{The research of the first
two authors was supported in part by NSF grants DMS 9623284, DMS 9203981.}}
\maketitle

\begin{abstract}
This is the first of two papers devoted to connections between asymptotic
functions of groups and computational complexity.
One of the main results of this paper
states that if for every $m$ the first $m$ digits of a real number $\alpha\ge 4$
are computable in time $\le C2^{2^{Cm}}$ for some constant $C>0$
then $n^\alpha$ is equivalent (``big O'')
to the Dehn function of a finitely presented group. The smallest
isodiametric function
of this group is $n^{3/4\alpha}$. On the other hand
if $n^\alpha$ is equivalent to the Dehn function of a finitely presented group
then the first $m$ digits of $\alpha$ are computable in time $\le C2^{2^{2^{Cm}}}$
for some constant $C$. This implies that, say, functions $n^{\pi+1}$,
$n^{e^2}$ and $n^\alpha$ for all rational numbers $\alpha\ge 4$
are equivalent to the Dehn functions of some finitely presented group and that
$n^\pi$ and $n^\alpha$ for all rational numbers $\alpha\ge 3$ are equivalent to
the smallest isodiametric functions of finitely presented groups.

Moreover we describe all Dehn functions of finitely presented groups $\succ n^4$
as time functions of Turing machines modulo two conjectures:
\begin{enumerate}
\item Every Dehn function is equivalent to a superadditive function.
\item The square root of the time function of a Turing machine is equivalent to
the time function of a Turing machine.
\end{enumerate}

\end{abstract}

\tableofcontents

\newtheorem{theo}{\quad Theorem}[section]
\newtheorem{lm}{\quad Lemma}[section]
\newtheorem{cy}{\quad Corollary}[section]
\newtheorem{df}{\quad Definition}[section]
\newtheorem{rk}{\quad Remark}[section]
\newtheorem{prop}{\quad Proposition}[section]
\newtheorem{prob}{\quad Problem}

\section{Preliminaries}

A function $f: \bn\to\bn$ is called an {\em isoperimetric function}
of a finite
presentation $\pp=\langle  X\ |\ R\rangle$ of a group $G$
if for every number $n$ and
every word $w$ over $X$ which is equal to 1 in $G$, $|w|\le n$, there exists a
\vk diagram over $\pp$ whose boundary label is $w$ and area $\le f(n)$;
in other words, $w$ is a product of at most $f(n)$ conjugates of
the relators from $R$ (see \cite{MadlenerOtto}, \cite{Gersten}, \cite{gromov1},
\cite{short}).

A function is called an isodiametric function of a finite presentation
$\pp=\langle  X\ |\ R\rangle$ if
for every number $n$ and
every word $w$ over $X$ which is equal to 1 in $G$, $|w|\le n$, there exists a
\vk diagram over $\pp$ whose boundary label is $w$ and diameter $\le f(n)$;
in other words, $w$ is a product of conjugates $x_i\iv r_ix_i$ of
the relators $r_i$ from $R$ and the lengths of the words $x_i$ are bounded
by $f(n)$ (see \cite{FHL}, \cite{Cohen}, \cite{Gersten2}).

The smallest isoperimetric function of a finite presentation $\pp$
is called the Dehn function of $\pp$.

Let $f, g:{\bf N}\to {\bf N}$ be two functions. We write $f \preceq g$ if there
exist non-negative constants $a, b, c, d$ such that $f(n)\leq ag(bn)+cn+d$.  All
functions $g(n)$ which we are considering in this paper grow at least as fast as
$n$. In this case $f(n)\preceq g(n)$ if and only if $f(n)\leq ag(bn)$ for some
positive constants $a, b$.  Two functions $f, g$ are called {\em equivalent} if
$f \preceq g$ and $g \preceq f$.

It is well known that Dehn functions and the smallest isodiametric functions
corresponding to different  finite presentations of the same group are
equivalent (see \cite{MadlenerOtto} or \cite{Alo}, \cite{Gersten2}).  This
allows us to speak about {\em the Dehn function} and the smallest isdiametric
function of a finitely presented group.  {\em Sometimes we shall also say that a
function $f(n)$ is the Dehn function (the smallest isodiametric function) of
group $G$ when in fact $f(x)$ is only equivalent to the Dehn function (the
smallest isodiametric function) of a presentation of this group.}

We shall not distinguish between equivalent functions in this
paper. Thus we do not distinguish, say, functions $2^n$ and $3^n$,
$n^\alpha$ and $[n^\alpha]$, $1$ and $2n$. But we distinguish $n^\alpha$
and $n^\beta$,
where $\alpha$ and $\beta$ are different numbers greater than 1. We also
distinguish, say $n^\alpha$ and $n^\alpha \log n$ for $\alpha\ge 1$.

It is known \cite{short} that every function $n^k$ where $k\ge 1$ is an integer
is the Dehn function of some finitely presented group. Some other functions like
$2^n$ are also Dehn functions. Some information about the class of Dehn
functions can be found in \cite{MadlenerOtto}, \cite{gromov}, \cite{Alo},
\cite{Gersten}, \cite{gromov1}, \cite{Gersten1} and other papers.  It is also
known that not every reasonable function is the Dehn function of some finitely
presented group. For example, as was shown by Gromov \cite{gromov}, there are no
Dehn functions which are strictly between $n$ and $n^2$.  Thus there exists a
gap between $n$ and $n^2$.
Information about isodiametric functions can be
found in \cite{FHL,Gersten2,Cohen}.

Several people including Hamish Short asked whether there are other gaps.
Answering another question of Short, M. Bridson proved that there are
finitely presented groups with Dehn functions equivalent $n^\alpha$ for some
non-integral rational numbers $\alpha$ \cite{Bridson}.

It is worth mentioning that a related class of {\em growth functions} of
finitely generated groups contains many gaps. For example, from a well known
theorem by Gromov \cite{gromov1} if the growth function of a finitely
generated group $G$ is bounded from above by a polynomial then $G$ has a
nilpotent subgroup of finite index. This easily implies that its growth
function is equivalent to $n^k$ for some natural number $k$. Thus there does
not exist a finitely generated group with growth function equivalent to,
say, $n^{4.5}$ or $n^{\pi+1}$.

In this paper, we shall prove that for every relatively fast computable
number $\alpha\ge 4$ there exists a finitely presented group with Dehn
function equivalent to $n^\alpha$ and the smallest isodiametric function
equivalent to $n^{3\alpha/4}$. On the other hand we shall prove that
if $n^\alpha$ is equivalent to the Dehn function of a finitely presented
group then $\alpha$ is a relatively fast computable number (see Corollary
\ref{cy304}). For example, $n^\alpha$ is a Dehn function of a finitely presented
group if $\alpha\ge 4$ is rational or if $\alpha=\pi+1$, etc.

We show that there exists a close relationship between
Dehn functions and complexity functions of Turing machines. This paper
is the first of two papers where we explore this relationship,
the next paper will be joint with A. Yu. Ol'shanskii.  One of the
main results of the present paper (Theorem \ref{th111}) says that every Dehn
function is a time function of some (not necessarily deterministic) Turing
machine.
The class of time functions of Turing machines is very large.  It
contains all relatively fast computable functions (see below).  So it is
very surprising that, as the other result of this paper states, there is
virtually no difference between the class of Dehn functions $\succ n^4$ and
the class of time functions $\succ n^4$.  As an immediate corollary we get
that the class of Dehn functions is very large.
For example all functions of the
form $n^\alpha(\log n)^\beta$ for $\alpha>4$ are Dehn functions of finitely
presented groups.

There is a conjecture that every isodiametric function is the space function of
some Turing machine.  Notice that it is not known that the set of space
functions $> n$ differs from the set of time functions $> n$.

In order to formulate the main results of our paper, we need some more
definitions.

Let $M$ be an arbitrary (deterministic or nondeterministic) Turing machine.  For
every natural number $n$ let $T(n)$ be the smallest number such that for every
acceptable word $w$ with $|w|\le n$ there exists a computation of length $\le
T(n)$ which accepts $w$. Notice that since $M$ may be nondeterministic, there
could be several accepting computations for the same initial configuration. The
function $T(n)$ is called the {\em time function} of the Turing machine $M$.

We call a function $f$ {\em superadditive} if for all natural numbers $m,
n$ we have $f(m+n)\ge f(m)+f(n)$. We do not know any Dehn function of a
finitely presented group which is not equivalent to a superadditive function
and Sapir conjectures that there are no such Dehn functions.  In
\cite{freeproducts} Guba and Sapir proved that every free product $G*H$
where $G$ and $H$ are non-trivial groups, has a superadditive Dehn function.
For example, for every group $G$ the free product $G*{\bf Z}$ has a
superadditive Dehn function (here ${\bf Z}$ is the infinite cyclic group).

The main results of our paper are the following.

\begin{theo} Every Dehn function of a finitely presented group is equivalent to
the time function of some (not necessarily deterministic) Turing machine.
\label{th111}
\end{theo}

\begin{theo} \label{1} Let ${\cal D}_4$ be the set of all Dehn functions
$d(n)\ge n^4$ of finitely presented groups. Let ${\cal T}_4$ be the set of time
functions $t(n)\geq n^4$ of arbitrary Turing machines \footnote{Recall that we
do not distinguish equivalent functions so ${\cal T}_4$ is actually the set of
functions which are equivalent to time functions of Turing machines.}.  Let
${\cal T}^4$ be the set of superadditive functions which are fourth powers
of time functions. Then $${\cal T}^4\subseteq {\cal D}_4\subseteq {\cal T}_4.$$
\end{theo}

This theorem is a corollary of Theorem \ref{th111} and the following result.

\begin{theo} \label{2} Let $L\subseteq X^+$ be a language accepted by a Turing
machine $M$ with a time function $T(n)$ for which $T(n)^4$ is superadditive.
Then there exists a finitely presented group $G(M) = \langle A \rangle $ with
Dehn function equivalent to $T(n)^4$, the smallest isodiametric function
equivalent to $T^3(n)$, and there exists an injective map  $K: X^+\to (A \cup
A^{-1})^+$ such that

\begin{enumerate}
\item $|K(u)|=O(|u|)$ for every $u\in X^+$;
\item $u\in L$ if and only if $K(u)=1$ in $G$;
\item $K(u)$ is computable in time $O(|u|)$ by a deterministic Turing machine.
\end{enumerate}
\end{theo}

It is not known whether ${\cal T}^4$ coincides (up to equivalence) with the set
of superadditive functions in ${\cal T}_4$. It is quite possible that these sets
coincide.  If this were true and in addition all Dehn function were
superadditive (we have discussed this conjecture above) then ${\cal D}_4$ would
be equal to the set of all superadditive functions from ${\cal T}_4$. This is
why we said before that the class of Dehn functions $\succ n^4$ and the class of
time functions $\succ n^4$ virtually coincide.

Using Theorem \ref{2} and a variant of the Aanderaa construction (see
\cite{Aanderaa} or \cite{Rotman}) one can give a structural description of
groups with word problem solvable in polynomial time.  Let $G= \langle A \rangle
$ be a finitely generated subgroup of a finitely presented group $H= \langle B
\rangle $ and $A\subseteq B$.  Then we define the {\em Dehn function
$d_{G,H}(n)$ of $G$ in $H$} as the smallest function $f(n)$ with the following
property:  for every number $n$ and every word $w$ over $A$ which is equal to 1
in $G$, $|w|\le n$, there exists a \vk diagram over $H$ whose boundary label is
$w$ and area $\le f(n)$; in other words, $w$ is a product of at most $f(n)$
conjugates of the relators from $H$.

It is well known and easy to see (\cite{Gersten}) that the word problem of a
finitely presented group is solvable if and only if the Dehn function is
recursive. Nevertheless the difference between the computational complexity of
the word problem and the Dehn function can be quite large (see
\cite{MadlenerOtto}).

Notice that if $d(n)$ is the Dehn function of a finitely presented group $G$
then there exists a ``trivial" (non-deterministic) machine which checks if a
word is equal to 1 in $G$. This machine simply inserts relators of $G$ into the
word until the word becomes 1. The time function of this machine is equivalent
to $d(n)$ (for the precise description of the ``trivial" machine and the proof
that its time function is equivalent to $d(n)$ see in the proof of Theorem
\ref{th111}).

Thus the word problem in a group $G$ can be decidable by a very fast machine but
at the same time the ``trivial" machine can be very slow.  The following
corollary shows that we can always embed $G$ into another finitely presented
group $H$ such that the trivial machine for $H$ solves the word problem in $G$
almost as fast as the fastest Turing machine can.

\begin{cy} \label{cym} The word problem in a finitely generated group $G$ is
solvable in polynomial time by a non-deterministic Turing machine if and only if
$G$ is embeddable into a finitely presented group $H$ such that the Dehn
function of $G$ in $H$ is bounded by a polynomial. Moreover, for every function
$T(n)\ge n$, if the word problem in a group $G$ is solvable in time $\le T(n)$
by a non-deterministic Turing machine then $G$ is embeddable into a finitely
presented group $H$ such that the Dehn function of $G$ in $H$ is bounded by
$O(T(n)^4)$.
\end{cy}

One can also view this corollary as a result about distortions of areas of words
in subgroups of finitely presented groups. By an {\em area} of a word $w=1$ in
$G$ we mean the smallest area of a \vk diagram with boundary label $w$. Then the
function $d_{G,H}(n)$ shows how distorted the areas of elements in the subgroup
$G$ of $H$ are.  Theorem \ref{th111} shows that the function $d_{G,H}(n)$ cannot
be lower than the computational complexity $T(n)$ of the word problem for $G$.
Corollary \ref{cym} shows that the distortion can reach as low as $T(n)^4$.

We are not going to prove Corollary \ref{cym} in this paper because in the next
paper (joint with A. Yu. Ol'shanskii) we shall prove the following more powerful
result.

\begin{theo} \label{olsh}
The word problem of a group is decidable in polynomial time
if and only if this group can be embedded into a group with polynomial Dehn
function. Moreover every group with the word problem solvable in time $T(n)$
can be embedded into a group with isoperimetric function $T(n)^5$.
\end{theo}

The following corollary of Theorem \ref{2} was the first result about Dehn
functions of groups proved by Rips. It was the start point of this work.

\begin{cy}\label{cy5} (Rips) There exists a finitely presented group with
undecidable conjugacy problem and Dehn function equivalent to $n^3$.
\end{cy}

The proof differs from the original proof of Rips and is given in the last
section of this paper.

The next corollary shows how to use Theorems \ref{1} and \ref{2} in order
to find many different Dehn functions and isodiametric functions.

With every function $f(n)$ one can associate two computational problems.
It is not very easy to define computability of functions by non-deterministic
machines. So when we talk about computability of functions,
we restrict ourselves to deterministic (multi-tape) Turing
machines (although one can replace them with
non-deterministic machines if one gives a ``right" definition of computability).

\vskip 0.1 in

{\it Problem A.  Given a natural number $n$ written in binary, compute $f(n)$
in binary. The size of $n$ is the number of digits of $n$, that is $[\log_2 n]+1$.}
\vskip 0.1 in

{\it Problem B. Given a natural number
$n$ written in unary (as a sequence of 1's), compute $f(n)$
in unary. The size of $n$ is $n$.}

It is easy to see that if $f(n)$ is a function such that Problem B is solvable
in time $O(f(n))$ then $f(n)$ is equivalent to the time function of a
deterministic Turing machine (a Turing machine which computes $f(n)$).

\begin{cy} Let $f(n)>n^4$ be a superadditive
function such that the binary representation of $f(n)$ is computable in
time $O(\sqrt[4]{f(n)})$ by a Turing machine (i.e., Problem
A is solvable in time $O(\sqrt[4]{f(n)})$).
Then $f(n)$ is equivalent to
the Dehn function of a finitely presented group and the smallest
isodiametric function of this group is equivalent to $T^{3/4}(n)$.
\label{cy1} \end{cy}

{\bf Remark.} Note that the size of $n$
in Problem A is $O(\log n)$. Thus if, say, Problem A for a function $f(n)$
is solvable in polynomial
time then the condition of the corollary holds and $f(n)$ is the Dehn function
of some finitely presented group and the smallest isodiametric function of this
group is equivalent to $T^{3/4}(n)$.

{\bf Proof.} By Theorem \ref{1} it is enough to show that
$g(n)=[\sqrt[4]{f(n)}]$ is equivalent to the time function of a Turing machine.

Consider the following Turing machine $M$. The input of $M$ is a natural
number $n$ written in unary. The machine has three tapes, with  $n$ (in unary)
initially written on tape 1, and tapes 2 and 3 are empty. Machine $M$ first
computes the binary representation of $n$.

The algorithm is well known:

\begin{center}
\begin{tabular}{l}
Until the number $n$ written on tape 1 is 0,\\
\qquad write the remainder of $n$ modulo $2$ on tape 2,
$[n/2]$ in unary on tape 3,\\
\qquad erase tape 1,\\
\qquad copy $[n/2]$ from tape 3 to tape 1 and erase tape 3\\
end of cycle.
\end{tabular}
\end{center}

This algorithm takes $O(n)$ steps. After this,  $M$ computes
$f(n)$ in binary. By assumption, this takes at most
$O(\sqrt[4]{f(n)})=O(g(n))$ steps.

Then $M$ computes $g(n)=[\sqrt[4]{f(n)}]$ in binary.
Note that for every number $m$ written in binary the number $[\sqrt[4]{m}]$
(in binary) can be computed in $O((\log m)^3)$ steps by the ordinary
bisection algorithm. So $M$ can compute the binary representation of
$g(n)$ from the binary representation of  $f(n)$ in $O((\log f(n))^3)$
steps.

Finally $M$ computes $g(n)$ in unary. The algorithm is opposite to the one
presented above. Again we use three tapes with $g(n)$ (in binary) initially
written on tape 1, and with tapes 2 and 3
empty.

\begin{center}
\begin{tabular}{l}
Write 1 on tape 2;\\
move the head one step right on tape 1;\\
until an empty square on tape 1 is to the right of the head, repeat:\\
\qquad copy the content of tape 2 twice on tape 3 and erase tape 2\\
\qquad copy the content of tape 3 on tape 2 and erase tape 3\\
\qquad if 1 is written to the right of the head on tape 1 then
write one more 1 on tape 2\\
\qquad move the head one step right on tape 1\\
end of cycle\\
\end{tabular}
\end{center}

This takes at most $O(g(n))$ steps. Thus the whole algorithm is executed in
time at most $O(n)+O(g(n))+O((\log f(n))^3) +O(g(n))=O(g(n))$.
This algorithm solves Problem B for $g(n)$. Therefore $g(n)$ is the
time function of
a Turing machine. By Theorem \ref{1},  $g(n)^4$ is equivalent to
the Dehn function of a finitely presented group.
Since $f(n)$ is equivalent to $g(n)^4$, we conclude that
$f(n)$ is equivalent to the Dehn function of a finitely presented group.
$\Box$

\vskip 0.1 in

This corollary allows one to construct finitely presented groups with
``arbitrary weird" Dehn functions and smallest isodiametric functions.
The following general fact follows from Corollary \ref{cy1}.

\begin{cy} \label{cy304}  For every real number $\alpha\ge 4$ such that
the first $m$ digits of $\alpha$ can be computed in time $\preceq 2^{2^m}$
(for every $m$)
the function $[n^\alpha]$ is equivalent to
the Dehn function of a finitely presented group and the smallest isodiametric
function of this group is $n^{3/4\alpha}$. On the other hand if
$n^\alpha$ is the Dehn function of a finitely presented group then
for every $m$ the first $m$ digits of $\alpha$ can
be computed deterministically
in time
 $\preceq 2^{2^{2^m}}$.
\end{cy}

{\bf Proof.} Notice that the function $[n^\alpha]$ is equivalent to the
function $2^{[\alpha[\log_2 n]]}$. The function $[\log_2 n]$ (i.e. the binary
expression of the number
of binary digits in $n$) is computable in time $\le O((\log_2 n)^2)$ by an
obvious algorithm: scan the number $n$ from left to right on one tape
and after each step add 1 to the number on the other tape.

Since the first $[\log_2(\log_2 n)]+1$ digits of $\alpha$ are computable in time
$O(n)$, the function $[\alpha[\log_2 n]]$ is computable in
time $O(n)\le O(n^{\alpha/4})$.

Notice also that
Problem A for a function equivalent to $2^m$ is solvable in time
$O(m)$. Indeed, we
can consider the unary expression of $m$ as the binary expression
of $2^{m+1}-1$
and use the second algorithm in the proof of Corollary \ref{cy1}.
It remains to apply Corollary \ref{cy1}.

Now suppose that $n^\alpha$ is the Dehn function of a finitely presented group. Then
by Theorem \ref{th111}, $n^\alpha$ is equivalent to the time function $T(n)$ of
some (non-deterministic) Turing machine $M$.
We can assume that when $M$
accepts, all its tapes are empty and that $M$ has only
one accept configuration
$c_0$. We can also assume that $M$ has an input tape which has number 1.

By Gromov's theorem, we can assume
that $\alpha\ge 2$.
Then for any $n>0$
$$\epsilon_1 n^\alpha\le T(n) \le \epsilon_2n^\alpha$$
for some positive constants $\epsilon_1$ and $\epsilon_2$.
Let number $n_0$ be such that $2^n>\log_2 \epsilon_2$ for every $n\ge n_0$.
Let $q=[\alpha]+1$.

Consider the following deterministic Turing machine
$M'$ which computes $T(n)$.
Let $k$ be the number of tapes of $M$. Then $M'$ will have $k+3$ tapes.
Tape $k+1$ will contain an input number $n$ in binary,
tape $k+2$ will contain $T(n)$ when $M'$ stops, tape $k+3$ is auxiliary.
In the input configuration,
$M'$ has number $n$ written on tape $k+1$ and all other tapes empty.
First of all it types the accept configuration of the machine $M$
on the first $k$ tapes and computes $n^q$ (in binary) on tape $k+3$.
Then it
considers all possible computations of length $\le n^q$ of the
machine $M\iv$ (commands of $M$ are applied backwards)
starting with the accept configuration $c_0$ of $M$. The number of such
computations is at most $r^{n^q}$ where $r$ is the number of commands of $M$.
For each of these computations $C$, it calculates the length $T(C)$
and the final word $W(C)$ written on the first tape.
It writes this information on tape $k+2$
provided $|W(C)|\le n$ and the final configuration of $C$ is an input configuration
of $M$. Thus it produces a sequence $\sss$, consisting of numbers $T(C)$ and words
$W(C)$. The words $W(C)$ in this sequence are all words
of length $\le n^q$
 accepted by $M$.
It is
clear that the time to produce this sequence and the length of it
does not exceed
$Dn^q r^{n^q}\le 2^{n^d}$ for some positive constants $D, d,r$.
Then for every word $W$ occurring in $\sss$, $M'$ calculates
the minimal $t(W)$ of the corresponding numbers $T(C)$ (for all computations $C$
such that $W=W(C)$). Then it finds the maximal number among all $t(W)$'s.
This number is obviously equal to $T(n)$. After that
$M'$ writes $T(n)$ on tape $k+2$, erases all other tapes and stops.
It is clear that in order to compute
$T(n)$ given the sequence $\sss$, $M'$ needs time at most
$|\sss|^2\le  2^{n^{d_1}}$ for some constant $d_1$. Thus $M'$ is a deterministic
machine which computes $T(n)$ in time $\le 2^{n^{d_2}}$ for some constant $d_2$.

Given this machine $M'$, consider the following Turing machine $M''$ which will
calculate the first $m$ digits of $\alpha$ (for every $m$). This machine has $k+4$ tapes
with tape $k+4$ being the input tape. It starts with number $m$ in binary
written on tape
$k+4$ and all other tapes empty. Then it calculates the number $n=2^{2^{m+n_0}}$
and
writes it on tape $k+1$ (using tape $k+2$ as an auxiliary tape and
cleaninig it after $n$ is computed). Then $M''$ turns on the machine $M'$
and produces $T(n)$ on tape $k+2$. Then it calculates $p=[(\log_2T(n)-\log_2\epsilon_1)/2^{n_0}]$
and writes it on tape $k+2$. Notice that $\alpha\log_2 n+\log_2\epsilon_1\le
\log_2 T(n)\le  \alpha\log_2 n+\log_2 \epsilon_2$. Therefore
$$[\alpha 2^m]\le p\le
\alpha 2^m+(\log_2\epsilon_2)/2^{n_0}.$$
Hence $p=[\alpha 2^m]$, so $p$ is the number formed by the first $m+k$ binary
digits of $\alpha$ where $k$ is the number of digits in $\alpha$ before
the period. From the construction of $M''$, it is clear that the time complexity
of $M''$ does not exceed $2^{2^{d_22^{m+n_0}}}\le 2^{2^{2^{d_3m}}}$
for some constant $d_3$.
$\Box$

Notice that rational numbers $\ge 4$ and fast computable irrational numbers
like $\pi+1$, $e^2$, $e\pi$ and others satisfy the conditions of Corollary
\ref{cy304} (see \cite{Bre}). It would be interesting to find out whether the
number $2^{2^{2^m}}$ in the second statement of Corollary \ref{cy304}
can be decreased to $2^{2^m}$.
This would give a necessary and sufficient
condition for $n^\alpha$ to be the Dehn function of a finitely presented group.

Corollary \ref{cy304} implies that not every polynomially bounded increasing
computable function $>n^2$ is the Dehn function of a finitely presented group.
 (just take a very slowly computable number
$\alpha$ and consider the function $n^\alpha$).
This answers a question of Gersten.

\bigskip

The plan of the paper is the following.
First we prove Theorem \ref{th111}: every Dehn function
of a finitely presented group is a time function of some Turing
machine. This implies one inclusion in Theorem \ref{1}.
In order to prove Theorem \ref{2} and the other inclusion in Theorem
\ref{1} we start with
an arbitrary multitape
nondeterministic Turing
machine with time function $T(n)$.
Our goal is to simulate a Turing machine in
a group. The mere fact that a Turing machine can be simulated by a group
is not new. Papers by Novikov and Boone showed it long ago. The problem that we
face in this paper is that we need the Dehn function of the group to be not
much bigger than the time function of the Turing machine.
In fact all known simulations of Turing machines in
groups (see \cite{KharSap})
lead to groups with exponential Dehn functions. The same is true if we use
Minsky machines instead of arbitrary Turing machines \cite{KharSap}.

In order to overcome this difficulty,  we introduce the concept of an
$S$-machine. An $S$-machine is a kind of Turing machine with one tape and many
heads (several heads can move at once). Heads are placed between cells on the
tape, so we consider their states as letters when we talk about the word written
on the tape.  A typical command of an $S$-machine is $q\to apb\iv$ where $q$ and
$p$ are states of some head and $a$ is a tape symbol. Notice that this command
does not depend on the content of the tape, only on the state of the head.  It
can be applied any time when the state of the head is $q$.  Thus usually
$S$-machines are extremely non-deterministic.  The alphabet of $S$-machine is
divided into subsets $Y$ and $Y\iv$ of the same size, so the words written on
the tapes can be considered as group words. After every step of the machine, the
word on the tape is automatically reduced.  We do not consider reducing the word
a separate step, it is a part of execution of a command.

For every Turing machine $M$ with time function $T(n)$ we find an $S$-machine
$\sss=\sss(M)$ which works in time $T(n)^3$ and recognizes almost the same
language.

Then we define a finitely presented
group $G_{N}(\sss)$ where $N$ is a natural number.
For $N=1$ this presentation resembles the presentation of Boone's group
(see \cite{Rotman}) but it does not have the Baumslag-Solitar type relations
$xyx\iv=y^2$ which make the Dehn function exponential.

In fact the $S$-machines are needed precisely to avoid such relations.

We shall prove that (for $N$ large enough), the Dehn function
of $G_{N}(\sss)$ is $T(n)^4$,
the smallest isodiametric function is $T^3(n)$
and that $G_{N}(\sss)$ satisfies the conditions of Theorem \ref{2}.

In order to do that, with every admissible word $W$ of $\sss$ we associate a
word $K(W)$ over the generators of $W$. It is relatively easy to show that if
$W$ is accepted by $\sss$ then there exists a ``standard" \vk diagram ({\em
disc}) with boundary label $K(W)$, so $K(W)=1$ in $G_N(\sss)$.  In order to
compute the Dehn function and the smallest isodiametric function
of $G_N(\sss)$ we consider an arbitrary \vk diagram
$\Delta$ over the presentation of $G_N(\sss)$. We define a planar graph
associated with this diagram, whose vertices are subdiscs of $\Delta$. For $N$
big enough the degree of each vertex of this graph is greater than 8, so  small
cancellation theory is applicable to this graph (see \cite{LS}). Using  small
cancellation theory, we decompose $\Delta$ in some natural way into simpler
subdiagrams which either do not have subdiscs or are discs themselves. This is
like decomposing a snowman into small snow balls. Then we show that the total
area of the diagrams without discs is bounded by $|W|^4$ and the sum of the
lengths of boundaries of discs is bounded roughly by $O(|W|)$. The area of a
disc with perimeter $n$ is bounded by the area of the corresponding computation
of $\sss$, that is by a function equivalent to $T(n)^4$ and its diameter
is equivalent to $T^3(n)$. Using the fact that
$T^4$ is superadditive, we conclude that the total area of the maximal discs
inside $\Delta$ is bounded from above by a function equivalent to $T(n)^4$.
This gives us the upper bound of the Dehn function.

To get the lower bound
we consider an arbitrary \vk diagram with boundary label $K(W)$
where $W$ is an admissible  word for $\sss$. We show that in this case $W$
is accepted by $\sss$ and that
the area (diameter) of this diagram cannot be smaller than a constant
times the
area (diameter) of a disc corresponding to an accepting computation for $W$.
This gives us the lower bound for the Dehn function and also proves that $K(W)=1$
implies that $W$ is accepted by $\sss$.

In the last section of the paper we prove Corollary \ref{cy5}.

\section{History}

This paper has a long and complex history. Several versions of this paper
containing gaps appeared as preprints during the last 2 years,
and many people have
these wrong versions. In this section, we try to explain what were the gaps
and how these mistakes were fixed. Without it several main ideas
of this paper would be unjustified. Another reason
for writing this section was that different authors had different input in
this paper.

We began working together on Dehn functions of groups during Rips'
visit to Lincoln in Summer 1994.  About a year before that Birget
formulated the ``embedding conjecture" that {\em every finitely presented group
with word problem solvable in polynomial time can be embedded into a finitely
presented group with polynomial Dehn function}. By that time the semigroup
analog of this conjecture was proved by Birget. It appeared later in preprint
\cite{Bir}. In this preprint, Birget also proved semigroup versions
of the main results of this paper about isoperimetric functions:
analogs of Theorems \ref{th111}, \ref{1}, \ref{2}.  He also presented
a sketch of a proof of our Theorem \ref{th111} for groups:
every Dehn function of a finitely presented group is equivalent to the time
function of some Turing machine (the complete proof presented here is due
to Sapir).

Rips' idea of proving the ``embedding conjecture" was
the following. Consider the classic Novikov-Boone-Higman-Aanderaa embedding
of a finitely generated group $G$ into a finitely
presented group $H$ \cite{Rotman}. The goal was to show that the
Dehn function of $H$ is not much bigger than the time complexity of a machine
$M$ solving the word problem in $G$. This idea was explored before in
Madlener and Otto \cite{MadlenerOtto} but the estimates of
the Dehn functions there
were too rough, they dealt only with Dehn functions which are well above
exponential. In fact the main problem is that the
``classic" Novikov-Boone-Higman-Aanderaa
construction produces groups with at least exponential Dehn function because
of the relations of the Baumslag-Solitar form $a\iv xa=x^2$.

Recall that the Novikov-Boone-Higman-Aanderaa embedding has two steps
\cite{Rotman}.
In the first step one takes a Turing machine $M$ and {\em simulates}
it in a group $B(M)$. This means that there exists a map $K$ from
the set of configurations of $M$ to the set of words in $B(M)$ such that
$K(c)=1$ in $B(M)$ if and only if the configuration
$c$ is accepted by $M$. In the second step
one uses the Turing machine $M$ solving the word problem of the group $G$
and several relatively simple HNN-extensions of $B(M)$ to produce a finitely
presented group $H$ containing $G$.

Rips' program was, first, by modifying
the Novikov-Boone construction get a group $B'(M)$ also simulating
$M$ with Dehn function
polynomially related to the time function of $M$ and, second,
show that  HNN-extensions of Aanderaa
does not change the Dehn function too much.

In order to compute the Dehn function of a group Birget and Rips wanted to use
geometric methods used by Rips in his proof of Corollary \ref{cy5} and in his
``geometrization" of the proof of Novikov-Boone theorem. This method counts the
number of bands and annuli of cells in a \vk diagram, see Section
\ref{forbidden} and Lemma \ref{nohub} of this paper. Notice that a similar
method was used in studying \vk diagrams over HNN-extensions by Miller and
Schupp in 1973 \cite{MS73}.

In 1994, Birget and Rips decided that by removing the Baumslag-Solitar relations
and modifying the Turing machine $M$ instead, they can still simulate the
machine $M$ in a group $B'(M)$ without exponential blow up of the Dehn function.
It also looked possible that the geometric methods will give not just an upper
bound of the Dehn function of the group $B'(M)$ but the Dehn function up to the
equivalence. This would give almost a complete description of all Dehn functions
of finitely presented groups as time functions of Turing machines with some
reasonable restrictions. This would solve a Short's problem of finding groups
with exotic Dehn function, which Sapir was working on at that time. So Sapir
joined Birget and Rips.

Trying to implement our program we constructed a modification $B'(M)$
of $B(M)$ and
proved several geometric lemmas similar to those in Section \ref{forbidden} were
proved.

Unfortunately these lemmas gave only an upper bound of the Dehn function which
was way above the lower bound that we could get. It was not even clear that
the Dehn function of $B'(M)$ is polynomial provided $M$ has a polynomial time
function. To compute the Dehn function up to equivalence, further
modification of $B'(M)$ (with many sectors)
and more geometric ideas (the snowman decomposition) were introduced and
implemented mainly by Sapir. This gave us results similar to Theorems
\ref{1} and \ref{2} of this paper only instead of exponent 4, we had
exponent 3. The paper was written in February 1996 and Sapir put a link
to it on his Web page. Many people downloaded the paper.

Immediately after that we started working on the second step of our program,
the Higman-Aanderaa construction. Birget and Sapir introduced a many sectors
modification of the Aanderaa construction and proved that this construction
gives a Higman-type embedding of $G$ into a finitely presented group $H$ which
was constructed as a sequence of HNN-extensions of $B'(M)$. But when we tried to
prove that the Dehn function of $H$ is not much bigger than the Dehn function of
$B'(M)$, we discovered that one of the geometric lemmas (Lemma 5.11 in the 1996
version of our paper) proved among the first results in 1994 was wrong.  In the
``proof" of this lemma we used a property of the Turing machine $M$ which we did
not have.  As a result not only the computation of the Dehn function of the
group $B'(M)$ was wrong but even the fact that $B'(M)$ simulates $M$ in the
sense mentioned above proved to be incorrect.

At first this problem seemed easy to fix by modifying the Turing machine $M$.
During seven months after that at least 20 different modifications of $M$ were
introduced. Each one of them seemed very promising at the beginning but turned
out to be not better than the original modification later.  The hope of ever
fixing the mistake almost faded when Sapir realized that although the group
$B'(M)$ does not simulate the Turing machine $M$, it still simulates a
``machine" related to $M$. These machines were called $S$-machines. After that
we faced a relatively standard problem in the theory of algorithms. We needed to
prove that $S$-machines are polynomially equivalent to ordinary Turing machines.
This turned out to be very non-trivial.  Sapir's proof of this fact presented in
this paper is probably the longest proof of equivalence of two computing devices
in the theory of algorithms.  After this proof was obtained, the geometric part
of the paper, -- both the snowman construction and the proof of the lower bound
of the Dehn function, -- needed modifications also. It was done by Sapir. He also
proved Corollary \ref{cy304} about Dehn functions of the form $n^\alpha$.  All
results about isodiametric functions in this paper are formulated and proved by
Sapir as well. Since the proofs of these results depend heavily on the technique
used in dealing with Dehn functions, it was unreasonable to write a separate
paper on isodiametric functions, so we included these results here.

Several comments by Ol'shanskii simplified the geometric part of our paper.  We
are grateful to him for these comments. We are also grateful to Victor Guba
whose comments helped in writing the section about $S$-machines.

Ol'shanskii used several ideas of the 1996 version of our paper in his paper
about distortions of subgroups in finitely presented groups \cite{Ol97}. He also
used a many-sector variant of the Aanderaa embedding. During his visit to
Lincoln in January 1997, he and Sapir realized that a combination of
$S$-machines and geometric ideas of this paper and ideas of \cite{Ol97} gives a
proof of Birget's ``embedding conjecture".  This completed the work started in
1994.

\section{Turing Machines}

In this section we collect all information about Turing machines that we need
in the proof of our main results. We also prove Theorem \ref{th111}.

We shall use the following standard notation for Turing machines.
A (multi-tape)
Turing machine has $k$ tapes and $k$ heads.
One can view it as a six-tuple $$M= \langle  X, \Gamma, Q,
\Theta, \vec s_1, \vec s_0 \rangle$$
where $X$ is the input alphabet, $\Gamma$
is the tape alphabet ($X\subseteq \Gamma$),
$Q=\bigcup Q_i, i=1,...,k$ is the set of states of the heads of the machine,
$\Theta$ is a set of
transitions (commands), $\vec s_1$ is the $k$-vector of start states, $\vec s_0$
is the $k$-vector of accept states.

The number of tapes $k$ of the machine is determined by
$\Theta$. We  assume that in the normal situation the
machine starts working with states of the heads forming the vector
$\vec s_1$,  with the head
placed at the right end of each tape, and
accepts if it reaches the state vector $\vec s_0$.
In general, the machine can be turned on in any configuration and turned off
at any time.

The leftmost square on every tape is always marked by $\cee$,
the rightmost square
is always marked by $\dol$.

The head is placed between two consecutive squares on the tape. When we
talk about the word written on the tape, we do not include
$\alpha$, $\omega$, and the state of the head.
We assume that the machine can insert
or delete squares on the tape but only at the right end (just before the
$\omega$-sign).

At every moment the head observes two squares on each tape.

A {\em configuration} of a tape of a Turing machine is a word $\alpha u q v \omega$
where $q$ is the current state of the head, $u$ is the word to the left
of the head and $v$ is the word to the right of the head.

A {\em configuration} $U$ of a
Turing machine is a $k$-tuple
$$(U_1, U_2, ...,U_k)$$
where $U_i$ is the configuration of tape $i$.
The length $|U|$ of this configuration is the sum of lengths of the words
$U_i$.

An {\it input configuration} is a configuration where the word written
on the first tape is in $X^+$, all other tapes are empty, the head observes
the right marker $\omega$, and the states form the start vector $\vec s_1$.
An {\em accept configuration} is any configuration
where the state vector is $\vec s_0$, the accept vector of the machine.

A transition of a Turing machine is determined
by the states of the heads and some of the $2k$ letters
observed by the heads. As a result of a
transition we can replace some of these $2k$ letters by other letters,
insert new squares in some of the tapes and move the head one square to the
left (right) with respect to some of the tapes.\footnote{In
different books, one can find different definitions of Turing machines, but
all these definitions are equivalent: the machines recognize the same languages,
and have equivalent complexity functions.}.

For example in a one-tape machine every transition is of the following
form:
$$uqv\to u'q'v'$$
where $u, v, u', v'$ are letters or empty words.
The only constraints are that no letter can be inserted to the left of $\alpha$
or to the right of $\omega$ and the end markers cannot be deleted.
This command means that
if the state of the head is $q$, $u$ is written to the left of $q$ and
$v$ is written to the right of $q$ then the machine must replace $u$ by $u'$,
$q$ by $q'$ and $v$ by $v'$.

For a general  $k$-tape machine a command is a vector
$$\{U_1\to V_1,...,U_k\to V_k\}$$
where $U_i\to V_i$ is a command of a 1-tape machine, the elementary commands
are listed in the order of tape numbers. In order to execute this
command, the machine checks if $U_i$ is a subword of the configuration
of tape $i$ ($i=1,...,k$), and then replaces $U_i$ by $V_i$.

Notice that for every command
$\{U_1\to V_1,...,U_k\to V_k\}$, the vector $\{V_1\to U_1,...,V_k\to U_k\}$
is also a command of a Turing
machine. These two commands are called {\em mutually inverse}.

A {\em computation} is a sequence of configurations $w_1,...,w_n$ such that for
every $i=1,..., n-1$ the machine passes from $w_i$ to $w_{i+1}$ by applying one
of the transitions from $\Theta$.  A configuration $w$ is said to be {\em
accepted} by a machine $M$ if there exists at least one computation which starts
with $w$ and ends with an accept configuration.

A word $u\in X^*$ is said to be {\em accepted} by the machine if the
corresponding input configuration is accepted.  The set of all accepted words
over the alphabet $X$ is called the {\em language accepted by the machine}.

Let  $C = (c_1, ..., c_g)$ be a computation of a machine $M$ such that for
every $j=1,...,g-1$ the configuration $c_{j+1}$ is obtained from $c_j$ by a
transition $r_j$ from $\Theta$.  Then we call the word $r_1 \ldots r_{g-1}$
the {\em
history} of this computation.  The number $g$ will be called the {\em time}
({\em duration}) of the computation.  Let $p_i$ ($i=1,...,g$) be the sum of the
lengths of the configurations of the tapes in the configuration $c_i$.
Then the sum of all $p_i$
will be called the {\em area of computation} and will be denoted by $\area(C)$.

If the set of transitions $\Theta$ is divided into a number of disjoint parts
$\bigcup_{i=1}^n\Theta_i$ then every computation $C$ can be represented as a
``concatenation" of computations
\begin{equation}
\label{concat}
C_1C_2....
\end{equation}
where each $C_i$ is a computation
of one of the $\Theta_j$, the last configuration of $C_i$ is the first
configuration of $C_{i+1}$ and $C_i$ are maximal blocks of $C$ where rules from
exactly one part $C_j$ are applied. We shall use such representations often when
a program of a Turing machine (and, later, a program of an $S$-machine) is
divided into subroutines.

With every Turing machine one can associate five functions: the {\em time
function} $T(n)$, the {\em space function} $S(n)$, the {\em generalized time
function} $T'(n)$, the {\em generalized space function} $S'(n)$, and the {\em
area function} $A(n)$. These functions will be called the {\em complexity
functions} of the machine.

For every natural number $n$ the number $T(n)$ is the minimal number $p$ such
that for every accepted input configuration $w$ (with $|w| \leq n$) there exists
at least one accepting computation of length at most  $p$. The number $S(n)$ is
the minimal natural number $p$ such that for every accepted input configuration
$w$ (with $|w| \leq n$) there exists at least one accepting computation
which contains only configurations of length $\le p$. The definitions of
$T'(n)$ and $S'(n)$ are similar but we consider arbitrary accepted
configurations $w$ of length $n$, not just the input configurations as in the
definitions of $T(n)$ and $S(n)$.  It is clear that  $T(n)\leq T'(n)$ and
$S(n)\leq S'(n)$ and it is easy to give examples where the inequalities are
strict. The area function $A(n)$ is defined as the minimal number $p$ such that
for every accepted configuration $w$ (with $|w| \leq n$) there exists at least
one accepting computation of area at most $p$.

We do not only consider {\em deterministic} Turing machines, for example,
we allow several transitions with the same left side. The following
proof of Theorem
\ref{th111} shows that nondeterministic machines arise naturally when one works
with the word problem.  In this proof we associate  with every group $G$ a
nondeterministic Turing machine $M$ such that $M$ accepts those and only those
words which are equal to 1 in $G$, and the time function of $M$ is equivalent to
the Dehn function of $G$. Actually the machine simulates the standard way of
deriving relations of a group from defining relations. The nondeterminism of the
machine reflects the fact that for every word there are usually many relations
applicable to this word.  Computations of the
machine $M$ correspond to derivations in the group $G$.  Although an idea of the
proof of this statement can be found in Birget \cite{Bir}, here we present the
first complete proof.

Let us recall the statement of Theorem \ref{th111}.
{\it The Dehn function of a finitely presented group is
equivalent to the time function of a two-tape Turing machine. The language
accepted by this machine coincides with the set of words equal to 1 in the
group. }

{\bf Proof of Theorem \ref{th111}.} Let
$G=\langle X \ |\ R\rangle$ be a finitely presented group.  We
assume that $R$ is closed under taking inverses and cyclic shifts.  Let $U$ be
the finite set of all pairs of words $(u,v)$ from $(X\cup X\iv)^*$ such that
$uv\iv$ is a relation from $R$. Consider a Turing machine $M$ with alphabet
$X\cup X\iv$, two tapes and the following program.

The program starts in state $q_1$ with the head at the right end on each
tape, a word $w$ written on the first tape and empty second tape. The head
always stays next to the right marker during the work of our Turing machine.
The alphabet is $X\cup X\iv$.
The program consists
of the following operations.

1 (move). If the word written on one of the tapes has letter $a$ at the right
end, the machine can erase the square containing this letter and insert a
square
with $a\iv$ in it on the other tape. This operation requires $1$ transition:
$\{aq_1\omega\to q_1\omega, q_1\omega\to a\iv q_1\omega\}$

2 (substitution).
If the word written on the first tape has suffix $u$ and there exists
a pair $(u, v)$ from $U$ then the machine can replace $u$ by $v$,
by deleting the $|u|$ squares on the first tape and inserting $|v|$ squares
on the first tape containing the word $v$. It is easy to see that
this operation requires at most $|u|+|v|$ transitions
of the forms $\{aq_i\omega\to q_j\omega, q_i\omega\to q_j\omega\}$,
$\{q_i\omega\to aq_j\omega, q_i\omega\to q_j\omega\}$.

3 (reduction). If the word written on the first tape and the
word written on the second
tape have the same rightmost tape letter, then the
machine can erase the squares
containing these letters. This operation requires one transition
of the form $\{aq_1\omega\to q_1\omega, aq_1\omega\to q_1\omega\}$

4 (accept). The machine accepts if both tapes are empty.

It is easy to prove by induction that after any number of these operations
the product of the word written on the first tape and the inverse
of the
word written on the second tape is equal to $w$ in the group $G$. Therefore
a word $w$ is accepted by our machine only if it is equal to 1 in $G$.

Conversely, if a word $w$ is equal to 1 in $G$ then there exist a derivation
$w\to w_1\to...\to 1$ such that at each step we either insert/delete a word
$r\in R$ or insert/delete a subword of the form $aa\iv$ where $a\in X\cup
X\iv$.
It is clear that there exists a computation of our Turing machine which
simulates this derivation. Therefore $w$ is accepted by the machine $M$.
Thus the language accepted by $M$ consists of all words which are
equal to 1 in $G$.

Let us prove that the time function $t(n)$ of the machine $M$ is equivalent to
the Dehn function $d(n)$ of the group $G$.
Let $(\alpha u_1p_1\dol, \alpha v_1p_1\dol), \ldots,
(\alpha u_m p_m\dol, \alpha v_mp_m\dol)$ be an accepting computation
of our machine where $u_1=w$, $p_1 = q_1$, $v_1=u_m=v_m=1$, $p_m =q_0$, and for every
$i=1, \ldots , m-1$ the machine passes from the $i$-th configuration
to the $i+1$-st configuration by executing one of the operations 1--4.
Then there exists a derivation $w\to w_1\to...\to 1$ of length at most $m$
where at each step one of relations from $R$ is applied. This implies that
there exists a \vk diagram with at most $m$ cells and boundary label $w$.
Thus $t(n)\ge d(n)$ ($n=|w|$).

On the other hand suppose that there exists a \vk diagram $\Delta$ over the
presentation of $M$ with boundary label $w$. Using this diagram we shall
construct an accepting computation of the machine
$M$ for the word $w$ of length at most
$d(|w|)+|w|$.

We shall call an edge of $\Delta$ a {\em bridge} if after removing this
edge the graph $\Delta$ becomes a union of two connected components (it is clear that we
won't get more than two connected components).

We shall describe a process of
``deconstructing" $\Delta$. For every $i=1, 2,...,$ we define a diagram
$\Delta_i$, a path $t_i$ consisting of bridges in $\Delta_i$, a vertex $v_i$, a
subdiagram $\Delta_i'$ of $\Delta_i$, and two words $x_i$ and $y_i$ forming a
configuration $\gamma_i$ of the machine $M$.  By induction on $i$ we shall
prove that for every $i\ge 1$:

\begin{enumerate}
\item [a)] $v_i$ is the terminal vertex of $t_i$ and initial vertex of
$\Delta_i'$, the initial vertices of $t_i$ and $\Delta_i$ coincide,
$t_i\iv$ is a subpath of $\partial(\Delta_i)$.
\item [b)] The path $t_i$ has no  common edge
with $\Delta_i'$, $t_i$ consists of bridges of $\Delta_i$.
\item [c)] If $i>1$ then the diagram $\Delta_i$ is
a subdiagram of $\Delta_{i-1}$.
\item [d)] The path $t_i$ is simple.
\item [e)] If $t_i$ is not empty then $\Delta_i'$ is a connected component
of the graph obtained from $\Delta_i$ by removing the last edge of
$t_i$.
\item [f)] If $\Delta_i$ has more than one vertex then
either $t_i$ is not empty or
$\Delta_i'$ is not empty.
\item [g)] If $i>1$ then $\gamma_i$ is obtained from $\gamma_{i-1}$
by one of the 3 operations of the machine $M$ (move, substitution, reduction).
\item [h)] $y_i=\Lab(t_i)$, $x_iy_i\iv$ is the label of the
contour of $\Delta_i$.
\end{enumerate}

Let $\Delta_1=\Delta$, let $v_1$ be the initial vertex of $\Delta$, let $t_1$ be
an empty path with the initial vertex $\imath(t)$ and terminal vertex $\tau(t)$
equal to $v_1$. Let $\Delta'_1=\Delta_1$. Then $y_i=1$ and $x_i=w$.  This means
that $\gamma_1$ is an initial configuration of the machine $M$. It is clear that
properties a)--h) hold.

Assume that for some $i\ge 1$ we have defined
$\Delta_i, \Delta'_i$, $t_i$ ($v_i$, $x_i$ and $y_i$ are determined by
$\Delta_i$, $\Delta_i'$ and $t_i$). If $\Delta_i$ consists of one vertex
then the process stops. So suppose that $\Delta_i$ has more than one vertex.

1 (move). Suppose first that $v$ is an initial vertex of
a bridge $e$ in $\Delta'_i$ such that $(t_ie)$ is a subpath of
$\partial(\Delta_i)$.
Then we set $v_{i+1}=\tau(e)$,
$t_{i+1}=t_ie$. Since $e$ is a bridge, removing this edge from $\Delta_i'$
divides $\Delta_i'$ into two connected components. The connected
component containing $\tau(e)$ will be denoted by $\Delta_{i+1}'$.
We also let $\Delta_{i+1}=\Delta_i$. Notice that the properties
a)--f) hold for the quadruple $(\Delta_{i+1}, \Delta_{i+1}', v_{i+1},
t_{i+1})$. The only non-obvious property is d), but it follows from
the fact that $e\in \Delta_i$ and that $t_i$ has no common edges with $\Delta_i$
(by the induction hypothesis). The configuration
$\gamma_{i+1}$ is obtained from $\gamma_i$ by a move. Indeed, notice
that the subpath $p$ of $\partial(\Delta)$
labelled by $x_i$ is not empty (otherwise $t_i$ would be a closed path,
a contradiction with d) ). Therefore
$e\iv$ is the last edge of $p$, so the label of $e$ is the last letter
of $w$. As a result of the operation, edge $e\iv$ is subtracted from $p$
and edge $e$ is added to $t_i$. Thus $(x_{i+1}, y_{i+1})$
is obtained from $(x_i, y_i)$
by removing the last letter of $x_i$ and adding the inverse of this letter
to $y_i$. This is the move operation. Thus properties g) and h) also
hold.

2 (substitution). Now suppose that $v_i$ is not an initial vertex of any bridge
$e$ in $\Delta_i'$
such that $(t_ie)\iv$ is a subpath of $\partial(\Delta_i)$
and suppose that there is an edge $e$ in $\Delta_i'$ with initial vertex $v_i$
and a cell
$\pi$ containing $e$. Let $\Delta_{i+1}$ (resp. $\Delta_{i+1}'$ )
be the result of removing $e$ from $\Delta_i$ (resp. $\Delta_i'$).
We also let $v_{i+1}=v_i$, $t_{i+1}=t_i$. It is clear that properties
a)--f) hold for the new triple $(\Delta_{i+1}, \Delta_{i+1}', t_{i+1})$.
The path $\partial(\Delta_{i+1})$ is obtained from
$\partial(\Delta_i)$ by replacing $e$ with the rest of the boundary of the cell
$\pi$. Thus $x_{i+1}$ is obtained from $x_i$ by replacing the last letter
$\Lab(e)$ with a word $r$ such that $\Lab(e)r\iv$ is a relation in $R$.
Thus $\gamma_{i+1}$ is obtained from $\gamma_i$ by a substitution.
Therefore properties g) and h) hold.

3 (reduction). Finally there exists a possibility that $\Delta_i'=\{v_i\}$.
By e), since $\Delta_i\ne \{v_i\}$,
$t_i$ is not empty. Let $e$ be the last edge of $t_i$. By e) $v_i$ is a
connected component of $\Delta_i\backslash \{e\}$. This means that the boundary
of $\Delta_i$ has a subpath $ee\iv$ or $e\iv e$. Then by removing the edge $e$
and the vertex $v_{i}$
from $\Delta_i$ we obtain a new diagram which will be denoted by $\Delta_{i+1}$.
Let $v_{i+1}$ be the initial vertex of $e$ and let $t_{i+1}$ be the result of
removing $e$ from $t_i$. Now if $t_{i+1}$ is empty then let
$\Delta_{i+1}'=\Delta_{i+1}$. If $t_{i+1}$ is not empty then let $e'$ be the
last edge of $t_{i+1}$. Since $e'$ was a bridge in the diagram $\Delta_i$, it is
a bridge in the subdiagram $\Delta_{i+1}'$. Therefore by removing $e'$ from
$\Delta_{i+1}$ we get a graph with two connected components. The connected
component containing $v_{i+1}$ is denoted by $\Delta_{i+1}'$. Again it is clear
that all properties a)--f) hold for the new triple $(\Delta_{i+1},
\Delta_{i+1}', t_{i+1})$. By assumption $x_i$ and $y_i$ have a common last
letter $\Lab(e)$. The words $x_{i+1}$ and $y_{i+1}$ are obtained from
$x_i$ and $y_i$ by removing this letter. This is a reduction. Therefore
properties g) and h) hold.

Let $E$ be the number of edges in $\Delta$. Let us prove that by using at most
$2E$ moves, substitutions and reductions, we can reduce $\Delta$ to a one-vertex
diagram. Indeed, at every step we either increase the path $t_i$ by one edge
or delete at least one edge from the diagram. All these edges belong to
the original diagram $\Delta$. Since $t_i$ is a simple path
the number of edges of $\Delta$ which belong to at least one of the $t_i$'s
does not exceed $E$. The number of reductions and substitutions
also cannot exceed $E$. Therefore after at most $2E$ steps the process
must end. But the process ends only when $\Delta_i$ consists of one vertex.
Therefore after at most $2E$ steps we reduce $\Delta$ to a one-vertex
diagram.

Recall that the triple $(\Delta_1, \Delta_1', t_1)$ corresponds
to the initial configuration $\gamma_1$ where $w=\Lab(\partial(\Delta))$
is written on the first tape and the second tape
is empty. By property h) for every $i>1$ the configuration $\gamma_{i+1}$
is obtained from $\gamma_i$ by one of the three operations (move, substitution,
reduction). Since in the last triple $\Delta_i$ is a one vertex diagram,
the last configuration has empty tapes, i.e. it is the accept configuration.
The number of operations needed to transform $\gamma_1$ into the accept
configuration is the same as the number of operations needed to ``deconstruct"
$\Delta$, that is at most $2E$.
Since each operation of $M$ translates into a
constant number of elementary transitions, there exists an accepting
computation of $M$
starting with $\gamma_1$ whose length does not exceed a constant times
$E$.  It is well known \cite{LS} that in any diagram
$\Delta$ over a finite presentation the number of edges does not exceed
the length of $\partial(\Delta)$ plus a constant multiple of the number of
cells. Therefore the time function $t(n)$ of $M$ does not exceed a constant
times $n+d(n)$, where $d(n)$ is the Dehn function of $G$.
It remains to notice that  $d(n)$ is equivalent to $C(d(n)+n)$, for any constant $C$.
The theorem is proved. $\Box$
\vskip 0.1 in

Notice that the machine that we constructed in this proof has two tapes and
it has no commands which move the head away from the endmarker $\omega$.
As we shall see later
every Turing machine
can be converted into a Turing machine recognizing the same language,
having the same complexity functions, and having no head moving
commands. On the other hand, it is easy to construct a one tape machine
which is equivalent (in the above sense) to the machine constructed in this
proof. One has to just concatenate the first tape of our machine
with the symmetric image of the second tape (of course, one would also
to remove the right endmarker of the first tape, the left endmarker on the second
tape and identify the state letters on the two tapes).
Thus {\em every Dehn function of a finitely presented group is equivalent to the
 time function of a one
tape Turing machine}.

In the next section we will show how to convert any Turing machine to an $S$-machine
The following lemma is a small
preliminary step in this conversion.

\begin{lm} \label{mach23} For every Turing machine $M$ recognizing a language $L$
there exists a Turing machine $M'$ with the following properties.
\begin{enumerate}
\item The language recognized by $M'$ is $L$.
\item $M'$ is symmetric, that is, with every command $U\to V$ it contains the
{\em inverse} command $V\to U$.
\item The time, generalized time, space and generalized space functions of
$M'$ are equivalent
to the time function of $M$. The area function of $M'$ is equivalent to
the square of the time function of $M$.
\item The machine accepts only when all tapes are empty.
\item Every command of $M'$ or its inverse
has one of the following
forms for some $i$
\begin{equation}\label{eq57}
\{q_1\omega\to q_1'\omega, ...,q_{i-1}\omega\to q_{i-1}'\omega, aq_i\omega\to q_i'\omega, q_{i+1}\omega\to q_{i+1}'\omega,...\}
\end{equation}
\begin{equation}\label{eq58}
\{q_1\omega\to q_1'\omega,...,q_{i-1}\omega\to q_{i-1}'\omega, \alpha q_i\omega\to \alpha q_i'\omega, q_{i+1}\omega\to q_{i+1}'\omega,...\}
\end{equation}
where $a$ belongs to the tape alphabet of tape $i$,
and $q_j, q_j'$ are state
letters of tape $j$. Thus if the head observes the right markers
at the beginning of a computation, it will observe the right markers during
the whole computation. If the head does not observe the right markers at
the beginning, then no
command is applicable, and so the computation is trivial.
\item The letters used on different tapes are from disjoint alphabets.
This includes the state letters.

\end{enumerate}
\end{lm}

{\bf Proof.}
First we show how to construct a machine $M'$ satisfying
the first four properties. Then we shall show how to transform it into
a machine satisfying all six properties.

The symmetrization of an arbitrary Turing machine is essentially
well known (see \cite{Bir}). We modify the construction from \cite{Bir} here.

Let $M$ be a $k$-tape Turing machine. We can assume without loss of generality
that the first tape of $M$ is the {\em input} tape, that is it
can contain letters
only from the input alphabet. If such a tape does not exist, we can always
add it to the machine without changing the complexity functions or
the language accepted by the machine.

Recall that in an {\em input configuration} of $M$, a word
is written on the first tape, all other tapes are empty, the head is next
to the right endmarker on each tape.

We can also assume that the machine $M$ has only one {\em accept command},
that is a command in which the states in the right sides form the
vector $\vec s_0$.

The machine $M'$ has $k+1$ tapes and
is a composition of three machines which we call {\em phases 1, 2 and 3}.
The tape alphabet of the
$k+1$-st tape is the
set of commands $\Theta$ of the machine $M$. The tape alphabets of
the other tapes are
the same as in $M$. The input alphabet of $M'$ coincides with
the input alphabet  of $M$.

In the first phase $M'$ inserts a
sequence of squares containing commands of machine $M$ on tape $k+1$.
The only constrain is that the first inserted square must contain the
accept command.

Thus after the first phase, the machine $M'$ has a sequence of commands of
$M$ written on the $k+1$-st tape. Then the machine checks
if all tapes except tapes 1 and $k+1$ are empty, and proceeds to the second
phase. This can be done by one transition of the form
\begin{equation}\label{eqin}
\{q_1\omega\to q_1'\omega, \alpha q_2\omega\to \alpha q_2'\omega,...,\alpha q_k\omega\to \alpha q_k'\omega,
q\omega\to q'\omega\}.\end{equation}

In the second phase
$M'$ tries to execute on the first $k$ tapes
the sequence of commands written on tape $k+1$, reading them
from the right to the left.
In order to make the machine $M'$ do this, we replace every command
$$\tau=\{U_1\to V_1,...,U_k\to V_k\}$$
of the machine $M$ by the command
$$\tau'=\{U_1\to V_1,...,U_k\to V_k, \tau q \to q\tau\}$$
where $q$ is the state of $M'$ on tape $k+1$ (this state is not changed
by any of the commands of $M'$). This command executes the transition $\tau$
on the first $k$ tapes of $M'$, checks if $\tau$ is written next to the left
of the head on tape $k+1$
(if $\tau$ is not written there, the command $\tau'$ is not executed),
and moves the head on tape $k+1$ one square to the left.

The second phase of the machine $M'$ is deterministic.  It is also {\em
injective}, that is for every configuration of phase 2 there is at most one
command of $M'$ whose inverse is applicable to this configuration.
Indeed, this command is determined by the letter on tape $k+1$ which is written
next to $q$ on the right.

If (and only if) it turns out that the sequence of commands written on tape
$k+1$ during the first phase is a valid history of an accepting computation of
machine $M$, that is if the head moves next to the left endmarker on tape
$k+1$, the machine $M'$ returns the head to the right endmarker on tape
$k+1$ and passes to the third phase.

In the third phase the machine erases one by one all squares on all tapes
and accepts. The third phase of $M'$ is deterministic.

Let $\Sym(M')$ be the machine obtained by {\em symmetrizing} $M'$, that is by
adding the inverses to all commands of $M'$. Consider an arbitrary reduced
computation of $\Sym(M')$. Since the second and the third phases of $M'$ are
deterministic and the first and the second phases are injective, every reduced
computation of $\Sym(M')$ can be represented in the form (\ref{concat}):
$$C=C_{3,1}\iv C_{2,1}\iv C_{1,1}\iv C_{1,2} C_{2,2} C_{3,2} C_{3,3}\iv
C_{2,3}\iv ....$$ where $C_{i,j}$ is a computation of phase $i$. Some prefix and
some suffix of this sequence may be empty, the other blocks are not empty.  Now
suppose that this computation is accepting. Then it contains a block
$C_{3,j}^{\pm 1}$ for some $j$. Let $\ell$ be the length of the first
configuration $c$ in $C_{3,j}^{\pm 1}$. One can apply to the configuration $c$
commands from phase 3 of $M'$ and get to an accept configuration in fewer than $\ell$
steps because each step on phase 3 removes one square of the tapes of $M'$.  Let
$C_{3,j}'$ be the corresponding computation. It is clear that the suffix of the
computation $C$ after the configuration $c$ has length not smaller than the
length of $C_{3,j}'$.  The area and the space of this suffix are also not
smaller than the area and space of $C_{3,j}'$.  So we can replace this suffix of
the computation $C$ by $C_{3,j}'$.  As a result we shall get a computation with
no greater length, area and space.  Thus we can assume that $C$ has only one
block of the form $C_{3,j}^{\pm 1}$ and of course (since $C$ is accepting) it
should have the form $C_{3,j}$.  This means that $C$ has the form $$C_{2,1}\iv
C_{1,1}\iv C_{1,2} C_{2,2} C_{3,2}$$ where the first 3 blocks may be empty. Let
$\ell$ be the length of the first configuration $c$ in $C$.

If the blocks $C_{2,1}, C_{1,1}, C_{1,2}$ are empty then the length of the
computation cannot exceed $2\ell$ because the length of the block $C_{2,2}$ does
not exceed the length of the word written on tape $k+1$ in $c$, and the length
of $C_{3,2}$ does not exceed $\ell$. The space of this computation does not
exceed $\ell$ because the length of configurations does not get bigger during
this computation.  For the same reason the area does not exceed $O(\ell^2)$.

Suppose that $C_{1,2}$ is not empty and let $c$  be the first configuration of
$C_{1,2}$. Then we can apply inverses of the commands in phase 1 and get to the
input configuration. Indeed, tapes $2,3,...,k$ in $c$ must be empty, and the
head must observe the right endmarker on each tape:  otherwise the machine would
not be able to execute the command (\ref{eqin}) and pass from $C_{1,2}$ to
$C_{2,2}$.  Also since the machine passes from $C_{2,2}$ to the third phase, the
first (from the left) square on tape $k+1$ must contain the accept command
during the computation $C_{2,2}$. Since the content of tape $k+1$ does not
change during the work of the second phase, the accept command must occupy the
first square of tape $k+1$ at the end of the computation $C_{1,2}$. Therefore by
applying inverses of the commands in phase 1, $M'$ can erase all squares on tape
$k+1$ and get to the input configuration of $M'$.

Let $(C_{1,2}')\iv$ be the corresponding
computation. Then $$C_{2,1}\iv C_{1,1}\iv (C_{1,2}')\iv C_{1,2}'C_{1,2}C_{2,2}C_{3,2}$$
is an accepting computation. The length and the space of the prefix
$C_{2,1}\iv C_{1,1}\iv (C_{1,2}')\iv$ is $O(\ell)$ and the area of this prefix is
$O(\ell^2)$. The first configuration $c$ of $C_{1,2}'$ is an input configuration
of $M'$ whose length is $\ell'\le \ell$. Since $M'$ accepts this configuration,
the sequence $s$ of commands written on the $k+1$-st tape after $C_{1,2}$ is
a history of accepting computation of the machine $M$. Therefore
$M$ accepts the word written on the first tape in $c$. Let $m$ be the length
of this history word $s$. From the description of $M'$, it is clear that
then the length of $C$ is at most $O(\ell)$+$O(m)$, the space is at most
$m+\ell$ and the area at most $O(m^2+\ell^2)$. Therefore if we replace $s$
by the history word of a shortest computation of the machine $M$ accepting
the input word of configuration $c$, the length, space and area of the
corresponding computation of $M'$ cannot be greater then the length, space and
area of $C$. Thus we can assume that $m=T(\ell')$, where $T$ is the time
function of the machine $M$. This implies that
the length, space and area of $C$ is, respectively, $O(T(\ell')+\ell)$, $O(T(\ell')+\ell)$,
$O(T^2(\ell')+\ell^2)$. Since $\ell'\le \ell$, we can conclude that
the time function, the generalized time  function, the space function and the generalized
function of $M'$ are equivalent to the function $T$,
and the area function of $M'$ is
equivalent to $T^2$.

We have also proved that every word accepted by $\Sym(M')$ is also accepted
by $M$. The converse statement is obvious (in fact every word accepted by $M$ is also
accepted by $M'$). This proves
that $\Sym(M')$ satisfies properties 1, 2, 3, 4
of the lemma.

Thus we can assume that the original machine $M$ satisfies properties 1 -- 4.

In order to get property 5, we divide every tape $i$ of $M$ into two tapes, and
number them $i$ and $i+1/2$.
If a configuration of the old tape $i$ is $\alpha u q v \omega$ then the configurations of the new tapes $i$ and $i+1/2$ will be,
respectively:

$$\alpha u q\omega$$
and
$$\alpha \bar v q\omega$$
where $\bar v$ is the word $v$ written from right to left. Then we replace
every command
$$\{u_1q_1v_1\to u'_1q'_1v'_1,..., u_kq_kv_k\to u_k'q_k'v_k'\}$$
by the command
$$\{u_1q_1\omega\to u'_1q'_1\omega, \bar v_1q_1\omega\to \bar v'_1q'_1\omega,...,u_kq_k\omega\to u_k'q_k'\omega, \bar v_kq_k\omega\to \bar v_k'q_k'\omega\}$$
Here $\bar v$ is either $v$, if $v$ is not the right marker, or $\alpha$
if $v=\omega$.
This vector of commands has $2k$ components listed in the order of the
tape numbers (1, 3/2, 2, 5/2,...).
It is clear that this machine recognizes the same language as $M$ and has the
same complexity functions.

Now one can replace each command $$\{w_1q_1\omega\to w'_1q'_1\omega,....\}$$
by a sequence of $2k$ commands of the form
\begin{equation}\label{eq5777}
\{q_1\omega\to q_1'\omega,...,
w_iq_i\omega\to w_i'q_i\omega,...\}
\end{equation}
(only the $i$-th component of this command has a non-state letter different from $\omega$).
This, of course, requires increasing the the set of states (each of these
new commands uses new state letters).
The application of each of these commands may change the word written
on only one tape. Together
these commands do what the initial command does. This transformation does
not change the language accepted by the machine, it does not decrease the
complexity
functions and can increase them not more than by a constant factor ($\le 2k$).
Thus we can assume that every command of $M$ has the form (\ref{eq5777}).

Notice that if $w_i=\alpha$ then the command
(\ref{eq5777})
is of the form (\ref{eq58}) because $w_i'$ is also equal to $\alpha$.
Otherwise we can replace this command by two commands
$$
\{q_1\omega\to q_1''\omega,...,
w_iq_i\omega\to q''_i\omega,...\}
$$
and
$$
\{q_1''\omega\to q_1'\omega,...,
q_i''\omega\to w_i'q'_i\omega,...\}
$$
where $q_j''$ are new state letters (which we add to the set of state letters).
Notice that these commands have the desired form (\ref{eq57}).
It is clear that the new machine recognizes the same language and each
 complexity function can increase by a factor of at most 2. All
 of the commands of the new machine have the
desired form (\ref{eq57}) or (\ref{eq58})

The last property is easy to get.
To get property 6,
one just needs to replace the alphabet of each tape by a disjoint copy of this
alphabet, and the set of state letters on each tape by a disjoint copy of this set.
Then of course one needs to change the commands of the machine accordingly.
$\Box$.

\section{$S$-machines}

In this section, we define $S$-machines and prove that $S$-machines are
``polynomially" equivalent to multi-tape Turing machines.

We shall define $S$-machines as rewriting systems \cite{KharSap}.
Let $n$ be a natural number. A {\em hardware} of an $S$-machine
is a pair $(Y,Q)$ where
$Y$ is an $n$-vector of (not necessarily disjoint) sets
$Y_i$, $Q$ is a $(n+1)$-vector
of disjoint sets $Q_i$, $\bigcup Q_i$ and $\bigcup Y_i$ are also disjoint.
The elements of $\bigcup Y_i$ are called {\em tape letters}, the elements of
$\bigcup Q_i$ are called {\em state letters}.

With every hardware $\sss=(Y,Q)$ we associate the
{\em language of admissible words} $L(\sss)=Q_1F(Y_1)Q_2...F(Y_n)Q_{n+1}$
where $F(Y_j)$ is the language of reduced group words in the alphabet
$Y_j\cup Y_j\iv$.
This language completely determines
the hardware. So instead of describing the hardware $\sss$, one can
describe the language of admissible words (which is in many cases more
convenient).

If $0\le i\le j\le n$ and $W=
q_1u_1q_2...u_nq_{n+1}$ is an admissible word then
the subword $q_iu_i...q_j$ of $W$ is called the $(Q_i,Q_j)$-subword of $W$
($i<j$).

An $S$-machine with hardware $\sss$ is a rewriting systems. The objects
of this rewriting system are all admissible words.

The rewriting rules, or {\em $S$-rules}, have the following form:
$$[U_1\to V_1,...,U_m\to V_m]$$
where the following conditions hold:
\begin{description}
\item Each $U_i$ is a subword of an admissible word starting with
a $Q_\ell$-letter and ending with a $Q_r$-letter (where $\ell=\ell(i)$
must not exceed $r=r(i)$, of
course).
\item If $i<j$ then $r(i)<\ell(j)$.
\item Each $V_i$ is also a subword of an admissible word whose $Q$-letters
belong to $Q_{\ell(i)}\cup...\cup Q_{r(i)}$ and which contains a $Q_\ell$-letter
and a $Q_r$-letter.
\item If $\ell(1)=1$ then $V_1$ must start with a $Q_1$-letter and if
$r(m)=n+1$ then $V_n$ must end with a $Q_{n+1}$-letter (so tape letters
are not inserted to the left of $Q_1$-letters and to the right of $Q_{n+1}$-letters).
\end{description}

To apply an $S$-rule to a word $W$ means to replace simultaneously subwords
$U_i$ by subwords $V_i$, $i=1,...,m$. In particular, this means that our rule is
not applicable if one of the $U_i$'s is not a subword of $W$. The following convention
is important:

{\bf After every application of a rewriting rule,
the word is automatically reduced. We do not consider reducing of an admissible
word a separate step of an $S$-machine.}

For example, if a word is $$q_1aaq_2bq_3ccq_4$$ and
$q_i\in Q_i$, $a\in Y_1$, $b\in Y_2$, $c\in Y_3$ and the $S$-rule is
\begin{equation}\label{rule1}
[q_1\to p_1a\iv,  q_2bq_3\to a\iv p_2b'q_3c],
\end{equation}
where $p_1\in Q_1$, $p_2\in Q_2$, $b'\in Y_2$, then the result of the
application of this rule is $$p_1p_2b'q_3cccq_4.$$

With every $S$-rule $\tau$ we associate the inverse $S$-rule $\tau\iv$
in the following way: if $$\tau=[U_1\to x_1V_i'y_1,\ U_2\to x_2V_2'y_2,...,
U_n\to x_nV_n'y_n]$$ where $V_i'$ starts with a $Q_{\ell(i)}$-letter and
ends with a $Q_{r(i)}$-letter, then $$\tau\iv=[V_1'\to x_1\iv U_1y_1\iv,\
V_2'\to x_2\iv U_2y_2\iv,...,V_n'\to x_n\iv U_ny_n\iv].$$

For example, the inverse of the rule (\ref{rule1}) is

$$[p_1\to q_1a,  p_2b'q_3\to aq_2bq_3c\iv].$$

It is clear that $\tau\iv$ is an $S$-rule,  $(\tau\iv)\iv=\tau$, and that
rules $\tau$ and $\tau\iv$ cancel each other (meaning
that if we apply $\tau$ and then $\tau\iv$, we return to the
original word).

The following convention is also important:

{\bf We always assume that an $S$-machine is symmetric, that is if an
$S$-machine contains a rewriting rule $\tau$, it also contains the rule
$\tau\iv$.}

As in the case of Turing machines, we can define the history of
a computation of an $S$-machine as the sequence (word)
of rules used in this computation.
A computation is called {\em reduced} if the history of this computation is
reduced, that is if two mutually inverse rules are never
applied next to each other.

As usual the {\em length} of a computation $W_1,...,W_n$ is $n$, the
{\em space} is max$\{|W_i|\ |\ i=1,...,n\}$, and the {\em area} is
$\sum_i|W_i|$.

We say that a computation of an  $S$-machine is {\em proper} if no {\bf new}
negative letters appear in a word during this computation.  We say that a computation is
{\em semiproper} provided a new negative letter inserted during one step of the
computation never disappears during the rest of the computation.

It is easy to see that if a computation $$C=(W_1, W_2,..., W_n)$$ is
semiproper then the inverse computation
$$C\iv=(W_n,...,W_2,W_1)$$ is also semiproper.

For example, let $k=2$,
$Y_1=\{\delta\}=Y_2$, $Q_1=\{\alpha\}$, $Q_2=\{q_1,q_0\}$,
$Q_3=\{\omega\}$. Suppose further that the $S$-machine $\cal S$ has just four
rules:
\begin{itemize}
\item[(1)] $q_1\to a\iv q_1a$,
\item[(2)] $\alpha q_1\to \alpha q_0$,
\end{itemize}
and their inverses.

Here is an example of a proper computation of this $S$-machine:
$$\alpha aaaq_1aa\omega \to \alpha
aaq_1aaa\omega \to \alpha aq_1aaaa\omega\to \alpha q_1aaaaa\omega\to \alpha
q_0aaaaa\omega.$$
The history of this computation is $(1)(1)(1)(2)$.

Here is an example of a semiproper computation:  $$\alpha
q_1\omega\to \alpha a\iv q_1a\omega\to \alpha a\iv a\iv q_1aa\omega.$$
The history of this computation is $(1)\iv(1)\iv$.
It is
also possible to prove that every reduced computation of this machine
is semiproper.

There exists a ``natural"
way to convert a Turing machine $M$ into an
$S$-machine $\sss$. The idea is simple. Take a Turing machine satisfying all
conditions of Lemma \ref{mach23}. Concatenate all tapes of the machine $M$
together and replace
every command $aq\omega\to q'\omega$ by $q\omega\to a\iv q'\omega$,
so that commands of $M$ become $S$-rules.
Unfortunately the $S$-machine $\sss$ constructed this way does not inherit
most of the properties of the original machine $M$.

Consider the following typical example.
Let $M$ be a one tape Turing machine with $X=Y=\{a\}$,
$Q=\{q_1, q_2, q_0\}$ and the following commands:
\begin{enumerate}
\item $aq_1\omega\to q_2\omega$
\item $aq_2\omega\to q_2\omega$
\item $\alpha q_2\omega\to \alpha q_0\omega$
\end{enumerate}
plus all inverses of these commands.
It is obvious that this machine accepts the word $a^n$ if and only if $n$ is
strictly positive. The corresponding $S$-machine has the following commands:

\begin{enumerate}
\item $q_1\omega\to a\iv q_2\omega$
\item $q_2\omega\to a\iv q_2\omega$
\item $\alpha q_2\omega \to \alpha q_2\omega$
\end{enumerate}
and all inverse commands. This machine accepts the word $\alpha a^nq_1\omega$
for every integer $n$. For example, let $n=-1$. Then the computation is the
following (we write the number of a command in parentheses; $i\iv$ denotes the
inverse of command number $i$):

$$\begin{array}{l}
\alpha a\iv q_1\omega \to (1) \\
\alpha a\iv a\iv q_2\omega \to (2)\iv \\
\alpha a\iv q_2\omega \to (2\iv)\\
\alpha q_2\omega\to (3)\\
\alpha q_0\omega
\end{array}
$$

In fact it is nontrivial to construct an $S$-machine which recognizes only
positive powers of $a$. Such a machine will be the
key block in the main construction of this section.

Let $\sss$ be an $S$-machine, let $W$ be and admissible word and
let $w_1,...,w_k$ be some words.  Then $\sss(W,w_1,...,w_k)$ is the statement
saying that there exists a computation of $\sss$ starting with $W$ and
ending with a word containing $w_1,...,w_k$. The set of all such reduced
computations will be denoted by $C\sss(W,w_1,...,w_k)$.  If this set has
only one element, this element will also be denoted by
$C\sss(W,w_1,...,w_k)$.  If $C$ is a computation then $\tau C$ is the last
word of $C$.

We shall need the following two general results. They have a
``diagram group nature" \cite{GS}.

\begin{lm}\label{gen}
Let $\sss$ be an $S$-machine, $W$ and $W'$ be admissible words and
let $W$ be positive. Suppose that there exists a computation of
$\sss$ connecting $W$ and $W'$. Then if every reduced computation starting
at $W$ is semiproper then every reduced computation starting at $W'$ is
semiproper. In addition, suppose that every reduced computation $C$
starting at $W$ has length
$\le f(|W|+|\tau C|)$ and
area $\le g(|W|+|\tau C|)$ for some functions $f$ and $g$.
Then every reduced computation starting with $W'$ has length
$\le f(|W|+|\tau C|)+f(|W|+|W'|)$ and
area $\le g(|W|+|\tau C|)+g(|W|+|W'|)$.
\end{lm}

{\bf Proof.} Let $C$ be a reduced computation connecting $W$ and $W'$ and let
$C'$ be a reduced computation starting with $W'$. Consider the concatenation of
computations $CC'$. Let us represent $C$ in the form $C_1D$ and $C'$ in the
form $D\iv C'_1$ where the computation $C_1C'_1$ is reduced. Then
the computation $C_1C'_1$ is semiproper since it starts at $W$.

Suppose that $C'$ is not semiproper. Then a negative letter inserted during a
step of $C'$ is deleted during a later step of $C'$. These two steps
cannot both occur in $C_1'$ because then $C_1C_1'$ would not be semiproper.
Thus the insertion occurs in $D\iv$.
But then $D$ contains a step removing a negative letter. Since $D$ is a suffix
of the computation $C=C_1D$ and $W$ is a
positive word, this negative letter should have been inserted and removed
during
the computation $C$. Thus $C$ is not semiproper,
a contradiction.

The facts about the length and the area can be proved similarly. Indeed, the
length (area) of $C'$ does not exceed the sum of the lengths (areas) of
$C$ and the reduced form of the computation $CC'$.
$\Box$

\begin{lm}\label{gen1} Let $\sss$ be an $S$-machine and $W$ and $W'$ be
admissible words. Suppose that there exists a computation connecting
$W$ and $W'$, $W'$ contains subwords $w_1,...,w_k$ and
$C\sss(W,w_1,...,w_k)$ contains only one computation.
Then $C\sss(W',w_1,...,w_k)$ contains only one (trivial) computation.
\end{lm}

{\bf Proof.} Indeed, suppose that there exists a non-trivial reduced
computation $C'$ from $C\sss(W',w_1,...,w_k)$.
Let $C$ be the reduced computation
connecting $W$ and $W'$. Then the reduced form of the computation
$CC'$ is equal to $C$. Therefore the reduced form of $C'$ is
trivial, a contradiction (we assumed that $C'$ is reduced and non-trivial).
$\Box$
\vskip 0.1 in
We shall construct eleven $S$-machines $\sss_1$ -- $\sss_9$, $\sss_\alpha$,
$\sss_{\omega}$ one by one.
Then we'll use these machines in order to construct an $S$-machine $\sss(M)$
simulating a Turing machine $M$.

The machine $\sss_1$ has the following hardware:
$Y(1)=(\{\delta\}, \{\delta\}, \{\delta\}, \{\delta\})$,
$Q(1)=(\{p_1,p_2, p_3\},\{q_1,q_2,q_3\},\{r_1,r_2,r_3\},
\{s_1,s_2,s_3\},\{t_1,t_2,t_3\})$.
The admissible words have the following form:
$$p\delta^{n_1}q\delta^{n_2}r\delta^{n_3}s\delta^{n_4}t$$
where $p, q, r, s, t$ may have indices 1, 2 or 3.

The program $P(1)$ consists of the following rules and their inverses.
\begin{enumerate}
\item $[q_1\to \delta^{-2}q_1\delta^2, r_1\to \delta\iv r_1\delta]$
\item $[p_1q_1\to p_2q_2, r_1\to r_2, s_1\to s_2, t_1\to t_2]$
\item $[p_1\delta q_1\to p_3\delta q_3, r_1\to r_3, s_1\to s_3, t_1\to t_3]$
\end{enumerate}

\begin{lm}\label{sss1}
(Machine $\sss_1$ divides a number by 2 and tells even from odd.)
Let $n\ge 0$ and $m$ be integers,
$W=p_1\delta^nq_1r_1s_1\delta^mt_1$ and $k=n/2$
if $n$ is even and $k=(n-1)/2$
if $n$ is odd. Let $W_2=p_2q_2\delta^kr_2\delta^ks_2\delta^mt_2$,
$W_3=p_3\delta q_3\delta^kr_3\delta^ks_3\delta^mt_3$.
Then
\begin{enumerate}
\item[(i)] Every computation of $\sss_1$ starting with $W$ is semiproper.
\item[(ii)] Every word in any computation $W=U_1\to U_2\to...$ has the form
$p_i\delta^{\ell_1}q_i\delta^{\ell_2} r_i\delta^{\ell_2} s_i\delta^mt_i$ where
$\ell_1+2\ell_2=n$, $i=1,2,3$.
\item[(iii)] $C\sss_1(W,q_1r_1)$ consists of one (trivial) computation.

\item[(iv)] $\sss_1(W,p_2)$ if and only if $n$ is even. In this case
$C\sss_1(W,p_2)$ consists of one computation of length $n/2$ and $\tau C\sss_1(W,p_2)=W_2$.
\item[(v)] $\sss_1(W,p_3)$ if and only if $n$ is odd. In this case
$C\sss_1(W,p_3)$ consists of one computation of length $(n-1)/2$
and $\tau C\sss_1(W,p_3)=W_3$.
\end{enumerate}
\end{lm}

{\bf Proof.}
(i) Indeed, let  $U_1\to U_2\to ...$ be a computation.
Since rules 2 and 3 and their inverses do not insert (delete) $\delta$, we can
assume that in this computation, we apply only rule 1 and its inverse.  Since
our computation is reduced, an application of rule 1 (resp. its inverse) in our
computation cannot be followed by an application of the inverse of rule 1 (resp.
the rule 1 itself).  Therefore, in this computation we apply either only rule 1
or only the inverse of rule 1. But rule 1 cannot remove a negative power of
$\delta$ inserted by previous applications of rule 1 (because rule 1 always
inserts $\delta\iv$ only to the left of $q_1$). Similarly the inverse of rule 1
cannot remove a negative power of $\delta$ inserted by a previous  application
of the inverse of rule 1. This proves that a negative letter which appears
during our computation cannot disappear, so our computation is semiproper.

(ii) can be easily proved by induction: every rule of $\sss_1$ preserves the
property from (ii).

(iii) Suppose that a reduced computation $C$ of $\sss_1$ takes a word $W$ to a
word $W'$ containing the subword $q_1r_1$. Then by (ii) $W'$ must be equal to $W$.
Since there is only one rule of $\sss_1$ which can be applied to a word
containing $p_2$ (resp. $p_3$) and our computation is reduced, we can assume
that all of the words in the computation $C$ contain $p_1$. Thus all rules
applied in the computation are 1 or 1$\iv$. Since the computation is reduced, we
can assume that only rule 1 is applied (then in the inverse computation only
rule 1$\iv$ is applied).  But each application of rule 1 inserts a new $\delta$
between $q_1$ and $r_1$, so it is impossible that the last word in the
computation contains the subword $q_1r_1$, a contradiction.

(iv), (v). The ``if'' parts of the first statement  in (iv) and (v) are obvious.
The ``only if'' part follows from (ii). Indeed (ii)  shows that if $W=U_1\to
U_2,...$ is a computation of $\sss_1$ then $U_j$ has the form
$p_i\delta^{\ell_1}q_i\delta^{\ell_2} r_i\delta^{\ell_2} s_i\delta^m t_i$ where
$\ell_1+2\ell_2=n$, Now if $n$ is odd then $\ell_1$ stays always odd and no
computation can take $W$ to a word containing the subword $p_1q_1$. But if we
never get this subword, we can never apply rule 2, so we never get a word
containing $p_2$. Similarly if $n$ is even, $\ell_1$ stays always even, so we
can never apply rule 3 and we can never get a word containing $p_3$.

Also from (ii), it follows that if $C$ is a reduced computation of $\sss_1$
starting with $W$ and $\tau C$ contains subword  $p_2q_2$ then $\tau C\equiv
W_2$; if $\tau C$ contains subword $p_3q_3$ then $\tau C\equiv W_3$ .

The fact that $C\sss_1(W,p_2)$ and $C\sss_1(W,p_3)$ contains at most one
computation immediately follows from (iii). Indeed, suppose that, say,
$C\sss_1(W,p_2)$ contains two computations $C_1$ and $C_2$.  By (iii) the end
words in these computations must be the same. So we can consider the
superposition $C_1C_2\iv$. This computation takes $W$ to $W$, so by (iii) the
reduced form of this computation must be trivial, so $C_1=C_2$.  $\Box$

\vskip 0.1 in
The hardware of the machine $\sss_2$ is the following:
$$Y(2)=Y(1), Q(2)=(\{p_1, p_2\},\{q_1, q_2\},
\{r_1, r_2\},\{s_1, s_2\},\{t_1, t_2\}).$$ The program $P(2)$ consists
of the following rules and their inverses:

\begin{enumerate}
\item $[q_2\to \delta q_2\delta\iv, s_2\to \delta\iv s_2\delta]$
\item $[q_2r_2s_2\to q_1r_1s_1, p_2\to p_1, t_2\to t_1]$
\end{enumerate}

The following lemma can be proved similarly to Lemma \ref{sss1}, so we omit
the proof.

\begin{lm}\label{sss2}
(Machine $\sss_2$ moves $q$ and $s$ toward each other until they meet at $r$.)
Let $k,\ell,m\ge 0$, $n$ be integers,
$W=p_2\delta^k q_2\delta^\ell r_2\delta^m s_2\delta^n t_2$.
Let $W_1=p_1\delta^{k+\ell}q_1r_1s_1\delta^{m+n}t_1$.
Then
\begin{enumerate}
\item[(i)]  Every reduced computation of $\sss_2$ starting with $W$
is semiproper.
\item[(ii)]  $C\sss_2(W,W)$ consists of one (trivial) computation.
\item[(iii)] $\sss_2(W,p_1)$ if and only if $k=\ell$. In this case $C\sss_2(W,p_1)$
consists of one computation and $\tau C\sss_2(W,p_1)=W_1.$
\end{enumerate}
\end{lm}

The machine $\sss_3$ is a {\em cycle} of machines $\sss_1$ and $\sss_2$.  One
can describe $\sss_3$ as the machine obtained by taking the union of $\sss_1$
and $\sss_3$ and identifying  two state vectors of $\sss_1$ with two state
vectors of $\sss_2$. This machine will divide a number (the exponent of $\delta$)
until it gets an odd number. In case when this exponent is 0, it will never
stop. Thus $\sss_3$ will tell non-zero numbers from zero.

More precisely, the
hardware $(Y(3), Q(3))$ of the machine
$\sss_3$ is the same as the hardware of $\sss_1$.
The program $P(3)$
consists of the following rules and their inverses. We divide the program into
two subprograms which we call {\em steps}:

\begin{description}
\item Step 1.
 \begin{description}
\item[1.1] $[q_1\to \delta^{-2}q_1\delta^2, r_1\to \delta\iv r_1\delta]$
\item[1.2] $[p_1q_1\to p_2q_2, r_1\to r_2, s_1\to s_2, t_1\to t_2]$
\item[1.3] $[p_1\delta q_1\to p_3\delta q_3, r_1\to r_3, s_1\to s_3, t_1\to t_3]$
\end{description}
\item Step 2.
\begin{description}
\item[2.1] $[q_2\to \delta q_2\delta\iv, s_2\to \delta\iv s_2\delta]$
\item[2.2] $[q_2r_2s_2\to q_1r_1s_1, p_2\to p_1, t_2\to t_1]$
\end{description}
\end{description}

Notice that Step 1 has the same rules as the machine $\sss_1$ and Step 2
has the same rules as the machine $\sss_2$.

\begin{lm}\label{sss3}
(The machine $\sss_3$ tells zero from non-zero.)
Let $n\ge 0$ and either $n=0$ or $n=2^e(2m+1)$. Let
$W=p_1\delta^n q_1r_1s_1t_1$.
Let
$W_3=p_3\delta q_3\delta^m r_3\delta^m s_3\delta^{n-2m-1}t_3$.
Then
\begin{enumerate}
\item Every reduced computation $C$ of $\sss_3$ starting with the word $W$ or $W_3$
is semiproper.
The length of $C$ does not exceed $O(|W|+|V|)$ and
the area does not exceed $O((|W|+|V|)^2)$ where $V=\tau C$.\footnote{In this paper
we write that $f(n)=O(g(n))$ if $C_1g(n)\le f(n)\le C_2g(n)$ for some positive
constants $C_1, C_2$ and every $n$. We write $f(n)\le O(g(n))$ if $f(n)\le Cg(n)$
for some constant $C$ and every $n$.}
\item Each of the sets  $C\sss_3(W,s_1t_1)$ and $C\sss_3(W_3, p_3)$
consists of one (trivial) computations
\item $\sss_3(W,p_3)$ if and only if $n>0$. In this case $C\sss_3(W,p_3)$
consists of one computation and $\tau C\sss_3(W,p_3)=W_3$. The length
of $C\sss_3(W,p_3)$ is $O(n)$ and the area is $O(n^2)$.
\end{enumerate}
\end{lm}

{\bf Proof.} In order to prove the lemma, we shall describe all possible
computations starting with the word $W$.  Notice that every computation $C$ of
our machine has a representation in one of the following forms
\begin{equation}
\label{comp1}
C_{1,1}C_{2,1}C_{1,2}C_{2,2}...
\end{equation}
or
\begin{equation}\label{comp2}
C_{2,1}C_{1,2}C_{2,2}\ldots .
\end{equation}
where $C_{i,j}$ is a computation
of step $i$, $i=1,2$. Without loss of generality we can assume that each
$C_{i,j}$ is non-trivial. The last word in each of these blocks
(except the last one) is the first word of the next block.

If $C$ consists of one block then all statements of our lemma follow from
Lemmas \ref{sss1} and \ref{sss2}. So suppose that $C$ has more than 1 block.

Case 1. Suppose that $n>0$ and $C$ has the form (\ref{comp1}).  If the
last word $U$ in $C_{1,1}$ contains $p_1$ then the first rule applied in
$C_{2,1}$ is rule 2.2$\iv$. This means that $U$ contains the subword
$q_1r_1s_1$. By Lemma \ref{sss1} (iii) in this case $C_{1,1}$ is trivial,
a contradiction. Therefore $U$ contains $p_2$.  By Lemma \ref{sss1}(iv),
$U$ has the form $$p_2q_2\delta^{n/2}r_2\delta^{n/2}s_2t_2$$ and $n$
is even.
Since the computation $C_{1,1}$ is semiproper (Lemma \ref{sss1}(i) ) and
$U$ is a positive word, this computation is proper. By Lemma
\ref{sss1}(iv) there is only one computation in $C\sss_1(W,p_2)$. It is
clear that the history of this computation is $(1.1)^{n/2}(1.2)$ (meaning
that first we apply rule (1.1) $n/2$ times and then rule (1.2) ).

Similarly, if a part $C_{1,i}$ of $C$ starts with a word
$$p_1\delta^kq_1r_1s_1\delta^\ell t_1,$$ $k\ge 0, \ell\ge 0$ and is
followed by $C_{2,i}$ then $k$ is even, the computation $C_{1,i}$ is
proper, the last word in $C_{2,i}$ is
$$p_2q_2\delta^{k/2}r_2\delta^{k/2}s_2\delta^\ell t_2$$ and the history of
$C_{1,i}$ is $(1.1)^{k/2}(1.2)$.

Now suppose that a part $C_{2,i}$ starts with a word of the form
$$U=p_2q_2\delta^{k}r_2\delta^ks_2\delta^\ell t_2$$ and there exists
$C_{1,i+1}$. Let $V$ be the last word in $C_{2,i}$.

Suppose that $V$ contains $p_2$. Then the first rule in $C_{1,i+1}$ is
rule 1.2$\iv$. Therefore $V$ contains $p_2q_2$. Thus $V$ has the form
$$p_2q_2\delta^{k_1}r_2\delta^{k_2}s_2\delta^{\ell'}t_2.$$ Since $r_2$ and
$p_2$ do not move during the computation of step 2, $k_1+k_2=2k$.  Since
$t_2$ also does not move, we have: $\ell'=\ell$.  Finally, it is easy to
prove by induction that during a computation of step 2 the difference
between the distance from $q_2$ to $r_2$ and the distance from $r_2$ to
$s_2$ stays the same. Since this difference is $0$ in the word $U$, we
have that $k_1=k_2$. This implies that $k_1=k_2=k$, $\ell'=\ell$, so
$U=V$. By Lemma \ref{sss2}(ii) $C_{2,i}$ is trivial, a contradiction.

Thus $V$ contains $p_1$. By Lemma \ref{sss2}(iii) $V$ must have the form
$$p_1\delta^{k}q_1r_1s_1\delta^{k+\ell}t_1.$$
Since this word is positive
and the computation is semiproper (Lemma \ref{sss2}(i)) the computation
is proper. By Lemma \ref{sss3}(iii) there exists only one computation
which takes $U$ to $V$. It is clear that the history of this computation is
$(2.1)^k(2.2)$.

Let $N$ be the number of blocks in the representation (\ref{comp1}) of $C$.
Combining the information obtained so far, we can conclude that
the history of the first $N-1$ blocks has the form:
$$(1.1)^{n/2}(1.2)(2.1)^{n/2}(2.2)(1.1)^{n/4}(1.2)... .$$

Now suppose that the last block $C'$ in (\ref{comp1}) is a computation
of step 1. By what we have proved before, $C'$ starts with the word
$$p_1\delta^kq_1r_1s_1\delta^{n-k}t_1$$ where $k=n/(2^c)$, $c=(N-1)/2$
(in this case $N$ is odd).
By Lemma \ref{sss1}(iv) if $C'$ ends with a word containing $p_2$ then
$k$ is even and the history of $C'$ is $(1.1)^{k/2}(1.2)$. If $C'$ ends
with a word containing $p_3$ then $k$ is odd, that is $k=2m+1$, $c=e$
(recall that $n=2^e(2m+1)$) and the history of $C'$ is $(1.1)^m(1.3)$.
If this computation ends with a word containing $p_1$ then the computation
$C'$ is semiproper and the history is $(1.1)^d$ for some (not necessarily
positive) integer $d$.

Thus we proved that if $V=\tau C$, $C$ has the
form (\ref{comp1}), $n>0$, and the last block in $C$ is a computation
of step 1, then $C\sss_3(U,V)$ contains only one
reduced computation of the form (\ref{comp1})
which is semiproper if $V$ contains $p_1$ and
is proper if $V$ contains $p_2$ or $p_3$. The word $V$ contains $p_3$
if and only if $V=W_3$.

Suppose now that $C'$ is a computation of step 2. Then it starts with the
word $$p_2q_2\delta^kr_2\delta^ks_2\delta^{n-2k}t_2$$ where $k=n/(2^c)$, $c=N/2$
(in this case $N$ is even).
If $C'$ ends with a word containing $p_1$ then it's
history has the form $(2.1)^k(2.2)$. Otherwise the history of $C'$ has the form
$(2.1)^d$ for some integer $d$.

Let again $V=\tau C$. Then if $C$ has the form (\ref{comp1}),
$n>0$, and the last block of $C$ is a computation of step 2,
then $V$ can contain
$p_1$ or $p_2$ (it cannot contain $p_3$),  $C\sss_3(U,V)$
contains only one reduced computation of the form (\ref{comp1})
and $C$ is semiproper if
$V$ contains $p_1$, and it is proper if $V$ contains $p_2$.

We have also proved that in this case the computation $C$  has
one of the following histories:

\begin{equation}\label{comp11}
(1.1)^{n/2}(1.2)(2.1)^{n/2}(2.2)(1.1)^{n/4}(1.2)...(1.1)^{k/2}(1.2)
\end{equation}
where $k=n/2^{c}, c=(N-1)/2$;
\begin{equation}\label{comp12}
(1.1)^{n/2}(1.2)(2.1)^{n/2}(2.2)(1.1)^{n/4}(1.2)...(1.1)^{m}(1.3)
\end{equation}
where the number of blocks is $2e+1$;
\begin{equation}\label{comp121}
(1.1)^{n/2}(1.2)(2.1)^{n/2}(2.2)(1.1)^{n/4}(1.2)...(2.2)(1.1)^{d}
\end{equation}
where the number of blocks is $\le 2e+1$;
\begin{equation}\label{comp13}
(1.1)^{n/2}(1.2)(2.1)^{n/2}(2.2)(1.1)^{n/4}(1.2)...(1.2)(2.1)^{k}(2.2)
\end{equation}
where the number of the blocks is $\le 2e$;
\begin{equation}\label{comp14}
(1.1)^{n/2}(1.2)(2.1)^{n/2}(2.2)(1.1)^{n/4}(1.2)...(1.2)(2.1)^{d}
\end{equation}
where the number of the blocks is $\le 2e$.

It is easy to check that each of these computations has length
$O(|W|+|V|)$ and area $O((|W|+|V|)^2)$ where $V$ is the last word
in the computation.

Case 2. Suppose that $n>0$ and $C$ has the form (\ref{comp2}). Then the first
rule applied in $C_{2,1}$ is (2.2)$\iv$ because the other command of step 2
is not applicable.

Therefore the second word in $C_{2,1}$ is
$$U_1=p_2\delta^nq_2r_2s_2t_2.$$
If the block $C_{2,1}$ ends with a word containing $p_1$ then by
Lemma \ref{sss2} this word must have the form
$$p_1\delta^nq_1r_1s_1t_1$$
and the computation $C_{2,1}$ is proper. This means that rules (2.1), (2.1)$\iv$
are not applied in $C_{2,1}$, so
the second rule
applied in $C_{2,1}$ is rule 2.2, and this computation is
not reduced, a contradiction.

Thus the last word in $C_{2,1}$ must contain $p_2$. Then the first rule
applied in
$C_{1,2}$ should be the rule (1.2)$\iv$. Thus the last word in $C_{2,1}$ must
contain $p_2q_2$. By Lemma \ref{sss2} the computation $C$ must therefore start
with an application of rule (2.2)$\iv$ and after
that rule (2.1)$\iv$ is applied $n$
times (each application of rule (2.1)$\iv$ moves $q_2$
one step closer to $p_2$).
Hence the last word in $C_{2,1}$ and the first word in $C_{1,2}$ is
$$V_1=p_2q_2\delta^nr_2\delta^ns_2\delta^{-n}t_2.$$

Now the first rule applied in $C_{1,2}$ must be (1.2)$\iv$.
So the second word in
$C_{1,2}$ must be $$p_1q_1\delta^nr_2\delta^ns_2\delta^{-n}t_2.$$ By Lemma
\ref{sss1} the computation $C_{1,2}$ is semiproper, so in the rest of $C_{1,2}$
only one rule applies (several times): (1.1) or (1.1)$\iv$.  If the rest of
$C_{1,2}$ consists of applications of rule (1.1) then the rest of $C_{1,2}$ cannot
have a word containing either $p_1q_1$ or $q_1r_1s_1$. Therefore $C$ does not
have block $C_{2,2}$.

Thus if $C_{2.2}$ exists then $C_{2,2}$ must start with
$$p_1\delta^{2n}q_1r_1s_1\delta^{-n}t_1$$ and the history of $C_{1.2}$ has the
form $(1.2)\iv(1.1)^{-n}$.

Continuing in this manner (induction on the length of $C$) one can prove that
the history of computation $C$ has one of the following forms:

\begin{equation}\label{com21}
(2.2)\iv(2.1)^{-n}(1.2)\iv(1.1)^{-n}(2.2)\iv(2.1)^{-2n}(1.2)\iv(1.1)^{-2n}...
(2.2)\iv(2.1)^k
\end{equation}
or
\begin{equation}\label{com22}
(2.2)\iv(2.1)^{-n}(1.2)\iv(1.1)^{-n}(2.2)\iv(2.1)^{-2n}(1.2)\iv(1.1)^{-2n}...
(1.2)\iv(1.1)^k
\end{equation}
where $k$ is an integer.

This immediately implies that the length of $C$ is $O(|W|+|V|)$ and
the area of $C$ is $O((|W|+|V|)^2)$ where $V=\tau C$.

Also if $V=\tau C$, $n>0$, and $C$ has the form (\ref{comp2})
then $V$ cannot contain $p_3$ and $C\sss_3(U,V)$ consists of one semiproper
computation. It is clear that this $V$ cannot be reached from $W$
by any computation of the form (\ref{comp1}). It is also clear that any word
reachable by a computation of the form (\ref{comp2}) cannot be reached by
a computation of the form (\ref{comp1}).

Comparing our results in cases 1 and 2, we can conclude that if $V=\tau C$,
$n>0$ then $C\sss_3(W,V)$ contains only one reduced computation, this
computation is semiproper if $V$ contains $p_1$ or $p_2$, and it is proper
if $V$ contains $p_3$. Moreover there exists only one reduced computation
in $C\sss_3(W,p_3)$. This computation has the form (\ref{comp1}) and the last
word in this computation is $W_3$.

Case 3. Finally consider the case when $n=0$. Then $W=p_1q_1r_1s_1t_1$. Then it is
easy to see using Lemmas \ref{sss1} and \ref{sss2} that every block $C_{i,j}$ in
the computation $C$, except the last one, has only two words, so the history of
each block is of length one and the command used in each of the blocks except
for the last one is either $(1.2)^{\pm 1}$ or $(2.2)^{\pm 1}$. Only the indices of
letters change during the computation $C$ until the last block. Thus the last
block starts either with $p_1q_1r_1s_1t_1$ or with $p_2q_2r_2s_2t_2$. In the
first case the last block either has two words and the rule applied in the last
block is 1.2, 2.2$\iv$, or it contains more than 2 words and the rule
applied in the last block is 1.1$^{\pm 1}$. In the second case, either the last
block has 2 words and the rule is 2.2 or 1.2$\iv$, or it contains more than two
words and the rule applied in this block is 2.1$^{\pm 1}$. Thus the computation
$C$ has one of the following forms:

\begin{equation}\label{com31}
((1.2)(2.2))^k
\end{equation}
for some integer $k$;
\begin{equation}\label{com32}
((1.2)(2.2))^k(1.1)^m
\end{equation}
for some integers $k$ and $m$;
\begin{equation}\label{com33}
((1.2)(2.2))^k(1.2)^{\pm 1}
\end{equation}
for some integer $k$;
\begin{equation}\label{com331}
((1.2)(2.2))^k(1.2)^{\pm 1}
\end{equation}
for some integer $k$.

Since $k$ and $m$ are $O(|V|)$,
in each case the length of the computation is $O(|W|+|V|)$ and the area
is $O((|W|+|V|)^2)$.

This completes the description of $C$.

Thus if $n=0$ then the set of computations $C\sss_3(W,p_3)$ is empty,
and every computation starting with $W$ is semiproper.

Now all statements of the lemma can be proved easily.

(i) We have proved that every computation starting with the word $W$
is semiproper and the length and area of it satisfies condition 1 of the
lemma. The fact that every computation starting with $W_3$ satisfies the same
properties follows from Lemma \ref{gen}.

(ii) The first part follows from the description of all computations starting
at $W$. The second part follows from Lemma \ref{gen1}.

(iii) This has been established above. The length of $C\sss_3(W,p_3)$ is
$O(n/2+n/4+...)=O(n)$ and the area is $O(n^2/2+n^2/4+...)=O(n^2)$.
$\Box$
\vskip 0.1 in
The machine $\sss_4$ is a {\em concatenation} of two copies of $\sss_3$
with common states $p_3,q_3, r_3, s_3$ and $t_3$.
More precisely, take a copy $\sss'_3$
of the machine $\sss_3$.
The machine $\sss'_3$ is obtained
by adding $'$ to all state letters of $Q(3)$ except $p_3, q_3, r_3, s_3$ and $t_3$.
Thus $\sss'_3$ has states
$Q'(3)=\{p'_1, p'_2, p_3\}\cup\{q'_1, q'_2,q_3\}\cup\{r'_1, r'_2,
r_3\}\cup\{s'_1, s'_2, s_3\}\cup\{t'_1, t'_2, t_3\}$.

Now take the union $Q(4)$ of $Q(3)$ and $Q'(3)$
and consider the union $P(4)$ of programs $P(3)$ and $P'(3)$
of the machines $\sss_3$ and $\sss_3'$.
The new machine with the hardware
$(Y(1), Q(4))$ and program $P(4)$ will be denoted by $\sss_4$.

\begin{lm}
\label{sss4} (The machine $\sss_4$ tells zero from non-zero and returns all
state letters to their original positions.)
Let $W=p_1\delta^nq_1r_1s_1t_1$, $n\ge 0$,
$W'=p'_1\delta^nq'_1r'_1s'_1t'_1$.
\begin{enumerate}
\item[(i)] Every reduced computation $C$ of $\sss_4$ starting with $W$ or $W'$ is
semiproper.
The length of $C$ does not exceed $O(|W|+|V|)$ and
the area does not exceed $O((|W|+|V|)^2)$ where $V=\tau C$.
\item[(ii)] Each of the sets $C\sss_4(W,s_1t_1)$, $C\sss_4(W', s_1't_1')$
consists of one (trivial) computation.
\item[(iii)]
 $\sss_4(W,s_1't_1')$ if and only if $n>0$. In this case
$C\sss_4(W,s'_1t'_1)$ consists of one computation and
$\tau C\sss_4(W,s'_1t'_1)=W'$. The length of this computation is $O(n)$.
\item[(iv)]
 $\sss_4(W',s_1t_1)$ if and only if $n>0$. In this case
$C\sss_4(W',s_1t_1)$ consists of one computation and
$\tau C\sss_4(W',s_1t_1)=W$. The length of this computation is $O(n)$ and
the area is $O(n^2)$.
\end{enumerate}
\end{lm}

{\bf Proof.} All statements of the lemma will again follow from a
description of computations of $\sss_4$ starting with $W$.

Every computation $C$ of $\sss_4$ can be represented in the form
$C_1C_2C_3....$ where $C_i$ is a computation of one of the machines
$\sss_3$ or $\sss_3'$, neighbor blocks are computations of different
machines, and each block is non-trivial.

Since $C$ starts with the word $W$, $C_1$ must be a computation of the machine
$\sss_3$.  If $C=C_1$ then $C$ is a computation of $\sss_3$ and we can apply
Lemma \ref{sss3}.  It is clear that in this case $\tau C$ cannot contain $p_1'$.

If $C\ne C_1$ then $C_2$ is a computation of $\sss'_3$. Therefore $C_2$ starts
with a word $U$ containing $p_3$. This word is the end word of $C_1$. So
statement $\sss_3(W,p_3)$ holds and by Lemma \ref{sss3}(iii) $n>0$, $U=p_3\delta
q_3\delta^mr_3\delta^ms_3\delta^{n-2m-1}t_3$ and there exists a unique reduced
computation of $\sss_3'$ which takes $U$ to $W'$. Since $U$ is a positive word,
the computation $C_1$ is proper (Lemma \ref{sss3}(i) ).  By Lemma \ref{sss3}(i)
the computation $C_2$ is semiproper.

Suppose that $C_3$ exists. Then it must be a computation of $\sss_3$. Therefore
it must start with a word containing $p_3$. By Lemma \ref{sss3}(ii)
then $C_2$ is trivial, a contradiction.

Thus $C$ consists of one or two blocks: $C=C_1$ or $C=C_1C_2$. In both cases
$C$ is semiproper. The computation $C$ can contain two blocks if and only if
$n>0$.

Now all statements of the lemma follow immediately from this description and
Lemmas \ref{gen} and \ref{gen1}.
$\Box$

\vskip 0.1 in
We shall need a copy of $\sss_4$, $\bar\sss_4$ obtained from $\sss_4$ by
putting $\bar{ }$ over all its state letters.
\vskip 0.1 in
Now take an arbitrary alphabet $Y$ and consider the following machine
$\sss_5$. The admissible words for $\sss_5$ have the following form:
$$EuxvFE'p\Delta_1q\Delta_2r\Delta_3s\Delta_4t\bar p\Delta_5\bar q\Delta_6\bar r\Delta_7\bar s\Delta_8\bar tF'$$
where $E,x,F,E',p,q,r,s,t,\bar p, \bar q, \bar r, \bar s, \bar t, F'$ are
state letters from different $Q$-components, $u\in (Y\cup Y\iv)^*$,
$\Delta_j\in \langle\delta\rangle$. If $z$ is one of the letters $p, q, r, s, t, \bar p, \bar q, \bar r,
\bar s, \bar t$ then $z$ belongs to a $Q$-component
$\{z_0,z_1, z_2, z_3, z_4, z_0', z_1', z_2'\}$. The alphabet of tape letters also includes
letter $x_4$ from the same component as $x$.

The program of $\sss_5$ consists of the following rules (one for each $a\in Y$)
and their inverses:

$\begin{array}{ll}
R(a)& \qquad   [x\to
a\iv xa,\ p_1'\to p_0',\ q_1'r_1's_1't_1'\bar p_0\to \delta\iv
q_0'r_0's_0't_0'\bar p_1\delta,\ \bar q_0\bar r_0\bar s_0\bar t_0\to
\bar q_1\bar r_1\bar s_1\bar t_1]
\end{array}$

\bigskip
The machine $\sss_6$ has the same hardware as $\sss_5$ and just one rule:

$\begin{array}{ll}
& \qquad   [p_0'\to p_1,\ q_0'r_0's_0't_0'\bar p_1'\to
q_1r_1s_1t_1\bar p_0,\ \bar q_1'\bar r_1'\bar s_1'\bar t_1'\to
\bar q_0\bar r_0\bar s_0\bar t_0]
\end{array}$

\bigskip
The machine $\sss_7$ has the same hardware as $\sss_5$ and $\sss_6$
and one rule:

$\begin{array}{ll}
& \qquad   [Ex\to Ex_4,\ p_1q_1r_1s_1t_1\bar p_0\to p_4q_4r_4s_4t_4\bar p_4,\
\bar q_0\bar r_0\bar s_0\bar t_0\to
\bar q_4\bar r_4\bar s_4\bar t_4]
\end{array}$

It is clear that every nontrivial reduced computation of $\sss_6$ or $\sss_7$ has
length 2.

The computations of $\sss_5$ are more complicated but also easy to
describe.

The proof of the following lemma is similar to the proofs of
the previous lemmas and is left to the reader.

We shall need one more definition. Let $\langle Y\rangle$ be the free group with
free basis $Y$ and let $\langle 1\rangle$ be the additive group of integers. Consider the
homomorphism from $\langle Y\rangle$ to $\langle 1\rangle$ which takes every element from $Y$ to 1.
Then the image of a word $w$ from $\langle Y\rangle$  under this homomorphism
will be denoted by $||w||$. Thus $||w||$ is the algebraic sum of the degrees
of the letters in $w$. For example $||ab\iv c||=1$.

\begin{lm} \label{sss5} The following statements hold.
\begin{enumerate}
\item[(i)] Any computation of $\sss_5$ has a history of the form
$$R(a_1)R(a_2)\iv R(a_3)...$$
or
$$R(a_1)\iv R(a_2)R(a_3)\iv... $$
where $a_i\in Y$.
If $z=a_1a_2\iv...$ or $z=a_1\iv a_2...$, $a_i\in Y$, then the computation
starting with a word $W$ and having the history corresponding to $z$
will be denoted by $C(z, W)$.
\item[(ii)] Every computation of $\sss_5$ is semiproper.
\item[(iii)] Suppose that
$W=EuxvFE'p\delta^nqrst\bar p\delta^m\bar q\bar r\bar s \bar tF'$
where $u,v \in (Y\cup Y\iv)^*$ and either
\begin{enumerate}
\item[(1)] $(p,q,r,s,t)=(p_1',q_1',r_1',s_1',t_1')$ and
$(\bar p, \bar q, \bar r, \bar s, \bar t)=(\bar p_0, \bar q_0, \bar r_0, \bar s_0, \bar t_0)$
\end{enumerate}
or
\begin{enumerate}
\item[(0)] $(p,q,r,s,t)=(p_0',q_0',r_0',s_0',t_0')$ and
$(\bar p, \bar q, \bar r, \bar s, \bar t)=(\bar p_1,\bar q_1, \bar r_1, \bar s_1, \bar t_1).$
\end{enumerate}
Let $C=C(z,W)$ be a computation of $\sss_5$ and $W'=\tau C$.
Then $$W'=
Euz\iv xzvFE'p\delta^{n-\epsilon}qrst\bar p\delta^{m+\epsilon}\bar q\bar r\bar s \bar tF'$$
and
\begin{enumerate}
\item[(0,0)] if $W$ contains $p_0'$ and  $W'$ contains $p_0'$ then
$z=a_ka_{k-1}\iv a_{k-2}...a_1\iv$ for some  even number $k$, $a_i\in Y$,
$i=1,...,k$, $\epsilon=0$. Such $z$ will be called words of type (0,0). Notice that
$||z||=0$
\item[(1,0)] if $W$ contains $p_1'$ and  $W'$ contains $p_0'$ then
$z=a_ka_{k-1}\iv...a_1$ for some  odd number $k$, $a_i\in Y$, $i=1,...,k$,
$\epsilon=1$.
Such $z$ will be called words of type (1,0). Notice that $||z||=1$.
\item[(0,1)] if $W$ contains $p_0'$ and  $W'$ contains $p_1'$ then
$z=a_k\iv a_{k-1}...a_1\iv$ for some odd number $k$, $a_i\in Y$, $i=1,...,k$,
$\epsilon=-1$.
Such $z$ will be called words of type (0,1). Notice that $||z||=-1$.
\item[(1,1)] if $W$ contains $p_1'$ and  $W'$ contains $p_1'$ then
$z=a_k\iv a_{k-1}...a_1$ for some  even number $k$, $a_i\in Y$, $i=1,...,k$,
$\epsilon=0$.
Such $z$ will be called words of type (1,1). Notice that $||z||=0$.
\end{enumerate}
Notice that in each of the four cases if $z$ is a word of the type $(i,j)$
then $||z||=i-j=\epsilon$.
\end{enumerate}
\end{lm}

We shall need three simple facts about words of type $(\ell,\ell')$,
$\ell,\ell'\in\{0,1\}$ introduced in the previous lemma.

\begin{lm} \label{type1} Let $z_j$, $j=1,...,k$, be an irreducible
non-empty word of
type $(\ell_j,\ell_j')$
where $\ell_j, \ell_j'\in\{0,1\}$ such that $\ell_{j+1}=1-\ell'_j$
for every $j=1,...,k-1$. Then the product $z_kz_{k-1}...z_1$ (in the free
semigroup) is an irreducible word.
\end{lm}

{\bf Proof.} Indeed, an easy inspection of the definition of the words
of type $(\ell,\ell')$, $\ell,\ell'\in \{0,1\}$, gives that
if $z_1$ is of type $(\ell,0)$ and $z_2$ is of type $(1,\ell')$ then
$z_2$ ends with a positive letter and $z_1$ starts with a positive letter.
If, on the other hand, $z_1$ is of type $(\ell,1)$ and $z_2$ is of the
type $(0,\ell')$ then $z_2$ ends with a negative letter and $z_1$ starts
with a negative letter. In both cases there are no cancellations in the
product $z_2z_1$. This immediately implies the statement of the lemma.
$\Box$

\begin{lm}\label{type2} Let $v$ be a positive word which does not end
with letter $a\in Y$. Then $va\iv$ cannot be represented
in the form $z_kz_{k-1}...z_1$ where $z_j$ is a word of type $(\ell_j,\ell_j')$,
$\ell_j, \ell_j'\in\{0,1\}$ such that $\ell_{j+1}=1-\ell'_j$ for every
$j=1,...,k-1$, and  $z_1$ is of type $(1,0)$.
\end{lm}

{\bf Proof.}  Let $v=a_1a_2...a_m$ where $a_i\in Y$, $i=1,...,m$ and $a_m\ne a$.
Suppose that $va\iv=z_kz_{k-1}...z_1$ where the $z_j$ satisfy the conditions of
our lemma. Since $z_1$ must be of type $(1,0)$, it cannot end with a negative
letter. Therefore the product $z_kz_{k-1}...z_1$ must have cancellations.
But this contradicts Lemma \ref{type1}.
$\Box$

\begin{lm}\label{type3} Let $v=a_ma_{m-1}...a_1$ be a word from $Y^+$ with
$a_j\in Y$. Suppose that $v = z_kz_{k-1}...z_1$ where $z_j$ is a word of type
$(\ell_j,\ell_j')$, $\ell_j, \ell_j'\in\{0,1\}$ such that $\ell_{j+1}=1-\ell'_j$
for every $j=1,...,k-1$.  Then $k=m$ and $z_j=a_j$ for every $j=1,...,k$.
\end{lm}

{\bf Proof.} Indeed, since the product $z_kz_{k-1}...z_1$ is reduced by
Lemma \ref{type1},
$z_j$ contains no negative letters. But it immediately follows from the
definition of words of type $(\ell,\ell')$ that
if $z_j$ is not a one letter word, it contains a negative letter.
Thus each $z_j$ is a one letter word: $m=k$ and $z_j=a_j$ for every $j=1,...,m$.
$\Box$

\bigskip
The machine $\sss_8$ is a cycle of the machines
$\sss_4$, $\sss_5$, $\bar\sss_4$,
$\sss_6$ and $\sss_7$.

More precisely, $\sss_8$ has the same hardware as $\sss_5$, $\sss_6$ and
$\sss_7$ and its program is the union of the programs of
$\sss_4$, $\sss_5$, $\bar\sss_4$, $\sss_6$ and $\sss_7$.

An informal description of $\sss_8$ is the following. Suppose that it starts
with a word $$W=EuxFE'p_1\delta^nq_1r_1s_1t_1\bar p_0\bar q_0\bar r_0\bar
s_0\bar t_0F'.$$ where $n\ge 0$.  In general the work of the machine can be
complicated, so we describe the ``ideal" computation assuming that $u$ is a
positive word and $n=|u|$.  The next lemma will show that this is the only
computation from $C\sss_8(W,p_4)$.  First $\sss_8$ checks whether $n>0$. If
$n=0$ then by our assumption $u=\emptyset$ and $\sss_8$ uses rules of
$\sss_7$ and
finishes with the word $$W_4=Ex_4FE'p_4q_4r_4s_4t_4\bar p_4\bar q_4\bar
r_4\bar s_4\bar t_4F'.$$ If $n>0$ then $\sss_8$ uses rules of
$\sss_5$ to move $x$ one
letter toward $E$, replaces $n$ by $n-1$ and inserts $\delta$ between $\bar p$
and $\bar q$.  Then $\sss_8$ uses rules of $\bar \sss_4$ to check that the word between
$\bar p$ and $\bar q$ is not empty (which is true), then it uses $\sss_4$
again to check if $n-1$ is not zero. If $n-1=0$ then $|u|=1$, and  $\sss_8$
uses $\sss_7$ and gets the word $$W_4=Ex_4uFE'p_4q_4r_4s_4t_4\bar
p_4\delta\bar q_4\bar r_4\bar s_4\bar t_4F'.$$ If $n-1>0$ then $\sss_8$ uses
$\sss_5$ again and moves $x$ one letter closer to $E$, replaces $n-1$ by $n-2$
and inserts another $\delta$ between $\bar p$ and $\bar q$. This process
repeats until $x$ reaches $E$ (and all $n$ letters $\delta$
move to the right of $\bar p$). After that $\sss_8$ uses $\sss_7$ and we obtain
the word $$W_4=Ex_4uFE'p_4q_4r_4s_4t_4\bar p_4\delta^n\bar q_4\bar r_4\bar
s_4\bar t_4F'.$$

One can see that this is certainly not the only computation which can start
with $W$. For example, at the beginning $\sss_8$ can apply $\sss_6$ (backward)
instead of $\sss_4$. Also after $\sss_5$, the machine can proceed with
using $\sss_4$ instead of $\bar\sss_4$ (if the computation of $\sss_5$ is
of the type (1,1) ), etc. A description of all computations of $\sss_8$ starting
with $W$ is contained in the following lemma.

\begin{lm} \label{sss8} ($\sss_8$ tells positive words in the alphabet $Y$
from almost positive words.)
Let
$$W=
EuxFE'p_1\delta^nq_1r_1s_1t_1\bar p_0\bar q_0\bar r_0\bar s_0\bar t_0F',$$
$$W_4=Ex_4uFE'p_4q_4r_4s_4t_4\bar p_4\delta^n\bar q_4\bar r_4\bar s_4\bar t_4F',$$
$n\ge 0$.
Then
\begin{enumerate}
\item[(i)] Every computation $C$ of $\sss_8$ starting with $W$ or $W_4$
is semiproper. The length of $C$ does not exceed
$O((|W|+|V|)^2)$ and
the area does not exceed $O((|W|+|V|)^3)$ where $V$ is the last word in the
computation $C$.
\item[(ii)] The set $C\sss_8(W, xF, s_1t_1)$ consists of one
(trivial) computation.
\item[(iii)] If $u$ is a positive word
then the statement $\sss_8(W,p_4)$ is true if and only if $n=|u|$.
In this case $C\sss_8(W,p_4)$
consists of one computation and $\tau C\sss_8(W,p_4)=W_4$. The length of this
computation is $O(n^2)$ and the area is $O(n^3)$.
\item[(iv)] If $u$ has the form $u'a\iv$ where $a\in Y$ and
$u'$ is a positive words, and $u'$ does not end with $a$, then for every
$n$ the statement
$\sss_8(W,p_4)$ is false.
\end{enumerate}
\end{lm}

{\bf Proof.} As before, we shall describe all computations of our machine
which begin with $W$.

Let $C$ be a reduced computation of our machine starting with $W$. It is clear
that $C$ has the form

\begin{equation}\label{eqsss8}
C_1C_2...C_N.
\end{equation}
where $C_i$ is a non-trivial computation of one of the machines
$\sss_4$, $\sss_5$,
$\bar\sss_4$, $\sss_6$ and $\sss_7$. Consecutive blocks should correspond
to different machines.

It is easy to see that only the last block can be a computation of $\sss_7$.

Let $U_i$ be the first word in the block $C_i$, $i=1,...,N$. Then there should
be rules of at least two machines among
$\sss_4, \sss_5, \bar\sss_4, \sss_6, \sss_7$ which are applicable to each of
these words $U_i$. Easy inspection (checking all 10 possible pairs of machines)
shows that $U_i$ must have the form
$$Ev_ixw_iFE'p\delta^{n_i}qrst\bar p\delta^{m_i}\bar q\bar r\bar s\bar tF'$$
where $v_i, w_i\in (Y\cup Y\iv)^*$, $n_i, m_i$ are integers,
and the word $pqrst\bar p\bar q\bar r\bar s\bar t$, which will be denoted
by $\chi(U_i)$, is equal to one of the following
words:

\begin{description}
\item $P_1=p_1q_1r_1s_1t_1\bar p_0\bar q_0\bar r_0\bar s_0\bar t_0$.
\item $P_2=p'_1q'_1r'_1s'_1t'_1\bar p_0\bar q_0\bar r_0\bar s_0\bar t_0$.
\item $P_3=p'_0q'_0r'_0s'_0t'_0\bar p_1\bar q_1\bar r_1\bar s_1\bar t_1$.
\item $P_4=p'_0q'_0r'_0s'_0t'_0\bar p'_1\bar q'_1\bar r'_1\bar s'_1\bar t'_1$.
\item $P_5=p_4q_4r_4s_4t_4\bar p_4\bar q_4\bar r_4\bar s_4\bar t_4$.
\end{description}

We are going to prove the following statements by induction on $i=1,...,N$.
All of them are trivial for $i=1$.

{\bf Fact 1.} $n_i\ge 0$.

{\bf Fact 2.} $m_i\ge 0$.

{\bf Fact 3.} $n_i+m_i=n$.

{\bf Fact 4.} $w_i=z_kz_{k-1}...z_1$, $v_i=uw_i\iv$ where each word $z_j$ ($j=1,...,k$)
is of type $(\ell_j,\ell'_j)$ and these types satisfy the following property:
$\ell_{j+1}=1-\ell'_{j}$, $j=1,...,k-1$. The number $k=k(i)$ is the number of
blocks $C_j$, $j<i$, which are computations of $\sss_5$.
If $C_N$ is a computation of
$\sss_7$ then the words $z_1, z_k$ are of type $(1,0)$.

Suppose that all these fact have been proven for all numbers $\le i$.

{\bf Proof of Fact 1.} If $C_i$ is a computation of one of the machines
$\sss_4$, $\bar\sss_4$, $\sss_6$ then $n_{i+1}=n_i$ (for $\sss_4$
and $\bar\sss_4$ this follows from Lemma \ref{sss4}, for $\sss_6$ it is obvious),
so $n_{i+1}\ge 0$
as desired. Thus the only non-trivial case is when $C_i$ is a computation
of $\sss_5$. It is easy to see that in this case
$\chi(U_i)=P_2$ or $\chi(U_i)=P_3$.

Suppose first that $\chi(U_i)=P_2$. There are only two machines from the list
$\{\sss_4$, $\sss_5$, $\bar\sss_4$, $\sss_6\}$ which have rules applicable to
words with $\chi=P_2$. These are $\sss_5$ and $\sss_4$. Therefore
$C_{i-1}$ is a computation of $\sss_4$. By Lemma \ref{sss4} then
$\chi(U_{i-1})=P_1$ and $n_{i-1}=n_i>0$. By Lemma \ref{sss5},
$n_{i+1}=n_i-\epsilon$ where $\epsilon\in\{0,1\}$. Therefore
$n_{i+1}\ge 0$.

Suppose now that $\chi(U_i)=P_3$. By Lemma \ref{sss5} then $n_{i+1}=n_i-\epsilon$
where $\epsilon\in\{0,-1\}$. Since $n_i\ge 0$ by assumption, $n_{i+1}\ge 0$.
This proves Fact 1.
//

{\bf Proof of Fact 2.} As in the proof of Fact 1, we can assume that $C_i$
is a computation of $\sss_5$. Again $\chi(U_i)=P_2$ or $\chi(U_i)=P_3$.

Suppose that $\chi(U_i)=P_2$. Then by Lemma \ref{sss5}, $m_{i+1}=m_i+\epsilon$
where $\epsilon\in\{1,0\}$, so $m_{i+1}\ge 0$.

Now suppose that $\chi(U_i)=P_3$. There are only two machines in our list
which are applicable to words with $\chi=P_3$. These are $\bar\sss_4$ and $\sss_5$.
Therefore $C_{i-1}$ is a computation of $\bar\sss_4$. Then by Lemma \ref{sss4},
$\chi(U_{i-1})=P_4$ and $m_{i-1}=m_i>0$. Since by Lemma \ref{sss5},
$m_{i+1}=m_i+\epsilon$ where $\epsilon\in\{0,-1\}$, we can conclude that
$m_{i+1}\ge 0$. This proves Fact 2.

{\bf Fact 3} is obvious because during the computation of $\sss_8$ the total
number of occurrences of $\delta$ stays the same.

{\bf Proof of Fact 4.} Machines $\sss_4$, $\bar\sss_4$, $\sss_6$ do not
change $v_i$ and $w_i$. By Lemma \ref{sss5},  $\sss_5$ replaces $x$
by $z\iv xz$ for some word $z\in (Y\cup Y\iv)^*$ of type $(\ell,\ell')$
where $\ell,\ell'\in\{0,1\}$. This shows that for every $i$,
$w_i=z_kz_{k-1}...z_1$, $z_j$ is of type $(\ell_j, \ell'_j)$,
and $u=v_iw_i$, where $k=k(i)$ is the number of computations
of the machine $\sss_5$ among the blocks $C_1,...,C_{i-1}$.
Let us prove that $\ell_{j+1}=1-\ell'_j$
for $j=1,...,k(i)-1$.

Again we assume that this statement is proved for all numbers $\le i$.
Since only $\sss_5$ can change $v_i$ or $w_i$,
we can assume that $C_i$ is a computation of $\sss_5$.

Let $w_i=z_kz_{k-1}...z_1$ be the above representation of $w_i$.

Suppose that $\chi(U_i)=P_2$. Then by Lemma \ref{sss5}, $v_{i+1}=v_iz\iv$,
$w_{i+1}=zw_i$ where $z$ is of type $(1,h)$ for some $h\in \{0,1\}$.
We need to prove that $\ell_k'=0$, that is $z_k$ is of type $(\ell_k,0)$.
Without loss of generality we can assume that $k(i)>0$.

As we mentioned in the proof of Fact 2, in this case $C_{i-1}$ is a computation
of $\sss_4$. Since $n_i\ge 0$ (by Fact 1), we can apply Lemma \ref{sss4}
and conclude that $\chi(U_{i-1})=P_1$ and $n_i>0$.
Since $k(i)>0$, $i-1>1$. There are only two machines applicable to
$U_{i-1}$. These are $\sss_4$ and $\sss_6$. Therefore $C_{i-2}$ is a computation
of $\sss_6$. Therefore $\chi(U_{i-2})=P_4$. Since $P_4\ne P_1$, $i-2>1$.
The block $C_{i-3}$ must be a computation of $\bar\sss_4$ because $\bar\sss_4$
and $\sss_6$ are the only machines in our list which have rules applicable to
words with $\chi=P_4$. Since $m_i\ge 0$, we can apply Lemma \ref{sss4}
and conclude that $\chi(U_{i-3})=P_3$. This implies that $i-3>1$ and $C_{i-4}$
is a computation of $\sss_5$. Since $C_{i-4}$ is the last computation of
$\sss_5$ among the blocks $C_j$, $j<i$, by our assumption, $C_{i-4}$ replaced
$x$ with $z_k\iv xz_k$ where $k=k(i)$. Since $U_{i-3}=\tau C_{i-4}$
and $\chi(U_{i-3})=P_3$, by Lemma \ref{sss5}, $z_k$ has type $(\ell,0)$
for some $\ell$. Therefore $\ell_k'=0$ as desired.

Now suppose that $\chi(U_i)=P_3$. In this case by Lemma \ref{sss5},
$C_{i}$ replaces $x$ by $z\iv xz$ where $z$ has type $(0,\ell')$
for some $\ell'\in\{0,1\}$. We need to prove that $\ell_k'=1$, that is
$z_k$ is of type
$(\ell_k, 1)$.

Arguing in the same way as before, we
conclude that $C_{i-1}$ is a computation of
$\bar\sss_4$, $\chi(U_{i-1})=P_4$; $C_{i-2}$ is a computation of $\sss_6$,
$\chi(U_{i-2})=P_1$; $C_{i-3}$ is a computation of $\sss_1$,
$\chi(U_{i-3})=P_2$; $C_{i-4}$ is a computation of $\sss_5$.
By our assumptions $C_{i-4}$ replaces $x$ by $z_k\iv xz_k$.
Since $U_{i-3}=
\tau C_{i-4}$ and $\chi(U_{i-3})$, by Lemma \ref{sss5}, the type of
$z_k$ is $(\ell_k,1)$ for some $\ell_k$, as desired.

Now let us prove that if $C_N$ is a computation of $\sss_7$ then
$z_1$ is of type $(1,0)$.
We have established
that if $C_i$ is the first block in (\ref{eqsss8}) which is a computation
of $\sss_5$ then $C_i$ replaces $x$ with $z_1\iv xz_1$. We have also
proved that $\chi(U_i)=P_2$ or $\chi(U_i)=P_3$, so $i>1$.

If $\chi(U_i)=P_2$ then
$C_{i-1}$ is a computation of $\sss_4$, $\chi(U_{i-1})=P_1$.
If $i-1>1$ then $C_{i-2}$ is a
computation of $\sss_6$, $\chi(U_{i-2})=P_4$ so $i-2>1$. Then $C_{i-3}$ is a computation
of $\bar\sss_4$, so $U_{i-3}=P_3$ and $i-3>1$. Then $U_{i-4}$ must be a
computation of $\sss_5$. This contradicts the fact that $C_i$ is the first
block in (\ref{eqsss8}) which is a computation of $\sss_5$. Thus $i-1=1$ that
is $i=2$. Therefore $C_{i-1}=C_1$ is a computation of $\sss_4$.
Therefore by Lemma \ref{sss4},
$$U_2=EuxFE'p'_1\delta^nq'_1r'_1s'_1t'_1\bar p_0\bar q_0\bar r_0\bar s_0\bar t_0F'.$$
and $n>0$.

Suppose that $\chi(U_3)=P_2$. Then by Lemma \ref{sss5} $z_1$ has type (1,1)
and
$$U_3=Euz_1\iv xz_1FE'p'_1\delta^nq'_1r'_1s'_1t'_1\bar p_0\bar q_0\bar r_0\bar s_0\bar t_0F'.$$
Therefore $C_3$ is a computation of $\sss_4$ (it cannot be a computation of
$\sss_7$), so $N>3$. This implies that
$$U_4=Euz_1\iv xz_1FE'p_1\delta^nq_1r_1s_1t_1\bar p_0\bar q_0\bar r_0\bar s_0\bar t_0F'.$$
Since $n>0$, $C_4$ cannot be a computation of $\sss_7$. So $N>4$
and $C_4$ is a computation of $\sss_6$ and
$$U_5=Euz_1\iv xz_1FE'p_0'\delta^nq_0'r_0's_0't_0'\bar p_1'\bar q_1'\bar r_1'\bar s_1'\bar t_1'F'.$$
Hence $C_5$ should be a computation of $\bar\sss_4$ and
$U_6$ should be equal to $$Euz_1\iv xz_1FE'p_0'\delta^nq_0'r_0's_0't_0'\bar p_1\bar q_1\bar r_1\bar s_1\bar t_1F'.$$
But this is impossible by Lemma \ref{sss4} since $m_6=0$.
This means that $\chi(U_3)=P_3$, so by Lemma \ref{sss5}
$z_1$ has type (1,0) as desired.

Now suppose that $\chi(U_i)=P_3$. Then $C_{i-1}$ is a computation of
$\bar\sss_4$. Then $\chi(U_{i-1})=P_4$ and $m_i=m_{i-1}>0$ (by Lemma \ref{sss4}).
Then $i-1>1$ and $C_{i-2}$ is a computation of $\sss_6$. This implies that
$\chi(U_{i-2})=P_1$, $m_{i-2}=m_{i-1}>0$. Since $m_{i-2}>0$, $i-2$ cannot be equal
to $1$ (since $U_{i-2}\ne U_1=W$). Therefore $C_{i-3}$ exists and is a
computation of $\sss_4$. Then $\chi(U_{i-3})=P_2$ and $i-3>1$. Therefore
$C_{i-4}$ exists and is a computation of $\sss_5$. This contradicts the
assumption that $C_i$ is the first block in (\ref{eqsss8}) which is a
computation of $\sss_5$.

It remains to prove that if $C_N$ is a computation of $\sss_7$ then $z_k$
has type $(1,0)$. We know that if $C_i$ is the last block in (\ref{eqsss8})
which is a computation of $\sss_5$ then $C_i$ replaces $x$ by $z_k\iv xz_k$.
As before, $\chi(U_i)$ is either $P_2$ or $P_3$.

Suppose first that $\chi(U_i)=P_2$. Then $C_{i-1}$ is a computation of $\sss_4$.
This implies (by Lemma \ref{sss4}) that $n_i=n_{i-1}>0$.

Suppose that $\chi(U_{i+1})=P_2$. Then $z_k$ is of type (1,1)
and so $n_{i+1}=n_i>0$. Furthermore, $C_{i+1}$ is a computation of $\sss_4$,
so $\chi(U_{i+2})=P_1$ and $n_{i+2}=n_{i+1}>0$. Since $n_{i+2}>0$,
$C_{i+3}$ cannot be a computation of $\sss_7$. Therefore it must be a
computation of $\sss_6$. Then $C_{i+4}$ must be a computation of $\bar\sss_4$
and $C_{i+5}$ must be a computation of $\sss_5$ (both $i+4$ and $i+5$
are of course smaller than $N$, since $C_N$ is a computation of $\sss_7$).
But this contradicts the fact that $C_i$ is the last block in (\ref{eqsss8})
which is a computation of $\sss_5$.
Therefore $\chi(U_{i+1})=P_3$,
so by Lemma \ref{sss5}, $z_k$ has type $(1,0)$ as desired.

Now suppose that $\chi(U_i)=P_3$. Then again $\chi(U_{i+1})=P_2$ or
$\chi(U_{i+1})=P_3$. We shall show that both possibilities lead to
contradictions.

Suppose first that $\chi(U_{i+1})=P_2$. Then $C_{i-1}$ is a computation
of $\sss_4$, $\chi(U_{i+2})=P_1$ and $n_{i+2}>0$ (Lemma \ref{sss4}).
Therefore $C_{i+3}$ is not a computation of $\sss_7$. Hence $C_{i+3}$ is
a computation of $\sss_6$, $C_{i+4}$ is a computation of $\bar\sss_4$
and $C_{i+5}$ is a computation of $\sss_5$, a contradiction.

Finally let $\chi(U_{i+1})=P_3$. Then by Lemma \ref{sss5} the type
of $z_k$ is $(0, 0)$. Since $C_{i-1}$ is a computation of $\sss_4$,
$n_i>0$ by Lemma \ref{sss4}. Since the type of $z_k$ is (0,0) then by
Lemma \ref{sss5}, $n_{i+1}=n_i>0$. Then $C_{i+1}$, $C_{i+2}$
are computations of, respectively, $\bar\sss_4$ and $\sss_6$.
By Lemma \ref{sss4}, $n_{i+2}=n_{i+3}=n_{i+1}>0$ and
$\chi(U_{i+3})=P_1$. Since $n_{i+3}>0$, $C_{i+4}$ is not a computation
of $\sss_7$, so it is a computation of $\sss_4$. Then the next block,
$C_{i+5}$, is a computation of $\sss_5$, a contradiction.

This completes the proof of Fact 4.
\vskip 0.2 in

Now let us complete the proof of the lemma.

Let $C$ be a computation starting with $W$ which has the form (\ref{eqsss8}).
Let again $U_i$ be the first word in $C_i$. By Facts 1 -- 4 that we have
proved,
$$U_i=
Ev_ixw_iFE'p\delta^{n_i}qrst\bar p\delta^{m_i}\bar q\bar r\bar s\bar tF'$$
where $v_i, w_i\in (Y\cup Y\iv)^*$, $n_i, m_i$ are integers,
the word $\chi(U_i)=pqrst\bar p\bar q\bar r\bar s\bar t$
is equal to one of the
words $P_1, P_2, P_3, P_4, P_5$, $w_i=z_kz_{k-1}...z_1$ where each $z_j$
is
of type $(\ell_j,\ell'_j)$ for some $\ell_j, \ell'_j$, $v_i=uw_i\iv$,
$n_i\ge 0$, $m_i\ge 0$.

Then by Lemmas \ref{sss4} and \ref{sss5} all blocks $C_i$ are semiproper.
Moreover, all $C_i$, $i\le N-1$, which are
computations of the machines $\sss_4$, $\bar\sss_4$,
$\sss_6$ and $\sss_7$
are proper: no negative letters can be inserted during
these computations because in the initial and final words in these computations,
the parts of the words  which can be modified by the rules of these machines
are positive (since $n_i\ge 0, m_i\ge 0$).
All negative letters inserted during computations
of $\sss_5$ between $E$ and $F$ belong to the subwords $z_j$. But by Lemma \ref{type1}
the word $z_kz_{k-1}...z_1$ is reduced, so no negative letter inserted between
$E$ and $F$ during one of the computations of $\sss_5$ can be cancelled during
another computation of $\sss_5$. We also already know that the blocks $C_i$
corresponding to $\sss_5$ did not insert negative letters ($\delta$'s)
between $E'$ and $F'$
because these blocks are semiproper (Lemma \ref{sss5}) and $n_i$ and $m_i$ are
always non-negative.

This proves the first part of (i): every computation of
$\sss_8$ starting with $W$ is semiproper. Let $N_1$ be the number
of computations of $\sss_5$ among blocks in (\ref{eqsss8}).
Then our analysis shows that $N\le 4N_1$ and $N_1$ is equal to the
number of subwords $z_i$. Thus $N_1\le |\tau C|$. If $C_i$ is a computation
of $\sss_5$ then the length of $C_i$ is bounded by the length of the
corresponding word $z_i$. Therefore it is bounded by the length
of $|\tau C|$. The area of this block is bounded by $|\tau C|^2$.
Thus the sum of lengths of these blocks is bounded by $|\tau C|^2$ and
the total area of these blocks is bounded by $|\tau C|^3$. By Lemma \ref{sss4}
the length
of each of the other blocks is bounded by $O(|\tau C|+|W|)$ (notice that the
lengths of $U_i$ do not decrease when $i$ goes from 1 to $N$). And again
the area of each of the other blocks does not exceed $O((|\tau C|+|W|)^2)$.
This implies that  the length of $C$ does not exceed $O((|\tau C|+|W|)^2)$
and the area does not exceed $O((|\tau C|+|W|)^3)$. This proves the second
part of (i).


Let $C$ be a computation in $C\sss_8(W,xF, s_1t_1)$. Let $W'=\tau C$.
By Fact 4 and Lemma \ref{type1} there are no blocks in the representation (\ref{eqsss8}) of $C$
which are computations of $\sss_5$ (otherwise the subword in $W'$ between
$x$ and $F$ would not be empty).

This implies that either $C$ has just one block which is a computation
of $\sss_4$ or $\sss_6$, or it has two blocks $C_1C_2$ where
$C_1$ is a computation of $\sss_6$ and $C_2$ is a computation of $\bar\sss_4$.

The latter case is impossible because then $W'$ would not have the subword
$s_1t_1$ ($\sss_6$ changes $s_1$ to $s_0'$ and $\bar\sss_4$ does not touch
state letters without a $``\,\bar{   }\,")$.

Any nontrivial reduced computation of $\sss_6$ consists of two words,
and there are no nontrivial computations in $C\sss_6(W,s_1t_1)$. Therefore
$C$ cannot be a computation of $\sss_6$.

Finally by Lemma \ref{sss4}, every computation of $C\sss_4(W,s_1t_1)$ is
trivial. Thus $C$ is not a computation of $\sss_4$. All cases have been
considered, so $C\sss_8(W,xF, s_1t_1)$ consists of one (trivial) computation.
This proves (ii).


Suppose that $u$ is a positive word and $n=|u|$. Then there exists a
computation of $\sss_8$ which starts with $W$ and ends with $W_4$. This
computation has the following representation as a sequence of blocks:

\begin{equation}\label{upos}
C_{4,1}C_{5,1}C_{\bar 4, 1}C_{6,1}C_{4,2}... C_{4,n}C_{5,n}C_{\bar 4, n}C_{6,n}C_{7,1}
\end{equation}
where $C_{i,j}$ is a computation of $\sss_i$ ($i=4,5,6,7$, $j=1,...,n$)
and $C_{\bar 4,j}$ is a computation of $\bar\sss_4$. Every computation
of $\sss_5$ in this sequence consists of two words, so all $z_i$ are
one letter words.

Now suppose that there exists a computation $C$ in $C\sss_8(W,p_4)$.
Let $W' = \tau C\sss_8(W,p_4)$. By Facts 3 and 4,  $W'=W_4$.
By Fact 4, $u=z_kz_{k-1}...z_1$ where the $z_i$'s satisfy the properties
listed in Fact 4. By Lemma \ref{type3}, since $u$ is a positive word,
every $z_i$ is a one letter word. Therefore by Fact 4, every
block in (\ref{eqsss8}) which is a computation of $\sss_5$ consists of two
words.
Our description of computations of $\sss_8$
starting with $W$ (see proofs of Facts 1--4)
show that the computation must have the form
$$C_{4,1}C_{5,1}C_{\bar 4, 1}C_{6,1}C_{4,2}... C_{4,k}C_{5,k}C_{\bar 4, k}C_{6,k}C_{7,1}$$
where $k=|u|$. Each block $C_{5,j}$ consists of an application of
a rule of the form  $R(a)$ which subtracts $1$ from $n$, the number of $\delta$'s between
$p$ and $t$ (here we omit indices in $p$ and $t$).
Other blocks do not touch this number.
Since there is no $\delta$ between $p$ and $t$ in the last word
in $C_{6,k}$ (otherwise we would not be able to apply the rule of $\sss_7$
to this word), $k$ must be equal to $n$.

Thus we proved (iii): if $u$ is positive then $C\sss_8(W,p_4)$ is
not empty if and only if $n=|u|$; in this case $C\sss_8(W,p_4)$ consists
of one computation (\ref{upos}). The length of this computation can be easily computed
with the help of Lemma \ref{sss4} (iii). It is equal to
$O(n+(n-1)+...+1)=O(n^2)$.

Statement (iv) follows from Fact 4 and Lemma \ref{type2}. Indeed, by
Fact 4, if $C\sss_8(W,p_4)$ is not empty then $u=z_kz_{k-1}...z_1$
where the $z_j$'s satisfy the properties listed in  Fact 4. By Lemma \ref{type2}
such a representation of $u$ is impossible if $u=u' a\iv$
where $a\in Y$, and $u'$ is a positive word
which does not end with $a$.

This completes the proof of the Lemma.
$\Box$.

\bigskip

The machine $\sss_9$ is obtained from $\sss_8$ in the same way $\sss_4$
was
obtained from $\sss_3$.
Let $\sss_8'$ be a copy of $\sss_8$ which is obtained by adding $\hat{ }$ to
all state letters except $E, F, E', F'$ and those state letters
which have index 4. Let $\sss_9$ be the
concatenation of $\sss_8$ and $\sss_8'$. Namely, the tape alphabet vector
$Y(9)$ of $\sss_9$ is the same as the tape alphabet of $\sss_8$, the state
alphabet vector $Q(9)$ is the union of state alphabet vectors of $\sss_8$
and $\sss_8'$, and the program of $\sss_9$ is the union of the programs of
$\sss_8$ and $\sss_8'$.

The following Lemma is similar to Lemma \ref{sss4}.

\begin{lm}
\label{sss9} ($\sss_9$ tells positive words in the alphabet $Y$
from almost positive words and returns the state letters to their original
positions.)
Let
$$W=
EuxFE'p_1\delta^nq_1r_1s_1t_1\bar p_0\bar q_0\bar r_0\bar s_0\bar t_0F'$$
where $u$ is either a positive word
or a reduced word of the form $u'a\iv$ where $u'$ is a positive word and $a$ is a letter.
Let
$$\hat W=Eu\hat xFE'\hat p_1\delta^n\hat q_1\hat r_1\hat s_1\hat t_1\hat{\bar p}_1\hat{\bar q}_0\hat{\bar r}_0\hat{\bar s}_0\hat{\bar t}_0F'$$
$n\ge 0$.
Then
\begin{enumerate}
\item[(i)] Every computation $C$ of $\sss_9$ starting with $W$ or $\hat W$
is semiproper. The length of $C$ does not exceed
$O((|W|+|V|)^2)$ and
the area does not exceed $O((|W|+|V|)^3)$ where $V$ is the last word in the
computation.
\item[(ii)] Each of the sets $C\sss_9(W, xF, s_1t_1)$ and $C\sss_9(\hat W, \hat xF, \hat s_1\hat t_1)$
consists of one
(trivial) computation.
\item[(iii)] The statement $\sss_9(W,\hat x)$ (statement $\sss_9(\hat W, x)$)
is true if and only if $u$ is a positive word and $n=|u|$.
In this case $C\sss_9(W,\hat x)$ (resp. $C\sss_9(\hat W, x)$)
consists of one computation and $\tau C\sss_9(W,\hat x)=\hat W$
(resp. $\tau C\sss_9(\hat W,x)=W$). The length of this computation is $O(n^2)$.
\item[(iv)] If $u$ has the form $u'a\iv$ where $a\in Y$ and
$u'$ is a positive word, and $u'$ does not end with $a$, then for
every $n$ the statement
$\sss_9(W,\hat x)$ and the statement $\sss_9(\hat W, x)$ are false.
\end{enumerate}
\end{lm}

{\bf Proof.}
Every computation $C$ of $\sss_9$ starting with $W$
can be represented in the form
$C_1C_2C_3....$
where $C_i$ is a computation of one of the machines $\sss_8$ or $\sss_8'$,
neighboring blocks are computations of different machines, and each block is
non-trivial.

Since $C$ starts with the word $W$, $C_1$ must be a
computation of the machine $\sss_8$.
If $C=C_1$ then $C$ is a computation of $\sss_8$ and we can apply
Lemma \ref{sss8}.
It is clear that in this case $\tau C$ cannot contain ${\hat x}$.

If $C\ne C_1$ then $C_2$ is a computation of $\sss'_8$. Therefore
$C_2$ starts with a word $U$ containing $p_4$. This word is the end word
of $C_1$. So statement $\sss_8(W,p_4)$ holds and
by Lemma \ref{sss8} (iii), (iv)  $u$ is positive and $n=|u|$.
Now by Lemma \ref{sss8} (i) and Lemma \ref{gen1} the computation $C_1$
is proper (it is semiproper and the end word is positive) and $C_2$ is
semiproper.

Suppose that $C_3$ exists. Then it must be a computation of $\sss_8$.
Therefore
it must start with a word containing $p_4$. By Lemma \ref{sss8}(ii)
then $C_2$ is trivial, a contradiction.

Thus $C$ consists of one or two blocks: $C=C_1$ or $C=C_1C_2$. In both cases
$C$ is semiproper. The computation $C$ can contain two blocks if and only if
$u$ is positive and $n=|u|$.

Now all statements of the lemma follow immediately from this description and
Lemmas \ref{sss8}, \ref{gen} and \ref{gen1}.
$\Box$

\bigskip
The admissible words of the $S$-machine $\sss_\alpha$ have the following form:
$$E\alpha^kx\alpha^nF$$
where $x$ is one of the letters $x, x_1, x_2$. The rules  of $\sss_\alpha$
are the following:
\begin{enumerate}
\item[($\alpha1$)] $x\to \alpha\iv x\alpha$,
\item[($\alpha_2$)] $Ex\to Ex_1$,
\item[($\alpha_3$)] $x_1\to \alpha x_1\alpha\iv$,
\item[($\alpha_4$)] $x_1F\to x_2F$.
\end{enumerate}

\begin{lm} \label{sssalpha} ($\sss_\alpha$ moves $x$ from $F$ to $E$ and back).
Let $W=E\alpha^nxF$, $W_2=E\alpha^n x_2F$, $n\ge 0$. Then:
\begin{enumerate}
\item Every reduced computation $C$ starting with $W$ or $W_2$ is semiproper.
The length of $C$ does not exceed $O(|W|+|V|)$ and the area does not exceed
$O((|W|+|V|)^2)$ where $V=\tau C$.
\item Each of the sets $C\sss_\alpha(W,xF)$ and $C\sss_\alpha(W_2,x_2)$
consists of one (trivial) computation.
\item $C\sss_\alpha(W,x_2)$ consists of exactly one computation
with the history
word $$(\alpha_1)^n(\alpha_2)(\alpha_3)^n(\alpha_4)$$ and
$\tau C\sss_\alpha(W,x_2)=W_2$, and $Ex_1$ is a subword of a word in this
computation.
\end{enumerate}
\end{lm}

{\bf Proof.} Indeed it is easy to see that the history $h$ of
every computation $C$ starting with $W$
is a prefix of the following word
$$(\alpha1)^k(\alpha2)(\alpha3)^\ell(\alpha4)$$
for some numbers $k$ and $\ell$.
This and Lemma \ref{gen} immediately imply statements 1 and 2 of the lemma.

If the history $h$ contains $(\alpha2)$ then $k=n$ because every application of $(\alpha1)$
moves $x$ one letter closer to $E$. Similarly if $h$ contains
$(\alpha4)$ then it must contain $(\alpha2)$ and $\ell$ must be equal to $n$.
This implies statement 3.
$\Box$

\bigskip

The machine $\sss_\omega$ is similar to $\sss_\alpha$.

The admissible words of the $S$-machine $\sss_\omega$ have the following form:
$$E'\omega^kx'\omega^nF$$
where $x'$ is one of the letters $x', x'_1, x'_2$. The rules  of $\sss_\omega$
are the following:
\begin{enumerate}
\item[($\omega1$)] $x'\to \omega x\omega\iv$,
\item[($\omega_2$)] $x'F'\to x'_1F'$,
\item[($\omega_3$)] $x'_1\to \omega\iv x'_1\omega$,
\item[($\omega_4$)] $E'x'_1\to E'x'_2$.
\end{enumerate}

The following lemma is similar to Lemma \ref{sssalpha}, so its proof is omitted.

\begin{lm} \label{sssomega} ($\sss_\omega$ moves $x'$ from $E'$ to $F'$
and back).
Let $W=E'x\omega^nF'$, $W_2=E'x_2\omega^nF'$, $n\ge 0$. Then:
\begin{enumerate}
\item Every reduced computation $C$ starting with $W$ or $W_2$ is semiproper.
The length of $C$ does not exceed $O(|W|+|V|)$ and the area does not exceed
$O((|W|+|V|)^2$ where $V=\tau C$.
\item Each of the sets $C\sss_\omega(W,E'x')$ and $C\sss_\alpha(W_2,x'_2)$
consists of one (trivial) computation.
\item $C\sss_\omega(W,x'_2)$ consists of exactly one
computation with the history
word $$(\omega_1)^n(\omega_2)(\omega_3)^n(\omega_4)$$ and
$\tau C\sss_\omega(W,x'_2)=W_2$, and $x'_1F'$ is a subword of a word in this
computation.
\end{enumerate}
\end{lm}

\bigskip
Now we are ready to take any Turing machine
$$M=\langle X,Y, Q, \Theta, \vec s_1, \vec s_0\rangle$$
satisfying the conditions of Lemma \ref{mach23} and to construct
an $S$-machine $\sss(M)$ simulating $M$.

For every $q\in Q$
we denote the word $q\omega$ by $F_q$. Let us replace $q\omega$
in every command of $M$ by $F_q$ . After that we won't have $q$'s in any
command of $M$ because of part 5 of Lemma \ref{mach23}. In order to
have notation similar to that of the machine $\sss_9$, we also denote the
left marker on tape $\# i$ by $E_i$. This gives us a Turing machine $M'$ such that
the configurations of each tape of $M'$ have the form $E_iuF_q$ where
$u$ is a word in the alphabet of tape $i$, and every command or its
inverse has one of the forms:
\begin{equation}\label{eq4567}
\{F_{q_1}\to F_{q_1'},..., aF_{q_i}\to F_{q'_i},..., F_{q_k}\to F_{q_k'}\}
\end{equation}
where $a\in Y$ or
\begin{equation}\label{eq4568}
\{F_{q_1}\to F_{q_1'},..., E_iF_{q_i}\to E_iF_{q'_i},...,F_{q_k}\to F_{q_k'}\}
\end{equation}
This machine recognizes the same language and has the same complexity
functions as $M$, so we can assume that $M$ itself has this form.

An admissible word of the $S$-machine $\sss(M)$ is a product of
three parts. The first part has the form

$$E(0)\alpha^{n_1}x(0)\alpha^{n_2}F(0)$$
The second part is a product of $k$ words of the form
$$
\begin{array}{l}
E(i)u_ix(i)v_iF(i)E'(i)p(i)\Delta_{i,1}q(i)\Delta_{i,2}r(i)\Delta_{i,3}s(i)\Delta_{i,4}t(i)\Delta_{i,5}\\
\bar p(i)\Delta_{i,6}\bar q(i)\Delta_{i,7}\bar r(i)\Delta_{i,8}\bar s(i)\Delta_{i,9}\bar t(i)\Delta_{i,10}
F'(i)
\end{array}
$$
$i=1,...,k$.
The third part has the form
$$
E'(k+1)\omega^{n_1'}x'(k+1)\omega^{n_2'}F'(k+1)
$$

Here $u_i, v_i$ are group words in the alphabet $Y_i$ of tape $i$,
and
$\Delta_{i,j}$ is a power of $\delta$. The letters $$E(i), x(i), F(i), E'(i),
p(i), q(i), r(i), s(i), t(i), \bar p(i), \bar q(i), \bar r(i), \bar s(i), \bar t(i),
F'(i)$$ belong to disjoint sets of state letters
which we shall denote, respectively, by $${\bf E}(i), {\bf X}(i), {\bf F}(i),
{\bf E'}(i), {\bf P}(i), {\bf Q}(i), {\bf R}(i), {\bf S}(i),
{\bf \bar P}(i), {\bf \bar Q}(i), {\bf \bar R}(i), {\bf \bar S}(i),
{\bf \bar T}(i), {\bf F'}(i),$$
$i=0,...,k+1$. The description of these sets is below.

Since $M$ is symmetric, for every command from $\Theta$ its
inverse is in $\Theta$.
Let us call commands of the form (\ref{eq4567}) {\em positive},
the inverses of these commands will be called {\em negative}. We also choose
one from each pair of mutually inverse commands of the form (\ref{eq4568})
and call it {\em positive}; the other command in this pair
will be called {\em negative}. As we have seen before, every command in
$\Theta$ is either positive or negative.

Now let us describe the set of state letters of $\sss(M)$. First
for every $i=0,...,k+1$ we include the letter $E(i)$ into ${\bf E}(i)$ and
the letter $E'(i)$ into ${\bf E'}(i)$. Then for every state
letter $F$
on tape $i$ we include $F$ into ${\bf F}(i)$ and $F'$ into ${\bf F'}(i)$.
We also include the
letters $$x(i), p(i), q(i), r(i), s(i), t(i), \bar p(i), \bar q(i), \bar r(i), \bar s(i), \bar t(i)$$
into the corresponding sets ${\bf X}(i)$,
${\bf P}(i)$, ${\bf Q}(i)$, ${\bf R}(i)$, ${\bf S}(i)$, ${\bf T}(i)$,
${\bf \bar P}(i)$, ${\bf \bar Q}(i)$, ${\bf \bar R}(i)$, ${\bf \bar S}(i)$, ${\bf \bar T}(i)$,
$(i=1,...,k)$. The letters that we just described will be called {\em standard}.

Now let $\tau$ be a positive command in $\Theta$.
Then for every $\gamma\in \{4, 9, \alpha, \omega\}$ and for each component
${\bf U}(i)$ of the vector of sets of state letters, we
include the letter $U(i,\tau,\gamma)$ into ${\bf U}(i)$.

Suppose that $\tau$ has
the form (\ref{eq4567}):
$$\tau=\{F_{q_1}\to F_{q_1'},..., aF_{q_i}\to F_{q'_i},...,F_{q_k}\to F_{q_k'}\}.$$
for some $i$ from $1$ to $k$.
For each $S$-machine $\sss_\gamma$ where $\gamma\in \{4, 9, \alpha, \omega\}$
we consider a copy $\sss_\gamma(\tau)$ where every state letter $z$
is replaced by $z(j, \tau,\gamma)$ where $j=i$ if $\gamma=4,9$, $j=0$ if
$\gamma=\alpha$ and $j=k+1$ if $\gamma=\omega$.
We include these state letters into
the corresponding sets
${\bf P}(i)$, ${\bf Q}(i)$, ${\bf R}(i)$, ${\bf S}(i)$, ${\bf T}(i)$,
${\bf \bar P}(i)$, ${\bf \bar Q}(i)$, ${\bf \bar R}(i)$, ${\bf \bar S}(i)$, ${\bf \bar T}(i)$
if $\gamma\in \{4,9\}$, in ${\bf E}(0)$, ${\bf X}(0)$, ${\bf F}(0)$ if
$\gamma=\alpha$ and into ${\bf E}({k+1})$, ${\bf X}({k+1})$, ${\bf F}({k+1})$
if $\gamma=\omega$.

The set of state letters that we just described is the set of all state letters
of $\sss(M)$.

Now let us describe the set of rules of $\sss(M)$. It will consist
of the rules of $\sss_4(\tau)$,
$\sss_9(\tau)$, $\sss_\alpha(\tau)$, $\sss_\omega(\tau)$
for all $\tau$ of the form (\ref{eq4567})
plus the following connecting rules.
First let $\tau$ be a command of the form
(\ref{eq4567}):
$$\tau=
\{F_{q_1}\to F_{q_1'},..., aF_{q_i}\to F_{q'_i},...,F_{q_k}\to F_{q_k'}\}.$$
Notice that $\tau$ determines the letter $a$ and the index $i$.

The machine $\sss(M)$ simulates this command as follows. First, using $\sss_4(\tau)$, it
checks whether the word between $E'(i)$ and $F'_{q_i}(i)$ is not empty.
If it is empty, the execution cannot proceed to the next step.
Otherwise the machine changes $q_j$ to $q_j'$ in the indices of the $F$'s,
inserts $a\iv$ next to the left of $x(i)$, removes one $\delta$
in the word between $E'(i)$ and $F'_{q_i}(i)$, removes one $\alpha$ and removes one
$\omega$. Then it uses $\sss_\alpha(\tau)$ and $\sss_\omega(\tau)$ to move
$x(0)$ from $F(0)$ to $E(0)$ and back and to move $x(k+1)$ from $E'(k+1)$
to $F'(k+1)$ and back (the purpose of these seemingly needless moves will
be clear later).
Finally it checks (using $\sss_9(\tau)$)
if after we inserted $a\iv$, the word between $E(i)$ and $F_{q'_i}(i)$ is positive.
If it is positive, the machine gets ready to execute the next transition
of $\Theta$ (by returning all state letters to their standard forms).
Otherwise, the machine cannot proceed to the execution of the
next transition from
$\Theta$.

Here is the description of the connecting rules in the case when
$\tau$ has the form (\ref{eq4567}).

\begin{enumerate}
\item[$R_4(\tau)$.] This rule is applicable to an admissible word $W$
if all state letters of $W$
are standard, and $W$ contains the subwords
$$x(j)F_{q_j}(j)E'(j)p(j)$$
and $$q(j)r(j)s(j)t(j)F'_{q_j}(j)$$ for every $1\le j\le k$, and also the subwords
$x(0)F(0)$ and $E'(k+1)x'(k+1)$.
It changes each standard state letter $z(j)$ to $z(j,\tau,4)$.
The meaning of this rule
is simple: it turns on the machine $\sss_4(\tau)$.

\item[$R_{4,\alpha}(\tau)$.] This rule is applicable when all state letters
in an admissible word $W$ have $(\tau,4)$ in their labels, and in addition
$W$ contains the subword $s'_1(i,\tau,4)t'_1(i,\tau,4)$ (that is when $\sss_4(\tau)$
finishes its work). It
\begin{itemize}
\item changes $q_j$ to $q_j'$ in the indices of $F$ ($j=1,...,k$),

\item replaces $$x(0,\tau,4)F(0,\tau,4)$$
by $$\alpha\iv x(0,\tau,\alpha)F(0,\tau,\alpha),$$

\item replaces $$E'(k+1,\tau,4)x(k+1,\tau,4)$$
by $$E'(k+1,\tau,\alpha)x(k+1,\tau,\alpha)\omega\iv,$$

\item replaces $$x(i,\tau,4)F_{q_i}(i,\tau,4)E'(i,\tau,4)p'_1(i,\tau,4)$$
by $$a\iv x(i,\tau,\alpha)F_{q_i'}(i,\tau,\alpha)E'(i,\tau,\alpha)p_1(i, \tau,\alpha)\delta\iv,$$

\item replaces $$
q_1'(i,\tau,4)r_1'(i,\tau,4) s_1'(i,\tau,4) t_1'(i,\tau,4)
$$
by
$$
q_1(i,\tau,\alpha)r_1(i,\tau,\alpha) s_1(i,\tau,\alpha)t_1(i,\tau,\alpha)
$$

\item replaces $(\tau, 4)$ by $(\tau, \alpha)$ in all other state letters.
\end{itemize}
This rule simulates execution of the transition $\tau$ and turns on
the machine $\sss_\alpha(\tau)$.

\item[$R_{\alpha,\omega}(\tau)$.] This rule is applicable when all state letters
in an admissible word
have $(\tau,\alpha)$ in their labels and $x_2(0,\tau,\alpha)$ occurs in
this word. It changes $(\tau,\alpha)$ to $(\tau, \omega)$ in the labels of all
state letters, and replaces $x_2(0,\tau,\alpha)$ by $x(0,\tau,\omega)$.
This rule turns on the machine $\sss_\omega(\tau)$.

\item[$R_{\omega,9}(\tau)$] This rule is applicable when all state letters
in an admissible word
have $(\tau,\omega)$ in their labels, and
the letter $x_2(k+1,\tau,\omega)$ occurs
in this word. It changes $(\tau,\omega)$ to $(\tau,9)$ and replaces
$x_2(k+1,\tau,\omega)$ by $x(\tau,9)$. This rule turns on
the $S$-machine $\sss_9$.

\item[$R_{9}(\tau)$.] This rule applies to an admissible word $W$
when all state letters in $W$ have $(\tau,9)$ in their labels, and
$W$ contains
$\hat x(i,\tau,9)$,
that is
when $\sss_9(\tau)$ ends its work.
The rule removes $``\hat{\hbox{   }}"$ from all letters, removes $(\tau,9)$ and
indices from all
state letters, i.e. this rule returns the state letters to
their standard form.
The meaning of this rule is that it turns off $\sss_9(\tau)$ and gets our
machine ready to simulating the next transition from $\theta$.
\end{enumerate}

We shall consider the $S$-machines
$R_4(\tau)$, $R_{4,\alpha}(\tau)$,
$R_{\alpha,\omega}(\tau)$, $R_{\omega,9}(\tau)$, $R_9(\tau)$
whose hardware is the same as the hardware of $\sss(M)$ and whose only rules
are, respectively,
$R_4(\tau)$, $R_{4,\alpha}(\tau)$,
$R_{\alpha,\omega}(\tau)$, $R_{\omega,9}(\tau)$, $R_9(\tau)$.
We shall call these $S$-machines {\em transition machines}.

\bigskip

Now let $\tau$ have the form (\ref{eq4568}),
$$\tau=
\{F_{q_1}\to F_{q_1'},..., E_iF_{q_i}\to E_iF_{q'_i},...,F_{q_k}\to F_{q_k'}\}.$$

In this case the simulation is much easier. It consists of just one
$S$-rule:
$$
\begin{array}{ll}P(\tau)&=[F_{q_1}(1)\to F_{q_1'}(1),...,\\
& E(i)x(i)F_{q_i}(i)E'(i)p(i)q(i)r(i)s(i)t(i)\bar p(i)\bar q(i)\bar r(i)\bar s(i)\bar t(i)F'_{q_i}(i) \\
& \to E(i)x(i)F_{q'_i}(i)E'(i)p(i)q(i)r(i)s(i)t(i)\bar p(i)\bar q(i)\bar r(i)\bar s(i)\bar t(i)F'_{q'_i}(i),\\
& ...,F_{q_k}(k)\to F_{q_k'}(k)]
\end{array}
$$

\vskip 0.2 in
For every configuration $c=(E_1u_1F_{q_1},...,E_ku_kF_{q_k})$ of the
machine $M$
let $\sigma(c)$ be the following admissible word of $\sss(M)$:

$$
\begin{array}{l}
E(0)\alpha^nx(0)F(0)\\
E(1)u_1x(1)F_{q_1}(1)E'(1)p(1)\delta^{||u_1||}q(1)r(1)s(1)t(1)\bar p(1)\bar q(1)\bar r(1)\bar s(1)\bar t(1)F'_{q_1}(1) ...\\
E(k)u_kx(k)F_{q_k}(k)E'(k)p(k)\delta^{||u_k||}q(k)r(k)s(k)t(k)\bar p(k)\bar q(k)\bar r(k)\bar s(k)\bar t(k)F'_{q_k}(k)\\
E'(k+1)x(k+1)\omega^{n}F'(k+1)
\end{array}
$$
where $n=|u_1|+...+|u_k|$.
Notice that $|\sigma(c)|=4|c|+13k+6$.

Conversely, if $W$ is an admissible word for $\sss(M)$, let $\mu(W)$
be the word obtained by removing
$$x, p, q, r, s, t, \bar p, \bar q, \bar r, \bar s, \bar t,
\alpha, \omega, \delta, E(0), F(0)$$ and all $E', F'$
with indices, replacing
state letters by their standard
forms (that is removing $(\tau,j)$ from their labels), and reducing the
resulting word. It is clear that if $W$ is positive then $\mu(W)$
is a configuration of the machine $M$. It is clear also that
$\mu(\sigma(c))=c$.

We shall call an admissible word $W$ of $\sss(M)$
{\em normal} if (using the notation in the definition of admissible words
for $\sss(M)$) for every $i$ from $1$ to $k$ we have that
$||u_iv_i||=||\Delta_{i,1}...\Delta_{i,10}||$ and if
$n_1+n_2=n_1'+n_2'=\sum_{i=1}^k||u_iv_i||\ge 0$.
Recall that $||u||$ is the algebraic sum of
the degrees of all letters in $W$.

The following statement is obvious.

\begin{lm} \label{normal} Let $W$ be an admissible word for $\sss(M)$.
Then if
$W$ is positive and normal and one of the rules $R_4(\tau),
R_9(\tau)\iv, P(\tau)$,
$\tau\in \Theta$ is applicable to $W$ then $W=\sigma(c)$ for
some configuration $c$ of the machine $M$.
\end{lm}

If a transition $\tau$ has the form (\ref{eq4567}):
$\tau=\{F_{q_1}\to F_{q_1'},..., aF_{q_i}\to F_{q'_i},...,F_{q_k}\to F_{q_k'}\}$,
then the subword between an ${\bf E}(i)$-letter and an ${\bf F'}(i)$-letter
in an admissible word of $\sss(M)$ will be called the {\em $\tau$-part}
of this word. The subword between an ${\bf E'}(i)$-letter and
an ${\bf F'}(i)$-letter will be called the $(\tau,\delta)$-part, the
subword between an ${\bf E}(i)$-letter and an ${\bf F}(i)$-letter
will be called the $(\tau, M)$-part, the subword between an ${\bf E}(0)$-letter
 and an ${\bf F}(0)$-letter
will
be called the $\alpha$-part, and the subword between an ${\bf E}'(k+1)$-letter and
an ${\bf F'}(k+1)$-letter
will
be called the $\omega$-part. Notice that the work of
$\sss_4(\tau)$ affects only the $(\tau,\delta)$ parts of admissible words
and the work of $\sss_9(\tau)$ affects only
the $\tau$-parts, the work of $\sss_\alpha(\tau)$ affects only the $\alpha$-parts
and the work of $\sss_\omega$ affects only the $\omega$-parts. Therefore we can
apply Lemmas \ref{sss4} and \ref{sss9} when we talk about
$\sss_4(\tau)$ and $\sss_9(\tau)$ even though the hardware of, say,
$\sss_4(\tau)$ is ``bigger" than the hardware of $\sss_4$.

Let $c_0$ be the accept configuration of the machine $M$
(all tapes are empty,
the indices of the $F$'s form the accept vector $\vec s_0$).
Then we can define
the generalized time  function $d(n)$ and the area function $a(n)$ of
the $S$-machine $\sss(M)$. For every $n>1$ we denote by $d(n)$ (resp. $a(n)$)
the smallest
number such that every admissible word $W$ of length $\le n$ for which the set of
computations $C\sss(M)(W,\sigma(c_0))$
is not empty, this set contains a computation of length (resp. area) $\le d(n)$
(resp. $\le a(n)$). The function $d(n)$ is called the {\em generalized time  function}
and the function $a(n)$ is called the {\em area function} of the machine
$\sss(M)$.

The union of the $S$-machines $$R_{4,\alpha}(\tau), \sss_{\alpha}(\tau),
R_{\alpha,\omega}(\tau),
\sss_{\omega}(\tau), R_{\omega,9}(\tau)$$
will be denoted by $R_{4,9}(\tau)$. The hardware of $R_{4,9}(\tau)$ is the same as
for $\sss(M)$. The following lemma is a corollary of
Lemmas \ref{sssalpha} and \ref{sssomega}. The machine $R_{4,9}$ is defined
as a composition of
submachines in the same manner as machines $\sss_4$ and $\sss_9$, so this lemma
can be proved in a similar
way as Lemmas \ref{sss4}, \ref{sss9}, hence we omit the proof.

\begin{lm} \label{sssr49} Let $\tau$ be of the form (\ref{eq4567}).
Let $W$ be a positive word to which the rule $R_{4,\alpha}(\tau)$
is applicable. Let $W'$ be obtained from $W$ by first applying
$R_{4,\alpha}(\tau)$ and then replacing $(\tau, \alpha)$
in all state letters by $(\tau,9)$. Then
\begin{enumerate}
\item Every computation of $R_{4,9}(\tau)$ starting with $W$ or
$W'$ is semiproper. The length of any such computation $C$ is bounded
by $O(|W|+|V|)$ and the area is bounded by $O((|W|+|V|)^2)$ where $V$ is the
last word in the computation.
\item There is only one reduced computation $C$  of $R_{4,9}(\tau)$ starting with $W$
and such that the rule $R_{\omega,9}(\tau)\iv$ is applicable to the
last word in this computation. The last word in this computation
is $W'$. This computation contains two transitions
such that the rule applied in one of these
transitions contains a word $Ex$
as one of its left sides,
where $E\in {\bf E}(0)$ and $x\in {\bf X}(0)$, and the rule applied in the other transition
has a
word $x'F'$ as one of its
left sides, where $x'\in {\bf X}(k+1)$, $F'\in {\bf F'}(k+1)$.

\item Every computation of the machine $R_{4,9}(\tau)$ starting with $W$ and
ending with a word to which the rule $R_{4,\alpha}(\tau)$ is applicable, is
trivial. Every computation of $R_{4,9}(\tau)$ starting with $W'$ and ending
with a word to which $R_{\omega,9}(\tau)\iv$ is applicable, is trivial.
\end{enumerate}
\end{lm}

The next proposition is the main statement of this section. In
this Proposition, we included all the information about the $S$-machine $\sss(M)$
that we'll use later.

\begin{prop} \label{mach1}
($\sss(M)$ simulates $M$.)
Let $M=\langle X,Y,Q, \Theta, \vec s_1, \vec s_0\rangle$ be a Turing machine satisfying
the conditions of Lemma \ref{mach23}. Then its time function, generalized time  function
and  generalized space function are equivalent to each other, and its area
function is
equivalent to the square of the time function $T(n)$.
Let $W_0=\sigma(c_0)$ where $c_0$ is the accept configuration of the Turing machine $M$.
Then:

\begin{enumerate}
\item Every computation of $\sss(M)$ starting at $W_0$ consists
of normal words and is semiproper.
If a rule $R_{4,\alpha}(\tau)$ is applied in this
computation then the word $W$ to which it applies is positive, the degree
of $\alpha$ in this word is positive, and the result of the application of this
rule is shorter than $W$.
\item A configuration $c$ of the machine $M$ is acceptable by $M$ if and only
if $\sss(M)$ can take $\sigma(c)$ to $W_0$.
\item If $c$ is an acceptable configuration of $M$ then there exists a
one-to-one correspondence $\psi$ between all accepting computations of $M$
starting with $c$ and all computations of $\sss(M)$ connecting $\sigma(c)$
and $W_0$. This correspondence satisfies the following properties:
\begin{enumerate}
\item For every accepting computation $C'$ of the machine $M$, all words in the computation $\psi(C')$
are positive. In particular, $\psi(C')$ is proper.
\item If the accepting computation $C'$ has length $T$ and space $S$ then $\psi(C')$ has
length between $\epsilon_1S^3$ and $\epsilon_2TS^2$ and area between
$\epsilon_3S^4$ and $\epsilon_4TS^3$ for some positive constants $\epsilon_1,
\epsilon_2, \epsilon_3, \epsilon_4$.
\end{enumerate}
\item The generalized time  function of $\sss(M)$ is equivalent to $T(n)^3$ and
the area function of $\sss(M)$ is equivalent to $T(n)^4$.
If there exists a computation of $\sss(M)$ connecting $W$ and $W_0$
then the minimal area of such computations is
$\le \epsilon_5T^4(\epsilon_6||W||)+O(|W|^3)$
for some constants $\epsilon_5$ and $\epsilon_6$.
\item Let $C=(W_1,...,W_n)$, $n\ge 3$, be a reduced computation starting with
a positive word
$W_1$ such that there exists a computation of $\sss(M)$ connecting $W_0$
and $W_1$. Suppose a rule $R_{4,\alpha}(\tau)$ is applied in the transition
$W_1\to W_2$ and
a rule $R_{4,\alpha}(\tau')$ is applied in the transition $W_{n-1}\to W_n$.
Then there exist $1<i,j<n$ such that the rule applied in the transition
$W_i\to W_{i+1}$ contains a word $Ex$
as one of its left sides,
where $E\in {\bf E}(0)$ and $x\in {\bf X}(0)$
 and the rule applied in the transition
and $W_j\to W_{j+1}$ has a
word $x'F'$ as one of its
left sides, where $x'\in {\bf X}(k+1)$, $F'\in {\bf F'}(k+1)$ .
\item $R_{4,\alpha}(\tau)^{\pm 1}$
are the only rules in $\sss(M)$ which change the degrees of $\alpha$ and
$\omega$ in admissible words. The parts involving $\alpha$ and $\omega$ in
this rule have the form
$xF\to \alpha\iv x_1F$ and $E'x'\to E'x_1'\omega\iv$
where $x,x_1\in {\bf X}(0)$, $F\in {\bf F}(0)$, $E'\in {\bf E'}(k+1)$,
$x', x_1'\in {\bf X}(k+1)$.
\item The rules involving $\alpha$ (resp. $\omega$) are $R_{4,\alpha}(\tau)$
and the rules of the form $x\to \alpha^{\pm 1}x\alpha^{\mp 1}$ where
$x\in {\bf X}(0)$
(resp. $x'\to \omega^{\pm 1}x'\omega^{\mp 1}$
where $x'\in {\bf X'}(k+1)$). For every $x\in {\bf X}(0)\cup {\bf X}(k+1)$
there exist at most two mutually inverse rules involving $x$ and $\alpha$
or $x$ and $\omega$.
\end{enumerate}
\end{prop}

{\bf Proof.} Let $c$ be a configuration of the machine $M$. Suppose that
a computation $C$ of the machine $\sss(M)$ starts with $\sigma(c)$.
Notice that all rules of $\sss(M)$ take normal words to normal words. Since
$\sigma(c)$ is a normal word, every word in $C$ is normal.

The computation $C$ can be represented in the form (\ref{concat}):
\begin{equation}\label{eqccc}
C_1, C_2,..., C_N
\end{equation}
where each $C_i$ is a computation of one of the machines:
$$\sss_4(\tau), \sss_9(\tau), R_{4,9}(\tau),
R_9(\tau), P(\tau')$$
where $\tau$ is any transition of $M$ of the form (\ref{eq4567}),
$\tau'$ is any transition of $M$ of the form (\ref{eq4568}).

The submachine which executes the block $C_i$ will be denoted by $\chi(C_i)$.
Now we shall describe the sequence $\chi(C)=
(\chi(C_1),\chi(C_2),...,\chi(C_N))$.

By the definition of the transition machines
it is clear that every non-trivial reduced computation
of each of these machines contains only two words.

Also from the definition it follows that

a) the types of consecutive
blocks $C_i$ and $C_{i\pm 1}$ can be the following:
\begin{itemize}
\item $R_4(\tau)$ and $\sss_4(\tau)$,
\item $\sss_4(\tau)$ and $R_{4,9}(\tau)$,
\item $R_{4,9}(\tau)$ and $\sss_9(\tau)$,
\item $\sss_9(\tau)$ and $R_9(\tau)$,
\item $R_9(\tau)$ and $R_4(\tau_1)$,
\item $R_4(\tau)$ and $R_4(\tau_1)$,
\item $R_4(\tau)$ and $P(\tau_1)$,
\item $R_9(\tau)$ and $R_9(\tau_1)$,
\item $R_9(\tau)$ and $P(\tau_1)$.
\end{itemize}
for some (different) $\tau$ and $\tau_1$, and

b) $\chi(C)$ cannot contain subsequences $Z_1Z_2Z_3$ where $Z_1, Z_2, Z_3$
belong to the set $$\{R_4(\tau), R_9(\tau), P(\tau)\ |\ \tau\in \Theta\}.$$

Let us prove the following five statements by induction on $i$ from 1 to $N$.

\begin{enumerate}
\item[A1.] if $i<N$ and $\chi(C_i)=\sss_4(\tau)$ then the degrees
of $\delta$ in the $(\tau,\delta)$-parts of all words in $C_i$ are $>0$
and all words in $C_i$ are positive.

\item[A2.] if $i<N$, $\chi(C_i)=\sss_9(\tau)$ and $C_i$ starts with
a positive word then all words in $C_i$ are positive.

\item[A3.] if $i\le N-2$, or if $i=N-1$ and $\chi(C_i)\ne R_{4,9}(\tau)$,
then all words in $C_i$ are positive.

\item[A4.] for every $i$, if $\chi(C_i)=R_{4,9}(\tau)$ then the first word in
$C_i$ is positive; moreover the degree of $\alpha$ in the word $W$ from $C_i$,
 to which
$R_{4,\alpha}(\tau)$ is applicable, is positive, and the result of the application
of this rule is shorter than $W$.

\item[A5.] If $1<i<N$ and $\chi(C_i)\in\{\sss_4, \sss_9, R_{4,9}\}$ then
$\chi(C_i)\ne \chi(C_{i+1})$.
\end{enumerate}

All five statements A1 -- A5 hold for $i=1$ because of Lemmas \ref{sss4},
\ref{sss9}, \ref{sssr49}.

Suppose that $i<N$ and $\chi(C_i)=\sss_4$. Then
$$\chi(C_{i+1}), \chi(C_{i-1})\in \{R_4(\tau), R_{4,9}(\tau)\}.$$
If $\chi(C_{i-1})=\chi(C_{i+1})$ then Lemma \ref{sss4} (ii) applies and
by this lemma $C_i$
is trivial, a contradiction.
So $\chi(C_{i-1})\ne \chi(C_{i+1})$. Then all conditions of
Lemma \ref{sss4}
apply to the $(\tau,\delta)$-part of the first and the last
words from $C_i$. Therefore by this lemma the degree of $\delta$
in the first and the last words of $C_i$ are $>0$, $C_i$ is a proper
computation and the degree of $\delta$ in the $(\tau,\delta)$-part does
not change during $C_i$. Since the work of $\sss_4(\tau)$ affects
only the $(\tau,\delta)$-parts of admissible words, we  deduce that
all words in $C_i$ are positive. This proves statement A1.

Similarly suppose that $\chi(C_i)=\sss_9(\tau)$. Then by Lemma \ref{sss9} (ii),
$\chi(C_{i-1})\ne \chi(C_{i+1})$ and the $\tau$-parts of the words
from $C_i$ are positive. This proves A2.

Suppose that $i\le N-2$ and $\chi(C_i)=R_{4,9}(\tau)$.
Then by Lemma \ref{sssr49} either the first rule applied in $C_i$ is
$R_{4,\alpha}(\tau)$
or the last rule applied in $C_i$ is
$R_{4,\alpha}(\tau)\iv$
(depending on whether  $\chi(C_{i+1})$ is $\sss_9(\tau)$
or $\chi(C_{i+1})$ is $\sss_4(\tau)$).
In the second case all words in $C_i$ are positive by Lemma \ref{sssr49}.

{\bf Remark.} Notice that in this case the degree of $\alpha$ in
the last word in $C_i$ is positive because $R_{4,\alpha}(\tau)\iv$ inserts
a new $\alpha$ in the $\alpha$-part of the word. Also notice that
$R_{4,\alpha}(\tau)\iv$ applied to a positive word always makes
the word longer.

Suppose that the first rule applied in $C_i$ is
$R_{4,\alpha}(\tau)$.
We have that $i>1$ since
$\chi(C_1)\in \{R_4(\tau), R_9(\tau), P(\tau)\ |\ \tau\in \Theta\}$.
Clearly $\chi(C_{i-1})=\sss_4$.
If $\chi(C_{i-1})=\chi(C_{i+1})$
then by Lemma \ref{sssr49}, $C_i$ is trivial. Therefore $\chi(C_{i+1})=\sss_9$.
Since the statement A1 has been proved already,
the degree of $\delta$ in the $\tau$-part of the first word
of $C_i$ is strictly positive.
Therefore the
rule $R_{4,\alpha}(\tau)$ does not insert
a negative letter in the $(\tau,\delta)$-part.

We know that all words in $C$ are normal. Since the first word in
$C_i$ is positive and normal, and the degree
of $\delta$ is greater than 0, the degree of $\alpha$ and the degree
of $\omega$ in this word are also positive. Therefore
the rule $R_{4,\alpha}$ does not insert a negative letter into the
$\alpha$-part and $\omega$-part.
\vskip 0.1 in

{\bf Remark.} Notice that in this case we showed that
the degree of $\alpha$ in the
word to which $R_{4,\alpha}(\tau)$
applies cannot be 0. Also since the degree of $\delta$ in the $(\tau,\delta)$-part
of the first word of $C_i$ is positive, the application of the rule
$R_{4,\alpha}(\tau)$ makes the $(\tau,\delta)$-part shorter by one letter.
It also makes the $\alpha$-part and the $\omega$-part shorter by one letter.
It can make the $(\tau,\bar Y)$-part longer by at most one letter and it does
not touch other parts of the word. Therefore the resulting word is at least 2 letters
shorter than the word to which this rule was applied.
\vskip 0.1 in

By Lemma \ref{sssr49} the $(\tau,\delta)$-part
of the first word in $C_{i+1}$ is positive.
If the $(\tau,M)$-part of this
word is not positive then it has the form
$$E(i,\tau,9)ua\iv x(i,\tau,9)F(i,\tau,9)$$
where $u$
is a positive word (inherited from the first word in $C_i$) which does not
end with $a$.

But then Lemma \ref{sss9} implies that $C_{i+1}$ cannot end with a
word containing a subword $stF$ (with indices), so $C_{i+2}$ does not exist
which contradicts the inequality $i\le N-2$. Therefore the $\tau$-part
of the last word in $C_i$ is positive. Since $R_{4,9}(\tau)$ does not
affect tape letters in other parts of admissible words,
all words in $C_i$ are positive.

Suppose that $i<N$, $\chi(C_i)\ne R_{4,9}(\tau)$ for any
$\tau$, and $C_i$ starts with a positive word.

Then $\chi(C_i)\in \{R_4, R_9, \sss_4, \sss_9\}$.
If $$\chi(C_i)\in \{P(\tau), R_4(\tau), R_9(\tau)\}$$ then clearly
all words in $C_i$ are positive. If $\chi(C_i)=\sss_4$ or $\chi(C_i)=\sss_9$
then A1 and A2 imply that all words in $C_i$ are positive.
This proves A3.

Suppose that $\chi(C_i)=R_{4,9}(\tau)$. If $i\le N-2$ then all words in
$C_i$ are positive by A3 and the degree of $\alpha$ in this word is positive
by the two remarks in the proof of A3. If $i=N-1$ or $i=N$
then $\chi(C_{i-1})\ne R_{4,9}(\tau)$ and
by A3,  all words in $C_{i-1}$ are positive. In particular, the first
word of $C_i$ (which is the same as the last word in $C_{i-1}$) is positive.
If the rule $R_{4,\alpha}(\tau)$ is applied at the beginning of $C_i$
then $\chi(C_{i-1})=\sss_4(\tau)$ and the degree of $\alpha$ in the first word
of $C_i$ is positive. If $R_{4,\alpha}(\tau)\iv$ is applied at the end
of $C_i$ then all words in $C_i$ are positive and so the degree of $\alpha$
in the last word of $C_i$ is positive ($R_{4,\alpha}(\tau)\iv$ inserts a
new $\alpha$ in the alpha part of the word). This proves A4.

Property A5 follows from A3 and Lemmas \ref{sss4}, \ref{sss9}, \ref{sssr49}.
Indeed these lemmas imply that if $\chi(C_i)\in \{\sss_4, \sss_9, R_{4,9}\}$
and $\chi(C_{i-1})=\chi(C_{i+1})$ then $C_i$ is empty.

\bigskip

Now let us prove the statements of the proposition. Notice that statements
6 and 7 can be proved by a simple inspection of all rules of $\sss(M)$,
so we need to prove statements 1 -- 5 only.


The fact that every word in a computation $C$ starting at $W_0$ is normal has
been established before. If a rule $R_{4,\alpha}(\tau)^{\pm 1}$ is applied in
$C$ then it is either the first rule or the last rule
applied in a block $C_i$ with
$\chi(C_i)=R_{4,9}(\tau)$.  If this is the first rule then the word to which
this rule is applied is positive by property A4.  If this is the last rule then
this rule is $R_{4,9}(\tau)\iv$. In this case the first word in $C_i$ is
positive by property A4, and by Lemma \ref{sssr49} there exists exactly one
computation of $R_{4,9}(\tau)$ starting at the first word of $C_i$ and ending
with a word to which $R_{4,9}(\tau)$ is applicable. In this computation all
words are positive.  Let us prove that $C$ is semiproper.

By statement A3, if $\chi(C_{N-1})\ne R_{4,9}(\tau)$ then
all words in $\bigcup_{i=1}^{N-1} C_i$ are positive. In particular, the
first word in $C_N$ is positive. Then Lemmas \ref{sss4}, \ref{sss9}, \ref{sssr49}
imply that the computation $C$ is semiproper.

If $\chi(C_{N-1})=R_{4,9}(\tau)$ then by statement A4 the first word in
$C_{N-1}$ is positive. By statement A3, all words in $\bigcup_{i=1}^{N-2} C_i$
are positive. If the rule $R_{4,\alpha}(\tau)\iv$ is applied in $C_{N-1}$
then the first word in $C_N$ is positive and as before $C$ is a semiproper
computation. If the rule $R_{4,\tau}(\tau)$ is applied in $C_{N-1}$
and inserts a negative letter $a\iv$ in the $\tau$-part of
the admissible word then
$\chi(C_N)=\sss_9(\tau)$, and by Lemma \ref{sss9} the computation
$C_N$ is semiproper and does not remove the letter $a\iv$.
So the computation
$C$ is semiproper.

This proves statement 1 of the proposition.


Notice that property A5 and the list of possible pairs $(\chi(C_i), \chi(C_{i+1}))$
imply that $\chi(C)$ can be represented
in the form $D_1(\tau_1),D_2(\tau_2),...,D_\ell(\tau_\ell)$
where for each $i<\ell$, $\chi(D_i(\tau_i))$ is a sequence of one of the following
forms:
\begin{itemize}
\item[(S1)]
$R_4(\tau_i),\sss_4(\tau_i),R_{4,9}(\tau_i),\sss_9(\tau_i),R_9(\tau_i)$,
\item[(S2)]
$R_9(\tau_i),\sss_9(\tau_i),R_{4,9}(\tau_i),\sss_4(\tau_i),R_4(\tau_i)$,
\item[(S3)] $P(\tau_i),$
\end{itemize}
and $\chi(D_\ell)$ is a prefix of one of these sequences.

Since all words in $C$ are normal,
by Lemma \ref{normal}, there exists
a sequence of configurations $C'=(c_1,...,c_\ell)$ such that the first word
in $D_i(\tau_i)$ is $\sigma(c_i)$, $i=1,...,\ell$.
Notice that if we apply one of the rules of the
submachines $R_4(\tau), \sss_4(\tau), \sss_9(\tau)$
and $R_9(\tau)$ to a positive admissible word $W$ then the corresponding
configuration $\mu(W)$ of the machine $M$ will not change.
If $W$ is positive and the result $W'$ of an application of $R_{4,\alpha}(\tau)$
is also positive then $\mu(W')$ is obtained from $\mu(W)$ by applying the
transition $\tau$. Therefore in  $C'$ every $c_i$ is obtained from
$c_{i-1}$ by applying a transition from $\Theta$.
Thus $C'$ is a computation
of the Turing machine $M$.

If $c$ is an acceptable
configuration of the machine $M$ then there exists a computation $C'$
of $M$ from $c$ to $c_0$. Let $\tau_1\tau_2...\tau_\ell$ be the history
of this computation. Then the  corresponding computation
$D_1(\tau_1),D_2(\tau_2),...,
D_\ell(\tau_\ell)$ takes $\sigma(c)$ to $\sigma(c_0)$. Conversely, if
there exists a computation $C$ which takes $\sigma(c)$ to $\sigma(c_0)$
then as we proved before, $C=D_1(\tau_1)...D_k(\tau_\ell)$
where each $D_i$, $i<\ell$ has one of the forms (S1), (S2), (S3) and
$D_\ell$ is a prefix of such a sequence. Since the last word
in $D_\ell$ is $\sigma(c_0)$, $D_\ell$ cannot be a proper prefix.
As we proved before, the word $\tau_1\tau_2...\tau_\ell$ is the history
of a computation which takes $c$ to $c_0$. This proves
statement 2 of the proposition.

Let $\tau$ be a positive transition from $\Theta$ of the form (\ref{eq4567})
and suppose that a configuration $c'$ is obtained from $c$ by applying $\tau$.
Then Lemmas \ref{sss4},
\ref{sss9}, \ref{sssr49}
imply that there exists only one computation of $\sss(M)$
of the form (S1) which takes $\sigma(c)$ to $\sigma(c')$ and only one
computation of the form $(S2)$ which takes $\sigma(c')$ to $\sigma(c)$.
If $\tau$ has the form (\ref{eq4568}) and $c'$ is obtained from $c$
by applying $\tau$ then $\sigma(c')$ is obtained from $\sigma(c)$
by applying $P(\tau)$.
Therefore the correspondence between computations of $M$ which take
$c$ to $c_0$, and computations of $\sss(M)$ which take $\sigma(c)$
to $\sigma(c_0)$, is one-to-one. Let us call this correspondence $\psi$
as in the proposition.

The first property of $\psi$ has been established before.

In order to establish the second property, let $C'=(c_1,...,c_T)$ be
an accepting computation of the machine $M$ of length $T$ and space $S$.
Let $C=\psi(C')$. As before $\psi(C)$ can be represented as a sequence
of blocks $D_1(\tau_1),D_2(\tau_2),...,D_T(\tau_T)$.
By Lemmas \ref{sss4},
\ref{sss9} the length of each block $D_i(\tau_i)$, $i=1,...,T$,
in this representation is $O(n_i+n_i^2)=O(n_i^2)$ where $n_i$ is the length of
the $\tau_i$-part of the first word in $D_i(\tau_i)$, $i=1,...,T$.
Since all words in the computation $C$ are normal,
the numbers $n_i$ do not exceed $O(S)$. Therefore the length of
each block $D_i$ does
not exceed $O(S^2)$.
Thus the length of $C$ does not exceed $O(TS^2)$.
Similarly the area of $C$ does not exceed $O(TS^3)$ because the area
of each block does not exceed $S^3$.

Let $c'$ be a configuration in $C'$ with maximal possible length.
By definition of the space of a computation, this length must be equal
to $S$.

Then there exists a number $i$ such that the length of the word
written on tape $i$ in the configuration $c'$ is at least $S/k$ where
$k$ is the number of tapes of the machine $M$.

Let $\tau_1...\tau_T$ be the history of the computation $C'$.
Since the computation $C'$ is accepting, the length of the word written
on tape $i$ in the last configuration of $C'$ is 0 (see Lemma \ref{mach23},
statement 4). Since every transition of the machine $M$ can remove at
most one letter from the word written on tape $i$, there exists
at least $S/(2k)-1$ transitions $\tau_{j_1}, ..., \tau_{j_\ell}$ in the
history of computation $C'$ such that the length of the word
written on the tape $i$ in configuration $c_{j_n}$, $n=1,...,\ell$,
is at least $S/(2k)$, and the transition $\tau_{j_n}$ removes a letter
from this word. By Lemmas \ref{sss4}, \ref{sss9} and \ref{sssr49} the length
of each of the corresponding blocks $D_{j_n}(\tau_{j_n})$ exceeds
$O(S^2/(4k^2))$ and the area exceeds $O(S^3/(8k^3))$. Therefore the
total length of the computation $\psi(C')$ exceeds $O(S^3)$ and
the area exceeds $O(S^4)$. This gives us the second
property of $\psi$.

Let us prove the fourth statement of the proposition. Let $d(n)$ be the
generalized time  function of $\sss(M)$, and let
$a(n)$ be the area function of $\sss(M)$.
Let $S(n)$ be the generalized space function, and let $T(n)$ be the
time function of the machine $M$.
We know that $S(n)$ is equivalent to $T(n)$. Since every accepting
computation ends with all tapes empty, $T(n)\ge n$. It is also
clear that $S(n)\ge n$. Thus $T(n)^p$ is equivalent to
$S(n)^p$ for every natural number $p>0$.

Fix a number $n>0$. Let $c$ be an accepted configuration of the machine $M$
with $|c|\le n$, such that the smallest space of an accepting computation
for $c$ is $S(n)$. Then by the second property of the function
$\psi$, the shortest computation of $\sss(M)$ connecting $\sigma(c)$ with
$W_0$ has length $\ge \epsilon_1 S(n)^3$ and the smallest area computation
connecting $W_0$ and $\sigma(c)$ has area $\ge \epsilon_3 S^4(n)$
for some constants $\epsilon_1, \epsilon_3>0$.
Since $|\sigma(c)|=4|c|+13k+6$, we have
$$
d(4n+13k+6)\ge \epsilon_1 S^3(n), \ a(4n+13k+6)\ge \epsilon_3 S^4(n)
$$
Since the function $S(n)$ is equivalent to $T(n)$, we have
\begin{equation}
T(n)^3\preceq d(n), \ T(n)^4\preceq a(n)
\label{in1}
\end{equation}

On the other hand take any admissible word $W$, $||W||\le n$ such that
there exists
a computation of $\sss(M)$ connecting $W$ to $W_0$.
Take the shortest computation $C_1$ and the smallest area computation $C_2$
which connect $W_0$ and $W$.
We can represent  $C_i$, $i=1,2$, in the form $\psi(C_i')C_i''$.
Notice that the length of the last configuration $c_i$ in $C'_i$ is smaller
than $||W||$. Indeed, recall that $C_i''$ is a composition of computations
of $\sss_4(\tau)$, $\sss_9(\tau)$, $R_4(\tau)$, $R_{4,9}(\tau)$,
$R_9(\tau)$ for some transition $\tau$. Therefore either $\mu(W)$
coincides
with $c_i$, or $\mu(W)$ is a configuration of $M$ obtained from $c_i$
by applying a command $\tau$ from $\Theta$, or $\mu(M)$ is obtained from $c_i$
by inserting a negative letter in one of the tapes. Thus $|c_i|\le |\mu(W)|+1$.
On the other hand $|\mu(W)|\le (||W||-13k-6)/6<|W|-1$. So $|c_i|<||W||$.

The computation $C_i'$, $i=1,2$,  is an accepting computation
for the configuration
$c_i$. Since $|c_i|\le n$, there exists an accepting computation $\tilde C_i$
for $c_i$ of length $T\le T(n)$. The space of an accepting computation of length
$T$ cannot exceed $T$ (because in the last configuration all tapes are
empty and every step of a computation can remove at most one cell).
By the second part of our proposition, the length
of $\psi(\tilde C_i)$ does not exceed $\epsilon_2 T^3\le \epsilon_2 T(n)^3$
and the area of this computation does not exceed $\epsilon_4 T^4\le \epsilon_4 T(n)^4$.
for some constants
$\epsilon_2, \epsilon_4$. By Lemmas \ref{sss4},
\ref{sss9}, \ref{sssr49} the length of
$C_i''$, $i=1,2$, does not exceed $O(|W'|+|W|)^2)$
and the area does not exceed
$O((|W'|+|W|)^3)$ where $W'=\sigma(c_i)$. Therefore
the length of $C_i''$, $i=1,2$, does not exceed $O(|W|^2)$
and the area does not exceed
$O(|W|^3)$. The computation
$\psi({\tilde C_1})C_1''$ takes $W_0$ to $W$.
So it cannot be shorter
than the computation $C_1$. Therefore the length of $C_1$ does not exceed
$O(T^3(||W||))+O(|W|^2)$.
Similarly, the area of $C_2$ does not exceed the area of
$\psi({\tilde C_2})C_2''$. So the area of $C_2$ does not exceed $O(T(n)^4)+O(|W|^3)$.
This proves the second part of statement 4 of the proposition.
Since $||W||\le |W|$ and $T(n)\ge n$, we deduce
\begin{equation}
\label{in2}
d(n)\preceq T(n)^3, \ a(n)\preceq T(n)^4.
\end{equation}
The inequalities (\ref{in1}) and (\ref{in2})
imply that $d(n)$
is equivalent to $T(n)^3$ and $a(n)$ is equivalent to $T(n)^4$.
This completes the proof of statement 4 of the proposition.

Let us prove property 5 of the proposition.
Let $C=(W_1,...,W_n)$, $n\ge 3$, be a reduced computation of $\sss(M)$
starting with
a positive word
$W_1$ such that there exists a computation $\bar C$ of $\sss(M)$
connecting $W_0$
and $W_1$.
Suppose a rule $R_{4,\alpha}(\tau)$ is applied in the transition
$W_1\to W_2$ and
a rule $R_{4,\alpha}(\tau')$ is applied in the transition $W_{n-1}\to W_n$.
Let $C=\tilde C C'$ and $\bar C=C''\tilde C\iv$ for some computations
$\tilde C$, $C'$, $C''$ where $\tilde C$ is a maximal prefix
of the computation of $C$ which cancels when we reduce
the computation $\bar C C$.
Then $C''\tilde C$ and $C''C'$ are two reduced computations starting
with $W_0$.

Suppose first that $\tilde C$ is not empty.

Then as before we can decompose the computations $C''\tilde C\iv$ and $C''C'$
into a sequence of blocks (\ref{eqccc}). Since the rule $R_{4,\alpha}(\tau)$
is applied in the transition $W_1\to W_2$, the block $B$ of the computation
$C''C'$ containing $W_1$ is such that $\chi(B)=R_{4,9}(\tau)$. Similarly
$W_{n-1}, W_n$ are contained in a block $B'$ of  the
computation $C''\tilde C\iv$ such that
$\chi(B')=R_{4,9}(\tau')$. Since there is only one rule from $R_{4,9}(\tau')$
applicable to $W_{n-1}$, the block $B'$ has length 2.

If $B$ is contained in $\tilde C\iv$ then by
Lemma \ref{sssr49}, $C'$ contains two transitions $W_i\to W_{i+1}$ and
$W_j\to W_{j+1}$ such that
the corresponding rules have $Ex$ and $x'F'$
as their left sides
where $E\in {\bf E}(0)$ and $x\in {\bf X}(0)$, $x'\in {\bf X}(k+1)$,
$F'\in {\bf F'}(k+1)$
as desired.

Suppose that $\tilde C\iv$ is not empty and is contained in $B$.
Let $T_1, T_2$ be the first
two words of the block $B$.
Let $B''$ be the block of the computation $C''C'$ containing $T_1$ and $T_2$.
Since the rule applied in the transition $T_1\to T_2$ belongs to $R_{4,9}(\tau)$,
$\chi(B'')=R_{4,9}(\tau)$. Since the block $B'$
has length 2, the block $B''$ is not the last block in the computation
$C''C'$. By property A3, $T_1$ is a positive word.
Therefore $B''$ must end with an application of the rule
$R_{4,\alpha}(\tau)\iv$. By Lemma \ref{sssr49}, there is only one reduced
computation of the machine $R_{4,9}(\tau)$ which starts with $T_1$ and ends
with a word to which $R_{4,\alpha}(\tau)$ is applicable. Thus $B=B''$.
But this contradicts the assumption that $\tilde C$ is the longest part
of $C$ which cancels in the product $\bar C C$.

Finally suppose that $\tilde C$ is empty, that is, the computation $\bar CC$ is
reduced. Then the words $W_1$ and $W_2$
are contained in a block $B$ of the computation $\bar CC$
and $\chi(B)=R_{4,9}(\tau)$. As before this is not the last block in
the computation $\bar CC$ because $W_{n-1}$ must belong to another block.
Therefore we can apply Lemma \ref{sssr49} and conclude that $B$ contains two
desired transitions.

This completes the proof of the proposition. $\Box$.

\section{The Group Presentation}

In this section we convert the $S$-machine $\sss(M)$
into a group presentation.

Let $\sss=\sss(M)$ be the $S$-machine
described in Proposition \ref{mach1}. Let
$Y$ be the vector of sets of tape
letters, and let  $Q$
be the vector of sets of state letters of $\sss$.
The vector $Q$ has $15k+6$ components which we shall denote
by $Q_1,...,Q_{15k+6}$. Notice that
$Q_1=\ex$, $Q_2=\bx$, $Q_3=\fx$, $Q_{15k+4}={\bf E'}(k+1)$,
$Q_{15k+5}={\bf X}(k+1)$, $Q_{15k+6}={\bf F'}(k+1)$.

Let $W_0$ be the same
admissible word as in Proposition \ref{mach1}.
We shall not use the whole definition of $\sss$,
only Proposition \ref{mach1}.

Let $\Theta$ be the set of rules of $\sss$. Let us call one of each pair
of mutually inverse rules from $\Theta$ {\em positive} and the other
one {\em negative}. The set of all positive rules will be denoted
by $\Theta_+$ and the set of all negative rules will be denoted by
$\Theta_-$.

Let $N$ be any positive integer. Let
$$A=\bigcup_{i=1}^{15k+6} Q_i\cup \{\alpha, \omega, \delta\}\cup
\bigcup_{i=1}^{k} Y_i \cup \{\kappa_j\ |\ j=1,\ldots, 2N\}\cup
\Theta_+.$$

Our group $G_{N}(\sss)$ is generated by the set $A$ subject to the set
$\pp_N(\sss)$ of relations described below.

{\bf 1. Transition relations}. These relations correspond to
elements of $\Theta_+$.

Let $\tau\in \Theta_+$, $\tau=[U_1\to V_1,...,U_p\to V_p]$.
Then we include relations $U_1^\tau=V_1,...,U_p^\tau=V_p$
into $\pp_N(\sss)$. Here $x^y$ stands for $y\iv xy$.
If for some $j$ from 1 to $15k+6$
the letters from $Q_j$ do not appear in any of the $U_i$ then also include the relations
$q_j^\tau=q_j$ for every $q_j\in Q_j$.

{\bf 2. Auxiliary relations}.

These are all possible relations of the form $\tau x=x\tau$ where
$x\in \{\alpha, \omega, \delta\}\cup \bigcup_{i=1}^k Y_i$, $\tau\in \Theta_+$
and all relations of the form
$\tau \kappa_i=\kappa_i\tau$, $i=1,...,2N$, $\tau\in \Theta_+$.

{\bf 3. The hub relation.}

For every word $u$ let $K(u)$ denote the following word:
\begin{equation}
\label{hub}
K(u)\equiv
(u\iv \kappa_1u\kappa_2u\iv \kappa_3u \kappa_4\ldots u\iv \kappa_{2N-1}u \kappa_{2N})
(\kappa_{2N}u\iv \kappa_{2N-1}u \ldots \kappa_2 u\iv \kappa_1u)^{-1}.
\end{equation}

Then the hub relation is

$$K(W_0)=1.$$

\setcounter{pdfour}{\value{ppp}}
Figure \thepdfour\ shows the diagram corresponding to the hub relation if $N=1$.

\vskip 0.2 in
\unitlength=0.50mm
\special{em:linewidth 0.4pt}
\linethickness{0.4pt}
\begin{picture}(103.67,72.00)
\put(89.66,68.34){\vector(1,-1){9.67}}
\put(40.33,20.67){\vector(1,-1){9.67}}
\put(40.33,58.33){\vector(0,-1){37.67}}
\put(50.33,68.33){\vector(1,0){39.33}}
\put(89.67,11.00){\vector(-1,0){39.33}}
\put(37.00,40.33){\makebox(0,0)[rc]{$W_0$}}
\put(40.67,13.66){\makebox(0,0)[cc]{$\kappa_2$}}
\put(69.33,7.67){\makebox(0,0)[cc]{$W_0$}}
\put(99.67,13.00){\makebox(0,0)[cc]{$\kappa_1$}}
\put(103.67,40.67){\makebox(0,0)[lc]{$W_0$}}
\put(99.67,68.00){\makebox(0,0)[cc]{$\kappa_2$}}
\put(69.00,72.00){\makebox(0,0)[cc]{$W_0$}}
\put(40.67,68.00){\makebox(0,0)[cc]{$\kappa_1$}}
\put(99.33,20.67){\vector(0,1){38.00}}
\put(40.33,58.33){\vector(1,1){10.00}}
\put(89.67,11.00){\vector(1,1){9.67}}
\end{picture}
\begin{center}
\nopagebreak[4]
Fig. \theppp.

\end{center}
\addtocounter{ppp}{1}
\vskip 0.2 in

We define the presentation $\pp_{N}(\sss)$ by taking the set of all
the relations defined above, their cyclic shifts and inverses of these
cyclic shifts. We shall say that a relation is {\em determined} by some
characteristic (a subword, a pair of letters, etc.) if it is determined up to
taking cyclic shifts and inverses. Cyclic shifts and inverses of transition
(resp. auxiliary and hub) relations will also be called {\em transition} (resp.
{\em auxiliary} and {\em hub}) relations.

The group $G_{N}(\sss)$ with which we will be working in
this paper is given by the presentation $\pp_{N}(\sss)$.

The following lemma contains some properties of $\pp_{N}(\sss)$ which
either follow from the definition of $S$-machines or are easy to verify.

\begin{lm} \label{cancel} a) Every relation except the hub contains exactly
two $\Theta^{\pm 1}$-letters: a letter from $\Theta$ and a letter from $\Theta\iv$.

b) If transition relations $s_1, s_2\in \pp_{N}(\sss)$ have a common $\Theta$-letter
and a common $Q_j$-letter
then $s_1$ is a cyclic shift of $s_2^{\pm 1}$.

c) If two auxiliary relations $s_1, s_2\in \pp_{N}(\sss)$ have a common 2-letter
prefix then $s_1=s_2$.

d) If $s_1$ and $s_2\iv$ are two cyclic shifts of the hub or its inverse
and have a common prefix of the form $\kappa_jW_0 \kappa_s$ then $s_1=s_2$.
\end{lm}

\section{Some Basic Definitions}
\label{sim}

We mainly use notation and definitions from Ol'shanskii \cite{Ol89}. In particular,
we use the so called {\em 0-refinement} of a \vk diagram.
 We do not want to define
it here precisely recall only that a $0$-edge is an edge labelled by 1, and a zero-cell
corresponds to a relation of the form $a1^na\iv 1^m$ where $a$ is one of the generators,
$m$ and $n$ are integers. One can insert zero-edges and zero-cells in a diagram in
order to separate two paths which touch each other, or in order to make a path, which touches
itself, simple. When we compute an area or diameter of a diagram we never take
zero-cells and zero-edges into account.

If $\Delta$ is an (ordinary) \vk diagram then $\partial(\Delta)$ denotes its
boundary. If $\Delta$ is an annular diagram then $\partial_o(\Delta)$ and
$\partial_i(\Delta)$ denote the outer and the inner boundaries of $\Delta$.  The
cells of diagrams over ${\cal P}_{N}(\sss)$ will be called by the names of the
corresponding relations: transition cells, auxiliary cells, hub cells (or simply
hubs).

We always assume that boundaries of \vk diagrams and boundaries
of cells are oriented clockwise (this  is an insignificant difference
with \cite{Ol89}). The {\em contour} of a cell or a diagram is the union
of the boundary and and its inverse.

A path in a diagram $\Delta$ is called {\em simple} if it does not cross itself.
A path is called {\em reduced} if it does not contain consecutive mutually
inverse edges. Every path in a diagram can be reduced by removing subpaths
of the form $ee\iv$.
It is mentioned
in \cite{Ol89} that
by a 0-refinement one can turn every simple path into an {\em absolutely simple}
path which does not cross and does not touch itself.

The length of a path $p$, denoted by $|p|$
is the number of non-zero edges in it.

We shall use the following operations  which can be performed on
arbitrary diagrams.

\vskip 0.1 in

{\bf Taking the inverse (mirror image).}
Let $\Delta$ be a \vk diagram over a symmetric
set of defining relations $\pp$.
Consider the mirror image $\Delta\iv$
of the graph $\Delta$ with respect to some straight line.
Since $\pp$ is symmetric,
the graph $\Delta\iv$ is again a \vk diagram over $\pp$. We call $\Delta\iv$
the {\em inverse of $\Delta$}. It is easy
to see that every two inverses of the same diagram $\Delta$ can be transformed into each
other by a homotopy of a plane, so they are {\em homotopic}.

{\bf Composition.} Let $\Delta_1$ and $\Delta_2$ be \vk diagrams over
a presentation $\pp$. Let $\partial(\Delta_1)=p_1p_1'$,
$\partial(\Delta_2)=p_2\iv p_2'$. Suppose that $\Lab(p_1)=\Lab(p_2)$ in the free
group. Here $\Lab(p)$ is the label of the path $p$. Then
by a zero refinement we can make both paths $p_1$ and $p_2$
absolutely simple and the labels of $p_1$ and $p_2$ identical.
After that we can glue $\Delta_1$ and $\Delta_2$ by
identifying the corresponding edges of $p_1$ and $p_2$.
This operation will be called the {\em composition} and the resulting diagram
will be denoted by $\Delta_1\circ_{p_1=p_2}\Delta_2$.

\vskip 0.1 in

Let $\Delta$ be any van Kampen diagram. The following definition is similar to
the definition of a dual graph of a diagram \cite{LS}, \cite{Ol89}.  Fix a point
in each of the cells in $\Delta$ and a point in the inside of each of the edges
of $\Delta$.
For each cell $\pi$ and each edge $e$ on the boundary of this cell,
fix a simple polygonal line $\ell(\pi,e)$ inside $\pi_i$ which connects the
point inside $\pi_i$ and the point inside $e$. We can choose these
lines in such
a way that $\ell(\pi, e)$ and $\ell(\pi, e')$ do not have common points except
for the fixed point inside $\pi$.

\vskip 0.1 in
The next definition of a band in a diagram is crucial for our paper.
\vskip 0.1 in

Let $S$ be a set of letters and let $\Delta$ be a van Kampen diagram.  Fix
pairs of $S$-edges in some cells from $\Delta$ (we assume that each of these
cells contains at
least two $S$-edges).

Suppose that
$\Delta$ contains a sequence of cells ($\pi_1,\ldots,\pi_n$)
such that for each $i=1,\ldots,n-1$ the cells $\pi_i$ and $\pi_{i+1}$ have a
common $S$-edge and this edge belongs to the pair of
$S$-edges fixed in $\pi_i$ and $\pi_{i+1}$.
Consider the line which is the union of the lines $\ell(\pi_i, e_i)$ and
$\ell(\pi_{i+1}, e_i)$, $i=1,\ldots,n-1$. This polygonal line
is called the {\em median} of this
sequence of cells. Then our sequence of cells ($\pi_1,\ldots,\pi_n$) with
common edges $e_1, \ldots, e_{n-1}$ is called an {\em $S$-band}
if the median is a simple curve or a simple closed
curve.

We say that two bands {\em cross} if their medians cross. We say that two bands
{\em touch} each other if their medians touch each other. We say that a band is
an {\em annulus} if its median is a closed curve.

Let $\cal B$ be an $S$-band with common edges $e_1, e_2,\ldots, e_n$
which is not an annulus. Then the first cell
has an $S$-edge $e$ which forms a pair with $e_1$
and the last cell of $\cal B$ has an edge $f$ which forms a pair with
$e_n$. Then we shall say that $e$ is the {\em start edge} of $\cal B$ and
$f$ is the {\em end edge of $\cal B$}.  If $p$ is a path in $\Delta$ then
we shall say that a band
starts (ends) on the path $p$ if $e$ (resp. $f$) belongs to $p$.

If $S$ and $T$ are two disjoint sets of letters, ($\pi$, $\pi_1$, \ldots,
$\pi_n$, $\pi'$) is an $S$-band and ($\pi$, $\gamma_1$, \ldots, $\gamma_m$,
$\pi'$) is a $T$-band then we say that these two bands form an {\em
$(S,T)$-annulus} if the medians of these bands form a simple closed
curve called the {\em median of the $(S,T)$-annulus};
and in addition the start and end edges of these bands are not contained in the
polygon bounded by this median.

If $\ell$ is the median of an $S$-annulus or an $(S,T)$-annulus then the
maximal subdiagram of $\Delta$ contained in the area
bounded by $\ell$ is called
the {\em inside diagram} of the annulus.

The union of cells of $\cal B$ forms a
subdiagram (perhaps after some 0-refinement).
The reduced boundary of this diagram, which we shall call the {\em boundary
of the band},
has the form $e^{\pm 1}pf^{\pm 1}q\iv$ (recall that we trace boundaries
of diagrams clockwise).
Then we say that $p$ is the
{\em top path} of $\cal B$, denoted by $\topp({\cal B})$,
and $q$ is the {\em bottom
path} of $\cal B$, denoted by $\bott({\cal B})$.

Suppose that an $S$-band
${\cal B}$ is a union of two $S$-bands ${\cal B}_1$ and ${\cal B}_2$
and the end edge of ${\cal B}_1$ coincides with the
start edge of ${\cal B}_2$. Then the top (bottom) path of ${\cal B}$
is obtained by multiplying the top (bottom) paths of
${\cal B}_1$ and ${\cal B}_2$ and reducing consecutive mutually inverse
edges.

If ${\cal B}=(\pi_1,\ldots,\pi_n)$ is an $S$-band then
($\pi_n,\pi_{n-1},\ldots,\pi_1$)
is also an $S$-band which we shall call the {\em inverse of $\cal B$}.
We shall denote this band by ${\cal B}^{-1}$.

We shall call an $S$-band {\em maximal} if it is not contained in any other
$S$-band. If an $S$-band $\ww$ starts on the contour of a cell $\pi$,
does not contain
$\pi$  and is not contained in any other $S$-band with these properties
then we call $\ww$ a {\em  maximal $S$-band starting on the contour
of $\pi$}.

Now let us return to our presentation $\pp_{N}(\sss)$.
Consider the following partition of the generating set $A$:

\begin{equation}
\label{part}
 \Theta \cup Q_1\cup...\cup Q_{15k+6}\cup \bigcup Y_i\cup\{\kappa_1\}...
\cup\{\kappa_{2N}\}.
\end{equation}

Notice that these subsets are disjoint (the fact that $Q_1,...,Q_{15k+6}$ are
disjoint follows from the definition of $S$-machines).

Let $B$ be a block of this partition.
In order to consider $B$-bands and annuli, we need to
divide $B$-letters of some relation of $\pp_{N}(\sss)$ into pairs.

Every relation except for the hub contains a letter from
$\Theta$ and the inverse of this letter. These letters form a $\Theta$-pair.

Let us denote $\{\alpha, \omega, \delta\}\cup\bigcup Y_i$ by $\bar Y$.
Every auxiliary relation containing a $\bar Y$-letter,
contains a $\bar Y$-letter and the inverse of this letter. These
letters form a $\bar Y$-pair.

Every auxiliary relation containing $\kappa_i$ contains also $\kappa_i\iv$.
These two letters form a $\kappa_i$-pair.

If a transition relation has a letter from $Q_j$ then it has exactly two
letters from $Q_j^{\pm 1}$ (this follows from the definition of an $S$-rule).
They form a $Q_j$-pair.

The hub relation contains an occurrence of $\kappa_i$ and an occurrence of $\kappa_i\iv$.
These occurrences form a $\kappa_i$-pair in the hub.

A $\kappa_j$-band is called {\em even} (resp. {\em odd}) if
$j$ is even (resp. odd). The number $j$ will be called the {\em index} of this
$\kappa$-band.

\section{Forbidden Bands and Annuli}
\label{forbidden}

In this section we shall consider \vk diagrams with a fixed boundary
label over the presentation $\pp_N(\sss)$. Let $\Delta$ be such a diagram.

For every $j=1,\ldots,2N$ we define an automorphism of the group $G'_{N}(\sss)$
given by the presentation $\pp_N(\sss)$ without the hub.
Take a perfect $4N$-gon with boundary labelled by the word
$K(1)=(\kappa_1...\kappa_{2N})(\kappa_{2N}...\kappa_1)\iv$ (see
formula (\ref{hub}) for the definition of $K(u)$).
Consider the reflection
with respect to the axis passing through the mid
points of the opposite edges labelled by $\kappa_j^{\pm 1}$.
Then this reflection induces a permutation $\phi_j$ on the
set $$\{\kappa_1, \ldots, \kappa_{2N}\}\cup\{\kappa_1\iv,...,\kappa_{2N}\iv\}$$
which satisfies the property $\phi_j(\kappa_i\iv)=\phi_j(\kappa_i)\iv$.
We extend $\phi_j$ to a permutation of the set $A\cup A\iv$ by
fixing all other letters from $A\cup A\iv$. It is easy to see that
$\phi_j$ takes every relation from $\pp_{N}(\sss)$, except for the hub,
to a non-hub relation
from $\pp_{N}(\sss)$. Therefore these maps induce automorphisms
of $G'_{N}(\sss)$ which
will be also denoted by $\phi_j$. One notices also that these
automorphisms are involutions (each one is its own inverse).

\begin{lm}\label{cent} If $\Delta$ contains a $\kappa_j$-annulus
then there
exists a
\vk diagram with the same boundary label as $\Delta$, smaller area
and the same or smaller diameter.
\end{lm}
{\bf Proof.}
Consider each $\kappa_j$-cell as a $\kappa_j$-band.
Then it is easy to verify
(see the list of relations in the presentation) that the label of the
top and the label of the bottom paths of this band
are images of each other under
the automorphism  $\phi_j$
of the free group.  Therefore
the boundaries of the annular diagram formed by any $\kappa_j$-annulus
$\kkk$
have labels one of which is the image of another one under $\phi_j$
in the free group (when we formed the boundaries
of the annular diagram, we reduce pairs of consecutive mutually inverse
edges and this does not change the values of the labels of these boundaries
in the free group). We can assume that $\kkk$ is an innermost $\kappa_j$-annulus,
that is,
there are no $k_j$-annuli inside the subdiagram $\Delta'$ bounded by the median of $\kkk$.
Since the contour of $\kkk$ contains no $\kappa_j$-edges, the subdiagram $\Delta'$
does not contain $\kappa_j$-edges (otherwise the maximal $\kkk_j$-band containing
this edge would have to be an annulus). This implies that $\Delta'$
does not contain hubs.
\setcounter{pdsixteen}{\value{ppp}}
That is, $\Delta'$ is a \vk diagram over the presentation
of $G'_N(\sss)$ (see Figure \thepdsixteen).

\unitlength=1.00mm
\special{em:linewidth 0.4pt}
\linethickness{0.4pt}
\begin{picture}(156.06,129.28)
\bezier{192}(9.28,121.61)(38.28,129.28)(45.61,112.61)
\bezier{488}(46.17,112.61)(76.83,67.28)(10.17,71.28)
\put(38.17,96.17){\oval(21.56,24.44)[]}
\put(38.50,96.50){\oval(13.33,18.44)[]}
\put(38.06,105.50){\line(0,1){2.89}}
\put(41.17,104.83){\line(1,4){0.89}}
\put(42.06,108.39){\line(0,0){0.00}}
\put(43.17,103.72){\line(3,1){4.44}}
\put(45.17,97.94){\line(1,0){3.78}}
\put(48.94,91.06){\line(-1,0){4.44}}
\put(40.94,87.50){\line(2,-3){2.44}}
\put(32.06,83.94){\line(3,5){2.44}}
\put(31.83,92.83){\line(-1,0){4.44}}
\put(27.39,102.39){\line(3,-1){4.67}}
\put(7.17,93.28){\makebox(0,0)[cc]{$\Delta$}}
\put(38.06,95.94){\makebox(0,0)[cc]{$\Delta'$}}
\put(30.06,102.61){\makebox(0,0)[cc]{$\kappa_j$}}
\put(37.61,110.39){\makebox(0,0)[cc]{$p$}}
\put(38.28,101.72){\makebox(0,0)[cc]{$\phi_j(p)$}}
\bezier{192}(88.50,121.61)(117.50,129.28)(124.83,112.61)
\bezier{488}(125.39,112.61)(156.06,67.28)(89.39,71.28)
\put(117.39,96.17){\oval(21.56,24.44)[]}
\put(117.72,96.50){\oval(13.33,18.44)[]}
\put(117.28,105.50){\line(0,1){2.89}}
\put(120.39,104.83){\line(1,4){0.89}}
\put(120.61,104.39){\line(0,0){0.00}}
\put(122.39,103.72){\line(3,1){4.44}}
\put(124.39,97.94){\line(1,0){3.78}}
\put(128.17,91.06){\line(-1,0){4.44}}
\put(120.17,87.50){\line(2,-3){2.44}}
\put(111.28,83.94){\line(3,5){2.44}}
\put(111.06,92.83){\line(-1,0){4.44}}
\put(106.61,102.39){\line(3,-1){4.67}}
\put(86.39,93.28){\makebox(0,0)[cc]{$\Delta$}}
\put(117.28,95.94){\makebox(0,0)[cc]{$\phi_j(\Delta')$}}
\put(109.28,102.61){\makebox(0,0)[cc]{$\kappa_j$}}
\put(116.83,110.39){\makebox(0,0)[cc]{$p$}}
\put(117.83,101.39){\makebox(0,0)[cc]{$p$}}
\bezier{192}(46.17,52.61)(75.17,60.28)(82.50,43.61)
\bezier{488}(83.06,43.61)(113.72,-1.72)(47.06,2.28)
\put(75.06,27.17){\oval(21.56,24.44)[]}
\put(78.28,35.39){\line(0,0){0.00}}
\put(44.06,24.28){\makebox(0,0)[cc]{$\Delta$}}
\put(74.94,26.94){\makebox(0,0)[cc]{$\phi_j(\Delta')$}}
\put(74.50,41.39){\makebox(0,0)[cc]{$p$}}
\put(57.06,98.28){\vector(1,0){10.33}}
\put(89.39,70.61){\vector(-1,-1){12.33}}
\bezier{276}(46.06,52.28)(22.72,26.94)(46.72,2.28)
\bezier{276}(8.72,121.28)(-14.61,95.94)(9.39,71.28)
\bezier{276}(88.39,121.28)(65.06,95.94)(89.06,71.28)
\end{picture}
\begin{center}
\nopagebreak[4]
Fig. \theppp.

\end{center}
\addtocounter{ppp}{1}
\vskip 0.2 in

Let us apply $\phi_j$ to all labels of $\Delta'$.
Since $\phi_j$ takes relations of $G'_N(\sss)$
to relations of this group, the subdiagram $\Delta'$ turns into
another \vk diagram $\phi_j(\Delta')$ over $\pp_{N}(\sss)$.
We can also notice that the labels of the inner and the outer contours
of our annulus are now the same in the free group. By applying
a 0-refinement one
can make these labels graphically equal. After that one can identify these
contours and remove all the cells of the annulus from the diagram.
The result will be
a \vk diagram $\Sigma$ with the same boundary label as $\Delta$ and smaller
area.

It is also easy to show that the diameter of the diagram can only decrease when
we remove a $\kappa$-annulus. Indeed, with every path $p$ in $\Delta$ one
can associate a path $p'$ in $\Sigma$ conecting the same vertices. This path
is obtained by removing some edges (the common edges of the removed annulus)
and identifying the end vertices of these edges, and by relabeling the parts
of $p$ which are contained in the subdiagram bounded by the annulus.
The length of $p'$ can be only smaller than the length of $p$. Therefore
the diameter of $\Sigma$ can be only smaller than the diameter of $\Delta$.
$\Box$

Now we can define {\em reduced} van Kampen diagram over $\pp_{N}(\sss)$.

\begin{df} {\rm A van Kampen diagram over a presentation $\pp$
is {\em reduced} if it does not contain
any of the following:

\begin{enumerate}
\item A 0-edge (that is an edge labelled by 1).
\item
\setcounter{pdfour}{\value{ppp}}
A {\em reducible} pair of cells (see Figure \thepdfour):
\begin{center}
\unitlength=1mm
\special{em:linewidth 0.4pt}
\linethickness{0.4pt}
\begin{picture}(98.00,25.00)
\put(74.00,10.00){\vector(0,1){14.33}}
\bezier{200}(74.00,10.33)(50.00,19.00)(74.00,25.00)
\bezier{200}(74.00,10.33)(98.00,17.00)(74.00,25.00)
\put(71.00,17.00){\makebox(0,0)[cc]{$u$}}
\put(59.00,17.67){\makebox(0,0)[cc]{$v$}}
\put(88.67,17.33){\makebox(0,0)[cc]{$v$}}
\put(68.67,23.33){\vector(3,1){2.67}}
\put(79.33,23.00){\vector(-2,1){2.33}}
\end{picture}
\begin{center}
\nopagebreak[4]
Fig. \theppp.

\end{center}
\addtocounter{ppp}{1}
\end{center}
Here $v$ is  not empty.
\item A $\kappa_i$-annulus (for any $i$).
\end{enumerate}
}
\end{df}

Lemma \ref{cent} shows that if $\Delta$ has a $\kappa$-annulus then its
area can be decreased without increasing the diameter. It is easy to see
that if $\Delta$ contains a
reducible pair of cells, its area can decrease too:
both cells
can be removed by removing the path labelled by $u$ and identifying the
paths labelled by $v$. Unfortunately the diameter of the diagram can increase.
This does not cause any difficulty in computing an upper bound for the smallest
isodiametric function, but this makes it more difficult to find a good lower
bound (see Section \ref{lowbound}).

For each word $w$ which
is equal to 1 in our group $G_{N}(\sss)$ there exists a reduced diagram with
boundary label $w$. Indeed, it is enough to take the \vk diagram with minimal
number of cells and boundary label $w$ (such a diagram exists by van Kampen's
lemma).

\begin{df} \label{annularreduced}
{\rm An annular \vk diagram $\Psi$ over $\pp_{N}(\sss)$ will be called
{\em reduced} if any ordinary \vk diagram obtained by cutting $\Psi$
along a simple path connecting $\partial_0(\Psi)$ with $\partial_i(\Psi)$
is reduced.
}
\end{df}

Notice that if we obtain an annular diagram $\Psi$ by removing a subdiagram
from a reduced \vk Diagram $\Delta$ then $\Psi$ is reduced.

For the rest of this section
we shall assume that $\Delta$ is a {\it reduced diagram} over $\pp_{N}(\sss)$.

\begin{lm}\label{q} If $\Delta$ is reduced and does not have hubs then
it does not have $Q_j$-annuli for any $j=1,...,15k+6$.
\end{lm}

{\bf Proof.} Suppose there is a $Q_j$-annulus. Let $\qq$ be the $Q_i$-annulus
with the smallest inside diagram $\Delta'$. Since $Q$-annuli cannot cross
(this follows from the definition of an $S$-rule),
$\Delta'$ does not contain $Q$-cells (otherwise we would get
a $Q$-annulus with smaller inside diagram).
Therefore all cells in $\Delta'$ are auxiliary.
Since every cell in $\qq$
has exactly two $Q_i$-edges, the boundary $\partial(\Delta')$ contains no
$Q$-edges.  Thus $\Delta'$ is a diagram over the presentation
consisting of all auxiliary (commutativity) relations.
The group $H$ given by this presentation
has a homomorphism onto the direct product of the
free group on $\Theta$ and
the free group on $\gam\cup\{\cee, \dol, \kappa_1,\ldots,\kappa_{2N}\}$.
Therefore the
word $u=\Lab(\partial(\Delta'))$ is equal to 1 in this direct product.
Every $Q$-cell in $\qq$, considered as a $Q$-band, contains exactly one
$\Theta$-edge on the
top path (resp. on the bottom path).
If two of these $\Theta$-edges cancel
when we form the top/bottom path of $\qq$, then the corresponding cells of
$\qq$
cancel, by Lemma \ref{cancel} (b). This cannot happen since
$\Delta$ is reduced. Therefore $u$ contains
a $\Theta$-letter.

This implies that $u$ contains a subword of the form
$\tau^{\pm 1}v\tau^{\mp 1}$ for some letter $\tau \in \Theta_+$
and some word $v$ which does not
contain $\Theta$-letters. Then the $\Theta$-letters in this subword must label
edges of consecutive cells $\pi$ and $\pi'$
in $\qq$ (again we use the fact that $\Theta$-edges do
not cancel when we form the top/bottom path of $\qq$). Then the cells
$\pi$ and $\pi'$ cancel. Indeed by Lemma \ref{cancel} (b)
the relation $s_1$  corresponding to the cell $\pi_1$ is
a cyclic shift of $s_2^{\pm 1}$ where $s_2$ is the relation corresponding
to $\pi_2$. We can assume that we start reading $s_1$ and $s_2$ at the beginning of
their common $Q_j$-edge. Since the $\Theta_+$-edges on the bottoms of these
cells are oriented toward each other, we conclude that $s_1=s_2\iv$ and each of the
words $s_1$ and $s_2$ contains only two $\Theta$-edges, $s_1=s_2\iv$.
So the cells $\pi_1$ and $\pi_2$ cancel.

This again contradicts the assumption that
$\Delta$ is reduced. $\Box$

\begin{lm} \label{rpath} If $\Delta$ is reduced, does not have hubs and
$\partial(\Delta)$ consists of $\Theta$-edges then $\Delta$ does not have
cells.
\end{lm}

{\bf Proof.} Indeed, by Lemmas \ref{q} and \ref{cent},
$\Delta$ has no $Q_j$-annuli. Thus it
cannot contain $Q_j$-cells.
Therefore all cells in $\Delta$ are auxiliary.
Hence, as in the proof of Lemma \ref{q}, the boundary
label of $\Delta$ is equal to 1 in the free group generated by $\Theta$.
Since $\Delta$ is reduced, it does not contain cells. $\Box$

\begin{lm} \label{qr} If $\Delta$ is reduced and
does not contain hubs, then it does not have
$(Q_j,\Theta)$-annuli.
\end{lm}

{\bf Proof.} Suppose that $\Delta$ contains a $(Q_j,\Theta)$-annulus
${\cal W}_1\cup {\cal W}_2$
where ${\cal W}_1$ is a $Q_j$-band and ${\cal W}_2$ is a $\Theta$-band.
Let $\Delta'$ be the
inside diagram of this annulus. We can assume that $\Delta'$ has the smallest
possible area.

Then the contour of $\Delta'$ does not
contain $Q$-edges
(otherwise there would be a $(Q_i, R)$-annulus with a smaller
inside diagram).
Therefore $\Delta'$ does not contain any transition cells (by Lemma \ref{q}
there
are no $Q$-annuli). So it
can only contain auxiliary cells. Thus $\Delta'$ is a diagram over the group
$H$
given by the auxiliary relations only.
We know from the proof of Lemma \ref{q} that this group has a homomorphism onto
the direct product of the free group over $\Theta$ and the free group generated
by other letters (not from $\bigcup Q_i\cup \Theta$).  Thus if we remove all
non-$\Theta$-letters from the label $w'$ of the boundary of $\Delta'$, we get a
word which is equal to 1 in the free group generated by $\Theta$. Therefore $w'$
contains a subword $t=\tau^{\pm 1}v\tau^{\mp 1}$ as in Lemma \ref{q} or else the
contour of $\Delta'$ does not have $\Theta$-edges at all.

Suppose first that the contour of $\Delta'$ contains $\Theta$-edges.  The part
of the contour of $\Delta'$ which is contained in the contour of ${\cal W}_2$
cannot contain $\Theta$-edges because the contour of each $\Theta$-cell contains
exactly two $\Theta$-edges which form an opposing pair.  Thus both $\tau$'s in
the subword $t$ must be labels of some edges on the part of $\partial(\Delta')$
which is contained in ${\cal W}_1$.  These two $\tau$'s must belong to two
neighboring transition cells, so $\Delta$ is not reduced (we apply the same
argument as in Lemma \ref{q}).

\setcounter{pdten}{\value{ppp}}

Now suppose that $\partial(\Delta')$ does not contain $\Theta$-edges (see Figure
\thepdten).

\unitlength=1.00mm
\special{em:linewidth 0.4pt}
\linethickness{0.4pt}
\begin{picture}(113.00,45.00)
\put(98.67,38.67){\line(0,1){6.33}}
\put(98.67,45.00){\line(1,0){14.33}}
\put(113.00,45.00){\line(0,-1){16.67}}
\put(113.00,28.33){\line(-3,-4){3.67}}
\put(109.33,23.33){\line(-1,0){7.33}}
\put(101.67,23.00){\line(-3,-4){3.33}}
\put(98.33,18.67){\vector(0,-1){9.67}}
\put(98.33,9.00){\line(0,-1){7.67}}
\put(98.33,1.33){\line(1,0){14.67}}
\put(113.00,1.33){\line(0,1){17.33}}
\put(113.00,18.67){\line(-4,5){3.67}}
\put(105.67,26.67){\makebox(0,0)[cc]{$q$}}
\put(100.00,33.00){\makebox(0,0)[cc]{$\tau$}}
\put(100.67,14.33){\makebox(0,0)[cc]{$\tau'$}}
\bezier{800}(98.33,38.67)(-0.33,26.00)(98.67,9.00)
\bezier{568}(98.33,28.33)(27.67,26.00)(98.33,19.00)
\put(70.00,22.33){\vector(-3,-1){11.33}}
\put(82.67,20.33){\vector(-2,-3){4.67}}
\put(98.50,28.17){\vector(0,1){10.50}}
\put(100.83,21.83){\line(2,3){1.00}}
\put(90.83,28.00){\vector(-1,4){2.33}}
\put(82.00,27.67){\vector(-1,2){4.00}}
\put(73.33,27.00){\vector(-2,3){4.50}}
\put(64.50,25.83){\vector(-4,1){11.83}}
\put(101.67,23.33){\line(-2,3){3.33}}
\end{picture}

\begin{center}
\nopagebreak[4]
Fig. \theppp.

\end{center}
\addtocounter{ppp}{1}
\vskip 0.2 in

Then ${\cal W}_1$ contains just 2 cells, the intersection cells of ${\cal W}_1$
and ${\cal W}_2$.  Then the edges labelled by $\tau$ and $\tau'$ on
Fig. \thepdten\ belong to ${\cal W}_2$. The $\tau$'s in every $\Theta$-cell of
${\cal W}_2$ are the
same. Therefore $\tau=\tau'$.  This implies that the two $Q$-cells in ${\cal
W}_1$ form a reducible pair (Lemma \ref{cancel} b) ), a
contradiction. $\Box$

\begin{lm} \label{cr} Suppose that (not necessarily reduced)
$\Delta$ contains
a $(\kappa_j, \Theta)$-annulus $\aaa$ consisting of a $\kappa_j$-band ${\cal W}_1$
and a $\Theta$-band ${\cal W}_2$. Suppose that ${\cal W}_1$ does not have hubs.
Then $\ww_1$  contains a reducible pair of cells which belong to the same
$\Theta$-band.
\end{lm}

{\bf Proof.}
If there is a $\kappa_j$-cell in ${\cal W}_1$ between the intersection cells
$\rho_1$ and $\rho_2$ of $\ww_1$ and $\ww_2$ then this cell will belong
to a $\Theta$-band which forms
a $(\kappa_j, \Theta)$-band with a part of $\ww_1$ and the inside diagram of this
annulus would be smaller than the inside diagram of $\aaa$. Therefore we can assume
that the cells $\rho_1$ and
$\rho_2$ are consecutive cells on $\ww_1$. The start and the end edges of $\ww_2$
form an opposing pair with the same label.
Therefore the cells $\rho_1$
and $\rho_2$ correspond to mutually inverse relations and have a common
$\kappa_j$-edge. Thus these cells form a reducible pair and belong to the same
$\Theta$-band. $\Box$

\vskip 0.1 in

\begin{lm} \label{r} If $\Delta$ is reduced and
does not contain hubs then it does not
contain
$\Theta$-annuli.
\end{lm}

{\bf Proof.} Take an innermost $\Theta$-annulus $\bb$. Let $\Delta'$ be the
inside subdiagram formed by this annulus.
Then $\Delta'$ does not contain $\Theta$-annuli.
Since the contour of $\Delta'$ contains no $\Theta$-edges, the diagram
$\Delta'$ contains no $\Theta$-edges either. Since the only cells that
contain no $\Theta$-edges are the hubs (which were ruled out here), $\Delta'$
contains no cells. Therefore $\Lab(\partial(\Delta'))=1$ in the free group.
Since $\Delta$ contains no $(Q,\Theta)$-annuli (Lemma \ref{qr}),
$\partial(\Delta')$ contains
no $Q$-edges.
Therefore $\partial(\Delta)$ contains only ($\bar Y$)-edges
and all cells of $\bb$ are auxiliary.
Since $\Lab(\partial(\Delta'))=1$ there
are two neighboring edges in $\partial(\Delta')$ with mutually inverse edges.
The cells of $\bb$ containing these two edges form a reducible pair,
by Lemma \ref{cancel} (c).
This contradicts the
assumption that $\Delta$ is reduced. $\Box$

\begin{lm} \label{gammapath} If $\Delta$ is reduced,
does not contain hubs and
$\partial(\Delta)$
has no $\Theta$-edges then $\Delta$ does not have cells.
\end{lm}

{\bf Proof.} Since $\Delta$ contains no $\Theta$-annuli (Lemma \ref{r}),
the assumptions of the Lemma imply that $\Delta$ contains no $\Theta$-edges
at all.
Since every relation in $\pp_{N}(\sss)$ except for the hub contains letters from
$\Theta$, the diagram $\Delta$ has no cells. $\Box$

\begin{lm}\label{ar} If $\Delta$ is reduced and
contains no $Q$-edges then it does not
have $(\bar Y, \Theta)$-annuli.
\end{lm}

{\bf Proof.}
Since $\Delta$ contains no $Q$-edges, all cells in $\Delta$
are auxiliary. Take an innermost $(\bar Y, \Theta)$-annulus.
Let $\Delta'$ be the
inside diagram of this annulus. Then $\partial(\Delta')=pq$ where $\Lab(p)$
contains no $\Theta$-edges and $\Lab(q)$ contains no $\bar Y$-edges ($p$
is the bottom or the top path of a $\Theta$-band, $q$ is the bottom or the top
path of a $\bar Y$-band).
It  follows from the structure of auxiliary $\bar Y$-cells that the
top and the bottom paths of every non-empty $\bar Y$-band must
contain $\Theta$-edges.
Therefore either $q$ contains $\Theta$-edges or $q$ is empty. The first case is
impossible because our $(\bar Y, \Theta)$-annulus was innermost. Therefore
$q$ is empty. Then $\partial(\Delta')$ contains no $\Theta$-edges.
By Lemma \ref{gammapath}, $\Delta'$ contains no cells. Therefore
$\Lab(\partial(\Delta'))=1$ in the free group. This implies, as before,
that our $(\bar Y, \Theta)$-annulus
is not reduced.  $\Box$

\begin{lm}\label{a} If $\Delta$ is reduced and
does not have hubs then it does not
have $\bar Y$-annuli.
\end{lm}

{\bf Proof.} Suppose that $\Delta$ has a ${\bar Y}$-annulus. Take the
innermost ${\bar Y}$-annulus $\cal Y$. Then the outer contour $p$
of this annulus
cannot contain $Q$-edges since a ${\bar Y}$-band cannot contain transition cells
(by the definition of ${\bar Y}$-bands).
Lemma \ref{q} implies that the subdiagram bounded by $p$ does not contain
$Q$-edges. It also cannot contain $\kappa$- or $\Theta$-edges because of
Lemmas \ref{ar} and \ref{r}. Thus the diagram bounded by $p$ must be empty,
a contradiction (this diagram contains ${\cal Y}$).
$\Box$

\section{Diagrams Without Hubs}
\label{diagrams}

In this section, we estimate the area and diameter of a diagram without
hubs over the presentation $\pp_N(\sss)$ and also consider the process of reducing
a non-reduced diagram without hubs.

With every diagram $\Delta$ over $\pp_N(\sss)$ we associate two numbers:
$n(\Delta)$ is the length of the boundary of $\Delta$ and $h(\Delta)$ is
the number of hubs in $\Delta$.

\begin{lm} Let $\Delta$ be a reduced diagram without hubs
over the presentation $\pp_N(\sss)$, with perimeter $n=n(\Delta)$
and with $m$ maximal $Q$-bands.
Then $m\le n$ and the area of
$\Delta$ is at most $Cn^2m$ where $C$ is a constant. Therefore
the area of $\Delta$ does not exceed $Cn^3$. The diameter of $\Delta$ does
not exceed $C'n$ for some constant $C'$.
\label{nohub}
\end{lm}

{\bf Proof.} Every cell in $\Delta$ belongs to a maximal $\Theta$-band.
Since $\Delta$ contains no $\Theta$-annuli (Lemma \ref{r}), each of
these $\Theta$-bands starts and ends on the contour of $\Delta$.
Therefore there are at most $n$ maximal $\Theta$-bands. A similar reasoning
shows that $m\le n$.

Every transition cell is an intersection of a $\Theta$-band and a $Q_j$-band.
Since $\Delta$ does not contain $(Q_j, \Theta)$-annuli (Lemma \ref{qr}), each
$\Theta$-band intersects each $Q_j$-band at most once. Therefore $\Delta$
contains at most $mn/4$ transition cells (the total number of intersections is
bounded by $mn$ but each transition cell is  counted at least 4 times because it
belongs to two mutually inverse $Q_j$-bands and to two mutually inverse
$\Theta$-bands).

Let $\Phi$ be the set of all maximal $\bar Y$-bands.  Every non-transition cell
with $\bar Y$-edges in $\Delta$ belongs to two mutually inverse bands in $\Phi$.
Therefore we need to compute the number of bands in $\Phi$ and their lengths.

Since $\Delta$ does not contain $\bar Y$-annuli (Lemma \ref{a}),
each of the $\bar Y$-bands in $\Phi$ starts either on
the contour of  $\Delta$
or on the contour of one of the transition cells. The number of
$\bar Y$-edges on a transition cell
is bounded by a constant. Therefore the number of
$\bar Y$-bands in $\Phi$
is at most $Cmn$ for some constant $C$. Every cell in each of the
bands in $\Phi$ is an intersection of a $\bar Y$-band and a $\Theta$-band.
By Lemma \ref{ar}, each of the $\bar Y$-bands in $\Phi$ intersects each of the
$\Theta$-bands in $\Delta$ at most once (we use the fact that $\bar Y$-bands
in
$\Phi$ do not contain transition cells,
and Lemma \ref{qr}). So the length of each
of the bands in $\Phi$ does not exceed $n/2$.
Therefore the number of non-transition $\bar Y$-cells
in $\Delta$ does not exceed $C_1mn^2$ for some constant $C_1$.

Each non-transition cell which does not contain $\bar Y$-edges
is the intersection
of a $\Theta$-band and a $\kappa_j$-band. Since $\Delta$ does not contain
$\kappa$-annuli, the total number of
maximal $\kappa$-bands (including inverses) is at most $n$.
By Lemmas \ref{cr} each $\Theta$-band intersects each of the
$\kappa$-band
at most once. Therefore the number of non-transition cells without
$\bar Y$-edges in $\Delta$ does not exceed $n^2/4$.

Thus the total number of cells in $\Delta$ does not exceed
$C_1mn^2+n^2/4 < C_2mn^2$.

In order to estimate the diameter of $\Delta$ pick a vertex $v$ in $\Delta$.
Let $d$ be the distance from $v$ to the boundary of $\Delta$. It is clear
that $n/2$ plus the maximal among these $d$'s for all vertices $v$
is the upper bound for the diameter
of $\Delta$, so we need to estimate $d$. Without loss of generality
we can assume that $v$ is not on the boundary of $\Delta$.
The vertex $v$ belongs either to a $\kappa$-cell, or to a $\bar Y$-cell.
In the first case $v$ belongs to the boundary of a $\kappa$-band in $\Delta$.
Since the length of this $\kappa$-band is at most $n$, we conclude that
$d\le n$.

In the second case $v$ is within a constant distance $d_1$
from a $\bar Y$-edge $e$.
Consider the maximal $\bar Y$-band $\aaa$ containing $e$. We know that the length $d_2$
of this band is at most $n$. The length of the boundary of this band is $O(n)$.
The band $\aaa$  ends either on the contour of
the diagram $\Delta$ or on the boundary of a transition cell. In the first
of these cases $d\le d_1+O(n)\le O(n)$. In the second case we can conclude that
$v$ is within distance $O(n)$ from the boundary of a $Q$-band. As we know, the
length of every maximal $Q$-band in $\Delta$ is at most $n$. Therefore the
length of the boundary of any  $Q$-band is $\le O(n)$. Since every maximal
$Q$-band in $\Delta$ ends and starts on the contour of $\Delta$, we can conclude
that $d\le O(n)$. Therefore the diameter of $\Delta$ is $O(n)$.

The lemma is proved. $\Box$

Later we shall need to analyze the process of reducing a
non-reduced diagram without hubs over $\pp_{N}(\sss)$.
First we prove the following general results.

\begin{lm} \label{identification} Let $\Psi$ be a non-reduced diagram over any
presentation.  Consider a  process of cancelling a sequence of reducible pairs of cells in
$\Psi$.  If two edges $e$ and $f$ in $\Psi$ get identified in this process then
there exists a path $t$ in $\Psi$ connecting the initial vertices $\imath(e)$
and $\imath(f)$ and such that $\Lab(t)=1$ in the free group.  \end{lm}

{\bf Proof.} We shall use induction on the length of the cancellation process.
If $e=f$ in $\Psi$ (that is if the length of the process is 0) then
the statement is obvious. Suppose that $e\neq f$ in $\Psi$
and $e$ is identified with $f$ after $c$ cancellations.
Let $\Psi'$ be the diagram obtained from $\Psi$ after $c-1$ cancellations.
There exist a pair of cells $\pi$ and $\pi'$ in $\Psi'$
such that
$e\iv\in \partial(\pi)$, $f\in \partial(\pi')$, $\pi$ and $\pi'$ form a
reducible pair in $\Psi'$ and $e$ and $f$ are identified after we cancel
$\pi$ and $\pi'$. Since $\pi$ and $\pi'$ form a reducible pair in $\Psi'$,
the boundaries of these cells have a common edge $g$. Since $e$ and $f$
are corresponding edges of $\pi$ and $\pi'$ the path which
connects $\imath(g)$ with $\imath(e)$ on $\partial(\pi)\iv$ (oriented
counterclockwise) and the path connecting  $\imath(g)$ with $\imath(f)$
on $\partial(\pi')$
have equal labels. Let $\pi_1$ and $\pi'_1$ be the ``preimages"
of the cells $\pi$ and $\pi'$ in $\Psi$. Then these cells contain edges
$g_1$ and $g_2$ respectively such that
\begin{itemize}
\item the path $t$ connecting $\imath(g_1)$ with $\imath(e)$ on $\partial(\pi_1)\iv$
and the path $t'$ connecting $\imath(g_2)$ with $\imath(f)$ on
$\partial(\pi_1')$ have equal labels and
\item the edges $g_1$ and $g_2$ are identified (with the edge $g$) after
$c-1$ cancellations.
\end{itemize}
By the induction hypothesis there exists a path
$p$ in $\Psi$ connecting $\imath(g_1)$ with $\imath(g_2)$ whose label is equal to
1 in the free group. Then the path $t\iv pt'$ in $\Psi$
connects $\imath(e)$ with $\imath(f)$ and its label equals 1 in the free group.
The lemma is proved.  $\Box$

\begin{lm}\label{cells}
Let $\Psi$ be a non-reduced diagram over any presentation.
Consider a process of cancelling a sequence of reducible pairs of cells in $\Psi$.
If two cells $\pi$ and $\pi'$ in $\Psi$ cancel in this process
then these cells correspond to mutually inverse relations
and there exists a path $p$ in $\Psi$ connecting corresponding vertices
in $\partial(\pi)$ and $\partial(\pi')$ whose label is equal to 1 in the
free group.
\end{lm}

{\bf Proof.} Indeed it is clear that $\pi$ and $\pi'$ correspond to mutually
inverse relations. After a number of cancellations in $\Psi$ an edge $e$ on
$\partial(\pi)$ is identified with the corresponding edge $f$ on
$\partial(\pi')\iv$. It remains to apply Lemma \ref{identification}.  $\Box$

\begin{lm} \label{kbands}
Let  $\Psi$ be a non-reduced diagram over the presentation $\pp_{N}(\sss)$
without $\kappa$-annuli and without hubs. Suppose that each $\kappa$-band in
$\Psi$ is reduced.  Consider a process of cancelling a sequence of reducible pairs of cells in
$\Psi$.  Then no $\kappa$-annuli can appear as a result of this process and no
two $\kappa$-cells from the same $\kappa$-band cancel.
\end{lm}

{\bf Proof.} If we cancel two non-$\kappa$ cells then $\kappa$-bands
do not change.
Let $\pi$ and $\pi'$ be $\kappa$-cells which form a reducible pair in $\Psi$.
Let $\bb$ and $\bb'$ be the maximal $\kappa$-bands which contain $\pi$ and $\pi'$
respectively. Since $\pi$ and $\pi'$ are not hubs, they contain two $\kappa$-edges
each. Let $e$ and $f$ be the $\kappa$-edges in $\pi$ and $e'$ and $f'$ are the
$\kappa$-edges in $\pi'$.
Notice
that the start and end edges of $\bb$ and $\bb'$ belong to the contour
of $\Psi$. Finally assume without loss of generality that the median
of $\bb$ (resp. $\bb'$)
first crosses $e$ and then $f$ (resp. $e'$ and $f'$). When we cancel $\pi$ and
$\pi'$, the $\kappa$-edges in $\pi$ are identified with the corresponding
$\kappa$-edges in $\pi'$. By replacing $\bb$ and $\bb'$ with their inverses
and renaming the edges $e, f, e', f'$ if necessary, we can
assume that $e$ corresponds to $e'$ and $f$ corresponds to $f'$.
Let $\bb=(\pi_1,\ldots,\pi_i, \pi, \pi_{i+1}, \ldots, \pi_{t})$,
$\bb'=(\pi_1',\ldots,\pi_{j}, \pi', \pi_{j+1},\ldots,\pi_{s})$.
Then it is easy to check that after the cancellation we shall have two
$\kappa$-bands ($\pi_1,\ldots,\pi_i, \pi_j, \pi_{j-1},\ldots,\pi_1$) and
($\pi_{t},\ldots,\pi_{i+1}, \pi_{j+1},\ldots,\pi_s$) (and their inverses) instead of
$\bb$ and $\bb'$ (one or both of these bands may be empty).
Both these $\kappa$-bands start and end on the contour
of $\Psi$. All other maximal $\kappa$-bands do not change. Therefore
we do not get $\kappa$-annuli by reducing $\pi$ and $\pi'$.

Now suppose that two cells $\pi$ and $\pi'$ of $\Psi$
from the same $\kappa$-band $\bb$
eventually
cancel as a result of our process of cancellation. Then by Lemma \ref{cells}
these cells correspond to mutually inverse relations and
there exists a path $t$ in $\Psi$
connecting corresponding vertices $v$ and $v'$
of these cells
whose label is 1 in the free group. It is easy to see that $v$ and $v'$
both belong to either $\topp(\bb)$ or to $\bott(\bb)$. Therefore they
are connected by a path $t'$ which is contained in $\topp(\bb)$ or $\bott(\bb)$.
Since $t$ and $t'$ connect the same vertices in $\Psi$, their labels must be
the same modulo the presentation obtained from $\pp_{N}(\sss)$ by removing the hub.
Therefore $\Lab(t')=1$ modulo this presentation.
But $\Lab(t')$ is a word over $\Theta$. By Lemma \ref{rpath},
this word must be equal to 1 in the free group. Therefore $\topp(\bb)$ contains
two consecutive edges whose labels are mutually inverse.
Then the corresponding cells form a reducible pair.
This contradicts the assumption that every $\kappa$-band in $\Psi$ is reduced.
The lemma is proved. $\Box$

\section{Sectors}
\label{st12}

By a {\em sector} we mean a reduced diagram $\Delta$ over $\pp$
with boundary divided into four parts,
$\partial(\Delta)=p_1p_2p_3\iv p_4\iv$,
such that the following properties hold:
\begin{itemize}
\item $p_1$ is the top path of an odd $\kappa$-band $\kkk$;
\item $p_3$ is the bottom path of an even $\kappa$-band $\ttt$;
\item $p_2$ is the top path of a $\Theta$-band $\rr_t$, and the label $\Lab(p_2)$
is a reduced word;
\item $\Lab(p_4)=\kappa_j^{\pm 1}W\kappa_s^{\pm 1}$
where $W$ is an admissible word for which there exists a computation
connecting $W$ with $W_0$
\setcounter{pdeleven}{\value{ppp}}
(see Figure \thepdeleven).
\end{itemize}

\unitlength=1.00mm
\special{em:linewidth 0.4pt}
\linethickness{0.4pt}
\begin{picture}(105.72,56.78)
\put(53.72,5.94){\line(1,0){36.33}}
\put(90.06,5.94){\line(1,3){15.67}}
\put(105.72,52.94){\line(-1,0){56.00}}
\put(49.72,52.94){\line(-1,0){6.00}}
\put(43.72,52.94){\line(1,-5){9.33}}
\put(59.06,5.94){\line(-1,0){6.33}}
\put(57.39,5.94){\line(-1,5){9.33}}
\put(85.72,5.94){\line(1,3){15.67}}
\put(101.39,52.61){\line(1,0){4.33}}
\put(101.72,52.78){\line(1,0){4.00}}
\put(52.39,9.11){\line(1,0){4.33}}
\put(51.56,12.78){\line(1,0){4.50}}
\put(50.72,17.28){\line(1,0){4.33}}
\put(49.89,22.28){\line(1,0){4.33}}
\put(48.89,27.28){\line(1,0){4.33}}
\put(47.89,31.78){\line(1,0){4.33}}
\put(47.06,36.44){\line(1,0){4.33}}
\put(46.06,40.94){\line(1,0){4.33}}
\put(45.06,46.11){\rule{4.33\unitlength}{-0.17\unitlength}}
\put(86.72,9.11){\line(1,0){4.50}}
\put(87.89,12.78){\line(1,0){4.50}}
\put(89.56,17.28){\line(1,0){4.17}}
\put(91.06,22.28){\line(1,0){4.50}}
\put(92.89,27.28){\line(1,0){4.17}}
\put(92.89,27.28){\line(1,0){4.33}}
\put(94.39,31.78){\line(1,0){4.33}}
\put(95.89,36.44){\line(1,0){4.33}}
\put(97.39,40.94){\line(1,0){4.33}}
\put(99.06,46.11){\line(1,0){4.50}}
\put(45.72,56.28){\makebox(0,0)[cc]{$\kappa_j^{\pm 1}$}}
\put(103.89,56.28){\makebox(0,0)[cc]{$\kappa_s^{\pm 1}$}}
\put(55.56,3.61){\makebox(0,0)[cc]{$\kappa_j^{\pm 1}$}}
\put(87.89,3.44){\makebox(0,0)[cc]{$\kappa_s^{\pm 1}$}}
\put(49.39,46.11){\line(1,0){50.00}}
\put(53.22,46.11){\line(0,1){6.83}}
\put(61.39,46.11){\line(0,1){6.83}}
\put(61.39,52.94){\line(0,0){0.00}}
\put(79.22,46.11){\line(0,1){6.83}}
\put(89.72,46.11){\line(0,1){6.83}}
\put(64.89,49.78){\makebox(0,0)[cc]{$\tau$}}
\put(45.72,27.94){\makebox(0,0)[cc]{$p_1$}}
\put(72.39,56.78){\makebox(0,0)[cc]{$p_2$}}
\put(100.22,28.11){\makebox(0,0)[cc]{$p_3$}}
\put(69.56,3.28){\makebox(0,0)[cc]{$p_4$}}
\put(71.22,24.61){\makebox(0,0)[cc]{$\Delta'$}}
\put(71.22,46.11){\line(0,1){6.83}}
\put(57.39,29.61){\makebox(0,0)[cc]{$\pleft$}}
\put(86.72,29.94){\makebox(0,0)[cc]{$\pright$}}
\put(70.06,9.94){\makebox(0,0)[cc]{$\pbot$}}
\end{picture}

\begin{center}
\nopagebreak[4]
Fig. \theppp.

\end{center}
\addtocounter{ppp}{1}
\vskip 0.2 in

Recall that a $\kappa_j$-band is called {\em even} (resp. {\em odd}) if
$j$ is even (resp. odd).

If $\Delta$ is a  sector then
by removing the bands $\kkk$, $\ttt$
we get the {\em inside} diagram $\Delta'$ which has a boundary divided into
4 parts $\pleft$, $\ptop$, $\pright\iv$, $\pbot\iv$ such that
$\pleft\subseteq \bott(\kkk)$, $\ptop\subseteq p_2$,
$\pright\subseteq \topp(\ttt)$, $\pbot\subseteq p_4$.

In this section, we shall provide a complete description of sectors.

We start with examples of  sectors which correspond to computations
of our machine $\sss$.

As in the definition of the hub, for every word $u$ let
$$K(u)\equiv
(u\iv\kappa_1u \kappa_2u\iv\kappa_3u\kappa_4\ldots u\iv\kappa_{2N-1}u\kappa_{2N})
(\kappa_{2N}u\iv\kappa_{2N-1}u\ldots \kappa_2 u\iv\kappa_1u)^{-1}.
$$
Notice that if $u=W_0$
then $K(u)$ is the boundary label of the hub.

With every computation $W_1,\ldots, W_g$ we associate a van Kampen diagram over
the presentation $\pp_{N}(\sss)$ in the following way.

If $W_1, W_2$ is a one-step computation which uses a
transition $\tau\in \Theta$
then there exists an annular diagram with boundary labels $K(W_1)$ and $K(W_2)$
of the  form shown in \setcounter{pdfive}{\value{ppp}}
 Figure \thepdfive (in the case when $N=1$):
\vskip 0.2 in
\unitlength=1.00mm
\special{em:linewidth 0.4pt}
\linethickness{0.4pt}
\begin{picture}(133.61,88.06)
\put(80.00,34.83){\line(-1,1){10.00}}
\put(70.00,44.50){\line(0,1){30.33}}
\put(70.00,74.83){\line(1,1){10.00}}
\put(80.00,84.83){\line(1,0){40.00}}
\put(120.00,84.83){\line(1,-1){10.00}}
\put(130.00,74.83){\line(0,-1){30.00}}
\put(130.00,44.83){\line(-1,-1){10.00}}
\put(120.00,34.83){\line(-1,0){40.00}}
\put(80.00,39.83){\line(-1,1){5.00}}
\put(75.00,44.83){\line(0,1){30.00}}
\put(75.00,74.83){\line(1,1){5.00}}
\put(80.00,79.83){\line(1,0){40.00}}
\put(120.00,79.83){\line(1,-1){5.33}}
\put(125.33,74.50){\line(0,-1){29.67}}
\put(125.33,44.83){\line(-1,-1){5.00}}
\put(120.33,39.83){\line(-1,0){40.33}}
\put(70.00,45.00){\line(1,0){5.00}}
\put(80.00,39.83){\line(0,-1){5.00}}
\put(120.00,34.83){\line(0,1){5.00}}
\put(125.33,45.00){\line(1,0){4.67}}
\put(125.33,74.67){\line(1,0){4.67}}
\put(120.00,79.83){\line(0,1){5.00}}
\put(80.00,84.83){\line(0,-1){5.00}}
\put(75.00,75.00){\line(-1,0){5.00}}
\put(73.61,41.39){\vector(1,-1){3.33}}
\put(76.61,43.17){\vector(1,-1){2.00}}
\put(76.50,76.28){\vector(1,1){1.89}}
\put(75.94,80.72){\vector(1,1){1.33}}
\put(122.17,77.72){\vector(1,-1){1.44}}
\put(124.94,79.83){\vector(1,-1){1.33}}
\put(122.17,41.50){\vector(1,1){1.44}}
\put(125.61,40.39){\vector(1,1){1.33}}
\put(102.61,84.72){\vector(1,0){4.33}}
\put(129.94,62.72){\vector(0,1){3.67}}
\put(101.61,79.72){\vector(1,0){5.33}}
\put(119.94,38.39){\vector(0,-1){1.67}}
\put(126.94,45.06){\vector(1,0){1.00}}
\put(127.94,45.06){\vector(0,0){0.00}}
\put(126.28,74.72){\vector(1,0){1.67}}
\put(119.94,81.06){\vector(0,1){1.33}}
\put(119.94,82.39){\vector(0,0){0.00}}
\put(79.94,80.72){\vector(0,1){1.67}}
\put(73.94,75.06){\vector(-1,0){1.67}}
\put(73.94,45.06){\vector(-1,0){2.00}}
\put(79.94,38.39){\vector(0,-1){1.67}}
\put(88.28,79.72){\line(0,1){5.00}}
\put(98.94,79.72){\line(0,1){5.00}}
\put(114.61,79.72){\line(0,1){5.00}}
\put(88.28,34.72){\line(0,1){5.00}}
\put(98.94,39.72){\line(0,-1){5.00}}
\put(114.61,34.72){\line(0,1){5.00}}
\put(106.94,39.72){\vector(-1,0){4.00}}
\put(106.94,34.72){\vector(-1,0){4.00}}
\put(74.94,66.39){\vector(0,-1){7.00}}
\put(69.94,66.39){\vector(0,-1){7.00}}
\put(125.28,59.39){\vector(0,1){7.00}}
\put(125.28,51.72){\line(1,0){4.67}}
\put(129.94,61.72){\line(-1,0){4.67}}
\put(125.28,69.06){\line(1,0){4.67}}
\put(69.94,66.72){\line(1,0){5.00}}
\put(74.94,57.39){\line(-1,0){5.00}}
\put(69.94,49.72){\line(1,0){5.00}}
\put(96.61,88.06){\makebox(0,0)[cc]{$W_1$}}
\put(96.61,76.06){\makebox(0,0)[cc]{$W_2$}}
\put(133.61,60.39){\makebox(0,0)[cc]{$W_1$}}
\put(121.61,60.39){\makebox(0,0)[cc]{$W_2$}}
\put(99.61,42.72){\makebox(0,0)[cc]{$W_2$}}
\put(99.61,29.72){\makebox(0,0)[cc]{$W_1$}}
\put(66.61,61.39){\makebox(0,0)[cc]{$W_1$}}
\put(77.61,61.39){\makebox(0,0)[cc]{$W_2$}}
\put(72.61,82.06){\makebox(0,0)[cc]{$\kappa_1$}}
\put(127.61,81.39){\makebox(0,0)[cc]{$\kappa_2$}}
\put(128.28,37.72){\makebox(0,0)[cc]{$\kappa_1$}}
\put(72.61,37.39){\makebox(0,0)[cc]{$\kappa_2$}}
\put(78.61,74.39){\makebox(0,0)[cc]{$\kappa_1$}}
\put(121.28,75.06){\makebox(0,0)[cc]{$\kappa_2$}}
\put(120.28,44.06){\makebox(0,0)[cc]{$\kappa_1$}}
\put(79.61,43.72){\makebox(0,0)[cc]{$\kappa_2$}}
\put(97.28,82.72){\makebox(0,0)[cc]{$\tau$}}
\put(117.61,82.06){\makebox(0,0)[cc]{$\tau$}}
\put(127.61,47.06){\makebox(0,0)[cc]{$\tau$}}
\put(82.61,37.06){\makebox(0,0)[cc]{$\tau$}}
\put(72.28,72.39){\makebox(0,0)[cc]{$\tau$}}
\end{picture}

\begin{center}
\nopagebreak[4]
Fig. \theppp.
\end{center}
\addtocounter{ppp}{1}
\vskip 0.2 in

This $\Theta$-annulus is composed of cells corresponding to relations
involving $\tau$. If we remove the $4N$ corners ($\kappa$-cells) then the
annulus breaks into $4N$ equal diagrams. Each of them is obtained in the
following way. Consider a path $p$ labelled by $W_1$. Let
$\tau=[U_1\to V_1, U_2\to V_2,...]$. Since $\tau$ is applicable to $W_1$,
this word contains all $U_i$'s. Attach transition cells corresponding
to the relations $U_i^\tau=V_i$ to the subpaths of $p$ labelled by $U_i$.
Attach auxiliary $\tau$-cells to the remaining edges of $p$. Finally
fold all pairs of edges with the same label and the same initial vertex.
It is easy to see that the resulting diagram is a $\Theta$-band whose common
edges are labelled by $\tau$, the top path is labelled by $W_1$ and the bottom
path is labelled by $W_2$.

Consider now an arbitrary computation $C$ = ($W_1, W_2, \ldots, W_g$)
with a history word $h$. Suppose that $C\sss(W_g,W_0)$ is not empty.
With every $i=1,\ldots,g-1$
we associate the annular diagram with boundaries $K(W_i)$ and
$K(W_{i+1})$ as in Fig. \thepdfive.
Since for every $i=1,\ldots, g-1$ the inner contour
of the $i$-th annular diagram is labelled by the same word as the outer contour of
the $(i+1)$-st diagram, we can ``concatenate" all these diagrams by gluing the next
diagram in the hole of the previous diagram. As a result
we obtain an annular diagram with outer contour labelled by
$K(W_1)$ and the
inner contour labelled by $K(W_g)$.
This diagram is called a {\em computational annulus} corresponding
to the computation $W_1,...,W_g$.
This annular diagram is a union of
$4N$ sectors and inverses of sectors which will be called
{\em computational
sectors corresponding to the computation $C$}. Notice that the history word
$h$ of the computation $C$ is equal to the label of the top (and bottom) path
of any maximal $\kappa$-band starting on the outer contour of the annulus.

If $C$ is an accepting computation then the inner boundary of the annular
diagram corresponding to $C$ is the boundary of the hub. Thus we can glue in the
hub and obtain an (ordinary) van Kampen diagram with boundary label $K(W_1)$.
We call it
the {\em computational disk corresponding
to
the computation} $C$.

It is easy to observe that the area of the computational
sector corresponding to
the computation $C$ is $O(\area(C))$ where $\area(C)$ is
the area of the computation $C$.  The area of a computational disc corresponding to
$C$ is thus also equal to $O(\area(C))$. We shall give an ``abstract"
definition of computational discs in the next section.

In this section we shall prove that every sector can be transformed
into a computational sector without changing the boundary label or increasing
the area.

Consider an arbitrary sector $\Delta$ with the inside diagram
$\Delta'$, with boundary  $$\partial(\Delta')=\pleft\ptop\pright\iv\pbot\iv,$$ with
an odd
$\kappa$-band
$\kkk$ and an even $\kappa$-band $\ttt$.
We assume that $\Delta$ has the smallest area among all sectors
with the same boundary label.

Notice that $\pleft$ and $\pright$ consist of $\Theta$-edges.
The maximal $\Theta$-bands starting on $\pleft$ must end on $\pright$
since $\Delta$ has no $(k,\Theta)$-annuli. For the same
reason, maximal $\Theta$-bands
starting on $\pright$ must end on $\pleft$. Therefore the maximal $\Theta$-bands
establish a one-to-one correspondence between edges on $\pleft$ and
edges on $\pright$. Since all $\Theta$-edges in any $\Theta$-band have
the same label, the labels of $\pleft$ and $\pright$ are the same.
Let us denote the maximal $\Theta$-bands in $\Delta$ starting on $\pleft$
by $\rr_1$,...,$\rr_\ell$ counting
from bottom, $\pbot$, up to $\ptop$.

For every $1<i\le \ell$ the path $\bott(\rr_i)\topp(\rr_{i-1})\iv$ bounds a subdiagram of
$\Delta$. The label of this path does not contain $\Theta$-edges. Therefore
by Lemma \ref{gammapath} this subdiagram contains no cells. Hence
$\bott(\rr_i)=
\topp(\rr_{i-1})$.

By the definition of a sector, the label of the bottom path of $\rr_1$, that is
the label of $\pbot$, is an admissible word $W$ such that $C\sss(W,W_0)$
is not empty. Therefore the label of $\bott(\rr_1)$ is a reduced word. The label
of $\topp(\rr_\ell)$ is a reduced word  by  the definition of a sector. We can
also assume that the labels of $\bott(\rr_i)$, $i=1,...,\ell$,
are reduced words. Indeed, if the
label of the path $\bott(\rr_i)$ is not reduced then it contains a subpath
$e_1e_2$ where $e_1$, $e_2$ are edges
labelled by $a$ and $a\iv$ respectively,
where $a$ is a letter.
\setcounter{pdtwo}{\value{ppp}}
These two edges must belong to two different cells
$\pi_1$ and $\pi_2$ of $\rr_i$ as in the diagram on the
left in the Figure \thepdtwo.

\unitlength=0.70mm
\special{em:linewidth 0.4pt}
\linethickness{0.4pt}
\begin{picture}(205.66,92.33)
\bezier{684}(0.33,20.67)(47.00,92.33)(92.66,20.67)
\bezier{528}(0.33,20.67)(80.00,-20.00)(92.66,20.67)
\put(37.33,23.33){\line(1,0){12.00}}
\put(49.33,23.33){\line(2,3){5.67}}
\bezier{160}(55.00,31.33)(38.00,42.67)(37.66,23.33)
\put(49.33,23.00){\line(5,-6){5.67}}
\bezier{160}(37.66,23.33)(44.33,2.00)(54.66,16.33)
\put(44.66,28.67){\makebox(0,0)[cc]{$\pi_1$}}
\put(44.66,18.00){\makebox(0,0)[cc]{$\pi_2$}}
\put(54.33,26.67){\makebox(0,0)[cc]{$a$}}
\put(54.33,21.33){\makebox(0,0)[cc]{$a$}}
\put(52.00,27.33){\vector(2,3){2.33}}
\put(52.00,20.00){\vector(1,-1){2.33}}
\put(97.33,20.67){\vector(1,0){10.00}}
\bezier{684}(113.33,20.34)(160.00,92.00)(205.66,20.34)
\bezier{528}(113.33,20.34)(193.00,-20.33)(205.66,20.34)
\put(150.33,23.00){\line(1,0){12.00}}
\bezier{160}(168.00,31.00)(151.00,42.34)(150.66,23.00)
\bezier{160}(150.66,23.00)(157.33,1.67)(167.66,16.00)
\put(157.66,28.34){\makebox(0,0)[cc]{$\pi_1$}}
\put(157.66,17.67){\makebox(0,0)[cc]{$\pi_2$}}
\put(162.00,23.00){\vector(1,0){11.67}}
\bezier{44}(174.67,23.00)(172.33,29.33)(168.00,31.00)
\bezier{40}(174.33,23.00)(167.33,16.67)(167.33,16.00)
\put(186.33,23.00){\vector(-1,0){12.67}}
\put(167.67,25.67){\makebox(0,0)[cc]{$a$}}
\put(179.33,25.67){\makebox(0,0)[cc]{$a$}}
\put(163.67,23.00){\circle*{1.33}}
\end{picture}

\begin{center}
\nopagebreak[4]
Figure \theppp.
\end{center}
\addtocounter{ppp}{1}

Then we fold $e_1$ and $e_2$, producing a new edge $e$ labelled
by $a$ and we introduce a new edge $f$ with label $a$ which has a common
end vertex with $e$ so that cells $\pi_1$ and $\pi_2$
have a common edge $e$. The other cells which were attached to $e_1$
and $e_2$ will be attached to the edge $f$.
This operation does not change the area or the boundary
of the diagram $\Delta$. The new diagram $\bar \Delta$ must
be reduced, otherwise after
reducing it, we shall get a sector with a smaller area and the same boundary
label as $\Delta$ (but we assumed that $\Delta$ has the smallest possible area).
Thus the diagram $\bar\Delta$ is again a sector
with the same number of maximal $\Theta$-bands starting on $\pleft$.
The bottom path of the $\Theta$-band $\rr_i$ in $\bar\Delta$ is shorter (by two edges)
than the bottom path of $\rr_i$ in $\Delta$. The bottom paths of the other
$\Theta$-bands are not affected by this operation.
Thus after
a finite number of such operation we get a sector
with the same boundary label as $\Delta$, in which all bands $\rr_i$
have top and bottom paths with reduced labels.

The contour of $\Delta'$ does not have $\kappa$-letters. Since $\Delta$
does not have $\kappa$-annuli (this diagram is reduced), $\Delta'$ does not
have $\kappa$-edges. In particular, it has no hubs.

For every component $Q_i$ of the vector of sets of state letters $Q$
there exists exactly one $Q_i$-edge on $\pbot$
(here we use the form of the word $W_0$).
Therefore $\Delta$ contains
exactly one maximal $Q_i$-band starting on $\pbot$ for every $i$.
These $Q_i$-bands do not intersect and cannot end on $\pbot$, $\pleft$
or $\pright$ ($\pleft$ and $\pright$ contain no $Q$-letters). They cannot
end on the contour of a hub because $\Delta$ has no hubs. Therefore
these $Q_i$-bands end on $\ptop$. Since $\Delta$ contains no
$(Q_j, \Theta)$-annuli (Lemma \ref{qr}), each $Q_i$-band starting on
$\pbot$ intersects each $\Theta$-band $\rr_j$ ($j=1,...,\ell$)
exactly once. Also this
implies that there are no $Q_j$-bands which start and end on $\ptop$.
Therefore for every component $Q_i$ of the vector $Q$,
$\Delta$ has exactly two mutually inverse maximal
$Q_i$-bands. These $Q_i$-bands connect $\ptop$ and $\pbot$ and intersect each
of the $\Theta$-bands $\rr_j$ exactly once ($j=1,...,\ell$).

Therefore each of the $\Theta$-bands $\rr_j$ contains exactly one $Q_i$-cell
for every component $Q_i$.

Consider one of the $\Theta$-bands $\rr=\rr_j$, $j=1,...,\ell$.
Let $W'$ be the label of the path
$\bott(\rr_j)$ and let $\tau\in \Theta$ be the label of the start edge of this
$\Theta$-band. Suppose that $W'$ is an admissible word.
Since $\rr$ starts on $\pleft$, the start edge of $\rr$ belongs to $\pleft$,
so it is oriented toward $\ptop$.

For convenience suppose that $\tau\in \Theta_+$ (the case when $\tau\iv \in\Theta_+$
is similar). Let $$\tau=[U_1\to V_1, U_2\to V_2,..., U_m\to V_m]$$
be the form of the $S$-rule $\tau$. By the definition of $S$-rules,
 $U_i$ starts and ends with state letters, and
 the state letters from different $U_i$'s belong to different
components of $Q$, so these state letters are different.

Let $q$ and $q'$ be the first and the last state letters of $U_i$, $q\in Q_a$,
$q'\in Q_b$ for some numbers $a$ and $b$, $a<b$. There exists only one (up to cyclic
shifts and inverses) relation in $\pp_N(\sss)$ containing a letter from
$Q_a$ and $\tau$. This is the relation $U_i^\tau V_i\iv=1$. The same relation is the
only relation containing a letter from $Q_b$. We have proved that
the $Q_a$-band and  the $Q_b$-band starting on $\pbot$ intersect with $\rr$. Therefore
their medians intersect with $\bott(\rr)$. Thus $\bott(\rr)$ contains
a $Q_a$-edge $e$ and a $Q_b$-edge $e'$. These edges are oriented
from $\pleft$ to $\pright$, that is both of these edges belong
to the path $\pbot(\rr)$.
The labels of
$e$ and $e'$ are $q$ and $q'$, respectively because these are the only $Q_a$- and
$Q_b$-letters appearing in relations of $\pp_N(\sss)$ together with $\tau$.
Let $\pi$ and $\pi'$
be cells in $\rr$ containing $e$ and $e'$. The boundary of $\pi$ must have the
form $e\iv tw'(t')\iv f\iv w\iv$ where $\Lab(t)=\Lab(t')=\tau$,
$\Lab(w')=V_i$, $\Lab(f)=q'$, $\Lab(ewf)=U_i$. If $f$ is not equal to $e'$
then $f$ must belong to the contour of the a cell $\pi''$ of $\rr$
which goes after
$\pi$.  Then $f$,  $e'$, the second $Q_b$-edge of $\pi$ and
the second $Q_b$-edge
of $\pi''$  belong to the same $Q_b$-band (as we know
there are only two mutually inverse $Q_b$-bands in $\Delta$).
Therefore this $Q_b$-band
intersects $\rr$ twice which is impossible by Lemma \ref{qr}. Therefore
$f=e'$. Thus $W'$ contains subword $U_i$ for every $i=1,...,m$,
$\rr$ has cells with boundary words $\tau V_i\tau\iv U_i\iv$ for every
$i=1,...,m$, and the subpaths of these boundaries labelled by
$U_i\iv$ are subpaths of $\pbot\iv$. Notice that the subpaths labelled by $V_i$
are the top paths of these cells considered as $\Theta$-bands.

If $W'$ has a $Q_j$-edge $e$ which does not belong to any of the $U_i$ then
the cells involving $e$ must correspond to the commutativity transition relation
$q^\tau=q$ where $q=\Lab(e)$.

Removing from $\rr$ all these transition cells which were discussed in the last two paragraphs,
we obtain several $\Theta$-bands
without $Q$-edges on their contours.
By Lemma \ref{q} these $\Theta$-bands
must consist of auxiliary cells. It is clear that such a $\Theta$-band is
uniquely determined by the label of its bottom path (which is equal to the label
of its top path because the top and the bottom paths of an auxiliary cell
considered as a $\Theta$-band, coincide).

Summarizing the information that we have obtained, we can conclude that
$W'$ has the form $$g_1P_1g_2P_2...g_cP_cg_{c+1}$$ where
$g_i$ are words without state letters, and either $P_i=U_j$ or $P_i\in Q_j$
such that none of the $U$'s contain a $Q_j$-letter. The label
of the top path of $\rr$ is obtained by reducing the word
$$g_1R_1g_2R_2...g_cR_cg_{c+1}$$ where $R_i=V_i$ if $P_i=U_i$ and
$R_i=P_i$ if $P_i$ is none of the $U$'s.
Thus
the label of $\topp(\rr)$ is obtained from $W'$ by applying the $S$-rule
$\tau$.

Notice that if the start and end edges of $\rr$
are labelled by $\tau\iv$ where
$\tau\in\Theta_+$ then the same argument shows that $\Lab(\topp(\rr))$ is
obtained by applying the rule $\tau\iv$ to $W'$.

We have proved that the sequence $C$ of
words $$\Lab(\pbot), \Lab(\topp(\rr_1)),\Lab(\topp(\rr_2)),...,
\Lab(\topp(\rr_\ell))$$  is a computation of the $S$-machine $\sss$.

Notice also that our argument shows that there is a unique $\Theta$-band with
top path labelled by $\Lab(\topp(\rr))$, the bottom path labelled by $W'$
and the $\Theta$-edges labelled by $\tau$.  In order
to create this band one needs to
draw a path $p$ labelled by $W'$, glue in cells with boundary labels
$\tau V_1\tau\iv U_1\iv$, $\tau V_2\tau\iv U_2\iv$, ...,
$\tau V_m\tau\iv U_m\iv$ to the
parts of the path $p$ labelled by $U_1$, ..., $U_m$ respectively.
Then attach auxiliary
$\tau$-cells to the edges of $p$ which are
not in the parts of $p$ which we used
before. If after that we get two edges with the same
labels and the same initial vertex,
we fold them. After all these
foldings we get a $\Theta$-band with the top path labelled
by $\Lab(\topp(\rr)$ and the bottom path labelled by $W'$.

Therefore the sector $\Delta$ coincides with the computational sector
corresponding to the computation $C$. Thus we have proved the following result.

\begin{prop} Every sector can be transformed into a computational sector
without increasing the area or the boundary label. The area of a computational
sector is ``big O'' of the area of the corresponding computation
\label{prop12}\end{prop}

\section{Computational Discs}
\label{compdiscs}

Let $\Delta$ be a reduced diagram over $\pp_{N}(\sss)$ with reduced boundary
label and exactly one hub. Suppose further that $\Delta$ does not have
$\Theta$-edges on the contour. In this case we call $\Delta$ a {\em disc}. The 
goal of this section is to provide a description of discs.

\begin{lm} \label{introd} If $\Delta$ is a  disc then every maximal
$\kappa$-band in $\Delta$ contains the hub, and every maximal
$Q$-band starts or ends on the contour of the hub.
\end{lm}

{\bf Proof.} Suppose that there exists a maximal $\kappa$-band $\kkk$
which does not contain the hub. By Lemma \ref{cent}, $\kkk$ starts and
ends on the contour of $\Delta$. Cutting $\Delta$ along the path
$\topp(\kkk)$ we obtain two diagrams $\Delta_1$ and $\Delta_2$. Since
$\Delta$ has only one hub, one of these diagrams does not contain a hub.
Without loss of generality assume that $\Delta_1$ is hub-free. The path
$\topp(\kkk)$ consists of $\Theta$-edges since $\kkk$ does not contain hubs.
These $\Theta$-edges cannot cancel when we form $\topp(\kkk)$, otherwise
$\kkk$
would not be reduced.  Consider the maximal $\Theta$-bands in $\Delta_1$
starting on $\topp(\kkk)$.  Since $\Delta_1$ does not contain hubs, these
$\Theta$-bands cannot end on $\topp(\kkk)$ (Lemma \ref{cr}).  Therefore
they
end on $\partial(\Delta)$.  This contradicts the assumption that
$\partial(\Delta)$ does not contain $\Theta$-edges.

If $\Delta$ has a maximal $Q$-band $\qq$ which does not start or end on
the contour of the hub then $\qq$ starts and ends on the contour of
$\Delta$ (Lemma \ref{q}).  Cutting $\Delta$ along $\topp(\qq)$ we obtain
two diagrams $\Delta_1$ and $\Delta_2$, one of which, say $\Delta_1$, does
not contain hubs.  The $\Theta$-edges do not cancel when we form the path
$\topp(\qq)$ because otherwise two cells in $\qq$ form a reducible pair.
Therefore the path $\topp(\qq)$ contains $\Theta$-edges. By Lemma \ref{qr}
the
maximal $\Theta$-band starting on one of these edges must end on
$\partial(\Delta)$. This again contradicts the assumption that
$\partial(\Delta)$ does not have $\Theta$-edges. The lemma is proved.
$\Box$

\vskip 0.1 in

Suppose that $\Delta$ is a disc.
Since every cell in $\Delta$ (except the hub) contains a $\Theta$-edge,
$\Delta$ is covered by $\Theta$-bands (except for the hub).
Since there are no $\Theta$-edges on the contour of $\Delta$,
each maximal $\Theta$-band in $\Delta$ is an annulus. Every
$\Theta$-annulus
must contain a hub in its inside diagram (Lemma \ref{r}). Therefore
every $\Theta$-annulus in $\Delta$ goes around the hub, so the
$\Theta$-annuli
in $\Delta$ form concentric annuli surrounding
the hub. We shall consider only the
$\Theta$-annuli that go clockwise around the hub (other $\Theta$-annuli
are the inverses
of these $\Theta$-annuli). Let $\rr_o$ be the outermost
$\Theta$-annulus and $\rr_i$ be the innermost $\Theta$-annulus.
Every edge of the contour of $\Delta$ belongs to the contour of a
cell. Indeed, otherwise
by removing this edge one gets two disjoint diagrams one
of which does not contain
hubs. Since this hub-free diagram does not have $\Theta$-edges on its
contour,
it must be empty (Lemma \ref{gammapath}). Since the $\Theta$-bands in
$\Delta$
form concentric rings, all contour
edges must belong to contours of cells from $\topp(\rr_o)$.
It is clear that $\bott(\rr_i)\cap \bott(\rr_i)\iv$ is the contour
of the hub (there cannot be any cells between $\bott(\rr_i)$ and the contour
of the hub).

Let us remove the hub $\pi$ from $\Delta$ and consider the maximal
$\kappa$-bands in the resulting annular diagram $\Delta_1$ starting on
$\partial(\pi)$.  Let us enumerate these bands clockwise by
$\bb_1,...,\bb_{4N}$.  Since $\Delta$ contains no maximal $\kappa$-bands
which do not contain the hub (Lemma \ref{introd}), and since there are no
$\kappa$-annuli in $\Delta$ (Lemma \ref{cent}), $\Delta_1$ is covered by
$4N$ subdiagrams $\Sigma_1$,..., $\Sigma_{4N}$ such that $\Sigma_i$ is
bounded by $\topp(\bb_i)$, a part of $\partial(\Delta)$,
$\bott(\bb_{i+1})$ (addition modulo $4N$) and a part of $\partial(\pi)$.
By definition each of $\Sigma_i$ is either a sector or the inverse of a
sector, depending on whether the $\kappa$-band $\bb_i$ is even or odd;
$\rr_o\cap \Sigma_i$ is the top $\Theta$-band in this sector.
Using Proposition \ref{prop12}, we can transform each sector
$\Sigma_i^{\pm 1}$ into a computational sector of the same or smaller area
and  with the same boundary.
Thus we can assume that each $\Sigma_i$ is a computational sector.

The label of $\pbot$ in each of these sectors is $W_0$.
Let $w$ be the history word of the sector $\Sigma_1$.
This word is equal to $\Lab(\topp(\bb_1))$.
By the definition of computational sectors,
$$\Lab(\topp(\bb_i))=\Lab(\bott(\bb_{i+1})=\Lab(\topp(\bb_{i+1}))$$
for every $i=1,...,4N$. Therefore the history words of all these sectors are
the same. By the definition of a computational sector, computational
sectors with the same label of the bottom
paths and the same history words are homotopic.
This implies that if we let $\Lab(\topp(\rr_o\cap\Sigma_1))=u$,
then
$$\Lab(\partial(\Delta))=K(u).$$

It is clear that the area of the disc is ``big O" of the  area of the 
corresponding computation. Let us compute the diameter of the disc. 
Take any vertex $v$ inside the disc. It belongs to one of the $\Theta$-annuli
of the disc. Let is be the $\Theta$-band number $i$ counted from the hub to the 
boundary of the disc. It is clear that this 
vertex is within a constant distance (less than the maximal length of a relator 
of $G_N(\sss)$)  from
the $\Theta$-band $i-1$. Therefore $v$ is within distance $Ci$ from the boundary
of the disc. Since the number of $\Theta$-annuli in the disc is equal to the
length of the corresponding computation, the diameter of the disc is bounded
above by a constant multiple of this length. On the other hand pick a vertex 
on the hub and consider any path $p$ connecting this vertex with a vertex on 
the boundary of the disc. By Jordan's lemma, this path must cross medians of all
$\Theta$-annuli in the disc. Therefore the length of $p$ cannot be smaller 
than the length  of the computation corresponding to the disc.

This proves the following statement.

\begin{prop} \label{prop13}
Every disc can be transformed without changing the boundary and without
increasing the area into a computational disc. The area of this disc is
``big O" of the area of the corresponding computation. The diameter of the
disc is ``big O" of the time of the corresponding computation.
\end{prop}

\section{The Upper Bound}
\label{upperbound}

In this section, we shall prove that for every word $w$ which is equal to 1
in $G_{N}(\sss)$ there exists a \vk diagram
$\Delta$ over $\pp_{N}(\sss)$ with $w=\Lab(\partial(\Delta))$ such that $\Delta$
can be decomposed in some standard way into
diagrams with at most one hub (like a snowman can be decomposed into small
snow balls). This decomposition will then help us find an upper bound for the
area and diameter of $\Delta$, as functions of $|w|$.
We assume that the function $T^4(|w|)$
(the area function of the machine $\sss=\sss(M)$ is superadditive.

Consider the following graph $G(\Delta)$ associated with the diagram $\Delta$.
The vertices of $\Delta$ are the hubs and the edges are maximal hub-free
$\kappa$-bands connecting the hubs. Since $\kappa$-bands cannot intersect except
in the hubs, $G(\Delta)$ is a plane graph (map).  Every internal vertex of this
graph has degree $4N$.  We use the theory of plane graphs (maps) similar to
small cancellation theory (see \cite{LS}). In order to do that we need the
following fact.

\begin{lm} \label{bigon} If $\Delta$ has the smallest possible number of hubs
among all reduced diagrams with the same boundary label, then $G(\Delta)$ does
not contain polygons with only one vertex and does not have bigons (polygons with only two
vertices).  \end{lm}

{\bf Proof.} The first part of the statement follows immediately from the fact
that $\Delta$ does not have $\kappa$-annuli.

Suppose that $\Delta$ has a bigon with two vertices $\pi_1$
and $\pi_2$ and two edges $\bb_1$ and $\bb_2$. We can assume that there
are no $\kappa$-bands inside the bigon between $\bb_1$ and $\bb_2$,
so (considered as edges of $G(\Delta)$), $\bb_1$ and $\bb_2$ are consecutive edges
of $\pi_1$.

This implies that one of the
bands $\bb_1$ or $\bb_2$ is a $\kappa_j$-band and another one is a $\kappa_{j+1}$-band.
Let $p_2$ (resp. $p_4$) be the shortest subpath of $\partial(\pi_1)$
(resp. $\partial(\pi_2)$) containing the start (resp. the end)
edges of $\bb_1$ and $\bb_2$. Then
$\Lab(p_2)=\Lab(p_4)\iv=(\kappa_jW_0 \kappa_{j+1})^{\pm 1}$.
Without loss of generality we can assume that $\Lab(p_2)=\Lab(p_4)\iv=
\kappa_jW_0 \kappa_{j+1}$. By renaming $\bb_1$ and $\bb_2$ and by
replacing $\Delta$ by $\Delta\iv$ if necessary we can also make
the path $\topp(\bb_1)p_2\bott(\bb_2)\iv p_4\iv$ bound a subdiagram $\Sigma$
containing the bands $\bb_1$ and $\bb_2$.
Every $\Theta$-band which crosses $\bb_1$ must cross $\bb_2$ (Lemma
\ref{cr}).
Therefore $p_2$ is the top path of a $\Theta$-band in $\Sigma$ (the one
which
is the closest to $\pi_1$).  Thus $\Sigma$ is a  sector.
By Proposition \ref{prop12}, we can assume that $\Sigma$ is a
computational sector.
Let $\partial(\pi_2)=p_4p_5$.

Let $w$ be the history word for $\Sigma$. Let $D$ be the
computational disc with a hub $\pi$ corresponding to the same history
word.
Then $D$ contains a copy $\Sigma_1$
of $\Sigma$. The boundary of the diagram
$D_1=D\backslash (\Sigma_1\cup \pi)$ has
the form $s_1s_2$ where $$\Lab(s_2)=\Lab(\topp(\bb_1)p_5\bott(\bb_2)\iv).$$
Thus we can cut $\Delta$ along the
path $\topp(\bb_1)p_5\bott(\bb_2)\iv$ and insert the diagram
$$D_1\circ_{s_1=s_1} D_1\iv$$ in the hole (\setcounter{pdthree}{\value{ppp}}see
Fig. \thepdthree).

\unitlength=1.00mm
\special{em:linewidth 0.4pt}
\linethickness{0.4pt}
\begin{picture}(153.94,105.72)
\put(21.28,32.39){\circle{10.85}}
\put(50.28,32.06){\circle{11.18}}
\bezier{116}(25.61,34.39)(37.61,45.39)(45.61,35.39)
\bezier{100}(26.61,33.39)(37.61,42.06)(44.61,34.06)
\bezier{92}(25.94,29.72)(34.61,23.06)(45.28,29.06)
\bezier{108}(24.28,27.72)(36.28,19.72)(46.28,27.72)
\put(112.94,50.39){\circle{10.85}}
\put(141.94,50.06){\circle{11.18}}
\bezier{116}(117.28,52.39)(129.28,63.39)(137.28,53.39)
\bezier{100}(118.28,51.39)(129.28,60.06)(136.28,52.06)
\bezier{92}(117.61,47.72)(126.28,41.06)(136.94,47.06)
\bezier{108}(115.94,45.72)(127.94,37.72)(137.94,45.72)
\bezier{196}(137.28,54.06)(136.61,79.06)(112.61,81.06)
\bezier{220}(112.61,81.06)(87.28,78.72)(87.61,49.39)
\bezier{192}(87.61,49.39)(88.28,21.72)(108.94,20.39)
\bezier{204}(108.94,20.39)(135.28,20.72)(137.94,45.39)
\bezier{224}(137.28,53.72)(153.94,86.39)(136.61,94.06)
\bezier{252}(136.61,94.06)(107.94,105.72)(81.28,88.06)
\bezier{180}(81.28,88.06)(70.94,80.06)(73.61,48.06)
\bezier{272}(73.61,48.06)(74.61,-0.61)(93.61,-2.61)
\bezier{400}(93.28,-2.61)(142.28,-5.28)(137.94,45.72)
\put(115.28,55.39){\line(4,5){27.33}}
\put(113.94,55.72){\line(4,5){28.00}}
\put(108.94,54.06){\line(0,1){44.33}}
\put(108.94,98.39){\line(-1,0){0.33}}
\put(110.94,55.06){\line(0,1){43.67}}
\put(107.61,51.06){\line(-1,0){34.33}}
\put(107.94,49.39){\line(-1,0){34.33}}
\put(61.94,31.72){\vector(1,0){6.33}}
\put(35.61,31.72){\makebox(0,0)[cc]{$\Sigma$}}
\put(127.94,49.39){\makebox(0,0)[cc]{$\Sigma$}}
\put(110.94,33.72){\makebox(0,0)[cc]{$D_1$}}
\put(94.61,9.72){\makebox(0,0)[cc]{$D_1\iv$}}
\put(112.94,49.72){\makebox(0,0)[cc]{$\pi_2$}}
\put(142.28,49.72){\makebox(0,0)[cc]{$\pi_1$}}
\put(50.28,31.72){\makebox(0,0)[cc]{$\pi_1$}}
\put(21.61,31.72){\makebox(0,0)[cc]{$\pi_2$}}
\end{picture}
\begin{center}
\nopagebreak[4]
Fig. \theppp.
\end{center}
\addtocounter{ppp}{1}

The resulting diagram will have the same boundary label as $\Delta$.
The maximal disc with the hub $\pi_2$ in this diagram corresponds to the
history word $w$. Thus the
boundary label of this disc will be $K(W_0)$,
the same as the boundary label of the hub $\pi_2$.
Therefore this disc without the hub is an annulus whose inner and outer
boundaries have the same labels. Thus we can remove this annulus and obtain
a diagram $\Delta'$ with the same boundary label and the same number of hubs
as $\Delta$ (see the proof of Lemma \ref{cent}). But in $\Delta'$, the
hubs $\pi_1$ and $\pi_2$ have a common path labelled by
$\kappa_jW_0 \kappa_s$.
By Lemma \ref{cancel} (d), these two cells cancel.
By cancelling these hubs we obtain a diagram with the same
boundary label and a smaller number
of hubs. This contradicts the assumption that $\Delta$ contains
the minimal possible
number of hubs. The lemma is proved.  $\Box$

\begin{lm} \label{triangle} $G(\Delta)$ does not have triangular faces.
\end{lm}

{\bf Proof.} Indeed, suppose that there are three hubs $\pi_1$, $\pi_2$,
$\pi_3$ connected by three bands $\bb_1, \bb_2, \bb_3$
such that $\bb_2$
and $\bb_3$ start on the contour of $\pi_1$ and end on the boundaries
of $\pi_3$ and $\pi_2$ respectively, $\bb_1$ connects boundaries of
$\pi_2$ and
$\pi_3$, none of the bands $\bb_1, \bb_2, \bb_3$ contains hubs, and there
are no hubs in the triangle bounded by these bands.
Since $\bb_2$ and $\bb_3$ are two
consecutive edges of $\pi_1$ in $G(\Delta)$,
the $\kappa$-edges in one of these bands have even index
and the $\kappa$-edges in the other band have odd index.
Let, say, $\bb_2$ be a $\kappa_i$-band with even $i$ and $\bb_3$ be a $\kappa_j$-band
with odd $j$.
The band $\bb_1$ connects $\pi_2$ and $\pi_3$. Since $\bb_1$ and
$\bb_3$ are consecutive edges of $\pi_2$, the index of $\kappa$-edges in
$\bb_1$ must be a even.
But a similar argument applied to the hub $\pi_3$ gives that
this index must be odd, a contradiction.  $\Box$

\vskip 0.1 in

Lemmas \ref{bigon} and \ref{triangle} shows that all faces
in the planar graph $G(\Delta)$ have degrees at least 4.

The following two lemmas are probably well known but we were not able to
find them in the literature.

\begin{lm} \label{lambda} Let $G$ be a non-empty plane graph in which every
face has degree at least 4. Then $G$ either consists of one vertex
or has two vertices of degree at most 3.
\end{lm}

{\bf Proof.} Suppose that $G$ contains more than one vertex.  We shall use the
notation from \cite{LS} (Chapter 5).  Let $V, E, F$ be the numbers of vertices,
edges and faces in $G$.  Let $Q$ be the number of connected components in $G$.
$\Sigma$ refers to the summation over all vertices $v$ or faces $D$ in $G$,
$E^{\cdot}$ denotes the number of all boundary edges of $G$ (the length of the
boundary of $G$). By $d(v)$ and $d(D)$ we denote the degree of a vertex $v$ or a
face $D$.  Then by Theorem 3.1 in Chapter 5 of \cite{LS} the following formula
holds (we take $p=q=4, h=0$ in the formula in this theorem):

\begin{equation}
4Q=\Sigma[4-d(v)]+\Sigma[4-d(D)]-E^\cdot.
\label{eq}
\end{equation}

But by assumption, $d(D)\geq 4$ for every face $D$.
Therefore $4Q < \Sigma[4-d(v)]$. The left hand side is at least $4$.
Hence either there exists a vertex with $d(v)=0$, or there are at least
two vertices of degree at most $3$. In the first case the vertex of degree
0 forms a connected component of the graph. Since $G$ has more than one vertex,
$Q\ge 2$. Then the left hand side of (\ref{eq}) is at least 8. This implies that
the right hand side contains at least two positive terms, that is $G$
contains two vertices of degree at most 3.  $\Box$

\vskip 0.1 in

Now we are in a position to define the snowman decomposition of $\Delta$.  For
the rest of this section we fix a large enough number $N$. Later
we'll show that we can take $N\ge 6$
where $k$ is the number of tapes of the original Turing machine.
Recall
that the number of sectors in every
computational disc is $4N$.

With every diagram $\Delta$ we associate two  numbers:

\begin{itemize}
\item $n(\Delta)=|\partial(\Delta)|$ --- the length of the boundary of $\Delta$;
\item $h(\Delta)$ --- the number of hubs in $\Delta$.
\end{itemize}

Let $\kkk$ be a $\kappa$-band in $\Delta$. Let $\Delta', \Psi$ be the
subdiagrams obtained by dividing $\Delta$ along $\topp(\kkk)$ such that $\Psi$
contains $\kkk$.  If $\Psi$ does not contain hubs then we call $\kkk$ a {\em
dividing $\kappa$-band}.  In particular a dividing $\kappa$-band is hub-free.

Suppose that $\Delta$ contains a dividing $\kappa$-band $\kkk$.  Clearly both
$\Delta'$ and $\Psi$ are reduced, $h(\Delta)=h(\Delta')$.  Suppose that
$\partial(\Psi)=\topp(\kkk)p_1$, $\partial(\Delta')=p_2\topp(\kkk)\iv$.  Since
$\kkk$ does not contain hubs, the word $\Lab(\topp(\kkk))$ consists of
$\Theta$-letters. The maximal $\Theta$-bands in $\Psi$ which start on $\topp(\kkk)$ cannot
end on $\topp(\kkk)$ because of Lemma \ref{cr}. Therefore all these $\Theta$-bands
end on $p_1$
\setcounter{pdtwelve}{\value{ppp}}
(see Fig. \thepdtwelve).
Since $p_1$ contains at least 2 $\kappa$-edges, we have $|p_1|\ge
|\topp(\kkk)|+2$. This implies that $|\partial(\Psi)|\le 2|\partial(\Delta)|$
and $|\partial(\Delta')|\le |\partial(\Delta)|-2$.

Thus we have proved
the following properties of the pair $(\Psi, \Delta')$.

\begin{enumerate}
\item[(D1)] $\Psi$ and $\Delta'$ are reduced.
\item[(D2)] By gluing $\Psi$ and $\Delta'$ we obtain a diagram with the
same boundary label as $\Delta$.
\item[(D3)] $\Psi$ does not contain hubs.
\item[(D4)] $\Delta'$ contains a smaller number of dividing $\kappa$-bands
than $\Delta$.
\item[(D5)] $n(\Psi)\le 2n(\Delta)$,
$n(\Delta')\le n(\Delta)-2$,
$h(\Delta')=h(\Delta)$.
\end{enumerate}

\unitlength=1.00mm
\special{em:linewidth 0.4pt}
\linethickness{0.4pt}
\begin{picture}(119.06,41.94)
\bezier{244}(57.39,38.28)(76.72,14.28)(95.06,38.28)
\bezier{352}(53.06,38.28)(73.72,2.28)(103.06,38.28)
\put(68.72,21.28){\line(2,3){3.67}}
\put(72.39,26.94){\line(0,1){11.00}}
\put(72.39,37.94){\line(0,1){0.33}}
\put(77.06,20.28){\line(0,1){18.00}}
\put(85.06,22.94){\line(-3,5){2.67}}
\put(82.39,27.28){\line(1,6){1.67}}
\put(92.06,27.28){\line(-4,3){4.33}}
\put(87.72,30.61){\line(1,6){1.33}}
\put(63.39,25.61){\line(1,1){3.33}}
\put(66.72,28.94){\line(0,1){9.33}}
\put(59.39,29.28){\line(3,2){4.00}}
\put(63.39,31.94){\line(0,1){6.33}}
\put(56.06,41.94){\makebox(0,0)[cc]{$\kappa$}}
\put(99.39,41.94){\makebox(0,0)[cc]{$\kappa$}}
\put(73.06,14.61){\makebox(0,0)[cc]{$w(\Theta)$}}
\put(80.39,20.61){\oval(77.33,35.33)[]}
\end{picture}

\begin{center}
\nopagebreak[4]
Figure \theppp.

\end{center}
\addtocounter{ppp}{1}
\vskip 0.2 in

In Figure \thepdtwelve\, $w(\Theta)$ denotes any group word over the alphabet $\Theta$.

Since $\Psi$ has no hubs, its area can be computed with the help of
Lemma \ref{nohub}. Thus in order to estimate
the area of $\Delta$ we can replace $\Delta$ by $\Delta'$.
We shall call the pair $(\Psi, \Delta')$ the {\em decomposition
of type 1 of $\Delta$}.

\bigskip

Suppose that $\Delta$ contains no dividing $\kappa$-bands.
Then by Lemma \ref{lambda} there
exists a vertex $\pi$ in $G(\Delta)$ of degree at most $3$.
Let $\bb_1,\ldots,\bb_{4N}$ be the maximal $\kappa$-bands starting on the
boundary of $\pi$, counted clockwise.
Since the degree of $\pi$ is at most 3, at least
${4N-3}$ of these bands do not contain hubs
(so they are not edges in $G(\Delta)$).

\begin{lm}
There exist at least two hubs $\pi$  and $\pi'$ in  $G(\Delta)$
of degree at most $3$ in $G($\dti$)$ such that
${4N-3}$ {\em consecutive} $\kappa$-bands starting on $\partial(\pi)$
(resp. $\partial(\pi')$)
do not contain hubs.
\label{consec}
\end{lm}

{\bf Proof. } Induction on the number of hubs.
By Lemma \ref{lambda}, $G($\dti$)$ has at least two vertices $\pi$ and $\pi'$
of degree at most 3.
Suppose that $\pi$ does not have ${4N-3}$ consecutive $\kappa$-bands without hubs.
Then there exist three maximal $\kappa$-bands
$\bb_i$, $\bb_j$, $\bb_m$, $i<j-1$, $j<m-1$ starting on
$\partial(\pi)$ such that $\bb_i$, $\bb_j$ and $\bb_m$ do not have hubs,
but at least one band among $\bb_{i+1},\ldots,\bb_{j-1}$ and at least one
band among $\bb_{j+1},\ldots,\bb_{m-1}$ contains a hub.

Consider the subdiagram $\Psi$
of $\Delta$ containing $\bb_{i},\ldots,\bb_{j-1}$ and $\pi$,
and bounded by $\topp(\bb_i)$, $\bott(\bb_{j-1})$, a part of
$\partial(\pi)$ and a part of $\partial($\dti$)$. The subdiagram $\Psi$
contains fewer hubs than $\Delta$ because among the
$\kappa$-bands $\bb_{j+1},\ldots,\bb_{m-1}$ there is at least one hub, which
is not contained in $\Psi$.
By the induction hypothesis, $\Psi$ contains two hubs $\pi_1$ and $\pi_2$
of degree at most 3 such that there are ${4N-3}$ consecutive maximal $\kappa$-bands starting
on $\pi_1$ (resp. $\pi_2$) which do not contain hubs. One of these
hubs, say, $\pi_1$ differs from $\pi$. It is easy to see that $\pi_1$
considered as a vertex in $G($\dti$)$ also has degree at most 3 and ${4N-3}$
consecutive $\kappa$-bands without hubs.

Now consider the subdiagram $\Psi'$ of $\Delta$ containing $\bb_{j+1},\ldots,\bb_{m-1}$
and $\pi$
bounded by $\topp(\bb_j)$, $\bott(\bb_m)$, a part of $\partial(\pi)$ and
a part of $\partial($\dti$)$. The same argument as in the case of the subdiagram
$\Psi$ shows that $\Psi'$ contains a hub $\pi_1'$
of degree at most 3 which has {4N-3} consecutive $\kappa$-bands without hubs.
Thus we found two hubs $\pi_1$ and $\pi_1'$ which satisfy the conditions of the lemma.
$\Box$

\vskip 0.1 in

\begin{lm} For every hub $\pi$ in $\Delta$ with hub-free consecutive
maximal
$\kappa$-bands
$\bb_1$,..., $\bb_{4N-3}$
starting on $\partial(\pi)$, let $\Psi_\Delta(\pi)$
be the subdiagram of $\Delta$ bounded by $\topp(\bb_1)$, $\bott(\bb_{{4N-3}})$,
$\partial(\Delta)$
and $\partial(\pi)$, which contains $\bb_1$, $\bb_{{4N-3}}$ and does not contain
$\pi$ (there is only one subdiagram in $\Delta$ satisfying these conditions).
Then there exists a hub $\pi$ such that $\Psi_\Delta(\pi)$
does not contain hubs
\setcounter{psi}{\value{ppp}}
(see Fig. \thepsi).
\label{psi}
\end{lm}

\vskip 1 in

\begin{center}
Fig. \thepsi.
\addtocounter{ppp}{1}

\end{center}

{\bf Proof.} Induction on the number of hubs in $\Delta$.
By Lemma \ref{consec}, there exists a hub $\pi$ in $\Delta$
with hub-free consecutive maximal $\kappa$-bands $\bb_1,\ldots,\bb_{{4N-3}}$ starting on
$\partial(\pi)$. Suppose that $\Delta'=\Psi_\Delta(\pi)$ contains hubs.
Notice that $\Delta'$
contains fewer hubs than $\Delta$ since $\Delta'$ does
not contain $\pi$.
By the induction hypothesis $\Delta'$ contains a hub $\pi'$ such that
$\Psi_{\Delta'}(\pi')$ does not contain hubs. Notice that none
of the $\kappa$-bands in $\Delta$
starting on the contour of $\pi'$ can intersect any of the $\kappa$-bands $\bb_1$,\ldots
$\bb_{{4N-3}}$ because these $\kappa$-bands do not contain hubs. Therefore
$\Psi_{\Delta'}(\pi')=\Psi_\Delta(\pi')$, so $\pi'$ is a hub for which
$\Psi_{\Delta}(\pi')$ is hub-free. The lemma is proved. $\Box$

\vskip 0.1 in

By Lemma \ref{psi} we can find a hub $\pi$ and ${4N-3}$ consecutive
$\kappa$-bands starting in $\partial(\pi)$ which do not have hubs
and moreover $\Psi_\Delta(\pi)$ is hub-free. We can assume that
these $\kappa$-bands have indices $1,..., 2N, 1, ...,2N-3$ (otherwise we can rename
$\kappa$'s).
We denote these ${4N-3}$
$\kappa$-bands by $\bb'_1,\ldots,\bb'_{{4N-3}}$ (as usual
we enumerate the $\kappa$-bands clockwise).

Let $\Pi$ be the maximal disc in $\Delta$ with the hub $\pi$.
By Proposition \ref{prop13}, we can assume that $\Pi$ is a computational disc.
We can also assume that $\Pi$ corresponds to a minimal (with respect
to area)
accepting
computation. By removing the interior of $\Pi$
from $\Delta$ we get an annular diagram $\Psi$.
Denote the parts of $\bb'_i$ ($i=1,...,4N$) in $\Psi$ by $\bb_i$.
Cutting $\Psi$ along
$\topp(\bb_1)$ and $\bott(\bb_{4N-3})$ we obtain two
ordinary \vk diagrams $\Psi_1$ and $\Psi_2$ such that
$\bb_1, \bb_2,\ldots,\bb_{4N-3}\subset \Psi_1$, $\bb_{4N-2},
\bb_{4N-1}, \bb_{4N}\subset \Psi_2$. Notice that $\Psi_1=\Psi_\Delta(\pi) - \Pi$.

\setcounter{triple}{\value{ppp}}
\addtocounter{ppp}{1}
\setcounter{sigma}{\value{ppp}}

We shall fix the following notation associated with the
triple $(\Psi_1, \Psi_2, \Pi)$  (see Figures \thetriple\ and \thesigma):
\begin{itemize}
\item $\partial(\Psi_1)=\topp(\bb_1)p_1\bott(\bb_{4N-3})\iv s_1 $;
\item $\partial(\Psi_2)=\topp(\bb_1)\iv s_2\iv \bott(\bb_{4N-3})p_2$;
\item $\partial(\Pi)=s_1\iv s_2$;
\item $p_1=e(\bb_1)t_1e(\bb_2)t_2\ldots t_{4N-4}e(\bb_{4N-3})$
where $e(\bb_i)$ is the end edge of $\bb_i$, $t_j$ is a subpath of $p_1$,
$i=1,\ldots,4N-3, j=1,\ldots,4N-4$;
\item $$s_1\iv =i(\bb_1)y_1i(\bb_2)y_2\ldots y_{4N-4}i(\bb_{4N-3}),$$
$$s_2=y_{4N-3}i(\bb_{4N-2})y_{4N-2}i(\bb_{4N-1})y_{4N-1}i(\bb_{4N})y_{4N}$$
where $i(\bb_j)$ is the start edge of $\bb_j$, for $j=1,...,4N-4$,
$y_j$ is a subpath of $s_1\iv$, for $j=4N-3,...,4N$, $y_i$ is a subpath of $s_2$;
\item $\Sigma_i$ ($i=1,\ldots, 4N-4$) is a subdiagram of $\Psi_1$
bounded by $\topp(\bb_i)$ on the left, $\bott(\bb_{i+1})$ on the right,
$p_1$ on the top and $s_1\iv$ on the bottom;
\item if $\Psi_1$ contains a $\Theta$-band $\rr$
which crosses $\bb_i$ and $\bb_{i+1}$, and such that the intersection
with $\bb_i$ precedes in $\rr$ the intersection with $\bb_{i+1}$ then
$t(\rr,i)$ is the portion of $\topp(\rr)$ which is contained in
$\Sigma_i$;  $t(\rr,i)= f(\rr, i)y(\rr,i)f(\rr,i+1)$ where
$f(\rr,i)$ and $f(\rr,i+1)$ are $\kappa$-edges;
\item $b=|\topp(\bb_1)|+|\topp(\bb_2)|+\ldots+|\topp(\bb_{4N-3})|$;
\item $c=|y_1|$.
\end{itemize}

\vskip 1 in
\begin{center}
Fig. \thetriple.
\end{center}

\vskip 1 in
\begin{center}
Fig. \thesigma.
\end{center}
\addtocounter{ppp}{1}

When we need to specify the triple ${\cal T}=(\Psi_1, \Psi_2, \Pi)$ to which
this notation is associated, we shall write $\bb_1({\cal T})$ instead of
$\bb_1$, etc.

The triple $(\Psi_1, \Psi_2, \Pi)$ satisfies the following properties.

\begin{enumerate}
\item[(P1)] By gluing together $\Psi_1$, $\Psi_2$ and $\Pi$ one can get a
diagram with the same boundary label as $\Delta$.
\item[(P2)] $\Psi_1$ is reduced and contains no hubs.
$\Pi$ is a computational disc  corresponding
to a minimal area accepting computation.
\item[(P3)] The $\kappa$-bands $\bb_i$, $i=1,\ldots,4N-3$, in $\Psi_1$
start on $s_1\iv$, and end on $p_1$. Every $\kappa$-band in $\Psi_1$
starts or ends on $s_1$.
\item[(P4)] The word $\Lab(f_1y_1\ldots f_{4N}y_{4N})$
is a cyclic shift of the word $K(W)$
for $W=\Lab(y_1)$. For all $i=1,\ldots, 4N$, $c=|y_i|$.
\item[(P5)] $h(\Psi_2)+1=h(\Delta)$.
\end{enumerate}

We assume from now on that we have a triple of diagrams $(\Psi_1, \Psi_2, \Pi)$
which satisfies the properties (P1)--(P5). We shall transform this triple into
another triple which satisfies the same properties plus some additional
conditions.

Suppose that there exists a $\Theta$-band in $\Psi_1$
which crosses the bands $\bb_1,\ldots,\bb_{4N-3}$. Then let $\rr$
be the $\Theta$-band
with this property whose intersection with $\bb_1$ is as far away from
$s_1$ as possible (the distance is counted along $\bb_1$).
We assume that $\rr$ starts on $\topp(\bb_1)$.
Consider the diagram $E$ bounded by $\topp(\bb_1)$, $\topp(\rr)$,
$\bott(\bb_{4N-3})$ and $s_1$. Then
$E'=E\circ_{s_1=s_1} \cup_{i=1}^{4N-4}\Sigma_i'$
is a union
of sectors and inverses of  sectors which correspond
to the same history word $w'$. Let us remove this subdiagram
from $\Psi_1$,
and consider the resulting diagram as a new $\Psi_1$. Then extend the
disc $\Pi$ by an annulus with the history word $w'$ and replace this
disc by a computational disc with the same boundary label which corresponds to
a computation of $\sss$ with the smallest possible area. Let it be
the new $\Pi$. After that, glue $\Psi_2$ to the union of the four diagrams
(sectors and inverse sectors
with the history word $w'$)
$\Sigma_{i}^{-1}$, $i = 4N-3, 4N-2, 4N-1, 4N$, then reduce the resulting
diagram. Let this diagram
be the new $\Psi_2$  \setcounter{move}{\value{ppp}} (see Fig. \themove).
The new triple $(\Psi_1, \Psi_2, \Pi)$ obviously satisfies properties
(P1)-(P5).
We shall call this
operation {\em moving $\Theta$-bands}. In addition, this triple satisfies the
following property.

\begin{enumerate}
\item[(P6)] There are no $\Theta$-bands in $\Psi_1$ which cross all the bands
$\bb_1,\ldots,\bb_{4N-3}$.
\end{enumerate}

\vskip 1 in
\begin{center}

Fig. \themove.
\addtocounter{ppp}{1}
\end{center}

\vskip 0.1 in
Now suppose that we have a triple ${\cal T}=(\Psi_1, \Psi_2, \Pi)$ which satisfies
conditions (P1)--(P6).

Take any number $i=1,\ldots,4N-4$. Suppose that the following property {\em
does not}
hold:

\begin{enumerate}
\item[(P7)]
There is no $\Theta$-band $\rr$ in
$\Sigma_i$ such that $|y(\rr,i)| < |y_i|$.
\end{enumerate}

\noindent Then take
a $\Theta$-band $\rr$ such that $|y(\rr,i)|$ \ (which is less than $|y_i|$) is the smallest possible
(over all $i$ and all  $\Theta$-bands $\rr$).

Consider the  subsector $T$ of $\Sigma_i$ bounded by $\topp(\bb_i)$,
$t(\rr,i)$, $\bott(\bb_{i+1})$, $h_iy_ih_{i+1}$ (with
$\rr$ as the top $\Theta$-band). Let $w'$ be the history word
of this sector. We can assume
by Proposition \ref{prop12}
that $T$ is a computational
sector. Consider the  computational
annulus $\bar E$ corresponding to the history word
$w'$ and initial configuration $\Lab(y(\rr,i))$. Then the disc $\Pi$ fits
in the hole of this annulus. By gluing $\Pi$ inside $\bar E$, and
replacing this computational
disc by a disc corresponding to a minimal area computation, we get a
computational disc which we denote by $\Pi_1$  \setcounter{adjust}{\value{ppp}}
(see Fig. \theadjust).

\vskip 1 in
\begin{center}
Fig. \theadjust.
\addtocounter{ppp}{1}
\end{center}

Let $E_1$ be the union of $4N-4$ consecutive subsectors of $\bar E$ such that the
part $p$ of the inner boundary of $\bar E$ contained in $E_1$ has the same label as
$s_1$. Let $E_2=\bar E\backslash E_1$. Let us form a diagram $\Psi_{1,1}$ by the
following procedure. First ({\em step 1}) we glue $\Psi_1$ and $E_1\iv$ along
$p=s_1$. Notice that the $k$-bands of $E_1\iv$ will be glued to the $\kappa$-bands of
$\Psi_1$ along the edges $f_i$, so the $\kappa$-bands are getting longer. Some
$\kappa$-cells in $\Psi_1$ will cancel with $\kappa$-cell in $E_1\iv$.  In {\em step 2} we
cancel these pairs of $\kappa$-cells.  Notice that if $j=1,\ldots,4N-3$ and $w'$ has
a common suffix with $\Lab(\topp(\bb_j)))$ of length $d_j$ then $d_j$ cells of
$\bb_j$ will cancel. In particular, since $w'$ is a suffix of
$\Lab(\topp(\bb_i))$, $|w'|$ cells in $\bb_i$ cancel. Then in {\em step 3} we
reduce other reducible pairs of cells.  The resulting (reduced) diagram is
$\Psi_{1,1}$.  Finally let us glue $\Psi_2$ and $E_2\iv$, reduce this diagram
and denote the resulting diagram by $\Psi_{2,1}$. We call the process of
obtaining the triple $(\Psi_{1,1}, \Psi_{2,1}, \Pi_1)$ an {\em adjustment}.

\begin{lm}
The triple ${\cal T}_1=(\Psi_{1,1}, \Psi_{2,1}, \Pi_1)$ satisfies the
properties (P1)--(P5). For every $j=1,\ldots,{4N-3}$,
$\Lab(\topp(\bb_j({\cal T}_1)))=
(w')\iv\Lab(\topp(\bb_j({\cal T})))$ (equality in the free group).
\label{tau1}
\end{lm}

{\bf Proof.} Properties (P1) and (P2) are obvious.

After step 1 of building $\Psi_{1,1}$
the $\kappa$-bands of the resulting
diagram will be just unions of corresponding
pairs of $\kappa$-bands in $\Psi_1$ and $E_1\iv$, so they will start on
$s_1({\cal T}_1)\subseteq \partial(E_1\iv)$ and end on $p_1({\cal T}_1)=p_1({\cal T})$.
This property will still hold after step 2 because in this step we are
cancelling cells from the same $\kappa$-bands (notice that $s_1({\cal T}_1)$
and $p_1({\cal T}_1)$ do not change after steps 2 and 3 because reducing
a diagram does not affect its boundary).
By Lemma \ref{kbands}, no $\kappa$-annuli
are produced in step 3. So step 3 is
just the process of cancelling reducible pairs of cells.
By Lemma \ref{kbands}, no $\kappa$-cells from the same
$k$-band cancel during this
process. Cells from different $\kappa$-bands
$\bb_{j_1}$ and $\bb_{j_2}$ do not form
a reducible pair because these cells do not have common edges.
Let $\Psi$ be the diagram obtained after a number of reduction from step 3.
Then for every $j=1,\ldots, 4N-4$ there exists
an ${\bf E}_0$-band which starts on $s_1({\cal T}_1)$ between $\bb_j$ and
$\bb_{j+1}$ and ends on $p_1({\cal T}_1)$ between the same $\kappa$-bands.
Therefore no cell in $\bb_j$ can have a common edge with a cell in
$\bb_{j+1}$. Since $\Psi$ does not have any other $\kappa$-bands,
no $\kappa$-cells
participate in the reduction process of step 3. Therefore the $\kappa$-bands
$\bb_j$, $j=1,\ldots,{4N-3}$, do not change during step 3. This proves property
(P3) and the fact that in the free group,
$\Lab(\topp(\bb_j({\cal T}_1)))=
(w')\iv\Lab(\topp(\bb_j({\cal T})))$.

Properties (P4) and (P5) immediately follow from the definition of the triple
${\cal T}_1$. $\Box$

\vskip 0.1 in

Since ${\cal T}_1$ satisfies properties (P1)--(P5) we can use all
the notation associated with triples satisfying these properties.
Notice that by construction ${\cal T}_1$ satisfies the property
$c({\cal T}_1) < c({\cal T})$. The triple ${\cal T}_1$ may not satisfy property
(P6). We shall discuss this property in the next lemma. Consider the
following property:

\begin{enumerate}
\item[(P6')] The words $\Lab(\topp(\bb_j))$,
$j=1,\ldots, 4N-3$, do not all have the same first letter.
\end{enumerate}

The following lemma shows, in particular, that (P6') is stronger than
(P6).

\begin{lm}
\label{tau}
a) If ${\cal T}_1$ does not satisfy property (P6')
then $w'$ is a prefix of the label of
$\topp(\bb_j({\cal T}))$ for every $j=1, \ldots, 4N-3$.
In this case $$b({\cal T}_1)<b({\cal T}).$$

b) If ${\cal T}_1$ satisfies property (P6') then
${\cal T}_1$ satisfies property (P6) (that is (P6') implies (P6)).

c) If ${\cal T}$ satisfies (P6') then so does ${\cal T}_1$.
\end{lm}

{\bf Proof.} a) If $w'$ is not a prefix of $\Lab(\topp(\bb_j({\cal T})))$ for
some $j$ then by Lemma \ref{tau1} the first letter of $\Lab(\topp(\bb_j({\cal
T}_1)))$ is the first letter of $(w')\iv$. On the other hand, by construction,
$w'$ is a prefix of $\Lab(\topp(\bb_i({\cal T})))$ (for the particular $i$
considered above). Therefore the first letter in $\Lab(\topp(\bb_i({\cal
T}_1)))$ is the letter number $|w'|+1$ in $\Lab(\topp(\bb_i({\cal T})))$. If
these two first letters in $\Lab(\topp(\bb_j({\cal T}_1)))$, resp.
$\Lab(\topp(\bb_i({\cal T}_1)))$ coincide, then $\Lab(\topp(\bb_i({\cal T})))$
contains a pair of consecutive mutually inverse letters; but this is impossible
since $\bb_i({\cal T})$ is reduced. Therefore if the first letter of each word
$\Lab(\topp(\bb_j({\cal T}_1)))$ is the same, then $w'$ is a prefix of
$\Lab(\topp(\bb_j({\cal T})))$ for $j=1,\ldots, 4N-3$. It is easy to see that in
this case $b({\cal T}_1)=b({\cal T})-(4N-3)*|w'|$.

b) If there exists a $\Theta$-band in
$\Psi_1({\cal T}_1)=\Psi_{1,1}$ which crosses
all the $\kappa$-bands $\bb_i({\cal T}_1)$, $i=1,\ldots,4N-3$,
then there exists a $\Theta$-band
containing the first cells of
all these $\kappa$-bands (this is the maximal $\Theta$-band
containing the first cell
in $\bb_1({\cal T}_1)$), so the indices of the $\Theta$-letters
in the relations corresponding to these cells are the same. Therefore
the first letters in the words $\Lab(\topp(\bb_i({\cal T}_1)))$ are the same, a
contradiction.

c) If the the words $\Lab(\topp(\bb_j({\cal T})))$, $j=1,\ldots,4N-3$,
do not all have the same first letter,
then $w'$ cannot be a prefix of all these words, and we can apply part (a).
The lemma is proved. $\Box$

\vskip 0.1 in

If ${\cal T}_1$ does not satisfy property (P6') then by Lemma \ref{tau},
$b({\cal T}_1)<b({\cal T})$. If ${\cal T}_1$ does not satisfy (P6) then we can
apply the band moving construction and obtain a triple ${\cal T}_2$ which
satisfies (P1)--(P6). Notice that the band moving construction strictly
decreases the parameter $b$. If ${\cal T}_2$ does not satisfy property (P7) then
we can repeat the adjustment and band moving constructions and obtain a sequence
of triples ${\cal T}_2$, ${\cal T}_3,\ldots$ such that $b({\cal T}_1)>b({\cal
T}_2)>\ldots$. For some $s>0$ the triple ${\cal T}_s$ satisfies condition (P6').
By Lemma \ref{tau}, ${\cal T}_s$ satisfies condition (P6) and the triples
obtained from ${\cal T}_s$ by any number of adjustments satisfy (P6');
therefore we do not need further applications of the band moving construction to
obtain ${\cal T}_{s+1}, {\cal T}_{s+2},\ldots$.

As we noticed before, $c({\cal T}_1)<c({\cal T})$. Since we do not use the band
moving construction to obtain ${\cal T}_{s+1},...$, we have $c({\cal T}_{s+1}) >
c({\cal T}_{s+2}) >\ldots $ (notice that the band moving construction can
increase the parameter $c$).  This sequence cannot be infinite so there exists
$t$ such that the triple ${\cal T}_t$ satisfies all seven properties (P1)--(P7).
Without loss of generality we can assume now that ${\cal T}$ itself satisfies
these properties.

\bigskip
Now let ${\cal T}=(\Psi_1, \Psi_2, \Pi)$ be a decomposition of $\Delta$
satisfying properties (P1)--(P7).  We shall prove that the perimeters of
$\Psi_1, \Psi_2, \Pi$ relate nicely to the perimeter of $\Delta$.

From now on we shall assume that
\begin{equation}\label{numbern}
N\ge 6.
\end{equation}

Notice that because of Property (P6) every $\Theta$-band in $\Psi_1$
starting on $\topp(\bb_1)\cup \bott(\bb_{4N-3})$ ends on $p_1$.

If a $\Theta$-band $\rr$ of $\Psi_1$ starting on an edge $e\in t_i$, $i\le 4N-4$,
ends on the
contour of
$\bb_{4N-3}$ (resp. $\bb_1$) then we paint  the edge $e$ in
{\em red}
(resp. {\em orange}). We also shall call cells of these $\Theta$-bands
{\em red}
(resp. {\em orange}).

Since $\Theta$-bands cannot intersect, there could be at most one $i$
between $1$ and $4N-4$ such that $t_i$ contains both red and orange edges.
Also it is clear that if $t_i$ contains a red (resp. orange) edge then $t_{i+1}$
(resp. $t_{i-1}$) contains no orange (resp. red) edges. Therefore either
$t_1,...,t_{2N-3}$ contains no red edges or $t_{2N},...,t_{4N-4}$ contain no orange edges.
Without loss of generality we assume that the second possibility holds.

In this case let us erase the red paint from all red edges on the last 4
of the $t_j$'s ($j=4N-7, 4N-6, 4N-5, 4N-4$). We shall save the last four of $\Sigma_j$
for the later construction, so we shall not touch the corresponding $t_j$ and $y_j$
by a paint brush. Notice that by our assumption these $t_j$
do not have orange edges. Therefore every $\Theta$-band $\rr$ in $\Psi_1$
which crosses $\bb_{4N-7}$, must have one of its ends on a $t_\ell$ with
$\ell<4N-7$. We paint this end of $\rr$ in {\em red} if it is not red already.
Then we call the maximal $\Theta$-band starting on this edge {\em red} as well.

Among the $2N-7$ numbers $2N, ..., 4N-8$ exactly $N-4\ge 2$ numbers
are odd. Notice that if $i$ is odd then $y_i$ starts with
a ${\bf E}(0)$-edge.

Pick one of the odd numbers $i \in \{2N+1, 2N+3,..., 4N-11\}$.
It will be clear later why we do not consider $i=4N-9$.

Let $W'$ be the word written on $y_i$. We shall compare $||W'||$ and the
length of $t_i$. In order to do that we establish a correspondence between
$\alpha^{\pm 1}$-edges of $y_i$ and edges of $t_i$. Then we shall prove that
only a constant number of $\alpha^{\pm 1}$-edges on $y_i$ can correspond to
the same edge in $t_i$.

Let us paint all colorless $\Theta$-edges of $t_i$ in {\em green}. All cells
of the maximal $\Theta$-bands in $\Sigma_i$
starting on these edges will be called
{\em green}. Notice that since $t_i$ contains no orange edges,
all $\Theta$-edges
of $t_i$ are either green or red. Therefore the $\Theta$-band starting on
a green edge of $t_i$ ends on $p_1$.
Since by definition of $y_i$
no $\Theta$-band
in $\Sigma_i$ which starts on $\bb_i$ can end on $\bb_{i+1}$, every transition
cell in $\Sigma_i$ is either red or green.

We shall paint an $\alpha^{\pm 1}$-edge of $y_i$ in {\em green} if the maximal
${\bar Y}$-band in $\Sigma_i$ starting on this edge ends
on a green transition cell.

\begin{lm} \label{green} No more than 4 ${\bar Y}$-bands starting on
$\alpha^{\pm 1}$-edges
of $y_i$
can end on a green cell belonging to the same $\Theta$-band.
\end{lm}

{\bf Proof.} Indeed, suppose that 5 ${\bar Y}$-bands $\yyy_1$,...,$\yyy_5$
starting on
$\alpha$-edges $e_1, e_2, e_3, e_4, e_5$ end on the contour of the same
$\Theta$-band $\rr$. We can assume that $e_1$ precedes $e_2$ precedes ...
$e_5$ on $y_i$. Then first three of these edges are to the left of the ${\bf X}(0)$-edge
of $y_i$ or the last three of these edges are to the right of this
${\bf X}(0)$-edge. These two cases are similar so we shall consider only the
first case. Let $\rho_1$, $\rho_2$, $\rho_3$ be the green cells on the contour
of which the ${\bar Y}$-bands $\yyy_1, \yyy_2, \yyy_3$ end. These cells are
transition cells corresponding to the rules of $\sss$ involving $\alpha$.
Such transition cells must contain a ${\bf X}(0)$-edge on the contour, so
$\rho_i$, $\rho_2$ and $\rho_3$ belong to some maximal ${\bf X}(0)$-bands
$\xxx_1$, $\xxx_2$, $\xxx_3$ respectively
(\setcounter{pdeighteen}{\value{ppp}}
see Fig. \thepdeighteen).
Consider the subdiagram $\Sigma'$
of $\Sigma_i$ bounded by $\yyy_1$, $\yyy_3$, $\rr$ and $y_i$. The part
of $\xxx_2$ which is contained in $\Sigma'$ cannot cross $\rr$ (Lemma \ref{qr}).
It also cannot cross $\yyy_1$ and $\yyy_3$ because $\yyy_1$ or $\yyy_3$ would
contain a transition cells and ${\bar Y}$-bands cannot contain transition cells.
Therefore $\xxx_2$ must start or end on $y_i$ between $e_1$ and $e_2$.

\unitlength=1.00mm
\special{em:linewidth 0.4pt}
\linethickness{0.4pt}
\begin{picture}(85.33,66.33)
\put(14.00,8.33){\line(1,0){71.33}}
\put(85.33,8.33){\line(0,1){52.00}}
\put(85.33,60.33){\line(-1,0){71.00}}
\put(14.33,60.33){\line(0,-1){52.00}}
\put(78.67,60.33){\line(0,-1){52.00}}
\put(19.33,8.33){\line(0,1){52.00}}
\bezier{488}(35.00,60.33)(34.67,6.00)(75.33,60.33)
\bezier{156}(28.67,8.33)(24.33,32.33)(37.67,38.33)
\bezier{136}(32.33,8.33)(28.00,29.00)(39.33,35.67)
\bezier{104}(37.67,8.33)(40.00,26.33)(42.00,33.67)
\bezier{100}(45.67,33.33)(41.33,12.00)(41.67,8.33)
\bezier{144}(54.00,8.33)(64.67,23.33)(56.33,39.33)
\bezier{156}(58.67,41.33)(69.33,26.67)(58.67,8.33)
\bezier{428}(37.67,60.33)(36.00,13.33)(72.67,60.33)
\bezier{52}(39.67,48.00)(39.00,55.33)(41.33,60.33)
\bezier{52}(45.33,60.17)(42.00,55.17)(42.67,48.50)
\bezier{64}(48.67,20.33)(50.33,10.00)(45.33,8.33)
\bezier{60}(52.00,21.33)(53.00,10.67)(49.00,8.33)
\put(16.67,65.67){\makebox(0,0)[cc]{$\eee$}}
\put(81.33,66.33){\makebox(0,0)[cc]{$\fff$}}
\put(35.67,63.00){\makebox(0,0)[cc]{$\rr$}}
\put(43.67,63.00){\makebox(0,0)[cc]{$\xxx_2$}}
\put(48.00,5.00){\makebox(0,0)[cc]{$x$}}
\put(30.00,6.00){\makebox(0,0)[cc]{$\alpha$}}
\put(39.67,5.67){\makebox(0,0)[cc]{$\alpha$}}
\put(56.00,5.67){\makebox(0,0)[cc]{$\alpha$}}
\put(66.33,8.33){\vector(1,0){5.33}}
\put(68.67,5.33){\makebox(0,0)[cc]{$x$}}
\bezier{116}(39.67,47.83)(46.00,36.17)(48.50,20.33)
\bezier{116}(42.50,48.50)(46.00,41.00)(52.17,21.50)
\put(67.33,8.33){\circle*{0.75}}
\end{picture}

\begin{center}
\nopagebreak[4]
Figure \theppp.
\end{center}
\addtocounter{ppp}{1}
\vskip 0.2 in

But this contradicts the fact that there are no ${\bf X}(0)$-edges
between $e_1$ and $e_3$. This contradiction proves our lemma. $\Box$

Notice that for every $x\in {\bf X}(0)$ there exists at most one relation
in $\pp_N(\sss)$ which does not contain ${\bf E}(0)$-
or ${\bf F}(0)$-letters but contains $\alpha$.

This relation has one of the forms $x^\tau(\alpha x\alpha\iv)\iv=1$
or $x^\tau(\alpha\iv x\alpha)\iv=1$. Such relations will be called
{\em $x\alpha$-relations}. A cell corresponding to an $x\alpha$-relation
will be called
an $x\alpha$-cell. The contour of every $x\alpha$-cell
has two $\alpha$-edges with opposite orientation.
Let us mark these
edges and consider {\em $\alpha$-bands} consisting of auxiliary $\alpha$-cells
and these $x\alpha$-cells with common $\alpha$-edges.
Notice that every ${\bar Y}$-band starting on an $\alpha^{\pm 1}$-edge is an
$\alpha$-band, but $\alpha$-bands are in general longer, they can contain
some transition cells ($x\alpha$-cell).
An $\alpha$-band in $\Sigma_i$ starting on $y_i$
can end either on $t_i$ or on $y_i$ or on a contour of a transition
${\bf F}(0)$-cell.

We shall need two general lemmas about $\alpha$-bands.

\begin{lm} \label{511} If a $\bar Y$-band $\yyy$ in a reduced diagram
starts on the contour of an
$x$-cell $\pi_1$ and ends on the contour of an $x$-cell $\pi_2$ and
these two cells belong to the same ${\bf X}(0)$-band ${\cal B}$
then ${\cal B}$ contains an ${\bf E}(0)$-cell or an ${\bf F}$-cell between
$\pi_1$ and $\pi_2$.
\end{lm}

{\bf Proof.} We can assume that $\pi_1$ is the first cell in ${\cal B}$ and
$\pi_2$ is the last cell in this band. Suppose that ${\cal B}$ does not
have a ${\bf E}(0)$-cell
and does not have a ${\bf F}(0)$-cell between $\pi_1$ and $\pi_2$. Then
by Proposition \ref{mach1} all
${\bf X}(0)$-edges of ${\cal B}$ are the same. Therefore all these cells
correspond to the same $x\alpha$-relation $x^\tau(\alpha^{\pm 1}x\alpha^{\mp 1})\iv$
therefore the labels of $\topp(\bb)$ and $\bott(\bb)$
have the form $(\tau^a \alpha^b)^\ell$
where $a,b\in \{-1, 1\}$. Therefore all $\alpha$-edges on $\topp(\bb)$
(resp. $\bott(\bb)$) have the same orientation. Therefore our ${\bar Y}$-band
$\yyy$ cannot start and end on $\topp(\bb)$ and cannot start and end on
$\bott(\bb)$. But it also cannot start (end) on $\topp(\bb)$ and end (start)
on $\bott(\bb)$ because otherwise the maximal ${\bf X}(0)$-band containing
$\bb$ would intersect the ${\bar Y}$-band $\yyy$ and a ${\bar Y}$-band cannot have
transition cells.$\Box$

\begin{lm} \label{512} A reduced \vk diagram
over $\pp_N(\sss)$ without hubs cannot contain
$\alpha$-annuli.
\end{lm}

{\bf Proof.} Indeed, let $\aaa$ be an $\alpha$-annulus. If $\aaa$ does not
have ${\bf X}(0)$-cells then it is a ${\bar Y}$-annulus which is ruled out by Lemma
\ref{a}. So $\aaa$ must contain an ${\bf X}(0)$-cell. Thus it has at least 2
${\bf X}(0)$-edges on its inner contour. Therefore a $\bx$-band $\bb$ intersects
the annulus $\aaa$ twice. Then $\bb$ cuts $\aaa$ into two parts, $\aaa_1$
and $\aaa_2$. We can assume that $\aaa_1$ does not contain $\bx$-cells.
Then $\aaa_1$ is a ${\bar Y}$-band starting and ending on the contour of
an $\bx$-band $\bb$. By Lemma \ref{511} this band must contain an $\ex$ or
an $\fx$-cell $\pi$ between the two intersections with $\aaa$. Since $\aaa$
cannot have $\ex$-cells or $\fx$-cells, $\pi$ must be inside the
subdiagram bounded by the median of $\aaa$. But then the maximal
$\ex$ or $\fx$-band containing $\pi$ must intersect with $\aaa$. This is a
contradiction since $\aaa$ cannot contain $\ex$ and $\fx$-cells.
$\Box$.

Now we can continue painting edges of $y_i$  and $t_i$.

If an $\alpha$-band starts on a green edge of $y_i$ and ends on $y_i\iv$ then
we paint the inverse of the end edge (it belongs to $y_i$) green also.

\begin{lm}\label{green1} The number of green edges on $y_i$ cannot exceed
8 times the number of green edges on $t_i$.
\end{lm}

{\bf Proof.} Indeed, by the definition of green edges, at least
half of the green edges $e$ on $y_i$
have the property that the ${\bar Y}$-band  $\yyy$ starting on $e$ ends on a
green $\Theta$-cell $\rho$. If $e$ is such an edge then we associate with $e$
the green edge on $t_i$ which is the start or the end edge of the
maximal $\Theta$-band containing $\rho$ (whichever of
these two edges belongs to
$t_i$; if both start and end edges of this $\Theta$-band are on $t_i$
then we take the leftmost edge). By Lemma \ref{green} at most 4
green edges of $y_i$ are associated with the same green edge on $t_i$.
Therefore the total number of green edges on $y_i$ does not exceed 8 times
the number of green edges of $t_i$. $\Box$

Now we need more colors.

If an $\alpha$-band starting on $y_i$ ends on $t_i$ then we paint the start
and the end edges of this band in yellow. The following lemma is obvious.

\begin{lm} \label{yellow1}
The number of yellow edges on $y_i$ coincides with the number of yellow edges
on $t_i$.
\end{lm}

Let us denote the maximal
$\ex$-band starting on the $\ex$-edge of $y_i$ by
$\eee$, the maximal $\fx$-band starting on $y_i$ by
$\fff$ and the maximal $\bx$-band starting on $y_i$ by $\xxx$.

If an $\alpha$-band $\aaa$ starting on a non-green and non-yellow
edge $e$ of $y_i$ ends on the contour of a red $\fx$-cell such that
the maximal $\fx$-bands $\fff_e^{\pm 1}$
containing this red cell are not $\fff^{\pm 1}$
then both start and end edge of $\fff_e$ belong to $t_i$ (an $\fx$-band
cannot cross a $\kappa$-band). We paint $e$ and the start and end edges
of $\fff_e$ in {\em blue}.
\setcounter{pdtwentyone}{\value{ppp}}
See Fig. \thepdtwentyone.

\begin{lm}\label{blue} No more than two maximal
$\alpha$-bands starting on blue edges of
$y_i$ can end on the contour of the same $\fx$-band.
\end{lm}

{\bf Proof.} If there are three such maximal $\alpha$-bands then at least
two of them have start edges on the same side of the $\bx$-edge
of $y_i$.
Let $e_1$ and $e_2$ be blue edges of $(y_i)^{\pm 1}$
which are on the same
side of the $\bx$-edge of $(y_i)^{\pm 1}$, $e_1^{\pm 1}$
preceding $e_2^{\pm 1}$ on $y_i$,
let $\aaa_1$ and $\aaa_2$ be the maximal $\alpha$-bands starting
on $e_1$ and $e_2$ respectively.
Then the edges $e_1$ and $e_2$ either both belong to $y_i$ or both belong to
$(y_i)\iv$ (if these edges were on
opposite sides of the $\bx$-edge of $y_i$, one of them could belong to
$y_i$ and the other one could belong to $(y_i)\iv$).

Suppose that $\aaa_i$, $i=1,2$, ends
on the contour of a red cell $\rho_i$. Let $\rr_1$ and $\rr_2$
be the red $\Theta$-bands containing $\rho_1$ and $\rho_2$.
Let $\fff'$ be one of the two mutually inverse
maximal $\fx$-bands containing $\rho_1$ and $\rho_2$.

Notice that the $\fx$-band $\fff'$ cannot
intersect $\rr_1$ or $\rr_2$ twice by Lemma \ref{qr}. Therefore
the (red) start  edges of $\rr_1$ and $\rr_2$ are between the start and the end
edges of $\fff'$.
Notice also that between
the $\alpha$-edge and the $\fx$-edge on the
top or bottom of a $\fx$-cell considered as a $\Theta$-band, there always
exists a $\bx$-edge. This implies that if one of the maximal $\bx$-bands
containing
$\rho_1$ or $\rho_2$ is inside the subdiagram $\Sigma'$
bounded  by $\fff'$ and $t_i$
then the $\alpha$-band $\aaa_1$ must intersect $\fff$ which is impossible.
Therefore the maximal $\bx$-bands containing $\rho_1$ (resp. $\rho_2$) must be
outside $\Sigma'$ (recall that an $\bx$-band and an $\fx$-band cannot cross).

Since the $\Theta$-bands $\rr_1$ and $\rr_2$ cross
$\partial(\bb_{4N-7})^{\pm 1}$,
and since $i<4N-10$, $\rr_1$ and $\rr_2$ cross both
$\kappa$-bands $\bb_{4N-9}$
and $\bb_{4N-8}$ in $\Sigma_{4N-9}$.
It is clear that $\rr_2$ is higher than $\rr_1$ (that is the intersection
of $\rr_1$ $\bb_{i+1}$ has bigger number on $\bb_{i+1}$ than
the intersection of $\rr_2$ and $\bb_{i+1}$) .

The paths $\topp(\rr_i)$  and $\bott(\rr_i)$
start on $t_i$ and cross $\bb_{i+1}$.
On the boundary of the cell
$\rho_j$, $j=1,2$ we can read either the word
$$(x(0,\tau_j, 4)F(0,\tau_j,4))^{\theta} (x(0,\tau_j, \alpha)F(0,\tau_j,\alpha)\alpha\iv)\iv$$
or the inverse of this word. In the first case we shall call $\rho_j$
{\em positive}, in the other case -- {\em negative}.

Since $\Sigma'$ does not contain the
$\bx$-edge of $\rho_j$,  either both the edge labelled by
$F(0,\tau_j, 4)$,
and the edge labelled by $x(0,\tau_j,4)$ belong to the path $\topp(\rr_j)\iv$
or both of these edges belong to the path $\bott(\rr_j)\iv$, so they are
not oriented toward $\bb_{i+1}$. Similarly the edge
of $\rho_j$ labelled by
$F(0,\tau_j, \alpha)$,
and the edge labelled by $x(0,\tau_j,\alpha)$ belong both to the path
$\bott(\rr_j)\iv$ or both to the path $\topp(\rr_j)\iv$.

Let $e_j'$ be the $\alpha$-edge of $\rho_j$.
Then $e_j'$ is the end edge of the $\alpha$-band $\aaa_j$.
In every cell of $\aaa_j$ two $\alpha$-edges form an opposing pair.
Therefore
$e_j'$ and $e_j$ form an opposing pair of edges on
the boundary of $\aaa_j$. Since $e_1$ and $e_2$ either both point
toward $\bb_{i+1}$ or both point in the opposite direction,
the edges $e_1'$ and $e_2'$ must also have the same orientation meaning that
if we trace the boundary of $\rho_1$ and $\rho_2$ clockwise then we
either pass through $e_1'$ and $e_2'$ or we pass through
$(e_1')\iv$ and $(e_2')\iv$.

This implies that either both $\rho_1$ and $\rho_2$ are positive or both
are negative. Therefore the common $\Theta_+$-edges of $\rho_j$ with
$\bb_{4N-9}$ either both point toward $y_{4N-9}$ or both point in the
opposite direction.

The subdiagram $\Gamma$ of $\Sigma_{4N-9}$ bounded by $y_{4N-9}$ on the bottom,
$\rr_2$ on the top, $\bb_{4N-9}$ on the left and $\bb_{4N-8}$ of the right
is a sector. By Proposition \ref{prop12} we can assume that $\Gamma$
is a computational sector. Let $C=W_1,...,W_m$ be the
corresponding computation connecting $W_1=\Lab(\topp(\rr_2)$ and $\Lab(y_{4N-9})$.
Then $\Lab(\bott(\rr_1))$ is one of the words $W_n$. Since there exists a
computation of $\sss$ connecting $W_m$ and $W_0$, there also exist computations
connecting $W_1$ and $W_n$ with $W_0$.
The rules used in the transitions $W_1\to W_2$
and $W_{n-1}\to W_n$ are $R_{4,\alpha}(\tau_1)^{\pm 1}$ and
$R_{4,\alpha}(\tau_2)^{\pm 1}$. The fact that we proved in the previous
paragraph shows that either both of these rules are positive or both are
negative.  Therefore part 5 of Proposition \ref{mach1} applies either to the
computation $W_1,...,W_n$ or to the inverse computation. In any case, in this
computation, there exists a transition $W_\ell\to W_{\ell+1}$ such that the
rule $\tau$ applied in this transition contains a word $Ex$
as one of its left sides
where $E\in \ex$, $x\in \bx$. Let $\bar\rr_3$ be the maximal
$\Theta$-band corresponding to this transition in the computational sector
$\Gamma$.  Let $\rr_3$ be the maximal $\Theta$-band in $\Psi_1$ containing $\bar
\rr_3$.

\vskip 0.2 in
\unitlength=1.00mm
\special{em:linewidth 0.4pt}
\linethickness{0.4pt}
\begin{picture}(149.67,82.00)
\put(12.33,10.00){\line(1,0){47.67}}
\put(60.00,10.00){\line(0,1){69.67}}
\put(12.67,10.00){\line(0,1){69.67}}
\bezier{464}(18.33,79.33)(28.33,23.33)(47.33,79.67)
\bezier{192}(33.33,10.00)(40.00,25.33)(37.00,56.00)
\bezier{272}(49.67,10.00)(59.33,48.00)(44.67,72.33)
\put(21.33,10.00){\circle*{1.49}}
\put(22.00,10.00){\vector(1,0){5.33}}
\put(23.33,6.00){\makebox(0,0)[cc]{$x$}}
\put(95.33,79.67){\line(0,-1){69.67}}
\put(95.33,10.00){\line(1,0){54.33}}
\put(149.67,10.00){\line(0,1){69.67}}
\bezier{664}(24.00,78.00)(29.33,32.67)(149.33,35.67)
\bezier{500}(38.67,78.67)(42.00,62.33)(149.67,52.67)
\bezier{512}(35.00,70.00)(73.00,35.67)(149.67,44.00)
\bezier{112}(37.67,55.33)(48.33,57.33)(55.33,73.33)
\put(47.00,61.33){\circle*{3.40}}
\put(116.00,32.33){\makebox(0,0)[cc]{$\rr_1$}}
\put(115.00,58.67){\makebox(0,0)[cc]{$\rr_2$}}
\put(114.67,46.00){\makebox(0,0)[cc]{$\rr_3$}}
\put(17.33,82.00){\makebox(0,0)[cc]{$\fff'$}}
\put(34.33,30.67){\makebox(0,0)[cc]{$\aaa_1$}}
\put(50.00,31.00){\makebox(0,0)[cc]{$\aaa_2$}}
\put(32.33,10.00){\vector(1,0){3.67}}
\put(32.33,10.00){\circle*{1.33}}
\put(47.67,10.00){\vector(1,0){4.33}}
\put(47.67,10.00){\circle*{1.33}}
\put(60.00,10.00){\dashbox{1.00}(35.33,0.00)[cc]{}}
\put(33.33,5.67){\makebox(0,0)[cc]{$e_1$}}
\put(49.00,6.00){\makebox(0,0)[cc]{$e_2$}}
\put(55.33,76.67){\makebox(0,0)[cc]{$\xxx_1$}}
\put(46.67,54.33){\makebox(0,0)[cc]{$\Gamma_1$}}
\put(18.67,68.67){\vector(1,0){3.33}}
\put(56.00,70.67){\vector(-4,1){4.67}}
\put(37.00,55.33){\circle*{4.22}}
\put(44.33,72.00){\circle*{4.22}}
\put(47.00,61.33){\circle*{4.67}}
\end{picture}

\begin{center}
\nopagebreak[4]
Figure \theppp.
\end{center}
\addtocounter{ppp}{1}
\vskip 0.2 in

Let $\Gamma_1$ be the subdiagram of $\Psi_1$ bounded by $\bb_{i+1}$,
$\rr_1$, $\rr_2$, $\fff'$. Let $\xxx_1$ be the maximal $\bx$-band in
$\Gamma_1$ starting on the boundary of $\rho_1$. Then $\xxx_1$ cannot
have two common cells with $\rr_1$, it cannot cross $\bb_{i+1}$ and it cannot cross $\fff'$ (it can only
touch $\fff'$). Therefore $\xxx_1$ crosses $\rr_2$. Since $\rr_3$
is between $\rr_1$ and $\rr_2$, it must cross $\xxx_1$. Let $\rho$ be the
intersection cell. Then $\rho$ corresponds to the relation containing
both an $\bx$-letter and $\tau$.
By definition of $\pp_N(\sss)$ there is only one (modulo taking inverses and
cyclic shifts) relation which contains both an $\bx$-letter and $\tau$.
This relation or its inverse is written on the boundary of the cell
in $\bar\rr_3$ which contains an ${\bf X}(0)^{\pm 1}$-edge. By the choice
of $\tau$, this relation contains an $(\ex)^{\em }$-letter.
Therefore the boundary of $\rho$ contains an $\ex$-edge. Thus $\Gamma_1$
contains an $\ex$-edge. Let $\eee_1$ be the maximal $\ex$-band in $\Sigma_i$
containing this edge. The $\ex$-band $\eee_1$ cannot cross $\fff'$ or $\bb_{i+1}$.
By Lemma \ref{q}, $\eee_1$ is not an annulus. Therefore $\eee_1$ crosses
$\rr_1$ between $\rho_1$ and $\bb_{i+1}$. Since the $\ex$-band $\eee_1$
cannot cross $\aaa_1$, either the start or the end edge of $\eee_1$ must belong to
$y_i^{\pm 1}$. This edge must be between $e_1$ and the $(\fx)^{\pm 1}$-edge
of $(y_i)^{\pm 1}$.
But $(y_i)^{\pm 1}$ does not contain $(\ex)^{\pm 1}$-edges between the $\alpha^{\pm
1}$-edges and the $(\fx)^{\pm 1}$-edge. This contradiction completes the proof
of our lemma. $\Box$.

This lemma immediately implies the following result.

\begin{lm} \label{blue1}
The number of blue edges on $y_i$ does not exceed the number of blue edges
on $t_i$.
\end{lm}

We continue painting the edges of $y_i$. We shall paint an $\alpha^{\pm 1}$-edge
$e$ of $y_i$ in {\em pink} if the $\alpha$-band starting on this edge
ends on the contour of a green $\Theta$-cell.

\begin{lm}\label{pink} There are no more than two $\alpha$-bands starting
on $(y_i)^{\pm 1}$ and ending on the boundary of the same green $\Theta$-band.
\end{lm}

{\bf Proof.} The proof is similar to the proof of Lemma \ref{green}. Suppose that
there
are three such $\alpha$-bands $\aaa_1$, $\aaa_2$, $\aaa_3$ and the start edge
of $\aaa_1$ precedes the start edge of $\aaa_2$ which precedes the start edge
of $\aaa_3$. Let $\rho_1$, $\rho_2$ and $\rho_3$ be the green cells containing
the end edges of $\aaa_1$, $\aaa_2$ and $\aaa_3$ respectively and belonging
to the same $\Theta$-band $\rr$. Then
the maximal $\fx$-band $\fff_2$ containing $\rho_2$ cannot cross $\rr$
again and cannot cross $\aaa_1$ and $\aaa_3$. Thus $\fff_1$ must start or end
on $(y_i)^{\pm 1}$. This contradicts the fact that $y_i$ cannot contain
an $(\fx)^{\pm 1}$-edge between two $\alpha^{\pm 1}$-edges. $\Box$

This lemma implies the following fact.

\begin{lm} \label{pink1} The number of pink edges on $y_i$ does not exceed
twice the number of green edges on $t_i$.
\end{lm}

Our next goal will be to prove that every $\alpha^{\pm 1}$-edge on $y_i$
is either green, yellow, blue or pink.

Let us suppose that it is not so. Then there exists
a colorless edge $e\in (y_i)^{\pm 1}$ labelled by $\alpha$
such that the maximal $\alpha$-band
starting on $e$ ends either on an edge $e'$ of $(y_i)\iv$
or on the contour of a $\fx$-cell which belongs to the $\fx$-band $\fff$ and
to a red $\Theta$-band.
In the first case $e'$ is colorless (it could be only green but then $e$ would also
be green and we assumed that $e$ is colorless). In this case we
paint both $e$ and $(e')\iv$ in {\em brown}. In the second case we
paint $e$ in {\em black}.

Notice that if there are brown or black edges on $y_i$ then
$t_i$ contains red edges. Therefore there exists a $\Theta$-band in $\Psi_1$
which starts on $t_i$ and crosses $\bb_{4N-7}$.

\begin{lm} \label{brown} $y_i$ does not contain brown edges.
\end{lm}

{\bf Proof.} Suppose that $y_i$ contains a brown edge $e$. Let $e'$ be the
end edge of the maximal $\alpha$-band $\aaa$ in $\Sigma_i$ starting on $e$.
As we know, $e'$ is also brown. Since $e$ belongs to $y_i$, $e'$ belongs to
$y_i\iv$ and the labels of $e$ and $e'$ are the same.
Therefore these edges cannot be on the same side of the $\bx$-edge
of $y_i$ (since the label of $y_i$ is a reduced word).
Thus one of these edges is between the $\ex$-edge and
the $\bx$-edge and the other one is between the $\bx$-edge and the $\fx$-edge.
Suppose that $e$ is between the $\ex$-edge and the $\bx$-edge. The case
when $e'$ is between the $\ex$-edge and the $\bx$-edge is completely similar
(although not identical because we have the additional assumption that
$e\in y_i$).

Since $e\in y_i$ and the label of $e$
is $\alpha$, the part of the word $\Lab(y_i)$
between the $\ex$-letter and the $\fx$-letter has the form $\alpha^m x\alpha^n$
where $m>0, n<0$.

Since $e$ is not a green or yellow edge,
the ${\bar Y}$-band starting
on $e$ ends on a red $\Theta$-cell. Let $\rr'$ be the maximal red $\Theta$-band
which contains this cell.
Since $e$ is to the left of the $\bx$-band $\xxx$, $\rr'$ intersects
$\xxx$.

Let $\rho$ be the first cell in $\xxx$ (its contour has a common path with
$y_i\iv$). Let $\rr$ be the maximal $\Theta$-band in $\Psi_1$
containing $\rho$. Since $\rr$
cannot intersect $\rr'$, it must cross $\bott(\bb_{4N-7})$.
So $\rr'$ is a red $\Theta$-band (\setcounter{pdtwenty}{\value{ppp}}
see Fig. \thepdtwenty).

\vskip 0.1 in
\unitlength=1.00mm
\special{em:linewidth 0.4pt}
\linethickness{0.4pt}
\begin{picture}(134.00,84.00)
\put(9.33,10.00){\line(0,1){68.00}}
\put(55.67,10.00){\line(0,1){67.67}}
\bezier{124}(35.33,10.00)(28.00,20.67)(30.00,38.33)
\bezier{132}(30.00,38.33)(37.00,64.67)(32.67,69.00)
\bezier{600}(20.00,10.33)(30.67,84.00)(47.67,10.33)
\put(84.00,10.00){\line(0,1){67.67}}
\put(134.00,10.00){\line(0,1){67.33}}
\bezier{676}(13.67,61.00)(22.00,7.33)(134.00,31.33)
\put(31.67,10.00){\line(-1,1){5.67}}
\put(26.00,15.67){\line(1,0){17.00}}
\put(43.00,15.67){\line(-5,-6){4.67}}
\put(38.33,10.00){\line(-1,0){6.67}}
\put(27.33,15.67){\vector(1,0){3.33}}
\put(41.00,15.67){\vector(-1,0){3.67}}
\put(26.00,15.67){\circle*{0.67}}
\put(42.33,15.33){\circle*{1.33}}
\put(25.67,15.67){\circle*{1.49}}
\put(31.67,10.00){\circle*{0.94}}
\put(31.67,10.00){\vector(1,0){5.67}}
\put(94.33,10.00){\line(-1,1){5.67}}
\put(88.67,15.67){\line(1,0){17.00}}
\put(105.67,15.67){\line(-5,-6){4.67}}
\put(101.00,10.00){\line(-1,0){6.67}}
\put(90.00,15.67){\vector(1,0){3.33}}
\put(103.67,15.67){\vector(-1,0){3.67}}
\put(88.67,15.67){\circle*{0.67}}
\put(105.00,15.33){\circle*{1.33}}
\put(88.33,15.67){\circle*{1.49}}
\put(94.33,10.00){\circle*{0.94}}
\put(94.33,10.00){\vector(1,0){5.67}}
\put(9.33,10.00){\line(1,0){46.33}}
\put(84.00,10.00){\line(1,0){50.00}}
\put(55.67,10.00){\dashbox{1.00}(28.33,0.00)[cc]{}}
\put(42.67,15.67){\line(1,0){13.00}}
\put(84.00,15.67){\line(1,0){50.00}}
\put(55.67,15.67){\dashbox{1.37}(28.33,0.00)[cc]{}}
\bezier{124}(31.67,10.00)(12.33,12.33)(11.67,23.67)
\bezier{104}(25.67,15.67)(16.00,14.00)(12.33,29.67)
\bezier{124}(98.33,10.00)(91.00,20.67)(93.00,38.33)
\bezier{132}(93.00,38.33)(100.00,64.67)(95.67,69.00)
\bezier{60}(28.67,15.67)(31.33,26.33)(35.00,26.67)
\bezier{88}(40.33,15.67)(41.33,23.33)(27.67,25.67)
\end{picture}

\begin{center}
\nopagebreak[4]
Figure \theppp.
\end{center}
\addtocounter{ppp}{1}
\vskip 0.2 in

Let $\Sigma'$ be the diagram bounded by $\rr$, $\bb_{4N-9}$, $\bb_{4N-8}$,
$y_{4N-9}$.
Then $\Sigma'$ is a sector which we can assume
to be a computational sector. This sector has only two mutually inverse
maximal
$\Theta$-bands (one of which is a part of $\rr$). Indeed, if there were
a $\Theta$-band $\rr''$ in $\Sigma'$ below $\rr$ then the
maximal $\Theta$-band in $\Psi_1$ containing $\rr''$ would intersect $\xxx$ and the
intersection cell would be lower than $\rho$. Thus the computation
corresponding to the sector $\Sigma'$ consists of one transition $W_1\to W_2$
where $W_1=\Lab(y_{4N-9})$. Recall that the $\alpha$-part of $\Lab(y_{4N-9})$
is equal to $\alpha^m x\alpha^n$ where $n\ne 0$ and $m\ne 0$.
Hence the only rules applicable to $W_1$ are rules of the form
$x^\tau(\alpha^\epsilon x\alpha^{-\epsilon})\iv$.
If $\epsilon=-1$ then $|W_2|<|W_1|$, a contradiction with property P7.
Thus the boundary of the cell in $\Sigma'$ which
has an $\bx$-edge has label $\tau\alpha x\alpha\iv\tau\iv x\iv$. Denote this cell by
$\rho''$.
Since the cells $\rho$ and $\rho''$ belong to the same $\Theta$-band and their
$\bx$-edges have the same orientation on $s_1$,
on the boundary of $\rho$ one can read the same word.
Let $\aaa_1$ and $\aaa_2$
be the maximal $\alpha$-bands in $\Sigma_i$ which start on the two $\alpha$-edges
$e_1$ and $e_2$
of the contour of $\rho$
where $e_1$ belongs to the subdiagram bounded by $\eee$, $\xxx$, $y_i$
and $t_i$ and $e_2$ belongs to the complement of this subdiagram in $\Sigma_i$.
These $\alpha$-bands cannot intersect the $\alpha$-band $\aaa$. Therefore
they must end on $y_i^{\pm 1}$ ($\alpha$-annuli are ruled out by
Lemma \ref{512}). The $\alpha$-band $\alpha_1$ cannot end to the left of the
$\bx$-edge of $y_i$ because the start and the end edges of an $\alpha$-band
must form an opposing pair. Therefore the end edge of $\aaa_1$ is to the right
of the $\bx$-edge of $y_i$. Similarly the end edge of $\aaa_2$ must be to the
left of this $\bx$-edge. This implies that $\aaa_1$ and $\aaa_2$ must intersect,
a contradiction. $\Box$

\begin{lm} \label{black} $y_i$ does not contain black edges.
\end{lm}

{\bf Proof.} Let $e$ be a black $\alpha$-edge on
$y_i^{\pm 1}$. Let $\aaa$ be the maximal
$\alpha$-band starting on $e$, and let
$\rho$ be the cell from $\fff$ on whose boundary
$\aaa$ ends. Let $\rr$ be the maximal red $\Theta$-band containing $\rho$ (\setcounter{pdtwentyfive}{\value{ppp}}
see Fig. \thepdtwentyfive).

The subdiagram $\Gamma$ bounded by $\bb_{4N-9}$, $\bb_{4N-8}$, $\rr$ and
$y_{4N-9}$
is a sector and by Proposition \ref{prop12}
can be converted into a computational
sector without changing the boundary or increasing the area. So we can assume
that $\Gamma$ is a computational sector. Let $C=(W_1,...,W_\ell)$ be the
corresponding computation where $W_1$ (resp. $W_2$)
is the label of the part of $\topp(\rr)$ (resp. $\bott(\rr)$)
between $\bb_{4N-9}$ and $\bb_{4N-8}$, $W_\ell=\Lab(y_{4N-9})=\Lab(y_i)=W$.
The label $\theta$
of the $\Theta$-edges of $\rr$ is an $S$-rule which involves
$\alpha$ and an $\fx$-letter. Therefore this rule is $R_{4,\alpha}(\tau)$
for some $\tau$.
Since there exists a computation of $\sss$
connecting $W_\ell$ and $W_0$ (this is the computation
corresponding to the computational disc $\Pi$), there exists a reduced
computation connecting $W_0$ and $W_1$.
By Proposition \ref{mach1} either $W_1$ or $W_2$ is a positive word and the
degree of $\alpha$ in $W_1$ and $W_2$ is positive.

Suppose that $W_\ell$ contains $\alpha\iv$. Let $W_m$, $m<\ell$,
be the last word
in the computation $C$ which does not contain $\alpha\iv$. Then the
$\alpha$-part of $W_{m+1}$ has the form $\alpha^a x \alpha^b$ where $x\in \bx$,
$ab<0$.  Therefore only $x\alpha$-rules are applicable to $W_{m+1}$. Since the
computation $C$ is semiproper, all words in the computation $W_{m+1},...,W_\ell$
contain $\alpha\iv$, so in this computation only $x\alpha$-rules can apply and
each of them inserts a new $\alpha\iv$. Since $x\alpha$-rules do not affect
non-$\alpha$ parts of admissible words, the length of $W_{\ell}$ is greater than
the length of $W_{m+1}$ provided $m+1<\ell$. This contradicts Property P7.
Thus $m+1=\ell$. The only rules which can insert $\alpha\iv$ are
$R_{4,\alpha}(\tau)$ and $x\alpha$-rules. Thus one of these rules is applied
in the
transition $W_m\to W_\ell$. But if $R_{4,\alpha}(\tau)$ applies in this
transition then the degree of $\alpha$ in $W_m$ should be 0. Since $W_m$ is
normal and positive, it must contain no $\bar Y$-letters.  Then
$|W_{m+1}|=|W_m|+4$ which contradicts Property P7.

If an $x\alpha$-rule is applied in the transition $W_m\to W_{m+1}$ then
$|W_{m+1}|=|W_m|+2$ and we again get a contradiction.

Thus the $\alpha$-part of $W_\ell$ does not contain $\alpha\iv$, so all
$\alpha^{\pm 1}$-edges of $y_{4N-9}$ have label $\alpha$. Since the labels
of $y_i$ and $y_{4N-9}$ are the same, all $\alpha^{\pm 1}$-edges
of $y_i$ have label $\alpha$.

\vskip 0.2 in
\unitlength=1mm
\special{em:linewidth 0.4pt}
\linethickness{0.4pt}
\begin{picture}(152.00,99.33)
\put(12.67,95.33){\line(0,-1){90.00}}
\put(12.67,5.33){\line(1,0){47.00}}
\put(59.67,5.33){\line(0,1){90.00}}
\put(59.67,5.33){\dashbox{1.00}(34.00,0.00)[cc]{}}
\put(93.67,95.33){\line(0,-1){90.00}}
\put(93.67,5.33){\line(1,0){55.33}}
\put(149.00,5.33){\line(0,1){90.00}}
\bezier{292}(40.67,5.33)(40.67,51.33)(59.67,70.33)
\put(59.67,70.00){\circle*{4.22}}
\bezier{572}(33.67,93.67)(47.33,55.33)(149.00,64.67)
\put(149.00,64.33){\circle*{3.89}}
\bezier{308}(149.00,64.00)(118.67,48.33)(149.00,18.33)
\put(149.00,18.33){\circle*{4.71}}
\bezier{720}(148.67,18.33)(25.67,14.33)(19.33,71.00)
\put(59.67,27.67){\circle*{4.47}}
\bezier{224}(59.67,27.33)(37.33,46.33)(59.67,60.00)
\bezier{612}(25.33,89.00)(76.00,30.00)(149.00,49.33)
\put(59.67,60.00){\circle*{3.89}}
\put(149.33,49.00){\circle*{4.81}}
\bezier{152}(148.33,48.67)(132.67,38.00)(149.00,28.33)
\bezier{640}(149.00,28.33)(37.33,32.33)(20.67,77.67)
\put(57.67,86.33){\vector(1,0){5.00}}
\put(145.67,84.00){\vector(1,0){6.33}}
\put(79.00,70.33){\makebox(0,0)[cc]{$\rr$}}
\put(22.00,9.00){\makebox(0,0)[cc]{$y_i$}}
\put(104.00,9.67){\makebox(0,0)[cc]{$y_{4N-9}$}}
\put(39.00,14.33){\vector(1,0){4.00}}
\put(134.67,54.00){\vector(1,0){5.00}}
\put(48.67,36.67){\vector(1,0){4.33}}
\put(59.67,99.00){\makebox(0,0)[cc]{$\fff$}}
\put(149.00,99.33){\makebox(0,0)[cc]{$\fff'$}}
\put(37.33,21.00){\makebox(0,0)[cc]{$\aaa$}}
\end{picture}

\begin{center}
\nopagebreak[4]
Figure \theppp.
\end{center}
\addtocounter{ppp}{1}
\vskip 0.2 in

Since $e$ has label $\alpha$, it must belong to $y_i$. The start and end edges
of every $\alpha$-band form an opposing pair on the contour of this band.
Therefore the end edge of $\aaa$ must belong to the inverse of the boundary of
the cell $\rho$. Therefore on the boundary of $\rho$, one can read
$F\iv x\iv \gamma\iv \alpha\iv x'F'\gamma$ where $\gamma\in \Theta_+$,
$x,x'\in \bx$, $F,F'\in \fx$. Therefore the same word can be read on the
boundary of the intersection $\rho'$ of $\rr$ and the $\fx$-band $\fff'$
of $\Gamma$ starting on $y_{4N-9}$. Hence the rule applied in the
transition $W_1\to W_2$ in the computation $C$ is $R_{4, \alpha}(\tau)\iv$
for some $\tau$. Hence the length of $W_2$ is greater than the
length of $W_1$ (Proposition \ref{mach1}).

Let $\rr'_1, ..., \rr'_{\ell-1}$ be the maximal $\Theta$-bands of $\Gamma$
starting on the contour of $\bb_{4N-9}$ counted from top to bottom
(so $y_{4N-9}=\bott(\rr_{\ell-1})$\ ) and let $\rr_1, \rr_2,...,\rr_{\ell-1}$
be the red $\Theta$-bands of $\Psi_1$ containing $\rr_1',...,\rr_{\ell-1}'$
(so $\rr=\rr_1$).
Then $W_2$ is written on the bottom of $\rr_1'$. Let $e'_1$ be the last
$\alpha$-edge of $\bott(\rr_1')$. Let $\aaa'_1$ be the maximal $\alpha$-band
in $\Gamma$ starting on $e_1'$. If $\aaa'_1$ ends on $y_{4N-9}$ then
the degree of $\alpha$ in $W_\ell$ is greater than the degree of $\alpha$
in $W_1$. Since $W_1$ is a positive word (by Proposition \ref{mach1}),
$W_\ell$ is longer than $W_1$
which contradicts Property P7.

The $\alpha$-band $\aaa_1'$ cannot end on the contour of $\rr_1$ because
then $W_1$ would contain $\alpha\iv$.
Therefore the only possibility is that $\aaa_1'$ ends on the contour of
the $\fx$-band $\fff'$. Let the $\fx$-cell $\rho_1'$
on the contour of which this
band ends belong to the $\Theta$-band $\rr_{j_1}$. Considering the
orientation of $\alpha$-edges of $\bott(\rr_1)$ and
the orientation of $\rho_1'$ we conclude that the rule applied in
the transition $W_{j_1}\to W_{j_1+1}$ in the computation $C$ is
$R_{4,\alpha}(\tau_1)$ for some $\tau_1$.
Let $\rho_1$ be the intersection of $\rr_{j_1}$ and $\fff$. Then the maximal
$\alpha$-band $\aaa_1$ starting on the $\alpha$-edge of the contour of $\rho_1$
cannot cross $\aaa$ and cannot end on $y_i^{\pm 1}$ (because all $\alpha^{\pm 1}$-edges
of $y_i$ have label $\alpha$). Therefore $\aaa_1$ must end on the contour
of an $\fx$-cell $\rho_2$ from $\fff$. This cell must belong to a red $\Theta$-band
$\rr_{j_2}$ for some $j_2$ between $1$ and $j_1$. Let $\rho_2'$ be the
intersection of $\rr_{j_2}$ with $\fff'$. Then the $\alpha$-band $\aaa_2'$
starting on the $\alpha$-edge of the contour of $\rho_2'$ ends on the
contour of some cell $\rho_3'$ of $\fff$ (it cannot cross $\aaa_1'$),
and the red $\Theta$-band $\rr_{j_3}$ containing this cell intersects $\fff$ at
some cell $\rho_3$. Continuing in this manner we shall construct a sequence
of different cells $\rho_1,\rho_2,\rho_3,...$ of $\fff$. This process obviously
cannot stop which contradicts the fact that $\fff$ contains only
finitely many cells (see Fig. \thepdtwentyfive).
This contradiction completes the proof of the lemma.
$\Box$

Thus we have proved that every $\alpha^{\pm 1}$
edge of $y_i$ is either green or yellow or blue
or pink.  This allows us to prove the following statement.

Let $n_i$ be the number of non-red edges of $t_i$. Recall that $W$ denotes the
label of $y_i$.

\begin{lm} \label{count} $||W||=||\Lab(y_i)||\le 40n_i+15k+6<61k n_i$
(where $k$ is the number of tapes of the machine $M$).
\end{lm}

{\bf Proof.} Indeed, by Lemmas \ref{green1}, \ref{yellow1}, \ref{blue1},
\ref{pink1}, the number of green and blue edges on $y_i$ does not exceed
$8+2=10$ times the number of green edges of $t_i$, the number of yellow edges of
$y_i$ does not exceed the number of yellow edges of $t_i$ and the number of pink
edges of $y_i$ does not exceed the number of pink edges on $t_i$. Therefore the
degree $m_i$ of $\alpha$ in $W=\Lab(y_i)$ does not exceed $10n_i$. Since $W$ is
a normal word, $||W||=4m_i+15k+6$. This implies the inequalities of the lemma.
$\Box$

We shall also need to estimate the length of $W$.

\begin{lm} \label{psi1}
$|W|\le (30k+12)|p_1|$
\end{lm}

{\bf Proof.} Indeed, every ${\bar Y}$-band in $\Psi_1$ starting on $y_i$ ($i$
between $1$ and $4N-4$) ends
either on $p_1$ or on the contour of a maximal $\Theta$-band in $\Psi$. Since
every maximal $\Theta$-band either starts or ends on $p_1$, it is enough to show
that the number of maximal ${\bar Y}$-bands in $\Psi_1$ starting on $y_i$ and
ending on the contour of the same $\Theta$-band does not exceed $30k+12$.  The
proof is similar to the proof of Lemma \ref{green}.  Suppose that there exist
 $30k+13$ $\bar Y$-bands starting on $y_i$ and ending on the contour of the same
$\Theta$-band $\rr$. Since the number of $Q$-edges on $y_i$ is $15k+6$, there
exist three $\bar Y$-bands $\yyy_1, \yyy_2, \yyy_3$ ending on $\rr$ and starting
on edges $e_1, e_2, e_3$ of $y_i$ such that $e_1$ precedes $e_2$ which
precedes $e_3$,
and there are no $Q$-edges between $e_1$ and $e_3$. Then the maximal $Q$-band
containing the intersection of $\yyy_2$ and $\rr$ must end or start on $y_i$
between $e_1$ and $e_3$ because it cannot cross $\rr$ twice and cannot cross
$\yyy_1$ and $\yyy_3$. This contradicts the fact that there are no $Q$-edges
between $e_1$ and $e_3$ on $s_1\iv$.  $\Box$

Now given our decomposition $(\Psi_1, \Psi_2, \Pi)$ of $\Delta$
we define a new decomposition $(\Gamma, E, \Pi)$ in the following way.
Let $\Sigma'$ be inverse of the union of subdiagrams $\Sigma_j$ where
$j=4N-7, 4N-6, 4N-5, 4N-4$ with the $\kappa$-bands
$\bb_{4N-7}$ and $\bb_{4N-3}$
removed.
Then $$\partial(\Sigma')=u_1\iv q' u_2(p')\iv$$
where $u_1=\bott(\bb_{4N-7})$, $p'\subset p_1$, $u_2=\topp(\bb_{4N-3})$,
$q'\subset s_1\iv$. Notice that $\Lab(u_2)=\Lab(\bott(\bb_{4N-3}))$
and $\Lab(q')$ is equal to $\Lab(s_2\iv)$ if we change indices of the labels
of $\kappa$-edges. Let us make this change of indices in the whole $\Sigma'$
and denote the resulting diagram by $\Sigma$.

Then we can glue $\Sigma$ and $\Psi_1$ along $\bott(\bb_{4N-3})$.
The resulting diagram will be denoted by $\Gamma$.
We also can glue $\Sigma\iv$ and $\Psi_2$ along $s_2\iv\bott(\bb_{4N-3})$.
The resulting diagram will be denoted by $E$.

Notice that $$\partial(\Gamma)=\bott(\bb_1)p_1qus's_1$$
where $|q|=\sum_{j=4N-7}^{4N-3} |t_i|+3<|p_1|$, $\Lab(u)=\Lab(\bott(\bb_{4N-7}))\iv$,
$\Lab(s')=\Lab(s_2\iv)$; $$\partial(E)=\bott(\bb_1)\iv u\iv q\iv p_2.$$

It is clear that by gluing $E$, $\Gamma$ and $\Pi$ we can get a diagram
with the same boundary label as $\Delta$. It is also clear that
the number of hubs in $E$ is $h(\Delta)-1$ and $\Gamma$ does not have hubs.

We shall call the triple $(\Gamma, E, \Pi)$ a decomposition of type 2 of $\Delta$.
We can now estimate the perimeters of $\Gamma, E, \Pi$.

\begin{lm}\label{nice} Let the boundary label of
the computational disc $\Pi$ be
$K(W)$ for some word $W$.
Then there exist positive constants $\epsilon_1$ and $\epsilon_2$
such that
\begin{enumerate}
\item $n(\Gamma)\le \epsilon_1 n(\Delta)$ and
\item $n(E)+\epsilon_2 ||W||\le n(\Delta)$.
\end{enumerate}
\end{lm}

{\bf Proof.}  We have that
$$n(\Gamma)=|\topp(\bb_1)|+|s_1|+|s_2|+|p_1|+|q|+|u|,$$

$$n(E)=|\topp(\bb_1)+|u|+|q|+|p_2|.$$

The number $|\topp(\bb_1)|$ is equal to the number of orange edges on $p_1$.

The path $u$ consists of $\Theta$-edges (recall that $\Lab(u)=\Lab(\topp(\bb_{4N-7}))$).
No maximal $\Theta$-band $\rr$ in $\Gamma$ starting on $u$
can end on $t_j$ where $j\ge 4N-7$: this follows from the fact that
$\Gamma-\Psi_1$ is symmetric to
the union of $\Sigma_j$, $j=4N-7, 4N-6, 4N-5, 4N-4$ (with $\bb_{4N-3}$ removed).
Thus it ends
either on
$t_j$ for some $j<4N-7$ or on $q$. In the first case the end of $\rr$ is a
red edge. In the second case let $e$ be the end edge of $\rr$. Since $q$ with
$\kappa$-letters removed is symmetric to the union of the $t_j$'s,
$j=4N-7,...,4N-4$,
there exists a symmetric edge $e'$ on one of these $t_j$. Then the maximal
$\Theta$-band $\rr'$ starting on $e'$ must cross $\bb_{4N-7}$, therefore it ends on a
red edge of one of the $t_j$'s, $j<4N-7$. The end edge of $\rr'$ is red. Thus we
have a one-to-one correspondence between edges of $u$ and some
red edges on $p_1$.
Therefore the length of $u$ does not exceed (in fact one can prove that
it equals) the number of red edges of $p_1$.

The length of $s_1s_2\iv$ is $4N|W|+4N$. Therefore by Lemma \ref{psi1}
$$|s_1|+|s_2|\le 4N(30k+12)|p_1|+4N<4N(30k+13)|p_1|.$$

Thus $$n(\Gamma)\le |p_1|+4N(30k+13)|p_1|+2|p_1|<\epsilon_1n(\Delta)$$
for some
$\epsilon_1$.


Notice that $n(\Delta)-n(E)=|p_1|-|\topp(\bb_1)|-|u|-|q|$.
Since the $|\topp(\bb_1)|$ is equal to the number of orange edges on $p_1$,
$|u|$ does not exceed the number of red edges
on $p_1$, $|p_1|-|\topp(\bb_1)|-|u|$
exceeds the number of green, yellow, blue and pink edges of $p_1$ plus the number
of colorless edges. By construction $|q|$ is equal to the sum of
the lengths
of the $t_j$'s, $j=4N-7, 4N-6, 4N-5, 4N-4$ plus 3 (the $\kappa$-edges).
Notice that the edges of these $t_j$ are colorless (because we removed the
red paint from these $t_j$ and picked $i<4N-7$
when we starting to paint $t_i$ in different colors). Therefore
$n(\Delta)-n(E)=|p_1|-|\topp(\bb_1)|-|u|-|q|$
exceeds the sum of the $n_i$'s for all odd
$i$ from $2N+1$ to $4N-11$, where $n_i$ is the number of non-red edges of $t_i$
plus the number of $\kappa$-edges between $t_1$
and $t_{4N-7}$ minus $3$. By Lemma \ref{count} each of
these $n_i$ is greater than $||W||/(61k)$.
The number of $\kappa$-edges between
$t_1$ and $t_{4N-7}$ is greater than 3, so
$$n(\Delta)-n(E)\ge ||W||/(61k).$$ This completes the proof of our lemma. $\Box$

The following proposition gives the upper bound for the Dehn function
of the group $G_{N}(\sss)$.

\begin{prop} \label{upper} If $w=1$ in $G_{N}(\sss)$ then there exists a \vk
diagram $\Delta(w)$ over $\pp_{N}(\sss)$ with boundary label $w$, area
$< c_1 \cdot T(c_2 \cdot |w|)^4$ and diameter $< c_1'\cdot T(c_2'\cdot |w|)^3$,
for some positive constants $c_1, c_2$ and $c_1', c_2'$.
\end{prop}

{\bf Proof.} By the \vk Lemma there exists a reduced diagram $\Delta$ over
$\pp_{N}(\sss)$ with boundary label $w$ and minimal possible number of hubs.
Then we can apply the snowman decomposition and obtain three sequences of
diagrams $E_1$,...,$E_s$, $\Gamma_1$,...,$\Gamma_{s-1}$,
$\Pi_1$,...,$\Pi_{s-1}$,
with the following properties:

\begin{enumerate}
\item[(S1)] $E_s = \Delta$.
\item[(S2)] $E_1$ contains no hubs.
\item[(S3)] For every $i=1,...,s-1$, we have:
either $(\Gamma_i, E_i)$ is a decomposition
of type 1 of $E_{i+1}$ (and $\Pi_i$ is empty);
or $(\Gamma_i, E_i, \Pi_i)$ is a
decomposition of type 2 of $E_{i+1}$ (in particular, $\Pi_i$ is a
computational disc).
\end{enumerate}

This implies that we can glue together the diagrams
$E_1, \Gamma_1,...,\Gamma_{s-1}, \Pi_1,...,\Pi_{s-1}$
so as to obtain
a diagram with boundary label $w$. Therefore in order to find
the upper bound  for the area of $\Delta$ we need to estimate
the sum of the areas of
$E_1, \Gamma_1,...,\Gamma_{s-1}, \Pi_1,...,\Pi_{s-1}$.

Property (D5) of decompositions of type 1, and part 2 of Lemma \ref{nice}
show that for every $i=1,...,s-1$,
$$n(E_i)<n(E_{i+1}).$$ Therefore $$n(E_i)\le |w|, \ i=1,...,s,$$
and $s<|w|$.

The diagram $E_1$ and for every $i$ the diagram $\Gamma_i$
contain no hubs.
So by Lemma \ref{nohub},
$$\area(\Gamma_i)\le
C_1n(\Gamma_i)^3$$
for some constant $C_1$
and $$\area(E_1)\le
2C_1n(E_1)^3.$$
We also know that $n(E_1)\le |w|$ and by Lemma \ref{nice}
and property (D5), \  $n(\Gamma_i)\le C_2n(E_{i+1})\le C_2|w|$,
\ for $i=1,...,s-1$ and some constant $C_2$.
Therefore
\begin{equation}
\begin{array}{l}
\area(E_1)+\area(\Gamma_1)+ \ldots +\area(\Gamma_{s-1}) \ \le \\
C_1n(E_1)^3+
C_1n(\Gamma_1)^3 +
C_1n(\Gamma_2)^2+ \ldots \\
+C_1(\Gamma_{s-1})^3 \ \le\\
C_1C_2s |w|^3\le C_3|w|^4.
\end{array}
\label{eqq1}
\end{equation}
for some constant $C_3$.

It remains to estimate the sum of the areas of discs $\Pi_i$.
We can assume that each $\Pi_1$ corresponds to a minimal area computation
of $\sss$ connecting some word $W_i$ and $W_0$.
Therefore the area of each $\Pi_i$
does not exceed ``big O'' of the area of the corresponding computation. By
Proposition \ref{mach1} the area of each $\Pi_i$
does not exceed $C_4T(||W_i||)^4 +C_5n(\Pi_i)^3$
for some constants $C_4$ and $C_5$.

$$\begin{array}{l}
\area(\Pi_1)+\area(\Pi_2)+ \ldots +\area(\Pi_{s-1})\ \le\\
C_4(T(||W_1||)^4+ \ldots +T(||W_{s-1}||)^4))+\\
C_5(n(\Pi_1)^3+n(\Pi_2)^3+...+n(\Pi_{s-1})^3) \ \le  \\
 C_4T(||W_1||+ \ldots +||W_{s-1}||)^4+ C_5(n(\Pi_1)^3+...+n(\Pi_{s-1})^3).
 \end{array}$$

Here we used the superadditivity of the function $T(n)^4$.
For every $i$ from $1$ to $s-1$, if $(\Gamma_i, E_i, \Pi_i)$ is a decomposition
of type 2 of $E_{i+1}$ then by Lemma \ref{nice} (2) we have
$$||W_i||\le C_6(n(E_{i+1})-n(E_i))$$
for some constant $C_6$. Therefore
$$||W_1||+...+||W_{s-1}|| \le C_6n(\Delta)=C_6|w|.$$

By Lemma \ref{psi1} $n(\Pi_i)\le C_8n(E_{i+1})\le C_8|w|$
for some constant $C_8$.
Therefore $$n(\Pi_1)^3+...+n(\Pi_{s-1})^3\le C_9|w|\cdot |w|^3=C_9 |w|^4.$$
Thus

\begin{equation}\label{eqq2}
\begin{array}{l}
\area(\Pi_1)+\area(\Pi_2)+ \ldots +\area(\Pi_{s-1})\ \le\\
C_4 T(C_6|w|)^4+C_{10}|w|^4.
\end{array}
\end{equation}
for some constant $C_{10}$.

Combining (\ref{eqq1}) and (\ref{eqq2}) we get

$$\area(\Delta)\le C_3|w|^4+C_4 T(C_6|w|)^4+C_{10}|w|^4\le C_{11}T(C_{12}|w|)^4.$$
for some constants $C_{11}$ and $C_{12}$.
This completes the area part of the proposition.

Now let us estimate the diameter $d(\Delta)$. If $\Delta=E_1$ that is if $\Delta$
has no hubs then we can apply Lemma \ref{nohub} and conclude that
$d(\Delta)\le O(|w|)$. If $\Delta$ contains hubs, we can use the
snowman decomposition and obtain three sequences of diagrams $E_1$,..., $E_s$,
 $\Gamma_1$, ...,$\Gamma_{s-1}$,
$\Pi_1$,...,$\Pi_{s-1}$ as above.  Let $d(\Gamma_i)$, $d(E_i)$, $d(\Pi_i)$
 be the diameters of the correspomding diagrams.
It is easy to see that
\begin{equation}\label{diameter}
\begin{array}{ll}
d(\Delta)\le & \hbox{max}_{i=1,...,s-1}\{d(E_{1})+\sum_{j=1}^{s-1}|\partial(\Gamma_{j})|+
\sum_{j=1}^{s-1}|\partial(\Pi_{s-1})|,\\
& d(\Gamma_i)+|\partial(E_1)|+
\sum_{j=1}^{s-1}|\partial(\Gamma_{j})|+
\sum_{j=1}^{s-1}|\partial(\Pi_j)|,\\
& |d(\Pi_i)|+|\partial(E_1)|+\sum_{j=1}^{s-1} |\partial(\Gamma_{j})|+
\sum_{j=1}^{s-1}|\partial(\Pi_j)|\}.
\end{array}
\end{equation}
Indeed, take a vertex inside the diagram $\Delta'$ obtained by gluing together
$E_1$, all $\Gamma_i$'s and all $\Pi_i$'s. It belongs to one of these
blocks. Take the shortest path to the boundary of this block, then we can get
to the boundary of $\Delta'$ going along the boundaries of the blocks.
The resulting path will have length bounded by the right part of
(\ref{diameter}).

We already know the upper bounds for
$|\partial(\Gamma_{j})|$ and $|\partial(\Pi_j)|$.
By Lemma \ref{nice} these do not exceed $O(n(\Delta))$.
Since $s\le n(\Delta)$, the sums
$\sum_{j=1}^{s-1}|\partial(\Gamma_{j})|$ and
$\sum_{j=1}^{s-1}|\partial(\Pi_j)|$
do not exceed $O(n(\Delta)^2)$.
By Lemma \ref{nohub}, the diameters of diagrams $E_1$ and $\Gamma_i$\ ($i=1,...,s-1$)
do not exceed ``big O" of their perimeters. By Lemma \ref{nice} these
perimeters are smaller than $O(n(\Delta))$. Thus the diameters of
$E_1$ and $\Gamma_i$ do not exceed $O(n(\Delta))$. By Proposition
\ref{prop13}, the diameter of the disc $\Pi_i$ does not exceed
$O(T(|\partial(\Pi_i)|)^3)$. We also know that
$|\partial(\Pi_i)|\le O(n(\Delta))$.
Therefore the diameter of
$\Pi_i$ does not exceed
$C_1'T(C_2'|w|)^3$ for some constants $C_1'$ and $C_2'$. Combining this
information with (\ref{diameter}) we conclude that
$$d(\Delta)\le  C_1'T(C_2'|w|)^3+C_3'n^2$$
for some constant $C_3$. Since $T(n)\ge O(n)$, we deduce that
$$d(\Delta)\le  c_1'T(c_2'|w|)^3$$
for some constants $c_1'$ and $c_2'$ as desired.

The proposition is proved.

$\Box$

\section{The Lower Bound}
\label{lowbound}

The goal of this section is to prove that the group $G_N(\sss)$ simulates the
machine $\sss$, that the Dehn function of $G_N(\sss)$ is bounded below
by a function equivalent to $T(n)^4$, and that the smallest isodiametric
function is bounded below by a function equivalent to $T(n)^3$.

We start with the following lemma.

\begin{lm} \label{lb1} Fix an admissible word $W$ for $\sss$.
We call a diagram $\Delta$ over the presentation $\pp_N(\sss)$
{\em quasi-sector
for $W$}
if the boundary of $\Delta$ can be represented in the  form $b_1s_1b_2\iv s_2\iv$
where $b_1$ is the top path of a $\kappa_1$-band $\bb_1$, $b_2$ is the bottom
path of a $\kappa_2$-band $\bb_2$, and both bands start on $s_2$;
$\Lab(s_1)=\kappa_1W\kappa_2$, $\Lab(s_2)=\kappa_1W_0\kappa_2$;
all hubs of $\Delta$ are on $\bb_1$.

1. Suppose that $\Delta$ is a quasi-sector for $W$ which either has the smallest
possible area among all quasisectors for $W$ or has the smallest length
of $\bb_1$ among all quasi-sectors for $W$.
Then $\Delta$ is a sector.

2. Suppose that $\Delta$ is a quasi-sector for $W$. Then the label of $\bott(\bb_2)$
is the history word of a computation of $\sss$ connecting $W$ and $W_0$.
\end{lm}

{\bf Proof.} Suppose that $\bb_1$ does not have hubs.
Then since $s_1$ and $s_2$ do not have $\Theta$-edges, every maximal
$\Theta$-band
which starts on $b_1$  ends on $b_2$ and every maximal $\Theta$-band
which starts on $b_2$ ends on $b_1$. This implies that $s_1$ is the top path
of a $\Theta$-band starting on $b_1$ and ending on $b_2$. Thus by definition,
$\Delta$ is a sector.

Now suppose that $\bb_1$ contains hubs. Let us number these hubs
$\pi_1,...,\pi_n$ along $\bb_1$. Let $\pi_i$ be one of these hubs.
The path $b_1$ contains exactly one $\kappa_{N+1}$-edge from
$\partial(\pi_i)$.
Let $\kkk_i$ be the maximal $\kappa_{N+1}$-band of $\Delta$
starting on this edge. Since $\Delta$ does not have hubs outside $\bb_1$,
each $\kkk_i$ ends on the boundary of some hub $\pi_j$.
Since $\kappa_{N+1}$-bands
cannot intersect, it is easy to see that there must exist a number $i$ from
$1$ to $n-1$ such that $\kkk_i$ ends on $\pi_{i+1}$.

\vskip 0.2 in
\unitlength=1.00mm
\special{em:linewidth 0.4pt}
\linethickness{0.4pt}
\begin{picture}(76.33,77.33)
\put(26.67,70.67){\line(1,-6){11.33}}
\put(38.00,2.67){\line(1,0){21.33}}
\put(59.33,2.67){\line(1,4){17.00}}
\put(76.33,70.67){\line(-1,0){49.67}}
\put(40.33,2.67){\line(-1,6){11.33}}
\put(57.00,2.67){\line(1,4){17.00}}
\put(51.00,75.00){\makebox(0,0)[cc]{$W$}}
\put(31.00,51.33){\circle*{4.06}}
\put(33.33,38.00){\circle*{4.22}}
\put(35.67,23.00){\circle*{3.89}}
\put(38.00,9.00){\circle*{4.06}}
\bezier{100}(31.00,51.33)(42.00,48.67)(33.33,37.67)
\bezier{220}(31.00,52.00)(57.00,54.00)(33.67,37.33)
\bezier{344}(31.00,52.33)(71.00,62.00)(34.00,36.00)
\put(27.33,77.33){\makebox(0,0)[cc]{$\bb_1$}}
\put(75.33,75.33){\makebox(0,0)[cc]{$\bb_2$}}
\put(47.67,6.67){\makebox(0,0)[cc]{$W_0$}}
\put(51.00,57.67){\makebox(0,0)[cc]{$\kkk_i$}}
\end{picture}
\begin{center}
\nopagebreak[4]
Figure \theppp.
\setcounter{pdtwentysix}{\value{ppp}}

\end{center}
\addtocounter{ppp}{1}
\vskip 0.2 in

Consider the subdiagram $\Sigma$ of $\Delta$ bounded by $\bb_1$,
$\kkk_i$, and the contours of the two hubs $\pi_i$ and $\pi_{i+1}$ such that
$\Sigma$ does not contain $\pi_i$ and $\pi_{i+1}$ (see Fig. \thepdtwentysix).
Let $\partial(\Sigma)=b_1's_2'(b_2')\iv (s_1')\iv$ be the corresponding
division of the boundary of $\Sigma$ where $b_1'\subset b_1$, $s'_1\subset
\partial(\pi_i)$, $s'_2\subset \partial(\pi_{i+1})\iv$, $b'_2 \subset
\bott(\kkk_i)$.
The diagram $\Sigma$ does not contain hubs and has $N+1$ $\kappa$-bands
$\ttt_1, ..., \ttt_{N+1}$, connecting edges on $s_1$ with edges on $s_2$
where $\ttt_1\subset \bb_1$, $\ttt_{N+1}=\kkk_i$.
An argument similar to that used in the proof of Lemma \ref{bigon} shows
that $\Sigma$ is a union of $N$ sectors and inverses of sectors
corresponding to the same computation $C$. By Proposition \ref{prop12}
we can assume that these are computational sectors.
As in the proof of Lemma \ref{bigon} consider the computational disc $\Pi$ corresponding
to the computation $C$. We can identify $\Sigma$ with a part of $\Pi$
and $\pi_i$ with the hub of $\Pi$.
 Let $\Gamma=\Pi - \Sigma - \pi_i$. Cut $\Delta$ along the
path $(b_1')\iv\bar s_1'b_2'$ where $\bar s_1\iv s_1=\partial(\pi_i)$,
and glue in the resulting hole the diagrams $\Gamma$ and $\Gamma\iv$.
Then $\Pi$ becomes a subdiagram of the resulting diagram $\Delta'$.
Then again as in the proof of Lemma \ref{bigon} we notice that the boundary
label
of $\Pi$ is the same as the boundary label of
the hub, so we can replace the subdiagram
$\Pi$ in $\Delta'$ by the hub $\pi_i$. After that $\pi_i$ and $\pi_{i+1}$
will form a cancellable pair of cells, so we can cancel them.
Let $\bar\Delta$ be the resulting diagram. Let $\bar\bb_1$ be the $\kappa_1$-band
in $\bar\Delta$ connecting $s_1$ and $s_2$. Finally let $\tilde \Delta$ be the
subdiagram of $\bar\Delta$ bounded by $\bar\bb_1$, $s_1$, $s_2$, $\bb_2$.
It is easy to see that the area of $\tilde\Delta$ is smaller than the area of
$\Delta$ and the length of $\bar\bb_1$ is smaller than the length of $\bb_1$.
Indeed, $\tilde\Delta$ is obtained from $\Delta$ by removing
$N$ sectors of $\Pi$ and two hubs and inserting $N$ sectors of $\Pi$.
Therefore $\area(\tilde\Delta)\le \area(\Delta)-2$; the band $\bar\bb_1$
is obtained from the band $\bb_1$ by removing two hubs, so the $|\bar\bb_1|=
|\bb_1|-2$.  This contradicts the assumption that $\Delta$
is a quasi-sector for $W$ which either has the smallest
possible area among all quasisectors for $W$ or has the smallest length
of $\bb_1$ among all quasi-sectors for $W$.
This comtradiction proves part 1 of the lemma.

In order to prove part 2 let $\Delta$ be a quasi-sector for $W$. Then
reducing the number of hubs of $\bb_1$ as above we can transform $\Delta$
into a sector $\tilde\Delta$ with $W$ as the label of top path.
By Proposition \ref{prop12} we can assume that $\tilde\Delta$ is a computational
sector.
Notice that our transformations
do not affect $\bb_2$. Therefore the label of $\bott(\bb_2)$ is the history
word of the computation corresponding to the sector $\tilde\Delta$.
$\Box$

\bigskip

The next proposition shows that the group $G_N(\sss)$ simulates the machine
$\sss$ and provides a lower bound for the Dehn function of this group.

\begin{prop} \label{low} 1. For every admissible word $W$ of the machine $\sss$
there exists a computation of $\sss$ connecting $W$ and $W_0$ if and only if
$K(W)=1$ in the group $G_N(\sss)$.

2. The Dehn function of $G_N(\sss)$ is
bounded below by a function equivalent to $T(n)^4$.
\end{prop}

{\bf Proof.}  If there exists a computation of $\sss$ connecting $W$
and $W_0$ then the computational disc corresponding to this computation
has boundary label $K(W)$, so $K(W)=1$ in the group $G_N(\sss)$.

Conversely suppose that $K(W)=1$ in $G_N(\sss)$. Then there exists a
\vk diagram $\Delta$ over the presentation $\pp_N(\sss)$ with boundary label
$K(W)$. Our goal is to prove that there exists a computation connecting
$W$ and $W_0$ and to estimate the area and diameter of $\Delta$ from below.
Thus we
assume that the area of $\Delta$ is minimal. Then $\Delta$ is reduced.
Let $e_i$, $i=1,...,4N$ be $\kappa_i^{\pm 1}$-edges on $\partial(\Delta)$, we
assume that $e_1$ precedes $e_2$,..., precedes $e_{4N}$
on $\partial(\Delta$). Let $\kkk_i$, $i=1,...,2N$ be the maximal
$\kappa_i$-band
starting on $e_i$. Since the contour of $\Delta$ contains two $\kappa_1$
edges and two $\kappa_2$-edges, $\kkk_1$ and $\kkk_2$ intersect. Therefore
$\Delta$ contains hubs.

\vskip 0.1 in
{\bf Remark.} Now it is tempting to use the graph theoretic
lemmas from the previous section, to
deduce that $\Delta$ contains just one hub, and then prove that
$\Delta$ is a disc.
But unfortunately we cannot use the lemmas from the previous section because
we cannot assume that $\Delta$ has minimal number of hubs:
this assumption may contradict the minimal area assumption or the minimal diameter
assumption. Thus our proof
proceeds in a different direction.
\vskip 0.1 in

Notice that for every $i$ from 1 to $2N$, $\Delta$ contains only two mutually
inverse $\kappa_i$-bands. Therefore each hub in $\Delta$ belongs to all
$\kkk_i$.

Let $\pi_1$,...,$\pi_m$ be all the hubs in $\Delta$, numbered along $\kkk_1$.
Let $\pi_n$ be the first of these hubs appearing in $\kkk_2$. Consider the
subdiagram $\Gamma$ of $\Delta$ bounded by $\kkk_1$, $\partial(\Delta)$,
$\kkk_2$ and $\partial(\pi_n)$ such that $\pi_n\not\in\Gamma$. It is easy to
see that $\Gamma$ is a quasi-sector for $W$.

By Lemma \ref{lb1} there exists a computational
sector $\Sigma$ with top path labelled by
$\kappa_1W\kappa_2$ and the bottom path labelled by $\kappa_1W_0\kappa_2$
such that $\area(\Sigma)\le \area(\Gamma)$.
By Proposition \ref{prop12} there exists a computation of $\sss$
connecting $W$ and $W_0$. This proves statement 1 of our proposition.

Let $S(n)$ be the space function of the Turing machine $M$. We
know that $S(n)$
is equivalent to $T(n)$.
Using the notation from Proposition \ref{mach1} let us assume that
$W=\sigma(c)$ where $c$ is an accepted configuration of $M$
of length $\le n$
such that any smallest space
accepting computation of $M$ for the configuration
$c$ has space $S(n)$. Then by Proposition \ref{mach1}, any computation
of $\sss$ connecting $W$ and $W_0$ has area exceeding $CS^4(n)$ and time exceeding
$C'S^3(n)$ for some positive constants
$C$, $C'$. Therefore by Proposition \ref{prop12} the area of $\Sigma$ exceeds
$O(S^4(n))$ and the diameter exceeds $O(S^3(n))$.

Since $\area(\Delta)\ge \area(\Gamma)\ge \area(\Sigma)$,
the area of $\Delta$ exceeds $O(S(n)^4)$. Therefore $S(n)^4$ is equivalent
to a lower bound for
the Dehn function of $G_N(\sss)$. Since $S^4(n)$ is equivalent to $T(n)^4$,
$T(n)^4$ is also equivalent to a lower bound of the Dehn function
of $G_N(\sss)$. $\Box$

\bigskip

The next proposition provides a lower bound for isodiametric functions of
$G_N(\sss)$.

\begin{prop}\label{low1}
Every isodiametric function of $G_N(\sss)$ is bounded below by
a function equivalent to $T(n)^3$.
\end{prop}

{\bf Proof.} Let ${\cal T(n)}$ be the ``time function'' of the machine $\sss$, that
is $\ttt(n)$ is the smallest number such that for every admissible word $W$
of length $\le n$, for which there exists a computation connecting it with $W_0$,
the smallest length of this computation is $\le \ttt(n)$. By the choice of
$\sss$, $\ttt(n)$ is equivalent to $T(n)^3$.

Let us fix a number $n\ge 1$ and let $W$ be an adsmissible word
for the $S$-machine $\sss$, $|W|\le n$.
Suppose that there exists a computation connecting
$W$ and $W_0$ and the smallest length of this computation is $t(n)$.
Then there exist a diagram (for example, computational disc) over
$\pp_N(\sss)$ with
boundary label $K(W)$. Let $\Delta$ be any
(not necessarily reduced) diagram over
the presentation $\pp_N(\sss)$ with boundary label $K(W)$. 
Assume that $\Delta$ has the smallest diameter among all such diagrams.

By Lemma \ref{cent} the operation of removing $\kappa$-annuli does not
increase the  diameter of $\Delta$. So we can assume that $\Delta$ does not
contain $\kappa$-annuli. This implies as before that every hub in $\Delta$
belongs to every $\kappa$-band and that $\Delta$ contains 2 mutually inverse
maximal $\kappa_i$-bands for every $i=1,...,2N$.
As in the proof of Proposition \ref{low1},
let $e_i$, $i=1,...,4N$ be $\kappa_i^{\pm 1}$-edges on $\partial(\Delta)$, we
assume that $e_1$ precedes $e_2$,..., precedes $e_{4N}$
on $\partial(\Delta$). Let $\kkk_i$, $i=1,...,2N$ be the maximal
$\kappa_i$-band
starting on $e_i$. Since the contour of $\Delta$ contains two $\kappa_1$
edges and two $\kappa_2$-edges, $\kkk_1$ and $\kkk_2$ intersect. Therefore
$\Delta$ contains hubs.

Let $\pi_1$,...,$\pi_m$ be all the hubs in $\Delta$, numbered along $\kkk_1$.
Let $\pi_n$ be the first of these hubs appearing in $\kkk_2$. As before,
consider the
subdiagram $\Gamma$ of $\Delta$ bounded by $\kkk_1$, $\partial(\Delta)$,
$\kkk_2$ and $\partial(\pi_n)$ such that $\pi_n\not\in\Gamma$.

Let $\bb_1$ (resp. $\bb_2$) be the $\kappa_1$-band (resp. $\kappa_2$-band)
of $\Gamma$ starting on the contour of $\Delta$.
Let $f_i$ be the end edge of $\bb_i$ ($i=1,2$).
Let $v$ be the initial vertex of the $\kappa_2$-edge $f_2$
and let $p$ be the
shortest path connecting $v$ with the boundary of $\Delta$.

In every pair of mutually inverse
maximal $\Theta$-bands in $\Delta$, call one band {\em positive} and the other
one {\em negative}. Consider the set $\Omega$
of all positive maximal $\Theta$-bands of $\Delta$ intersecting $\bb_2$.
At most $|\partial(\Delta)|$ of these bands start and end on the boundary
of $\Delta$. Therefore $|\Omega|-|\partial(\Delta)|$ bands in $\Omega$ are
$\Theta$-annuli. Our goal is to calculate the number of $\Theta$-annuli in $\Omega$
which contain $v$ in their inside diagrams. It is clear that the length
of $p$ cannot be smaller than this number.

Notice that given $\Delta$, $\Gamma$, $\bb_2$, $v$ and $|p|$ are determined
uniquely. Without loss of generality
we assume that $\Delta$ has the smallest possible pair $(|p|, |\bb_2|)$ (in
the lexicographic order) among all diagrams with boundary label $K(W)$.

Suppose that $p$ contains at least two $\kappa_2^{\pm 1}$-edges of $\bb_2$.
We shall show that this leads to a contradiction.

Then $p$ contains two edges $e_1$ and $e_2$ which are common
edges of $\bb_2$ or their inverses and such that between $e_1$ and $e_2$,
$p$ does not intersect the median of $\bb_2$. Let $t_1$ be the initial vertex
of $e_1$ and let $t_2$ be one of the vertices of $e_2$ which belong to the
same side of $\bb_2$ (top ot bottom) as $t_1$. Then there exists a reduced
path $q$ which connects $t_1$ and $t_2$ and
is a subpath of $\topp(\bb_2)^{\pm 1}$ or $\bott(\bb_2)^{\pm 1}$.

{\bf Claim.} $\Lab(q)$ is a reduced word.

\bigskip

Indeed, suppose that $q$ contains a subpath $ee'$ where $e$ and $e'$
are edges with mutually inverse
$\Theta^{\pm 1}$-labels. Then the corresponding cells $\rho_1$ and $\rho_2$
of $\bb_2$ form a reducible
pair. The common edge of this pair does not belong to $p^{\pm 1}$ because
$p$ does not have $\kappa_2$-edges from
$\bb_2$ between $e_1$ and $e_2$. Therefore
reducing this pair does not affect the length of $p$ (if this  edge
belonged to $p$ then reducing this pair of cells would make the shortest
path connecting $v$ and $\partial(\Delta)$ longer). When
we reduce $\rho_1$ and $\rho_2$, we get a diagram with
smaller pair $(|p|, |\bb_2|)$
which contradicts our choice of $\Delta$.

\bigskip

Now let $p'$ be the portion of $p$
which connects $t_1$ and $t_2$. Then $(q\iv p')^{\pm 1}$ bounds a subdiagram of
$\Delta$. Therefore $\Lab(p')=\Lab(q)$ in $G_N(\sss)$. There exists a
homomorphism from $G_N(\sss)$ to the free group generated by $\Theta$ which
takes all non-$\Theta$-letters to 1. Therefore if we remove all
non-$\Theta$-letters
from $p'$, we get a word $U$ which is freely equal
to $\Lab(q)$. Since by Claim $\Lab(q)$ is reduced, $|U|\ge |q|$.
Notice that $\Lab(p')$ contains
a non-$\Theta$-edge $e_1$. Therefore $|p'|$ is strictly greater than $|q|$.
Therefore we can substitute $p'$ by $q$ in $p$ and get a shorter
path connecting $v$ with $\partial(\Delta)$, a contradiction with the
choice of $p$. This contradiction shows that $p$ contain at most one $\kappa_2^{\pm 1}$-edge
of $\bb_2$.

\bigskip

Suppose that $p$ contains one $\kappa_2^{\pm 1}$-edge $e$ of $\bb_2$. One of the
vertices $v'$ of $e$ belongs to the same side (top or bottom) of $\bb_2$ as $v$.
Let $q$ be the subpath of $\bott(\bb_2)\iv$ connecting $v$ and $v'$. As in the
proof of Claim one can show that $\Lab(q)$ is a reduced word and that the part
$p'$ of the path $p$ connecting $v$ and $v'$ is not shorter than $q$.  Let $\bar
v$ be the end vertex of the $\kappa_2$ edge $f_2$ and let $\bar v'$ be the other
vertex of $e$. Then the vertices $\bar v$ and $\bar v$ are connected by a
subpath $\bar q$ of $\topp(\bb_2)\iv$. Notice that as before we can assume that
the word $\Lab(\bar q)$ is reduced. So the lengths of $q$ and $\bar q$ are
equal.  Thus we can substitute the subpath $p'$ of $p$ by $f_2\bar q$ and get
another path, say, $\bar p$, which connects $v$ with the boundary of $\Delta$,
has length not bigger than the length of $p$, and contains only one common
$\kappa_2$-edge with $\bb_2$, the edge $f_2$. Thus we can assume without loss of
generality that the path $p$ has these properties.

Now if $\bb_2$ contains reducible pairs of cells we can cancel them without
affecting the length of $p$.

In the case when $p$ does not have common $\kappa^{\pm 1}$-edges with $\bb_2$,
we can reduce $\bb_2$ without changing the length of $p$.

Thus we can assume that $\bb_2$ is reduced. Now we can apply Lemma \ref{cr}
and deduce that no $\Theta$-band from $\Omega$ can have two intersections
with $\bb_2$ (otherwise $\bb_2$ would not be reduced).
Thus the number of $\Theta$-bands in $\Omega$ is $|\bb_2|$.

Every $\Theta$-annulus
in $\Omega$ must contain $v$ in its inner subdiagram. Indeed,
the bottom path of $\bb_2$ crosses this annulus exactly once, so the index
of the median of this annulus (which is by definition a simple closed curve)
with respect to the point $v$ is 1.

Therefore the number
of annuli in $\Omega$ is at least $|\bb_2|-|\partial(\Delta)|$. Let us reduce
the diagram $\Delta$ by first reducing the band $\bb_1$ and then reducing other
reducible pairs of cells. Let $\bar\Gamma$
be the reduced diagram obtained as a result of this process.
Then $\bar\Gamma$ is a quasi-sector. By Lemma \ref{kbands} the band
$\bb_2$ does not change during the reduction process in $\Gamma$. Therefore by
Lemma \ref{lb1} $\Lab(\bb_2)$ is a history of a computation connecting
$W$ with $W_0$. Therefore by $|\bb_2|\ge \ttt(n)$. Therefore the
number of $\Theta$-annuli in $\Delta$ containing $v$ is not smaller than
$\ttt(n)-|K(W)|\ge \ttt(n)-4Nn-4N$. Hence the length of $p$ is at least
$\ttt(n)-4Nn-4N$. Since the length of $\partial(\Delta)$ does not exceed
$4Nn+4N$, we deduce that the isodiametric function of $G_N(\sss)$
is $\succ \ttt(n)\equiv T(n)^3$.
$\Box$

\section{Proof of Theorem 1.3}

The proof of Theorem \ref{2} can be completed as follows. Let $L\subseteq X^+$
be a language recognized by a Turing machine with time function $T(n)$.
Then there exists a symmetric nondeterministic
Turing machine $M$, with all the properties of Lemma \ref{mach23},
which recognizes the language $L$. Take any $N\ge 6$ and construct the
group $G_{N}(\sss)$. Let $D(n)$ be the Dehn function of this group. Then by
Proposition \ref{upper}, $D(n)\le C_1 T(C_2 n)^4$, for some constants $C_1$
and $C_2$.
By Proposition \ref{low}, $D(n)> c_1 T(c_2 n)^4$, for some positive
constants $c_1, c_2$.
Therefore $D(n)$ is equivalent to $T(n)^4$.

Similarly we can prove that the smallest isodiametric function of $\Delta$
is equivalent to $T(n)^3$.

Let $w$ be an input word for the machine $M$ and $c_w$ be the corresponding
input configuration of the machine $M$. Then
by Proposition \ref{low}
the correspondence
$w\to K(\sigma(c_w))$
satisfies the conditions of Theorem \ref{2}, in
particular $w\in L$ if and only if $K(\sigma(c_w))=1$ in $G_{N}(\sss)$.
Theorem \ref{2} is proved.

Theorem \ref{1} follows from Theorem \ref{2} and Theorem
\ref{th111}.

\section{Proof of Corollary 1.2}

This corollary almost immediately follows from Theorem \ref{2} and its proof.
Take any Turing machine $M$ which accepts a non-recursive language.
We can assume that $M$ satisfies the conditions
of Lemma \ref{mach23}.
Using Proposition \ref{mach1} we construct the $S$-machine $\sss=\sss(M)$
for which there is no algorithm deciding whether there exists a computation
of $\sss$ connecting a given admissible word $W$ with the fixed admissible
word $W_0$.
Now take $N\ge 6$ and consider
the group $G_N'(\sss)$ given by the presentation $\pp'_N(M)$ which is
obtained from
$\pp_N(M)$ by removing the hub relations.

In the group $G_N'(\sss)$
a word $K(W)$ where $W$ is an admissible word for $\sss$
is conjugate with $K(W_0)$ if and only if
$W$ is connected with $W_0$ by a computation of $\sss$.
Indeed, if such a computation exists then by Proposition \ref{prop13}
there exists a disc with boundary label $K(W)$.
By removing the hub from this disc we obtain an annular diagram
over the presentation $\pp'_9(M)$ which shows that $K(W)$
is conjugate to $K(W_0)$ (see \cite{LS}). On the other hand if
$K(W)$ is conjugate to $K(W_0)$ then there
exists an annular diagram over $\pp'_N(\sss)$ with boundary
labels $K(W)$ (the outer contour)
and $K(W_0)$ (the inner contour).
We can glue the hub in the hole of this annular diagram and obtain a diagram
with one hub over $\pp_N(\sss)$ which by definition  is a disc.
By Proposition \ref{prop13} then $W$ is connected with $W_0$ by a computation
of $\sss$.

Since the language accepted by $M$ is not recursive, $G'_9(\sss)$ has
undecidable
conjugacy problem. It remains to notice that by Lemma \ref{nohub}
the Dehn function of $G'_9(M)$ is at most $O(n^3)$ (indeed, the \vk diagrams
over $\pp'_N(\sss)$ are precisely the \vk diagrams over $\pp_N(\sss)$
without hubs).

Notice also that the group we constructed has linear isodiametric function
by Lemma \ref{nohub}.
$\Box$


\noindent Mark V. Sapir\\
Department of Mathematics and Statistics\\
University of Nebraska-Lincoln\\
http://www.math.unl.edu/$\sim$msapir\\

\noindent Jean-Camille Birget\\
Department of Computer Science\\
University of Nebraska-Lincoln\\
birget@cse.unl.edu\\

\noindent Eliyahu Rips\\
Department of Mathematics\\
Hebrew University of Jerusalem\\
rips@math.huji.ac.il

\end{document}